Національний технічний університет України
"Київський політехнічний інститут"

На правах рукопису

Блажієвська Ірина Петрівна

УДК 519.21

# Кумулянтні методи в задачах оцінювання перехідних функцій однорідних лінійних систем

01.01.05 - теорія ймовірностей і математична статистика

Дисертація на здобуття наукового ступеня кандидата фізико-математичних наук

Науковий керівник -
Булдигін Валерій Володимирович,
доктор фізико-математичних наук, професор;
Клесов Олег Іванович,
доктор фізико-математичних наук, професор

Київ - 2015

# Зміст

# ВСТУП

**Актуальність теми.** Предметом активних досліджень протягом останніх п'ятидесяти років є задача оцінювання або ідентифікації характеристик лінійних систем різноманітної фізичної природи, яка виникає у таких галузях, як радіофізика, гідролокація, сейсмографія, метеорологія, розвиток і взаємодія біологічних популяцій, теорія сигналів та автоматичного контролю, теорія фільтрації, економетрика, фінансова математика, тощо. В якості посилань, наприклад, можна звернутись до робіт Х. Акайка, Д. Бокса, Д. Бриллінджера, Е. Вен-тцеля, М. Дейстлера, Дж. Дженкінса, К. Жанга, В. Зайця, П. Кайнза, Л. Льюнга, Л.Овчарова, Дж. Сйоберга, П. Стойка, Е. Ханана, М. Шетцена та
Т. Шодерстрома.

Різноманітні статистичні методи оцінювання невідомих імпульсних перехідних функцій лінійних систем типу "чорної скриньки" й вивчення властивостей відповідних оцінок розглядались у роботах М. Арато, Х. Акайка, Е.-В. Баї, Дж. Бокса, Д. Бриллінджера, П. Броквелла, В. Булдигіна, Д. Ваттса, Е. Волтера, Р. Девіса, Г. Дженкінса, А. Дороговцева, В. Зайця, П. Кайнза, А. Кукуша, В. Курочки, Фу Лі, Г. Льюнга, М. Міланезе, Дж. Нортона, Ф. Уцета, М. Шетцена та інших.

У класі лінійних систем важливий підклас складають неперервні однорідні системи. Більшість методів дослідження властивостей статистичних оцінок в цьому випадку ґрунтуються на збуренні системи стаціонарними в широкому сенсі процесами. В якості відповідних оцінок використовують емпіричні періодограми (для стійких систем) або емпіричні корелограми (для нестійких систем) між процесами на вході та виході системи.

При збуренні однорідної лінійної системи гауссівським білим шумом, корелограмна оцінка для перехідної функції є незсуненою; дослідження її асимптотичних властивостей спирається на відомі факти про сумісні корелограми сумісно гауссівських та сумісно стаціонарних випадкових процесів. Взагалі, корелограми та сумісні корелограми для різноманітних класів випадкових процесів та полів детально вивчались у роботах Т. Андерсена, Р. Бенткуса, Д. Брилліджера, В. Булдигіна, О. Дем'яненко, О. Диховичного, В. Зайця, І. Ібрагімова, О. Іванова, Ю. Козаченка, М. Леоненка, А. Портнової, Ю. Розанова, Дж. Соле-Казалса, А. Стаднік, Ф. Уцета.

За умови, що на вхід однорідної лінійної системи подається сім'я сепарабельних центрованих стаціонарних гауссівських процесів, "близьких" до білого шуму, корелограмному оцінюванню імпульсної перехідної функції присвячені роботи В. Булдигіна, В. Зайця, В. Курочки, Фу Лі та Ф. Уцета.

Так, В. Булдигін та Фу Лі досліджують умови асимптотичної нормальності для відповідної інтегральної оцінки та похибки оцінювання, як у сенсі збіжності скінченновимірних розподілів, так і в сенсі збіжності розподілів у просторі неперервних функцій. Автори припускають, що перехідна функція $H \in L_2(\mathbb{R})$, а її частотна характеристика $H^*$ задовольняє двом умовам: $H^*$ неперервна майже всюди на $\mathbb{R}$ (відносно міри Лебега); $H^* \in L_1(\mathbb{R}) \bigcap L_\infty(\mathbb{R})$ або $H^* \in L_{2p}(\mathbb{R})$, якщо існує таке $p > 2$. Дослідження спирається на метод моментів та граничні властивості багатовимірних згорткових інтегралів.

Для корелограмної дискретної за часом оцінки у роботах В. Булдигіна, В. Зайця, В. Курочки та Ф. Уцета встановлюються умови асимптотичної незсуненості та конзистентності, а також асимптотичної нормальності як у сенсі збіжності скінченновимірних розподілів, так і в сенсі збіжності відповідних розподілів у просторі неперервних функцій. Результати цієї робо-ти вимагають мінімальне обмеження на порядок інтегрованості перехідної функції $H \in L_2(\mathbb{R})$. Застосування кумулянтного методу аналізу розподілу та властивостей багатовимірних інтегралів з циклічним зачепленням ядер сприяло отриманню результату, який в даній схемі поліпшенню не підлягає.

При даній постановці задачі, в дисертаційній роботі встановлюються: асимптотична незсуненість та консистентність інтегральної корелограмної оцінки; нові умови асимптотичної нормальності оцінки та похибки у сенсі збіжності скінченновимірних розподілів; умови асимптотичної нормальності оцінки та похибки в сенсі слабкої збіжності розподілів у просторі неперервних функцій. Всі результати отримані за умови $H \in L_2(\mathbb{R})$.

У дисертації для доведення слабкої збіжності оцінки та похибки використовуються кумулянти старших порядків, як, зокрема, і в роботах Д. Бриллінджера, Дж. Гріммета, Д. Нуаларта, Дж. Пекатті, Дж. Таггу, У. Хаберзеттла, тощо. Ці кумулянти зображаються у виді скінченної суми інтегралів, які містять циклічне зачеплення ядер та залежать від параметрів. Застосовуються верхні оцінки для таких інтегралів та властивості їх збіжності до нуля, запропоновані В. Булдигіним, В. Зайцем та Ф. Уцетом.

Відповідні покращення у схемі оцінювання ілюструються на численних прикладах. Так, для імпульсної перехідної функції

$$H(t) = \begin{cases} \dfrac{1}{(1+t)^\gamma}, & t \geq 0; \\ 0, & t < 0, \end{cases} \text{ де } \mu > 0 \text{ та } \gamma \in \left(\dfrac{1}{2}, \dfrac{3}{4}\right],$$

в околі особливих точок для $H^*$ має місце точна асимптотика:

$$|H^*(\lambda)| \sim \frac{C_1}{|\lambda|^{1-\gamma}}, \text{ при } |\lambda| \to 0;$$

$$|H^*(\lambda)| \sim \frac{C_2}{|\lambda|}, \text{ при } |\lambda| \to \infty,$$

де $C_k \notin \{0, \infty\}, k = 1, 2,$ - деякі додатні сталі.

Легко бачити, що $H^* \notin L_\infty(\mathbb{R}) \bigcap L_1(\mathbb{R})$, і $H^* \notin L_{2p}(\mathbb{R})$ при $p > 2$.

**Зв'язок роботи з науковими програмами, планами, темами.** Дисертаційна робота виконана на кафедрі математичного аналізу та теорії ймовірностей Національного технічного університету України "Київський політехнічний інститут" в рамках держбюджетних науково-дослідних робіт № 2200Ф "Дослідження актуальних проблем теорії випадкових процесів, математичного аналізу та крайових задач математичної фізики" (держреєстрація № 0109U001227) та № 2500Ф "Псевдорегулярні та спеціальні функції і їх застосування до задач стохастичного аналізу" (держреєстрація № 0112U001588). Тому напрям досліджень даної дисертаційної роботи є актуальним, своєчасним та перспективним.

**Мета та задачі дослідження.** Метою роботи є застосування кумулянтних методів при встановленні асимптотичної нормальності корелограмної оцінки та похибки оцінювання імпульсної перехідної функції однорідної лінійної нестійкої системи.

У роботі розглядаються такі задачі:

- отримання умов асимптотичної нормальності оцінки імпульсної перехідної функції у сенсі збіжності скінченновимірних розподілів;

- отримання умов асимптотичної нормальності оцінки імпульсної перехідної функції у сенсі слабкої збіжності розподілів у просторі неперервних функцій;

- отримання умови асимптотичної незсуненості оцінки імпульсної перехідної функції;

- отримання умов асимптотичної нормальності похибки оцінювання оцінки імпульсної перехідної функції у сенсі збіжності скінченновимірних розподілів;

- отримання умов асимптотичної нормальності похибки оцінювання імпульсної перехідної функції у сенсі слабкої збіжності розподілів у просторі неперервних функцій;

- отримання прикладів вхідних сигналів, внутрішнього шуму та імпульсних перехідних функцій.

**Об'єктом дослідження** є невідома дійснозначна імпульсна перехідна функція $H = (H(\tau), \tau \in \mathbb{R})$ неперервної однорідної лінійної системи. Розглядається ситуація, коли на вхід системи подається сім'я сепарабельних стаціонарних центрованих дійснозначних гауссівських процесів $(X_\Delta(t), t \in \mathbb{R})$, спектральні щільності яких збігаються на обмежених інтервалах до сталої $\frac{c}{2\pi}, c > 0,$ при $\Delta \to \infty$; при цьому відгук системи має вигляд:

$$Y_\Delta(t) = \int\limits_{-\infty}^{\infty} H(t-s) X_\Delta(s) ds,$$

або

$$Y_\Delta(t) = \int\limits_{-\infty}^{\infty} H(t-s) X_\Delta(s) ds + U(t),$$

якщо система містить внутрішній шум. При цьому припускається, що внутрішній шум описується сепарабельним стаціонарним центрованим дійснозначним гауссівським процесом $(U(t), t \in \mathbb{R})$, який є ортогональним до вхідного сигналу $X_\Delta$.

У роботі розглядаються нестійкі неперервні однорідні лінійні системи - - з внутрішнім шумом та без нього – за умови $H \in L_2(\mathbb{R})$.

**Предметом дослідження** є умови асимптотичної нормальності інтегральної корелограмної оцінки

$$\hat{H}_{T,\Delta}(\tau) = \frac{1}{cT} \int_0^T Y_\Delta(t+\tau) X_\Delta(t) dt, \tau \in \mathbb{R},$$

та похибки оцінювання

$$\hat{W}_{T,\Delta}(\tau) = \sqrt{T}[\hat{H}_{T,\Delta}(\tau) - H(\tau)], \tau \in \mathbb{R},$$

для імпульсної перехідної функції неперервної однорідної лінійної системи у випадках відсутності внутрішнього шуму та його наявності при $T, \Delta \to \infty$.

**Методика дослідження.** У роботі використовується кумулянтний метод аналізу розподілу; інтегральні представлення кумулянтів від статистик другого порядку; властивості багатовимірних інтегралів з циклічним зачепленням ядер, залежних від параметра; а також класичні результати про неперервність гауссівських процесів у просторі неперервних функцій.

**Наукова новизна одержаних результатів.** Усі отримані у дисертаційній роботі результати є новими. Основні з них такі:

- Отримано умови асимптотичної нормальності оцінки імпульсної перехідної функції у сенсі збіжності скінченновимірних розподілів;

- Отримано умови асимптотичної нормальності оцінки імпульсної перехідної функції у сенсі слабкої збіжності розподілів у просторі неперервних функцій;

- Отримано умови асимптотичної незсуненості оцінки імпульсної перехідної функції;

- Отримано умови асимптотичної нормальності похибки оцінювання оцінки імпульсної перехідної функції у сенсі збіжності скінченновимірних розподілів;

- Отримано умови асимптотичної нормальності похибки оцінювання імпульсної перехідної функції у сенсі слабкої збіжності розподілів у просторі неперервних функцій;

- Отримано приклади вхідних сигналів, внутрішнього шуму та імпульсних перехідних функцій.

**Практичне значення одержаних результатів.** Дисертаційна робота має теоретичний характер. Отримані результати можуть застосовуватись при оцінюванні або ідентифікації імпульсних перехідних функцій неперервних однорідних лінійних систем, зокрема, типу Вольтерра. Припущення інтегрованості у квадраті імпульсної перехідної функції дозволяє розглядати нестійкі неперервні однорідні лінійні системи з резонансними особливостями (кількість яких може бути зліченна).

Крім того, зважаючи на те, що неперервні однорідні лінійні системи часто виникають у задачах радіофізики, гідролокації, сейсмографії, метеорології, теорії сигналів та автоматичного контролю, теорії фільтрації, економетрики, тощо, слід очікувати на прикладні застосування результатів дисертації.

**Особистий внесок здобувача.** Всі результати дисертаційної роботи отримані здобувачем самостійно. За результатами дисертації здобувач опублікував п'ять робіт, з них три у співавторстві з науковим керівником проф. Булдигіним В.В., в яких Булдигіну В.В. належить постановка задач та загальне керівництво роботою. Дві роботи є авторськими.

**Апробація результатів.** езультати дисертації доповідались та обговорювалися на:

- міжуніверситетській науковій конференції з математики та фізики для студентів та молодих вчених (м. Київ, 2009 р.);

- Міжнародній конференції "Stochastic analysis and random dynamics" (м. Львів, 2009 р.);

- конференції "Фрактали і сучасна математика" (м. Київ, 2009р.);

- XIII Міжнародній науковій конференції імені академіка М. Кравчука (м. Київ, 2010 р.);

- Міжнародній конференції "Modern Stochastics: Theory and Applications II" (м. Київ, 2010 р.);

- XV міжнародній науковій конференції імені академіка М. Кравчука (м. Київ, 2014 р.);

- засіданнях наукового семінару "Теорія випадкових процесів" при кафедрі математичного аналізу і теорії ймовірностей НТУУ "КПІ" (керівник - доктор фіз.-мат. наук, проф. В.В. Булдигін), м. Київ, 2010, 2011, рр;

- засіданні наукового семінару при кафедрі електротехніки, інформатики та математики Університету Падерборну (керівник - доктор фіз.-мат. наук, проф. К.-Х. Індлекофер), м. Падерборн, 2011 р.;

- засіданнях наукового семінару "Статистичні проблеми для випадкових процесів та полів" при кафедрі математичного аналізу і теорії ймовірностей НТУУ "КПІ" (керівники - доктор фіз.-мат. наук, проф. О.І. Клесов, доктор фіз.-мат. наук, проф. О.В. Іванов), м. Київ, 2012, 2013 рр.;

• засіданні наукового семінару "Стохастика та її застосування" при кафедрі дослідження операцій на факультеті кібернетики Київського національного університету імені Тараса Шевченка (керівник - доктор фіз.-мат. наук, проф. О.М. Іксанов), м. Київ, 2014 р.;

• засіданні наукового семінару з теорії ймовірностей та математичної статистики при кафедрі при кафедрі теорії ймовірностей, статистики та актуарної математики механіко-математичного факультету Київського національного університету імені Тараса Шевченка (керівники - доктор фіз.-мат. наук, проф. Ю.В. Козаченко, доктор фіз.-мат. наук, проф. Ю.С. Мішура), м. Київ, 2014 р.;

• засіданні наукового семінару відділу математичних методів дослідження операцій Інституту кібернетики імені В.М. Глушкова НАН України (керівник - доктор фіз.-мат. наук, проф. П.С. Кнопов), м. Київ, 2015 р.

**Публікації.** За результатами дисертаційної роботи опубліковано 5 статтей у фахових виданнях ([123], [124], [125], [126], [127]) і 6 тез доповідей на конференціях ([128], [129], [130], [131], [132],[133]).

**Зміст дисертації.** Дисертація складається зі вступу, трьох розділів, розбитих на підрозділи, висновків, додатків та списку використаної літератури. Повний обсяг дисертації становить 191 сторінку, додатки займають 12 сторінок, список використаних джерел займає 18 сторінок і містить 133 найменування.

У **вступі** обґрунтовано актуальність теми дисертаційної роботи, визначено мету й задачі дослідження, виділено наукову новизну та практичну значущість отриманих результатів, охарактеризовано зміст кожного розділу дисертації (матеріал викладається зі збереженням оригінальної нумерації означень і теорем).

У **першому розділі** зроблено короткий історичний огляд літератури за тематикою дисертації, висвітлено сучасний стан вивчення проблем, близьких до тих, що розглядаються в дисертаційній роботі.

**Другий та третій розділи** містять основне дослідження, яке присвячене знаходженню умов асимптотичної нормальності інтегральної корелограмної оцінки, а також відповідної похибки оцінювання, для невідомої дійснозначної імпульсної перехідної функції $H = (H(\tau), \tau \in \mathbb{R})$ неперервної однорідної лінійної системи без внутрішнього шуму (другий розділ) та системи з внутрішнім шумом (третій розділ). Сформулюємо два базових припущення, притаманних всьому дослідженню.

*Перше припущення* стосується порядку інтегрованості імпульсної перехідної функції, а саме її інтегрованість у квадраті: $H \in L_2(\mathbb{R})$. Відмітимо, що зазвичай у літературі досліджуються стійкі лінійні системи, тобто, припускається, що функція $H$ абсолютно інтегрована ($H \in L_1(\mathbb{R})$). Таким чином, при вказаномій умові частотна характеристика системи $H^*$ (перетворення Фур'є - Планшереля функції $H$ у сенсі $L_2(\mathbb{R})$) може бути локально необмеженою. Цим самим в область дослідження потрапляють системи, що мають резонансні особливості. Деякі обмеження на порядок інтегрованості функції $H^*$ природньо з'являються як альтернатива збіжності ентропійних інтегралів.

*Друге припущення* стосується характеристик збурюючих процесів. Так,

в якості вхідних сигналів до неперервної однорідної лінійної системи розглядається сім'я вимірних сепарабельних стаціонарних центрованих дійснозначних гауссівських процесів $X_\Delta = (X_\Delta(t), t \in \mathbb{R})$, $\Delta > 0$, (див. Додаток А.4). Припустимо, що спектральні щільності $f_\Delta = (f_\Delta(\lambda), \lambda \in \mathbb{R})$, $\Delta > 0$, процесів $X_\Delta$ задовольняють умови:

$$f_\Delta(\lambda) = f_\Delta(-\lambda), \lambda \in \mathbb{R}; \tag{1а}$$

$$\sup_{\Delta > 0} \| f_\Delta \|_\infty < \infty; \tag{1б}$$

$$f_\Delta \in L_1(\mathbb{R}); \tag{1в}$$

$$\exists c \in (0,\infty) \ \forall a \in (0,\infty): \lim_{\Delta \to \infty} \sup_{-a \le \lambda \le a} \left| f_\Delta(\lambda) - \frac{c}{2\pi} \right| = 0; \tag{1г}$$

$$K_\Delta \in L_1(\mathbb{R}), \tag{1д}$$

де $K_\Delta(t) = \mathrm{E} X_\Delta(s+t) X_\Delta(s) = \int_{-\infty}^{\infty} e^{i\lambda t} f_\Delta(\lambda) d\lambda, \lambda \in \mathbb{R}$, - кореляційна функція процесу $X_\Delta$.

Умова (1г) показує, в якому саме сенсі слід розуміти "близькість" сім'ї процесів $X_\Delta, \Delta > 0$, до гауссівського білого шуму. (Відповідні приклади наведені у підрозділі 2.1).

Оцінку для імпульсної перехідної функції $H$ будемо шукати у виді інтегральної сумісної корелограми

$$H_{T,\Delta}(\tau) = \frac{1}{cT} \int_0^T Y_\Delta(t+\tau) X_\Delta(t) dt, \tau \in \mathbb{R}, \tag{2}$$

де $Y_\Delta(t), t \in \mathbb{R}$, - процес на виході системи; $c$ - стала з умови (1г), $\Delta$ - параметр схеми серій, й $T$ - довжина інтервалу усереднення $[0,T]$. Взагалі, для побудови оцінки процеси $X_\Delta$ та $Y_\Delta$ спостерігаються на всій дійсній осі.

Інтеграл, що фігурує у (2), розуміється як середньоквадратичний інтеграл Рімана.

Оцінка $H_{T,\Delta}$, по-перше, залежить від двох параметрів $T$ й $\Delta$, і, по-друге, є зсуненою, тобто $\mathrm{E} H_{T,\Delta}(\tau) \neq H(\tau)$. Ці факти ускладнюють асимптотичний аналіз оцінки і не дозволяють безпосередньо використовувати відомі результати для емпіричних сумісних корелограм.

Центральним місцем у дисертаційній роботі є дослідження асимптотичних властивостей нормованої похибки корелограмного оцінювання

$$W_{T,\Delta}(\tau) = \sqrt{T}[\hat{H}_{T,\Delta}(\tau) - H(\tau)], \tau \in \mathbb{R}, \qquad (3)$$

при прямуванні параметрів $T, \Delta$ до безмежності. Використовується наступний прийом - процес $W_{T,\Delta} = (W_{T,\Delta}(\tau), \tau \in \mathbb{R})$ розщеплюється на суму

$$W_{T,\Delta} = Z_{T,\Delta} + V_{T,\Delta},$$

де доданки справа визначаються далі, й для них спрацьовують такі факти:

- Вибір послідовності збурюючих процесів $(X_\Delta, \Delta > 0)$, апроксимуючої гауссівський білий шум, - є специфічним: *умови (1а) - (1д) є достатніми для нормалізації процесу $Z_{T,\Delta}$ при довільній збіжності параметрів $T, \Delta$ до безмежності.*

- Знищення вкладу невипадкової функції $V_{T,\Delta}$ у виразі для $W_{T,\Delta}$ можливе за рахунок використання $\delta$-видності сім'ї $K_\Delta$, підбору характеру сумісного прямування параметрів $T, \Delta$ до безмежності, та накладення умов на порядок локальної гладкості $H$.

Перейдемо до стислого огляду результатів другого та третього розділів дисертації. Повністю формулюються лише основні результати.

**Другий розділ** присвячено дослідженню умов асимптотичної нормальності інтегральної корелограмної оцінки (2), а також відповідної похибки оцінювання (3), для невідомої дійснозначної імпульсної перехідної функції $H = (H(\tau), \tau \in \mathbb{R})$ неперервної однорідної лінійної системи. Отримані результати далі застосовуються до задачі корелограмного оцінювання імпульсних перехідних функцій неперервних однорідних лінійних систем з внутрішнім шумом.

У *підрозділі 2.1* вводиться статистична оцінка $H_{T,\Delta}(\tau), \tau \in \mathbb{R}$, невідомої дійснозначної імпульсної перехідної функції $H(\tau), \tau \in \mathbb{R}$, неперервної однорідної лінійної системи (див. (2)). Вводяться загальні обмеження (1а)-(1д) на спектральні щільності $f_\Delta = (f_\Delta(\lambda), \lambda \in \mathbb{R})$ вхідних процесів $X_\Delta, \Delta > 0$. У лемі 2.1 встановлюються умови, при яких відгук системи, який описується процесом, коректно визначений:

$$Y_\Delta(t) = \int_{-\infty}^{\infty} H(s) X_\Delta(t-s) ds, \ t \in \mathbb{R},.$$

У *підрозділі 2.2* доводиться асимптотична нормальність емпіричного процесу

$$Z_{T,\Delta}(\tau) = \sqrt{T}[H_{T,\Delta}(\tau) - \mathrm{E}H(\tau)], \tau \in \mathbb{R},$$

як в сенсі збіжності скінченновимірних розподілів, так і в сенсі слабкої збіжності відповідних розподілів у просторі неперервних функцій при довільному прямуванні параметрів $T \to \infty, \Delta \to \infty$.

Даний підрозділ складається з чотирьох пунктів.

У *пункті 2.2.1* у лемі 2.2 встановлюється вигляд кореляційної функції процесу $Z_{T,\Delta}$ при всіх $T > 0, \Delta > 0$ та $\tau_1, \tau_2 \in \mathbb{R}$:

$$\mathrm{E}Z_{T,\Delta}(\tau_1)Z_{T,\Delta}(\tau_2) = C_{T,\Delta}(\tau_1, \tau_2) =$$

$$= \frac{2\pi}{c^2} \int_{-\infty}^{\infty}\int_{-\infty}^{\infty} [e^{i(\tau_1-\tau_2)\lambda_2} |H^*(\lambda_2)|^2 + e^{i(\tau_1\lambda_1+\tau_2\lambda_2)} H^*(\lambda_1)H^*(\lambda_2)] \times$$

$$\times \Phi_T(\lambda_2 - \lambda_1) f_\Delta(\lambda_1) f_\Delta(\lambda_2) d\lambda_1 d\lambda_2;$$

де $\Phi_T$ - ядро Фейєра:

$$\Phi_T(\lambda) = \frac{1}{2\pi T}\left(\frac{\sin(T\lambda/2)}{\lambda/2}\right)^2, \lambda \in \mathbb{R};$$

$c$ - стала з умови (1г), та $H^*$ - перетворення Фур'є-Планшереля функції $H$ у просторі $L_2(\mathbb{R})$.

Зокрема, в пункті встановлюється рівномірна оцінка по $T > 0, \Delta > 0$ та $\tau_1, \tau_2 \in \mathbb{R}$, для цієї кореляційної функції.

Далі запис $(T, \Delta) \to \infty$, означає, що одночасно $T \to \infty$ та $\Delta \to \infty$.

У *пункті 2.2.2* знайдено границю кореляційної функції $C_{T,\Delta}$ емпіричного процесу $Z_{T,\Delta}$ при $(T, \Delta) \to \infty$.

**Теорема 2.1.** *Нехай $H \in L_2(\mathbb{R})$, тоді для всіх $\tau_1, \tau_2 \in \mathbb{R}$, має місце рівність*

$$\lim_{(T,\Delta)\to\infty} \mathrm{E} Z_{T,\Delta}(\tau_1) Z_{T,\Delta}(\tau_2) = C_\infty(\tau_1, \tau_2),$$

*де*

$$C_\infty(\tau_1, \tau_2) = \frac{1}{2\pi} \int_{-\infty}^{\infty} \left[ e^{i(\tau_1-\tau_2)\lambda} |H^*(\lambda)|^2 + e^{i(\tau_1+\tau_2)\lambda} (H^*(\lambda))^2 \right] d\lambda.$$

У *пункті 2.2.3* доводиться асимптотична нормальність скінченновимірних розподілів процесу $Z_{T,\Delta} = (Z_{T,\Delta}(\tau), \tau \in \mathbb{R})$ при $(T, \Delta) \to \infty$, за допомогою методу кумулянтів та моментів.

Розглянемо центрований дійснозначний гауссівський процес $Z = (Z(\tau), \tau \in \mathbb{R})$ з кореляційною функцією $C_\infty$, визначеною у теоремі 2.1.

**Теорема 2.2.** *Нехай $H \in L_2(\mathbb{R})$, тоді для всіх $m \in \mathbb{N}$ та $\tau_1, \tau_2, ..., \tau_m \in \mathbb{R}$, має місце рівність*

$$\lim_{(T,\Delta)\to\infty} cum(Z_{T,\Delta}(\tau_j), j=1,...,m) = \begin{cases} 0, & m=1; \\ C_\infty(\tau_1,\tau_2), & m=2; \\ 0, & m\geq 3, \end{cases}$$

де $cum(Z_{T,\Delta}(\tau_j), j=1,...,m)$ - *сумісний кумулянт набору випадкових величин* $Z_{T,\Delta}(\tau_j), j=1,...,m.$

*Зокрема, всі скінченновимірні розподіли процесу $(Z_{T,\Delta}(\tau), \tau\in\mathbb{R})$ слабко збігаються до відповідних скінченновимірних розподілів центрованого гауссівського процесу $(Z(\tau), \tau\in\mathbb{R})$.*

Твердження теореми 2.2 можна сформулювати за допомогою моментів.

**Теорема 2.3.** *Нехай $H \in L_2(\mathbb{R})$, тоді для всіх $m \in \mathbb{N}$ та $\tau_1, \tau_2, ..., \tau_m \in \mathbb{R}$, має місце рівність*

$$\lim_{(T,\Delta)\to\infty} E[\prod_{j=1}^{m} Z_{T,\Delta}(\tau_j)] = E[\prod_{j=1}^{m} Z(\tau_j)].$$

*Зокрема, всі скінченновимірні розподіли процесу $(Z_{T,\Delta}(\tau), \tau\in\mathbb{R})$ слабко збігаються до відповідних скінченновимірних розподілів центрованого гауссівського процесу $(Z(\tau), \tau\in\mathbb{R})$.*

У *пункті 2.2.4* вивчаються умови слабкої збіжності процесу $Z_{T,\Delta}$ до процесу $Z$ при $(T,\Delta) \to \infty$, у просторі неперервних функцій. Вказані процеси вважаємо сепарабельними в силу їх стохастичної неперервності.

Нехай $C[a,b]$ - банахів простір неперервних дійснозначних функцій, визначених на сегменті $[a,b] \subset \mathbb{R}$, з рівномірною нормою. Далі запис

$$Z_{T,\Delta} \stackrel{C[a,b]}{\Rightarrow} Z$$

означає слабку збіжність процесу $Z_{T,\Delta}$ до процесу $Z$ у просторі неперервних функцій при $(T,\Delta)\to\infty$.

**Теорема 2.4.** *Нехай $H \in L_2(\mathbb{R})$ та виконується нерівність*

$$\int_{0+} \mathcal{H}_{\sqrt{\sigma}}(\varepsilon)d\varepsilon < \infty, \qquad (4)$$

*(де $\mathcal{H}_{\sqrt{\sigma}}(\varepsilon)$ - метрична ентропія відрізка $[0,1]$ відносно псевдометрики*

$$\sqrt{\sigma}(\tau_1,\tau_2) = \left[\int_{-\infty}^{\infty} \sin^2\frac{|\tau_1-\tau_2|\lambda}{2}|H^*(\lambda)|^2 d\lambda\right]^{\frac{1}{4}}, \tau_1,\tau_2 \in \mathbb{R},)$$

*тоді для будь-якого $[a,b] \subset \mathbb{R}$ мають місце наступні твердження*:

*(I)  $Z \in C[a,b]$ майже напевно;*

*(II)  $Z_{T,\Delta} \in C[a,b]$ майже напевно, $T>0, \Delta>0$;*

*(III)  $Z_{T,\Delta} \stackrel{C[a,b]}{\Rightarrow} Z$.*

*Зокрема, для будь-якого $x>0$*

$$\lim_{(T,\Delta)\to\infty} \mathrm{P}\left\{\sup_{\tau\in[a,b]}|Z_{T,\Delta}(\tau)|>x\right\} = \mathrm{P}\left\{\sup_{\tau\in[a,b]}|Z(\tau)|>x\right\}.$$

Відзначимо (зауваження 2.7), що твердження (I) теореми 2.4 виконується, якщо замість ентропійної умови (4) вимагати, щоб при деякому $\beta>0$

$$\int_{-\infty}^{\infty} |H^*(\lambda)|^2 \ln^{1+\beta}(1+|\lambda|)d\lambda < \infty.$$

У *підрозділі 2.3* вивчаються умови асимптотичної незсуненості та консистентності оцінки $H_{T,\Delta}$. Покладемо

$$\hat{v}_\Delta(\tau) = [\mathrm{E}H_{T,\Delta}(\tau) - H(\tau)], \tau \in \mathbb{R}.$$

Додатково до умов (1а)-(1д) відносно функцій $f_\Delta, K_\Delta$ виконуються умови:

(А)  Нехай при деякому $\alpha \in (0,1]$  $H \in Lip_\alpha(\mathbb{R})$;

(Б)  Припустимо, що при заданому $\alpha \in (0,1]$:

$$\forall \delta > 0 \quad \lim_{\Delta \to \infty} \int_\delta^\infty K_\Delta(t) dt = 0; \qquad (5а)$$

$$\forall \delta > 0 \quad \lim_{\Delta \to \infty} \int_\delta^\infty K_\Delta^2(t) dt = 0; \qquad (5б)$$

$$\exists \delta > 0 \quad \lim_{\Delta \to \infty} \int_{-\delta}^\delta |K_\Delta(t)| |t|^\alpha dt = 0. \qquad (5в)$$

Після введення умов (А) та (Б) мають місце наступні твердження.

**Теорема 2.5.** *Нехай задане $\alpha \in (0,1]$; $H \in Lip_\alpha(\mathbb{R}) \cap L_2(\mathbb{R})$ та виконуються умови (5а) - (5в), тоді:*

*(I)  для будь-якого $\tau \in \mathbb{R}$ $\quad \lim_{\Delta \to \infty} \hat{v}_\Delta(\tau) = 0$;*

*(II)  для будь-якого $[a,b] \subset \mathbb{R}$ $\lim_{\Delta \to \infty} \sup_{\tau \in [a,b]} |\hat{v}_\Delta(\tau)| = 0$.*

**Теорема 2.6.** *Нехай задане $\alpha \in (0,1]$; $H \in Lip_\alpha(\mathbb{R}) \cap L_2(\mathbb{R})$ та виконуються умови (5а) - (5в), тоді для будь-якого $\tau \in \mathbb{R}$*

$$\lim_{(T,\Delta) \to \infty} \mathrm{E} |H_{T,\Delta}(\tau) - H(\tau)|^2 = 0.$$

*Тобто, $H_{T,\Delta}$ є конзистентною у середньому квадратичному в точці $\tau$.*

У *підрозділі 2.4* досліджується асимптотична поведінка нормованої похибки корелограмного оцінювання виду (3):

$$W_{T,\Delta}(\tau) = \sqrt{T}[H_{T,\Delta}(\tau) - H(\tau)], \tau \in \mathbb{R},$$

при прямуванні параметрів $T, \Delta$ до нескінченності. Підрозділ складається з трьох пунктів.

У *пункті 2.4.1* вивчаються умови ануляції невипадкової функції

$$V_{T,\Delta}(\tau) = \sqrt{T}[\mathrm{E}H_{T,\Delta}(\tau) - H(\tau)], \tau \in \mathbb{R}.$$

Додатково до умов (1а)-(1д) відносно функцій $f_\Delta$, $K_\Delta$ виконуються наступні умови:

(А) Нехай при деякому $\alpha \in (0,1]$ $H \in Lip_\alpha(\mathbb{R})$;

(Б) При заданому $\alpha \in (0,1]$ параметри $T \to \infty, \Delta \to \infty$ так, що:

$$\sqrt{T}[1 - \frac{2\pi f_\Delta(0)}{c}] \to 0; \qquad (6а)$$

$$\forall \delta > 0 \quad \sqrt{T}\int_\delta^\infty K_\Delta(t)dt \to 0; \qquad (6б)$$

$$\forall \delta > 0 \quad T\int_\delta^\infty K_\Delta^2(t)dt \to 0; \qquad (6в)$$

$$\exists \delta > 0 \quad \sqrt{T}\int_{-\delta}^\delta |K_\Delta(t)||t|^\alpha dt \to 0. \qquad (6г)$$

У лемі 2.5 показано, що при деякому $\alpha \in (0,1]$, якщо $H \in Lip_\alpha(\mathbb{R}) \cap L_2(\mathbb{R})$ та $T \to \infty, \Delta \to \infty$ так, що виконуються умови (6а) - (6г), тоді:

(I) для будь-якого $\tau \in \mathbb{R}$ $\quad V_{T,\Delta}(\tau) \to 0$;

(II) для будь-якого $[a,b] \subset \mathbb{R}$ $\quad \sup_{\tau \in [a,b]} |V_{T,\Delta}(\tau)| \to 0$.

У *пунктах 2.4.2 - 2.4.3* доводиться асимптотична нормальність процесу $W_{T,\Delta} = (W_{T,\Delta}(\tau), \tau \in \mathbb{R})$. Усі твердження випливають з представлення

$$W_{T,\Delta} = Z_{T,\Delta} + V_{T,\Delta},$$

та відповідно теорем 2.1 - 2.4.

**Теорема 2.7.** *Нехай для деякого $\alpha \in (0,1]$ $H \in Lip_\alpha(\mathbb{R}) \cap L_2(\mathbb{R})$. Якщо $T \to \infty, \Delta \to \infty$ так, що виконуються умови (6а) - (6г), тоді для всіх $\tau_1, \tau_2 \in \mathbb{R}$, має місце співвідношення*

$$\mathrm{E} W_{T,\Delta}(\tau_1) W_{T,\Delta}(\tau_2) \to C_{\infty}(\tau_1, \tau_2) =$$

$$= \frac{1}{2\pi} \int_{-\infty}^{\infty} [e^{i(\tau_1 - \tau_2)\lambda} |H^*(\lambda)|^2 + e^{i(\tau_1 + \tau_2)\lambda} (H^*(\lambda))^2] d\lambda.$$

З теореми 2.7 випливає, що гранична кореляційна функція процесу $W_{T,\Delta}$ співпадає з кореляційною функцією гауссівського центрованого процесу $Z = (Z(\tau), \tau \in \mathbb{R})$. Зокрема, має місце твердження про асимптотичну поведінку скінченновимірних розподілів $W_{T,\Delta}$ при прямуванні $T, \Delta$ до безмежності.

**Теорема 2.8.** *Нехай для деякого $\alpha \in (0,1]$ $H \in Lip_{\alpha}(\mathbb{R}) \cap L_2(\mathbb{R})$. Якщо $T \to \infty, \Delta \to \infty$ так, що виконуються умови (6а) - (6г), тоді для всіх $m \in \mathbb{N}$ та $\tau_1, \tau_2, ..., \tau_m \in \mathbb{R}$, має місце співвідношення*

$$\mathrm{cum}(W_{T,\Delta}(\tau_j), j = 1,...,m) \to \begin{cases} 0, & m = 1; \\ C_{\infty}(\tau_1, \tau_2), & m = 2; \\ 0, & m \geq 3, \end{cases}$$

де $\mathrm{cum}(W_{T,\Delta}(\tau_j), j = 1,...,m)$ - сумісний кумулянт набору випадкових величин $W_{T,\Delta}(\tau_j), j = 1,...,m$.

*Зокрема, всі скінченновимірні розподіли процесу $(W_{T,\Delta}(\tau), \tau \in \mathbb{R})$ слабко збігаються до відповідних скінченновимірних розподілів центрованого гауссівського процесу $(Z(\tau), \tau \in \mathbb{R})$ при вказаному характері прямування $T$ і $\Delta$ до нескінченності.*

**Теорема 2.9.** *Нехай для деякого $\alpha \in (0,1]$ $H \in Lip_{\alpha}(\mathbb{R}) \cap L_2(\mathbb{R})$. Якщо $T \to \infty, \Delta \to \infty$ так, що виконуються умови (6а) - (6г), тоді для всіх $m \in \mathbb{N}$ та $\tau_1, \tau_2, ..., \tau_m \in \mathbb{R}$, має місце співвідношення*

$$\mathrm{E}[\prod_{j=1}^{m} W_{T,\Delta}(\tau_j)] \to \mathrm{E}[\prod_{j=1}^{m} Z(\tau_j)].$$

*Зокрема, всі скінченновимірні розподіли процесу $(W_{T,\Delta}(\tau), \tau \in \mathbb{R})$ слабко збігаються до відповідних скінченновимірних розподілів центрованого гауссівського процесу $(Z(\tau), \tau \in \mathbb{R})$ при вказаному характері прямування $T$ і $\Delta$ до нескінченності.*

Будемо вважати, що гауссівський процес $W_{T,\Delta} = (W_{T,\Delta}(\tau), \tau \in \mathbb{R})$ є сепарабельними. (Таке припущення не є обмежувальним в силу стохастичної неперервності цього процесу). Має місце наступне твердження про асимптотичну нормальність процесу $W_{T,\Delta}$ у просторі неперервних функцій.

**Теорема 2.10.** *Нехай для деякого $\alpha \in (0,1]$ $H \in Lip_\alpha(\mathbb{R}) \bigcap L_2(\mathbb{R})$, та виконується нерівність (4), тобто*

$$\int_{0+} \mathcal{H}_{\sqrt{\sigma}}(\varepsilon) d\varepsilon < \infty,$$

*тоді для будь-якого $[a,b] \subset \mathbb{R}$ мають місце наступні твердження:*

*(I) $Z \in C[a,b]$ майже напевно;*

*(II) $W_{T,\Delta} \in C[a,b]$ майже напевно, $T>0, \Delta>0$;*

*Крім того, якщо $T \to \infty, \Delta \to \infty$ так, що виконуються умови (6а) - (6г), тоді*

*(III) $W_{T,\Delta} \overset{C[a,b]}{\Rightarrow} Z$.*

*Зокрема, при вказаному характері прямування $T$ і $\Delta$ до нескінченності, для будь-якого $x>0$*

$$\mathrm{P}\left\{\sup_{\tau \in [a,b]} |W_{T,\Delta}(\tau)| > x\right\} \to \mathrm{P}\left\{\sup_{\tau \in [a,b]} |Z(\tau)| > x\right\}.$$

У *підрозділі 2.5* отримані результати ілюструються на конкретних прикладах імпульсних перехідних функцій (приклади 2.10 - 2.18) та вхідних сигналів (приклади 2.6 - 2.9), які задовольняють схемі оцінювання.

**Третій розділ** присвячено дослідженню умов асимптотичної нормальності інтегральної корелограмної оцінки (2), а також відповідної похибки оцінювання (3), для невідомої дійснозначної імпульсної перехідної функції $H = (H(\tau), \tau \in \mathbb{R})$ неперервної однорідної лінійної системи з внутрішнім шумом. Зазначимо, що результати даного розділу узагальнюють результати другого розділу.

У *підрозділі 3.1* вводиться статистична оцінка $H_{T,\Delta}(\tau), \tau \in \mathbb{R}$, виду (2) для невідомої дійснозначної імпульсної перехідної функції $H(\tau), \tau \in \mathbb{R}$, неперервної однорідної лінійної системи з внутрішнім шумом. Вводяться загальні обмеження (1а)-(1д) на спектральні щільності $f_\Delta = (f_\Delta(\lambda), \lambda \in \mathbb{R})$ вхідних процесів $X_\Delta, \Delta > 0$. Припускається, що внутрішній шум системи описується вимірним сепарабельним стаціонарним центрованим дійснозначним гауссівським процесом $U = (U(t), t \in \mathbb{R})$, який є ортогональним до $X_\Delta$; тобто $\mathrm{E} X_\Delta(s) U(t) = 0, s, t \in \mathbb{R}$. Припускається, що процес $U$ має спектральну щільність $g = (g(\lambda), \lambda \in \mathbb{R})$, яка є парною функцією. У лемі 3.1 встановлюються умови, при яких відгук системи, який описується процесом

$$Y_\Delta(t) = \int\limits_{-\infty}^{\infty} H(s) X_\Delta(t-s) ds + U(t), \ t \in \mathbb{R},$$

коректно визначений.

У *підрозділі 3.2* доводиться асимптотична нормальність процесу

$$Z_{T,\Delta}(\tau) = \sqrt{T}[H_{T,\Delta}(\tau) - \mathrm{E} H(\tau)], \tau \in \mathbb{R},$$

як в сенсі збіжності скінченновимірних розподілів, так і в сенсі слабкої збіжності відповідних розподілів у просторі неперервних функцій при довільному прямуванні параметрів $T \to \infty, \Delta \to \infty$. Даний підрозділ складається з чотирьох пунктів.

У *пункті 3.2.1* у лемі 3.2 встановлюється вигляд кореляційної функції процесу $Z_{T,\Delta}$ при всіх $T>0, \Delta>0$ та $\tau_1, \tau_2 \in \mathbb{R}$:

$$\mathrm{E} Z_{T,\Delta}(\tau_1) Z_{T,\Delta}(\tau_2) = C_{T,\Delta}(\tau_1, \tau_2) =$$

$$= \frac{2\pi}{c^2} \int_{-\infty}^{\infty}\int_{-\infty}^{\infty} [e^{i(\tau_1 - \tau_2)\lambda_2}(|H^*(\lambda_2)|^2 f_\Delta(\lambda_2) + g(\lambda_2)) +$$

$$+ e^{i(\tau_1\lambda_1 + \tau_2\lambda_2)} H^*(\lambda_1) H^*(\lambda_2) f_\Delta(\lambda_2)] \Phi_T(\lambda_2 - \lambda_1) f_\Delta(\lambda_1) d\lambda_1 d\lambda_2;$$

де $\Phi_T$ - ядро Фейєра:   $\Phi_T(\lambda) = \frac{1}{2\pi T}\left(\frac{\sin(T\lambda/2)}{\lambda/2}\right)^2$, $\lambda \in \mathbb{R}$;

$c$ - стала з умови (1г), та $H^*$ - перетворення Фур'є-Планшереля функції $H$ у просторі $L_2(\mathbb{R})$.

Зокрема, у пункті наводиться рівномірна по всіх $T>0, \Delta>0$ та $\tau_1, \tau_2 \in \mathbb{R}$ оцінка для кореляційної функції процесу $Z_{T,\Delta}$.

У *пункті 3.2.2* знайдено границю кореляційної функції $C_{T,\Delta}$ емпіричного процесу $Z_{T,\Delta}$ при $(T, \Delta) \to \infty$.

**Теорема 3.1.**  *Нехай $H \in L_2(\mathbb{R})$ та $g \in L_1(\mathbb{R})$, тоді для всіх $\tau_1, \tau_2 \in \mathbb{R}$, має місце рівність*

$$\lim_{(T,\Delta)\to\infty} \mathrm{E} Z_{T,\Delta}(\tau_1) Z_{T,\Delta}(\tau_2) = C_\infty(\tau_1, \tau_2),$$

*де*

$$C_\infty(\tau_1, \tau_2) = \frac{1}{2\pi} \int_{-\infty}^{\infty} \left[ e^{i(\tau_1 - \tau_2)\lambda}(|H^*(\lambda)|^2 + \frac{2\pi}{c} g(\lambda)) + e^{i(\tau_1 + \tau_2)\lambda} (H^*(\lambda))^2 \right] d\lambda.$$

У *пункті 3.2.3* доводиться асимптотична нормальність скінченновимірних розподілів процесу $Z_{T,\Delta} = (Z_{T,\Delta}(\tau), \tau \in \mathbb{R})$ при $(T, \Delta) \to \infty$, за допомогою методу кумулянтів та моментів.

Розглянемо центрований дійснозначний гауссівський процес $Z = (Z(\tau), \tau \in \mathbb{R})$ з кореляційною функцією $C_\infty$, визначеною у теоремі 3.1.

**Теорема 3.2.** *Нехай* $H \in L_2(\mathbb{R})$ *та* $g \in L_1(\mathbb{R})$, *тоді для всіх* $m \in \mathbb{N}$ *та* $\tau_1, \tau_2, ..., \tau_m \in \mathbb{R}$, *має місце рівність*

$$\lim_{(T,\Delta) \to \infty} cum(Z_{T,\Delta}(\tau_j), j = 1, ..., m) = \begin{cases} 0, & m = 1; \\ C_\infty(\tau_1, \tau_2), & m = 2; \\ 0, & m \geq 3, \end{cases}$$

*де* $cum(Z_{T,\Delta}(\tau_j), j = 1, ..., m)$ - *сумісний кумулянт набору випадкових величин* $Z_{T,\Delta}(\tau_j), j = 1, ..., m$.

*Зокрема, всі скінченновимірні розподіли процесу* $(Z_{T,\Delta}(\tau), \tau \in \mathbb{R})$ *слабко збігаються до відповідних скінченновимірних розподілів центрованого гауссівського процесу* $(Z(\tau), \tau \in \mathbb{R})$.

Твердження теореми 3.2 можна сформулювати за допомогою моментів.

**Теорема 3.3.** *Нехай* $H \in L_2(\mathbb{R})$ *та* $g \in L_1(\mathbb{R})$, *тоді для всіх* $m \in \mathbb{N}$ *та* $\tau_1, \tau_2, ..., \tau_m \in \mathbb{R}$, *має місце рівність*

$$\lim_{(T,\Delta) \to \infty} \mathrm{E}[\prod_{j=1}^{m} Z_{T,\Delta}(\tau_j)] = \mathrm{E}[\prod_{j=1}^{m} Z(\tau_j)].$$

*Зокрема, всі скінченновимірні розподіли процесу* $(Z_{T,\Delta}(\tau), \tau \in \mathbb{R})$ *слабко збігаються до відповідних скінченновимірних розподілів центрованого гауссівського процесу* $(Z(\tau), \tau \in \mathbb{R})$.

У *пункті 3.2.4* вивчаються умови слабкої збіжності процесу $Z_{T,\Delta}$ до процесу $Z$ при $(T,\Delta) \to \infty$, у просторі неперервних функцій. Вказані процеси вважаємо сепарабельними в силу їх стохастичної неперервності.

**Теорема 3.4.** *Нехай $H \in L_2(\mathbb{R})$, $g \in L_1(\mathbb{R})$ та виконується нерівність*

$$\int_{0+} \mathcal{H}_{\sqrt{\sigma}}(\varepsilon) d\varepsilon < \infty, \qquad (7)$$

*(де $\mathcal{H}_{\sqrt{\sigma}}(\varepsilon)$ - метрична ентропія відрізка [0,1] відносно псевдометрики*

$$\sqrt{\sigma}(\tau_1, \tau_2) = \left[ \int_{-\infty}^{\infty} \sin^2 \frac{|\tau_1 - \tau_2|\lambda}{2} (|H^*(\lambda)|^2 + g(\lambda)) d\lambda \right]^{\frac{1}{4}}, \tau_1, \tau_2 \in \mathbb{R},)$$

*тоді для будь-якого $[a,b] \subset \mathbb{R}$ мають місце наступні твердження:*

*(I)* $Z \in C[a,b]$ *майже напевно;*

*(II)* $Z_{T,\Delta} \in C[a,b]$ *майже напевно, $T > 0, \Delta > 0$;*

*(III)* $Z_{T,\Delta} \overset{C[a,b]}{\Rightarrow} Z$.

*Зокрема, для будь-якого $x > 0$*

$$\lim_{(T,\Delta) \to \infty} \mathrm{P}\left\{ \sup_{\tau \in [a,b]} |Z_{T,\Delta}(\tau)| > x \right\} = \mathrm{P}\left\{ \sup_{\tau \in [a,b]} |Z(\tau)| > x \right\}.$$

Відзначимо (зауваження 3.4), що твердження (I) теореми 3.4 виконується, якщо замість ентропійної умови (7) вимагати, щоб при деякому $\beta > 0$

$$\int_{-\infty}^{\infty} (|H^*(\lambda)|^2 + g(\lambda)) \ln^{1+\beta}(1+|\lambda|) d\lambda < \infty.$$

Зокрема (зауваження 3.5), додаткова ентропійна умова (7), що фігурує у теоремі 3.4, має місце, якщо при деякому $\beta > 0$

$$\int_{-\infty}^{\infty} (|H^*(\lambda)|^2 + g(\lambda)) \ln^{4+\beta}(1+|\lambda|) d\lambda < \infty.$$

У *підрозділі 3.3* вивчаються умови асимптотичної незсуненості та консистентності оцінки $H_{T,\Delta}$.

Наступне твердження про асимптотичну незсуненість дублює відповідне твердження у розділі 2.

**Теорема 3.5.** *Нехай задане $\alpha \in (0,1]$; $H \in Lip_\alpha(\mathbb{R}) \cap L_2(\mathbb{R})$ та виконуються умови (5а) - (5в), тоді:*

(I) *для будь-якого $\tau \in \mathbb{R}$*   $\lim\limits_{\Delta \to \infty} \hat{v}_\Delta(\tau) = 0$;

(II) *для будь-якого $[a,b] \subset \mathbb{R}$   $\lim\limits_{\Delta \to \infty} \sup\limits_{\tau \in [a,b]} |\hat{v}_\Delta(\tau)| = 0$.*

**Теорема 3.6.** *Нехай задане $\alpha \in (0,1]$; $H \in Lip_\alpha(\mathbb{R}) \cap L_2(\mathbb{R})$, $g \in L_1(\mathbb{R})$ та виконуються умови (5а) - (5в), тоді для будь-якого $\tau \in \mathbb{R}$*

$$\lim\limits_{(T,\Delta) \to \infty} \mathrm{E} |H_{T,\Delta}(\tau) - H(\tau)|^2 = 0.$$

*Тобто, $H_{T,\Delta}$ є конзистентною у середньому квадратичному в точці $\tau$.*

У *підрозділі 3.4* досліджується асимптотична поведінка нормованої похибки корелограмного оцінювання виду (3):

$$W_{T,\Delta}(\tau) = \sqrt{T}[H_{T,\Delta}(\tau) - H(\tau)], \tau \in \mathbb{R},$$

при прямуванні параметрів $T, \Delta$ до нескінченності. Підрозділ складається з трьох пунктів.

У *пункті 3.4.1* вивчаються умовиануляції невипадкової функції $V_{T,\Delta}$.

У лемі 3.5 "дублюються" твердження леми 2.5, тобто, якщо при деякому $\alpha \in (0,1]$, якщо $H \in Lip_\alpha(\mathbb{R}) \cap L_2(\mathbb{R})$ та $T \to \infty, \Delta \to \infty$ так, що виконуються умови (6а) - (6г), тоді:

(I) для будь-якого $\tau \in \mathbb{R}$   $V_{T,\Delta}(\tau) \to 0$;

(II) для будь-якого $[a,b] \subset \mathbb{R}$   $\sup\limits_{\tau \in [a,b]} |V_{T,\Delta}(\tau)| \to 0$.

У *пунктах 3.4.2 - 3.4.3* доводиться асимптотична нормальність процесу $W_{T,\Delta} = (W_{T,\Delta}(\tau), \tau \in \mathbb{R})$. Усі твердження випливають з представлення $W_{T,\Delta} = Z_{T,\Delta} + V_{T,\Delta}$, та відповідно теорем 3.1 - 3.4.

**Теорема 3.7.** *Нехай для деякого $\alpha \in (0,1]$ $H \in Lip_\alpha(\mathbb{R}) \cap L_2(\mathbb{R})$, та $g \in L_1(\mathbb{R})$. Якщо $T \to \infty, \Delta \to \infty$ так, що виконуються умови (6а) - (6г), тоді для всіх $\tau_1, \tau_2 \in \mathbb{R}$, має місце співвідношення*

$$\mathrm{E} W_{T,\Delta}(\tau_1) W_{T,\Delta}(\tau_2) \to C_\infty(\tau_1, \tau_2) =$$

$$= \frac{1}{2\pi} \int_{-\infty}^{\infty} \left[ e^{i(\tau_1 - \tau_2)\lambda}(|H^*(\lambda)|^2 + \frac{2\pi}{c} g(\lambda)) + e^{i(\tau_1 + \tau_2)\lambda}(H^*(\lambda))^2 \right] d\lambda, \ \tau_1, \tau_2 \in \mathbb{R}.$$

З теореми 3.7 випливає, що гранична кореляційна функція процесу $W_{T,\Delta}$ співпадає з кореляційною функцією гауссівського центрованого процесу $Z = (Z(\tau), \tau \in \mathbb{R})$. Зокрема, має місце таке твердження про асимптотичну поведінку скінченновимірних розподілів $W_{T,\Delta}$ при прямуванні параметрів $T, \Delta$ до безмежності.

**Теорема 3.8.** *Нехай для деякого $\alpha \in (0,1]$ $H \in Lip_\alpha(\mathbb{R}) \cap L_2(\mathbb{R})$, та $g \in L_1(\mathbb{R})$. Якщо $T \to \infty, \Delta \to \infty$ так, що виконуються умови (6а) - (6г), тоді для всіх $m \in \mathbb{N}$ та $\tau_1, \tau_2, ..., \tau_m \in \mathbb{R}$, має місце співвідношення*

$$\mathrm{cum}(W_{T,\Delta}(\tau_j), j=1,...,m) \to \begin{cases} 0, & m = 1; \\ C_\infty(\tau_1, \tau_2), & m = 2; \\ 0, & m \geq 3, \end{cases}$$

*де $\mathrm{cum}(W_{T,\Delta}(\tau_j), j=1,...,m)$ - сумісний кумулянт набору випадкових величин $W_{T,\Delta}(\tau_j), j=1,...,m.$*

*Зокрема, всі скінченновимірні розподіли процесу $(W_{T,\Delta}(\tau), \tau \in \mathbb{R})$ слабко збігаються до відповідних скінченновимірних розподілів центрованого гауссівського процесу $(Z(\tau), \tau \in \mathbb{R})$ при вказаному характері прямування $T$ і $\Delta$ до нескінченності.*

**Теорема 3.9.** *Нехай для деякого* $\alpha \in (0,1]$ $H \in Lip_\alpha(\mathbb{R}) \cap L_2(\mathbb{R})$, *та* $g \in L_1(\mathbb{R})$. *Якщо* $T \to \infty, \Delta \to \infty$ *так, що виконуються умови (6а) - (6г), тоді для всіх* $m \in \mathbb{N}$ *та* $\tau_1, \tau_2, ..., \tau_m \in \mathbb{R}$, *має місце співвідношення*

$$\mathrm{E}[\prod_{j=1}^{m} W_{T,\Delta}(\tau_j)] \to \mathrm{E}[\prod_{j=1}^{m} Z(\tau_j)].$$

*Зокрема, всі скінченновимірні розподіли процесу* $(W_{T,\Delta}(\tau), \tau \in \mathbb{R})$ *слабко збігаються до відповідних скінченновимірних розподілів центрованого гауссівського процесу* $(Z(\tau), \tau \in \mathbb{R})$ *при вказаному характері прямування* $T$ *і* $\Delta$ *до нескінченності.*

Будемо вважати, що гауссівський процес $W_{T,\Delta} = (W_{T,\Delta}(\tau), \tau \in \mathbb{R})$ є сепарабельними. (Таке припущення не є обмежувальним в силу стохастичної неперервності цього процесу).

Має місце наступне твердження про асимптотичну нормальність процесу $W_{T,\Delta}$ у просторі неперервних функцій.

**Теорема 3.10.** *Нехай для деякого* $\alpha \in (0,1]$ $H \in Lip_\alpha(\mathbb{R}) \bigcap L_2(\mathbb{R})$; $g \in L_1(\mathbb{R})$, *та виконується нерівність (7), тобто*

$$\int_{0+} \mathcal{H}_{\sqrt{\sigma}}(\varepsilon) d\varepsilon < \infty,$$

*тоді для будь-якого* $[a,b] \subset \mathbb{R}$ *мають місце наступні твердження:*

(I)     $Z \in C[a,b]$ *майже напевно;*

(II)    $W_{T,\Delta} \in C[a,b]$ *майже напевно,* $T > 0, \Delta > 0$;

*Крім того, якщо* $T \to \infty, \Delta \to \infty$ *так, що виконуються умови (6а) - (6г), тоді*

(III)    $W_{T,\Delta} \overset{C[a,b]}{\Rightarrow} Z$.

*Зокрема, при вказаному характері прямування T і Δ до нескінченності, для будь-якого $x > 0$*

$$P\left\{\sup_{\tau \in [a,b]} |W_{T,\Delta}(\tau)| > x\right\} \to P\left\{\sup_{\tau \in [a,b]} |Z(\tau)| > x\right\}.$$

У *пункті 3.4.3* отримані результати ілюструються на конкретних прикладах імпульсних перехідних функцій $H \in L_2(\mathbb{R})$, характеристик збурюючих процесів $X_\Delta$ й процесу внутрішнього шуму $U$ (приклад 3.5).

У **висновках** сформульовано основні результати дисертаційної роботи.

# РОЗДІЛ 1

# ОГЛЯД ЛІТЕРАТУРИ ЗА ТЕМОЮ ДИСЕРТАЦІЇ

Далі наведемо огляд робіт за тематикою дисертації за схемою: постановка задачі - обраний метод - існуючі розв'язки задачі - інструментарій.

**Постановка задачі.** Ідентифікація та оцінювання деяких характеристик стохастичних лінійних та нелінійних систем є предметом активних досліджень протягом останніх п'ятидесяти років, і вивчались в роботах Х.Акайка [42, 43, 44], І. Густавсона, Л. Льюнга та Т. Шодерстрома [80], Г. Гудвіна, П. Рамаджа та П. Кайнза [77], Л. Льюнга та Т. Шодерстрома [93], Р. Ісерманна [86], П. Кайнза та С. Лафортуна [64], П. Варайї та П. Кумара [88], П. Кайнза [63], М. Нієдзвеськи [97], Е. Ескіната, С. Джонсона та В. Лейбена [72], П. Янга, А. Чотаї та В. Тіча [118], В. Мармареліза [94], П. Ван дер Хофа та Р. Шрама [115], Т. Андерсена та П. Пукара [45], Л. Льюнга [91], М. Норгарда, Н. Польсена та О.Равна [99], К.-Г. Ванга та Я. Жанга [116], М. Флізза та Х. Сіра-Раміреса [75], Д. Ванга та Ф. Дінга [117], Дж. Красідіза та Дж. Дженкінса [66], Р. Пінтелона та Дж. Шоукенза [104], Е. Ханана та М. Дейстлера [83] і багатьох інших. Одна з найпростіших моделей, що описується такими системами, розглядає так звану "чорну скриньку," що приймає певний вхідний сигнал та видає відомий відгук на нього. При цьому вхідний та вихідний сигнали можуть бути простими або багатозначними, чистими або містити шуми. Різні варіанти постановки задачі оцінювання або ідентифікації характеристик "чорної скриньки" породжують ряд моделей, що можуть застосовуватись, зокрема, параметричні (див. роботи М. Гріна та Б. Андерсена [78], М. Дейстлера та Б. Андерсена [67], Я. Іноує та Т. Матсуї [85], М. Флізза та Х. Сіра-Раміреса [75]) або непараметричні (див. роботи

Х.Акайка [42, 43], Р. Бенткуса [4], М. Нієдзвеськи [97], Д. Боска та О. Лессі [54], В.Булдигіна, Ф.Уцета та В. Зайця [62], Л. Льюнга [91]). Галузі застосування цих моделей дуже широкі: радіофізика, гідролокація, сейсмографія, метеорологія, розвиток і взаємодія біологічних популяцій, теорія сигналів і автоматичного контролю, теорія фільтрації, економетрика та інші. В якості посилань, наприклад, можна звернутись до робіт Дж. Бокса і Г. Дженкінса [6], М. Шетцена [107, 108], Д. Брилліджера [7], Е. Вентцеля та Л. Овчарова [18], Дж. Сйоберга, К. Жанга, Л. Льюнга, А. Бенвеністе, Б.Дейона, П.-І. Глоренца, Г. Х`ялмарсона та А. Джудітські [110], Дж. Сйоберга та Л. Нгія [109], Т. Шодерстрома та П. Стойка [111], Е.-В. Баї [51], Е. Ханана та М. Дейстлера [83], тощо.

Одна із задач у теорії лінійних та нелінійних систем розглядає *оцінювання невідомої імпульсної перехідної функції* (або - часової перехідної функції) по спостереженнях за відгуками системи на відомі вхідні сигнали. Різноманітні статистичні методи для розв'язання цієї задачі й вивчення властивостей відповідних оцінок розглядались в роботах Г. Дженкінса та Д. Ваттса [21], Дж. Бокса і Г. Дженкінса [6], Г. Льюнга та Дж. Бокса [92], Д. Брилліджера [7, 55], М. Шетцена [108], М. Арато [46], Дж. Бендата та А. Пірсола [1], П. Кайнза [63], А. Дороговцева [23, 24], П. Броквелла та Р. Девіса [56], А. Кукуша [32], Е. Волтера, Дж. Нортона, Г. Пієт-Лаганьєра, М. Міланезе [100], Е.-В. Баї [51] та інших математиків. Так, Х. Акайк [42] вивчав стійкі лінійні системи з багатьма входами та одним виходом (MISO-system: multiple-input single-output system), що збурюються стаціонарними гауссівськими процесами, та отримав оцінки для перетворення Фур'є імпульсної перехідної функції по кожній компоненті. Пізніше, цей же автор у статті [43] вивчав оцінки перехідної функції при збуренні стійкої лінійної системи з одним входом та одним виходом (SISO-system: single-input single-output system) стаціонарним негауссіським процесом. Подібна задача розглядалася і у роботі М. Кьюмона [89], де вивчались векторнозначні стаціонарні часові ряди з дійснозначними компонентами, утвореними дискретним векторнозначним аналогом SISO-системи. Метод його дослідження спирався на внутрішньо - зовнішню факторизацію стійкої перехідної матриці. Деякі узагальнення для нелінійних систем були отримані, виходячи з властивостей ядер Вольтерра. Наприклад, у статті Дж. Каріоларо та Дж. ді Мазі [65] було знайдено явну формулу для обчислення математичного сподівання відгуку дискретної у часі нелінійної системи у термінах кумулянтів вхідного сигналу, які містять ядра Вольтерра. Застосування властивостей ядер Вольтерра також зустрічається у Д. Боска та О. Лессі [54], де запропонований метод рекурсивної ідентифікації невідомої перехідної функції нелінійної системи. У роботі В. Булдигіна, Ф. Уцета та В. Зайця [60] досліджуються асимптотичні властивості першого та другого моментів оцінки імпульсної перехідної функції (вибіркової моментної функції, побудованої за вхідними та вихідним сигналами) у білінійній системі Вольтерра за умови, що кожен з незалежних вхідних сигналів наближається до білого шуму, а кількість спостережень прямує до нескінченності. Для двовимірної лінійної системи, що збурюється білим шумом, у статті Фу Лі [38] вивчаються умови асимптотичної нормальності розподілів нормованої похибки оцінювання перехідної функції.

Далі зупинимось на найпростішому випадку з поставленого питання.

Розглянемо задачу оцінювання невідомої дійснозначної імпульсної перехідної функції $H = (H(\tau), \tau \in \mathbb{R})$ однорідної неперервної лінійної SISO-системи. Це означає, що відгук системи описується за допомогою згортки між вхідним процесом $x(t), t \in \mathbb{R},$ та перехідною функцією $H$, і має вигляд

$$y(t) = \int_{-\infty}^{\infty} H(t-s)x(s)ds. \qquad (1.1)$$

Зазначимо, за умови фізичної здійснимості: $H(\tau) = 0, \tau < 0$, система (1.1) називається системою типу Вольтерра.

**Обраний метод.** У дисертаційній роботі для оцінювання перехідної функції $H$ системи (1.1) застосовується метод емпіричних корелограм: в якості оцінки будується інтегральна (або дискретна за часом) сумісна корелограма між стаціонарним гауссівським вхідним процесом $X(t), t \in \mathbb{R}$, та процесом відгуку системи $Y(t), t \in \mathbb{R}$ (що також є гауссівським процесом). Мотивація такого вибору - у наступному [108]: якщо $K_{YX}(\tau) =$
$= \mathrm{E} Y(t+\tau)X(t), \tau \in \mathbb{R}$, - крос-корелограма процесів $Y$ та $X$, які припускаються центрованими, тоді підстановка виразу $Y(t) = \int_{-\infty}^{\infty} H(t-s)X(s)ds,$ у $K_{YX}$ дає

$$K_{YX}(\tau) = \int_{-\infty}^{\infty} H(s) K_{XX}(\tau - s) ds, \qquad (1.2)$$

де $K_{XX}(\tau) = \mathrm{E} X(t+\tau)X(t), \tau \in \mathbb{R},$ - кореляційна функція процесу $X$.

Якщо формально $X$ - гауссівський "білий шум", то $K_{XX}(\tau) = A\delta_0(\tau),$ $\tau \in \mathbb{R},$ де $\delta_0$ - дельта-функція Дірака у нулі. При підстановці цієї кореляційної функції у вираз (1.2), отримаємо

$$K_{YX}(\tau) = AH(\tau), \tau \in \mathbb{R}.$$

Тобто, сумісна крос-корелограма між процесами на вході та виході системи (1.1) пропорційна перехідній функції $H$ за умови, що вхідний сигнал - гауссівський білий шум. При цьому відповідна інтегральна (або дискретна за часом) корелограмна оцінка є незсуненою; її різноманітні асимптотичні властивості можна досліджувати, спираючись на відомі факти про сумісні корелограми сумісно гауссівських та сумісно стаціонарних випадкових процесів, див., наприклад, роботи І. Ібрагімова та Ю. Розанова [25], О. Іванова [26], Д. Бриллінджера [7], Ю. Козаченка [27, 28, 29], В. Булдигіна [9, 10], В. Булдигіна та Ю. Козаченка [14], тощо.

Взагалі, використання корелограмного методу при дослідженні властивостей статистик другого або старших порядків, поліноміальних статистик випадкових процесів та полів, зокрема гауссівських, зустрічається як у схемі однієї вибірки, див., монографії Г. Дженкінса та Д. Ваттса [21], Д. Бриллінджера [7], М. Ядренка [41], А. Яглома [40], Дж. Бендата та А. Пірсола [1] і статті В. Булдигіна [8, 10], М. Леоненка та А. Портнової [90], Є. Островського та С. Циркунової [35], В. Булдигіна та В. Зайця [17], В. Булдигіна та О. Диховичного [12], В. Булдигіна та Фу Лі [58, 59], В. Булдигіна та Ю. Козаченка [14], В. Булдигіна та В. Курочки [57], В. Булдигіна, Ф. Уцета та В. Зайця [60, 61, 62]; так і у схемі декількох незалежних виборок, див., роботи В. Булдигіна та Є. Іларіонова [13], М. Леоненка та О. Іванова [33], В. Булдигіна та В. Зайця [15, 16], В. Булдигіна та О. Дем'яненко [11], тощо.

**Існуючі розв'язки задачі.** Оскільки білий шум є математичною ідеалізацією фізичних процесів, то в реальних ситуаціях вводиться параметр схеми серій і використовується "наближення" до білого шуму. Так, в роботах В. Булдигіна та Фу Лі [58, 59], В. Булдигіна та В. Курочки [57], В. Булдигіна, Ф. Уцета та В. Зайця [60, 61, 62], Фу Лі [38], В. Зайця [119] прямування вхідного сигналу $X$ системи (1.1) до білого шуму відбувається за рахунок сім'ї центрованих стаціонарних гауссівських процесів $X_\Delta(t), t \in \mathbb{R}, \Delta > 0$, спектральні щільності яких збігаються до сталої на кожному обмеженому інтервалі, а кореляційні функції є дельта-видною послідовністю при $\Delta \to \infty$.

У статтях В. Булдигіна та Фу Лі [58, 59] для невідомої перехідної функції $H$ системи Вольтерра (1.1) розглядаються інтегральна корелограмна оцінка

$$\hat{H}_{T,\Delta}(\tau) = \frac{1}{cT}\int_0^T Y_\Delta(t+\tau)X_\Delta(t)dt, \tau > 0, \qquad (1.3)$$

та похибка оцінювання $\sqrt{T}[\hat{H}_{T,\Delta}(\tau) - H(\tau)], \tau > 0$, і встановлюються умови їх асимптотичної нормальності, як у сенсі збіжності скінченновимірних розподілів, так і в сенсі збіжності відповідних розподілів у просторі неперервних функцій. Оцінка (1.3) містить два параметри - довжину інтервалу спостереження $T$ та параметр схеми серій $\Delta$ - і є зсуненою ($c$ - деяка додатна стала). Крім *основного припущення* $H \in L_2(\mathbb{R})$ *на перехідну функцію*, також вимагаються дві умови на її перетворення Фур'є-Планшереля $H^*$:

(1) неперервність функції $H^*$ майже всюди на $\mathbb{R}$ (відносно міри Лебега);

(2) $H^* \in L_1(\mathbb{R}) \bigcap L_\infty(\mathbb{R})$ або $H^* \in L_{2p}(\mathbb{R})$, при деякому $p > 2$,

- і вводяться балансні умови між прямуванням параметрів $T, \Delta$ до нескінченності. Зауважимо, що для доведення асимптотичної нормальності оцінки та похибки використовується метод моментів, що в гауссівському випадку спирається на перевірку формули Іссерліса із застосуванням граничних властивостей багатовимірних згорткових інтегралів, залежних від параметра (див., Р.Бенткус [2, 3, 4], В. Булдигін та О. Диховичний [12]).

У роботах В. Булдигіна та В. Курочки [57], В. Булдигіна, Ф. Уцета та В. Зайця [60, 61, 62] для невідомої перехідної функції $H$ системи (1.1) розглядаються дискретна за часом корелограмна оцінка

$$\hat{H}_{N,h,\Delta}(\tau) = \frac{1}{cN}\sum_{n=1}^N Y_\Delta(nh+\tau)X_\Delta(nh), \tau \in \mathbb{R}, \qquad (1.4)$$

та похибка оцінювання $\sqrt{Nh}[\hat{H}_{N,h,\Delta}(\tau) - H(\tau)], \tau \in \mathbb{R}$ (тут $c$ - деяка додатна стала). Оцінка (1.4) містить три параметри: крок дискретизації $h$, довжину інтервалу спостереження $T = Nh + \tau$ та параметр схеми серій $\Delta$, - і є зсуненою. У роботах В. Булдигіна, Ф. Уцета та В. Зайця [60], В. Булдигіна та В. Курочки [57] за припущення $H \in L_2(\mathbb{R}) \bigcap L_4(\mathbb{R})$ та специфічної залежності між параметрами встановлюються умови асимптотичної незсуненості та конзистентності оцінки, а також асимптотична нормальність відповідних оцінки та похибки у сенсі збіжності скінченновимірних розподілів за допомогою методу моментів. Пізніше, у роботі В. Булдигіна, Ф. Уцета та В. Зайця [62] доведено асимптотичну нормальність оцінки (1.4) та похибки оцінювання, як у сенсі збіжності скінченновимірних розподілів, так і в сенсі збіжності відповідних розподілів у просторі неперервних функцій. Результати цієї роботи є остаточними, оскільки вимагають мінімальне обмеження на порядок інтегрованості перехідної функції $H \in L_2(\mathbb{R})$. Застосовано кумулянтний метод дослідження розподілу та інтегральні зображення кумулянтів від центрованих та нормованих випадкових величин $\hat{H}_{N,h,\Delta}(\tau_1), \hat{H}_{N,h,\Delta}(\tau_2), ..., \hat{H}_{N,h,\Delta}(\tau_m)$ у виді скінченної суми інтегралів, кожен з яких містить циклічне зачеплення ядер. Зауважимо, що асимптотична нормальність доводиться за допомогою нерівностей для багатовимірних інтегралів з циклічним зачепленням ядер, залежних від параметра, та властивості їх збіжності до нуля, встановлених у цій же роботі.

Спільне у дискретній та інтегральній схемах оцінювання - техніка встановлення слабкої збіжності розподілів у просторі неперервних функцій, яка спирається на ентропійні методи, а також класичні функціональні теореми про неперервність гауссівських випадкових процесів, див. роботи Р. Дадлі [68], [69], К. Ферніка [73, 74], Г. Крамера та Лідбеттера [31], І. Ібрагімова та Ю. Розанова [25], І. Гіхмана та А. Скорохода [19], Ю.

Прохорова та Ю. Розанова [36], П. Біллінгслі [5], Н. Джейна та М. Маркуса [87], М. Ліфшиця [34], В. Булдигіна та Ю. Козаченка [14]. Різне - техніка доведення збіжності скінченновимірних розподілів.

Природньо виникло запитання - чи можна, так само як і у дискретній, в інтегральній схемі оцінювання функції $H$ системи (1.1) обмежитись умовою $H \in L_2(\mathbb{R})$? Частина досліджень у дисертаційній роботі присвячена позитивній відповіді, яка є результатом поєднання кумулянтного методу аналізу розподілу з теорією інтегралів, що містять циклічне зачеплення ядер.

**Інструментарій.** Ключовим моментом у дисертаційній роботі є встановлення слабкої збіжності розподілів нормованої і центрованої корелограмної оцінки перехідної функції та її похибки до розподілів деякого гауссівського центрованого процесу. Це вдається зробити, показавши, що всі прості сумісні кумулянти порядку три і вище для сім'ї процесів збігаються до нуля; і спираючись на однозначну характеризацію гауссівського розподілу за допомогою простих сумісних кумулянтів (наслідок з теореми Маркова [36, 4]).

Основний прийом, що застосовується у роботі для дослідження статистичних властивостей відповідної оцінки, - кумулянтний метод Д. Бриллінджера [7]. В роботах В. Булдигіна, Ф. Уцета та В. Зайця [61], В. Зайця та Дж. Соле-Казалса [120], В. Зайця, Дж. Соле-Казалса, П. Марті-Піга та Р. Рейг-Болако [121] показано, що кумулянти від різних статистик другого порядку гауссівських та негауссівських випадкових векторів, часових рядів, випадкових процесів та полів мають специфічне інтегральне зображення. Зокрема, такий же клас інтегралів зустрічається й у інших непараметричних задачах із застосуванням кумулянтних статистик старших порядків, див., наприклад, роботи Ф. Аврама [47], Ф. Аврама та Р. Фокса [48], Ф. Аврама та М. Таггу [50], Ф. Аврама, М. Леоненка та Л. Сахно [49], де розглядаються граничні теореми для аддитивних функціоналів від стаціонарних випадкових полів. Інтеграли, що виникають при цьому, завжди мають однакову структуру - є багатовимірними згортковими інтегралами з циклічним зачепленням ядер. Деякі посилання на статті з таким типом інтегралів та їх застосуванням у статистичних задачах можна знайти також у В. Булдигіна [10] В. Булдигіна та О. Диховичного [12], В. Булдигіна та В. Курочки [57], В. Булдигіна, Ф. Уцета та В. Зайця [60, 62] та інших. Історія виникнення таких інтегралів -- розенблатівська апроксимація квадратичних форм, див. М. Розенблатт [105], де показано, що явна форма логарифма характеристичної функції розподілу Розенблата - це скінченна сума інтегралів з циклічним зачепленням ядер.

Оскільки слабка збіжність оцінок часто доводиться за допомогою кумулянтів старших порядків, див., Х.-Ч. Жанг та П. Шаман [122], Дж. Гріммет [79], Дж. Тагнайт [114], К. Нікіас та Дж. Мендель [98], Д. Брилінджер [55], У. Хаберзеттл [82], Д. Нуаларт та Дж. Пекатті [101], Дж. Пекатті та К. Тьюдор [102], Дж. Пекатті та М. Таггу [103], тощо, то виокремлене вивчення інтегралів з циклічним зачепленням ядер дає кращий результат, ніж вивчення їх як інтегралів взагалі. Так, інтерпретація інтегралів з циклічним зачепленням ядер як багатовимірної згортки із зачепленими аргументами допомагає отримати гарні верхні межі, див. статті В. Булдигіна, Ф. Уцета та В. Зайця [61, 62], Ф. Аврама, М. Леоненка та Л. Сахно [49]. Ранні роботи з інтегральних зображень кумулянтів статистик другого порядку стаціонарних випадкових процесів належать литовським статистикам Р. Бенткусу [2, 3] та В. Статулевичесу [112]. Пізніша робота Р. Бенткуса [53] пов'язана з кумулянтами поліноміальних статистик. Книга А. Мазаї та С. Провоста [95] - добре посилання на квадратичні форми, а роботи А. Мазаї, С. Провоста та Т. Найякави [96] та Б. Холмквеста [84] - сконцентровані на дослідженнях квадратичних та білінійних форм від нормальних випадкових величин. У роботі В. Булдигіна, Ф. Уцета та В. Зайця [61] показано, що зображення кумулянтів від різних білінійних форм гауссівських випадкових векторів, стаціонарних випадкових процесів та однорідних гауссівських випадкових полів містять скінченні суми інтегралів з циклічним зачепленням ядер. У пізнішій роботі тих же авторів [62] встановлено деякі нерівності для інтегралів з циклічним зачепленням ядер, залежних від декількох параметрів, та умови їх збіжності до нуля.

Інші інтегральні зображення кумулянтів періодограм однорідних випадкових полів розглядались у роботах Х. Гуйона [81], А. Бенна та Р. Кулпергера [52], М. Розенблатта [106]. Зокрема, цікаві класи багатовимірних випадкових інтегралів були запропоновані Д. Сургалісом [113], Д. Енджелом [71], Р. Фоксом та М. Таггу [76].

# РОЗДІЛ 2

# АСИМПТОТИЧНІ ВЛАСТИВОСТІ КОРЕЛОГРАМНИХ ОЦІНОК ПЕРЕХІДНИХ ФУНКЦІЙ ЛІНІЙНИХ ОДНОРІДНИХ СИСТЕМ

У розділі 2 розглядається метод емпіричних корелограм для оцінювання невідомої дійснозначної перехідної функції $H$ неперервної однорідної лінійної системи з умовою $H \in L_2(\mathbb{R})$. Таке припущення на порядок інтегрованості перехідної функції дозволяє досліджувати нестійкі системи з резонансними особливостями. Припускається, що на вхід однорідної лінійної системи подається сім'я вимірних сепарабельних стаціонарних центрованих гауссівських процесів близьких, в деякому сенсі, до білого шуму.

Використовуючи схему однієї вибірки, досліджуються умови асимптотичної незсуненості та консистентності корелограмних оцінок інтегрального типу, а також умови асимптотичної нормальності цієї оцінки та відповідної похибки оцінювання, як у сенсі збіжності скінченновимірних розподілів, так і в сенсі збіжності відповідних розподілів у просторі неперервних функцій. Зокрема, наводяться приклади збурюючих процесів та імпульсних перехідних функцій, які задовольняють розв'язуваній задачі.

Знайдені результати далі застосовуються до розв'язання більш загальної задачі - корелограмного оцінювання перехідної функції $H \in L_2(\mathbb{R})$ неперервної однорідної лінійної системи з внутрішнім шумом.

## 2.1 Оцінка імпульсної перехідної функції

У підрозділі визначається статистична оцінка імпульсної перехідної функції неперервної однорідної лінійної системи. Ця оцінка - відповідним чином нормована сумісна інтегральна корелограма між деяким стаціонарним гауссівським випадковим процесом, що збурює систему, та реакцією системи на це збурення.

**Процеси, що збурюють систему.** Нехай $X_\Delta = (X_\Delta(t), t \in \mathbb{R})$, $\Delta > 0$, - сім'я вимірних сепарабельних стаціонарних центрованих дійснозначних гауссівських процесів, що збурюють однорідну лінійну систему (див. додаток А.4). Сім'я невід'ємних неперервних функцій $f_\Delta = (f_\Delta(\lambda), \lambda \in \mathbb{R})$, $\Delta > 0$, є спектральними щільностями процесів $X_\Delta$ та задовольняє умовам:

$$f_\Delta(\lambda) = f_\Delta(-\lambda), \lambda \in \mathbb{R}; \tag{2.1а}$$

$$\sup_{\Delta>0} \| f_\Delta \|_\infty < \infty; \tag{2.1б}$$

$$f_\Delta \in L_1(\mathbb{R}); \tag{2.1в}$$

$$\exists c \in (0, \infty) \; \forall a \in (0, \infty) : \lim_{\Delta \to \infty} \sup_{-a \leq \lambda \leq a} \left| f_\Delta(\lambda) - \frac{c}{2\pi} \right| = 0; \tag{2.1г}$$

$$K_\Delta \in L_1(\mathbb{R}), \tag{2.1д}$$

де $K_\Delta(t) = \mathrm{E} X_\Delta(s+t) X_\Delta(s) = \int_{-\infty}^{\infty} e^{i\lambda t} f_\Delta(\lambda) d\lambda$, $\lambda \in \mathbb{R}$, - кореляційна функція процесу $X_\Delta$.

Далі в роботі припускається, що умови (2.1а) - (2.1д) завжди виконані. Умова 2.1г показує, в якому саме сенсі слід розуміти "близькість" сім'ї процесів $X_\Delta, \Delta > 0$, до гауссівського білого шуму при $\Delta \to \infty$.

Зазначимо, що парність функції $f_\Delta$ гарантує дійснозначність кореляційної функції $K_\Delta$. З умов (2.1б),(2.1в) випливає, що $f_\Delta \in L_2(\mathbb{R})$, тому з теореми

Фур'є - Планшереля (див. додаток А.1) маємо, що $K_\Delta \in L_2(\mathbb{R})$. Крім того, з умови (2.1в) слідує, що $K_\Delta(t), t \in \mathbb{R}$, є неперервною функцією на дійсній осі, тому стаціонарний процес $X_\Delta$ є неперервним у середньому квадратичному.

*Приклад 2.1.* Наступні спектральні щільності $f_\Delta$ та відповідні їм кореляційні функції $K_\Delta$ процесів $X_\Delta$:

$$f_\Delta = (\frac{c}{2\pi}\exp\left(-\frac{\lambda^2}{\Delta}\right), \lambda \in \mathbb{R}) \text{ та } K_\Delta = (\frac{c}{2}\sqrt{\frac{\Delta}{\pi}}\exp(-\frac{\Delta t^2}{4}), t \in \mathbb{R});$$

$$f_\Delta = (\frac{c}{2\pi}\frac{\Delta}{\Delta+\lambda^2}, \lambda \in \mathbb{R}) \text{ та } K_\Delta = (\frac{c\sqrt{\Delta}}{2}\exp(-\sqrt{\Delta}\,|t|), t \in \mathbb{R});$$

$$f_\Delta = (\frac{c}{2\pi}\exp\left(-\frac{|\lambda|}{\sqrt{\Delta}}\right), \lambda \in \mathbb{R}) \text{ та } K_\Delta = (\frac{c}{\pi}\frac{\sqrt{\Delta}}{1+\Delta t^2}), t \in \mathbb{R});$$

$$f_\Delta = (\frac{c}{2\pi}(1-\frac{|\lambda|}{\sqrt{\Delta}})\mathbb{I}_{[-\sqrt{\Delta},\sqrt{\Delta}]}(\lambda), \lambda \in \mathbb{R}) \text{ та } K_\Delta = (\frac{c\sqrt{\Delta}}{2\pi}\cdot(\frac{\sin(\frac{\sqrt{\Delta}t}{2})}{\sqrt{\Delta}t})^2, t \in \mathbb{R}),$$

задовольняють умовам (2.1а) - (2.1д).

**Вигляд оцінки для перехідної функції.** Нехай $H \in L_2(\mathbb{R})$ - перехідна функція однорідної лінійної системи (див. додаток А.4).

Розглянемо випадковий процес

$$Y_\Delta(t) = \int_{-\infty}^{\infty} H(s) X_\Delta(t-s) ds, \ t \in \mathbb{R},$$

який є відгуком системи на вхідний процес $X_\Delta$.

Відомо (див. теорему А.4), що інтеграл, який задає випадковий процес $Y_\Delta$, коректно визначений як відповідна границя у середньому квадратичному, тоді й лише тоді, коли існує інтеграл

$$\int_{-\infty}^{\infty}\int_{-\infty}^{\infty} K_\Delta(t-s) H(s) H(t) ds dt. \tag{2.2}$$

Відповідь про коректність визначення процесу $Y_\Delta$ дає таке твердження.

**Лема 2.1.** *Нехай $K_\Delta \in L_1(\mathbb{R})$ та $H \in L_2(\mathbb{R})$, тоді $Y_\Delta$ коректно визначений і є стаціонарним центрованим неперервним у середньому квадратичному гауссівським випадковим процесом зі спектральною щільністю*

$$\phi_\Delta(\lambda) = |H^*(\lambda)|^2 \, f_\Delta(\lambda), \lambda \in \mathbb{R}, \qquad (2.3)$$

*де $H^*$ - перетворення Фур'є - Планшереля функції $H$. Крім того, процеси $X_\Delta$ і $Y_\Delta$ є сумісно стаціонарними та сумісно гауссівськими випадковими процесами.*

*Доведення.* З нерівності Коші-Буняковського та теореми Фубіні - Тонеллі (див. теорему А.1) випливає, що

$$\int_{-\infty}^{\infty}\int_{-\infty}^{\infty} |K_\Delta(t-s)H(s)H(t)|\,dsdt = \int_{-\infty}^{\infty} |H(t)| \left(\int_{-\infty}^{\infty} |K_\Delta(t-s)||H(s)|\,ds\right)dt =$$

$$= \int_{-\infty}^{\infty} |H(t)| (|K_\Delta| * |H|)(t)\,dt \leq \|H\|_2 \| |K_\Delta| * |H| \|_2.$$

З нерівності Юнга (див. теорему А.2) випливає оцінка

$$\| |K_\Delta| * |H| \|_2 \leq \|H\|_2 \|K_\Delta\|_1.$$

Таким чином,

$$\int_{-\infty}^{\infty}\int_{-\infty}^{\infty} |K_\Delta(t-s)H(s)H(t)|\,dsdt \leq \|H\|_2^2 \|K_\Delta\|_1 < \infty.$$

Звідси випливає існування інтеграла (2.2), а разом з тим й процесу $Y_\Delta$. Інші властивості процесу $Y_\Delta$ мають місце в силу його означення (див. додаток А.4). Таким чином, лему 2.1 доведено. □

*Зауваження 2.1.* Інтеграл (2.2) існує, наприклад, якщо $H \in L_1(\mathbb{R})$, тобто якщо система є стійкою. У цьому випадку умова (2.1д) стає зайвою.

**Означення 2.1** Оцінку для $H$ будемо шукати у виді інтегральної сумісної корелограми

$$H_{T,\Delta}(\tau) = \frac{1}{cT}\int_0^T Y_\Delta(t+\tau)X_\Delta(t)dt, \tau \in \mathbb{R}, \qquad (2.4)$$

де $c$ - стала з умови (2.1г), $\Delta$ - параметр схеми серій, й $T$ - довжина інтервалу усереднення $[0,T]$. Зазначимо, що для побудови оцінки процеси $X_\Delta$ і $Y_\Delta$ мають спостерігатись на всій дійсній осі.

Інтеграл, що фігурує у визначенні оцінки $H_{T,\Delta}$, слід розуміти як середньоквадратичний інтеграл Рімана.

З формули (2.4), для довільних $T>0, \Delta>0$, при всіх $\tau \in \mathbb{R}$, маємо

$$\mathrm{E}H_{T,\Delta}(\tau) = \frac{1}{c}\mathrm{E}Y_\Delta(t+\tau)X_\Delta(t) = \frac{1}{c}\int_{-\infty}^{\infty} K_\Delta(\tau-s)H(s)ds. \qquad (2.5)$$

Взагалі кажучи, $\mathrm{E}H_{T,\Delta}(\tau) \neq H(\tau)$, тобто оцінка (2.4) є зсуненою.

**Вигляд похибки оцінювання.** Для подальших досліджень додатково розглянемо *нормовану похибку корелограмного оцінювання* виду

$$W_{T,\Delta}(\tau) = \sqrt{T}[H_{T,\Delta}(\tau) - H(\tau)], \tau \in \mathbb{R}, \qquad (2.6)$$

та поставимо питання про її асимптотичні властивості при прямуванні параметрів $T, \Delta$ до безмежності. Для цього процес $W_{T,\Delta} = (W_{T,\Delta}(\tau), \tau \in \mathbb{R})$ розщеплюється на суму

$$W_{T,\Delta} = Z_{T,\Delta} + V_{T,\Delta},$$

де доданки справа визначаються далі, й для них спрацьовують такі факти:

- Вибір послідовності збурюючих процесів $(X_\Delta, \Delta > 0)$, апроксимуючої гауссівський білий шум, - є специфічним: *умови (2.1а) - (2.1д) є достатніми для нормалізації процесу $Z_{T,\Delta}$ при довільній збіжності параметрів $T, \Delta$ до безмежності.*

- Знищення вкладу невипадкової функції $V_{T,\Delta}$ у виразі для $W_{T,\Delta}$ можливе за рахунок використання $\delta$-видності сім'ї $K_\Delta$, підбору характеру сумісного прямування параметрів $T, \Delta$ до безмежності, та накладення умов на порядок локальної гладкості $H$.

## 2.2 Асимптотична поведінка процесу $Z_{T,\Delta}$

Даний підрозділ присвячений дослідженню асимптотичних властивостей нормованої оцінки $H_{T,\Delta}$, центрованої своїм середнім значенням.

### 2.2.1 Кореляційна функція процесу $Z_{T,\Delta}$

У цьому пункті встановлюється вид кореляційної функції емпіричного процесу $Z_{T,\Delta}$, пов'язаного з оцінкою перехідної функції. Досліджуються властивості цієї кореляційної функції. Нагадаємо, що оцінка $H_{T,\Delta}$ та функції $f_\Delta$, визначені у підрозділі 2.1.

Покладемо

$$Z_{T,\Delta}(\tau) = \sqrt{T}[H_{T,\Delta}(\tau) - \mathrm{E}H_{T,\Delta}(\tau)], \tau \in \mathbb{R}, \qquad (2.7)$$

і зазначимо, що визначений процес є центрованим $L_2$-процесом.

**Лема 2.2.** *Нехай $H \in L_2(\mathbb{R})$, тоді для всіх $T > 0, \Delta > 0$, та $\tau_1, \tau_2 \in \mathbb{R}$, має місце*

*рівність*

$$\mathrm{E}Z_{T,\Delta}(\tau_1)Z_{T,\Delta}(\tau_2) = C_{T,\Delta}(\tau_1,\tau_2) =$$

$$= \frac{2\pi}{c^2}\int_{-\infty}^{\infty}\int_{-\infty}^{\infty}[e^{i(\tau_1-\tau_2)\lambda_2}\,|\,H^*(\lambda_2)|^2 + e^{i(\tau_1\lambda_1+\tau_2\lambda_2)}H^*(\lambda_1)H^*(\lambda_2)] \times \qquad (2.8)$$

$$\times \Phi_T(\lambda_2-\lambda_1)f_\Delta(\lambda_1)f_\Delta(\lambda_2)d\lambda_1 d\lambda_2;$$

*де $\Phi_T$ - ядро Фейєра:*

$$\Phi_T(\lambda) = \frac{1}{2\pi T}\left(\frac{\sin(T\lambda/2)}{\lambda/2}\right)^2, \lambda \in \mathbb{R};$$

*$c$ - стала з умови (2.1г), та $H^*$ - перетворення Фур'є-Планшереля функції $H$ у просторі $L_2(\mathbb{R})$.*

*Доведення.* Оскільки кумулянт другого порядку - лінійна функція своїх аргументів, інваріантна відносно зсувів [7], то

$$\mathrm{E}Z_{T,\Delta}(\tau_1)Z_{T,\Delta}(\tau_2) = (\sqrt{T})^2 cov(H_{T,\Delta}(\tau_1)H_{T,\Delta}(\tau_2)) =$$

$$= T\cdot(\frac{1}{cT})^2 cov(\int_0^T Y_\Delta(u+\tau_1)X_\Delta(u)du, \int_0^T Y_\Delta(v+\tau_2)X_\Delta(v)dv) =$$

$$= \frac{1}{c^2 T}\int_0^T\int_0^T cov(Y_\Delta(u+\tau_1)X_\Delta(u), Y_\Delta(v+\tau_2)X_\Delta(v))dudv. \qquad (2.9)$$

Для підінтегральної функції - кумулянта від набору, що є добутками центрованих випадкових величин, - застосуємо теорему Леонова - Ширяєва - Бриллінджера (див. теорему А.5) та лему 2.1:

$$cov(Y_\Delta(u+\tau_1)X_\Delta(u), Y_\Delta(v+\tau_2)X_\Delta(v)) =$$
$$= cov(Y_\Delta(u+\tau_1), Y_\Delta(v+\tau_2))\cdot cov(X_\Delta(v), X_\Delta(u)) +$$

$$+\mathrm{cov}(Y_\Delta(u+\tau_1), X_\Delta(v)) \cdot \mathrm{cov}(Y_\Delta(v+\tau_2), X_\Delta(u)) =$$

$$= K_{Y_\Delta}(u-v+\tau_1-\tau_2) \cdot K_{X_\Delta}(v-u) + K_{Y_\Delta,X_\Delta}(u-v+\tau_1) \cdot K_{Y_\Delta,X_\Delta}(v-u+\tau_2),$$

де $K_{Y_\Delta}$ та $K_{X_\Delta}$ - відповідно кореляційні функції процесів $Y_\Delta$ та $X_\Delta$, а $K_{Y_\Delta,X_\Delta}$ - сумісна кореляційна функція процесів $Y_\Delta$ та $X_\Delta$.

Зробивши заміну змінних $\begin{cases} t = u - v; \\ s = u + v, \end{cases}$ в останньому виразі, та врахувавши зміну області $\{(u,v) \in [0;T] \times [0;T]\}$ на $\{(t,s) \in [-T;T] \times [-T+|t|; T-|t|]\}$, й перехід $dudv = \frac{1}{2} dtds$, в результаті з формули (2.9) дістанемо:

$$\mathrm{E} Z_{T,\Delta}(\tau_1) Z_{T,\Delta}(\tau_2) = \qquad (2.10)$$

$$= \frac{1}{c^2} \int\limits_{-T}^{T} [K_{Y_\Delta}(t+\tau_1-\tau_2) \cdot K_{X_\Delta}(-t) + K_{Y_\Delta,X_\Delta}(t+\tau_1) \cdot K_{Y_\Delta,X_\Delta}(-t+\tau_2)](1-\frac{|t|}{T})dt.$$

Підставивши наступні спектральні представлення:

$$K_{Y_\Delta}(t) = \int\limits_{-\infty}^{\infty} e^{it\lambda} |H^*(\lambda)|^2 f_\Delta(\lambda) d\lambda;$$

$$K_{X_\Delta}(t) = \int\limits_{-\infty}^{\infty} e^{it\lambda} f_\Delta(\lambda) d\lambda;$$

$$K_{Y_\Delta,X_\Delta}(t) = \int\limits_{-\infty}^{\infty} e^{it\lambda} H^*(\lambda) f_\Delta(\lambda) d\lambda,$$

у формулу (2.10), отримаємо вираз:

$$\mathrm{E} Z_{T,\Delta}(\tau_1) Z_{T,\Delta}(\tau_2) =$$

$$= \frac{1}{c^2} \int\limits_{-T}^{T} (\int\limits_{-\infty}^{\infty} e^{i(-t)\lambda_1} f_\Delta(\lambda_1) d\lambda_1) \cdot (\int\limits_{-\infty}^{\infty} e^{i(t+\tau_1-\tau_2)\lambda_2} |H^*(\lambda_2)|^2 f_\Delta(\lambda_2) d\lambda_2)(1-\frac{|t|}{T})dt +$$

$$+\frac{1}{c^2}\int\limits_{-T}^{T}(\int\limits_{-\infty}^{\infty}e^{i(t+\tau_1)\lambda_1}H^*(\lambda_1)f_\Delta(\lambda_1)d\lambda_1)\times(\int\limits_{-\infty}^{\infty}e^{i(-t+\tau_2)\lambda_2}H^*(\lambda_2)f_\Delta(\lambda_2)d\lambda_2)(1-\frac{|t|}{T})dt=$$

$$=\frac{1}{c^2}\int\limits_{-\infty}^{\infty}\int\limits_{-\infty}^{\infty}e^{i(\tau_1-\tau_2)\lambda_2}|H^*(\lambda_2)|^2 f_\Delta(\lambda_1)f_\Delta(\lambda_2)[\int\limits_{-T}^{T}e^{it(\lambda_2-\lambda_1)}(1-\frac{|t|}{T})dt]d\lambda_1 d\lambda_2+$$

$$+\frac{1}{c^2}\int\limits_{-\infty}^{\infty}\int\limits_{-\infty}^{\infty}e^{i(\tau_1\lambda_1+\tau_2\lambda_2)}H^*(\lambda_1)H^*(\lambda_2)f_\Delta(\lambda_1)f_\Delta(\lambda_2)[\int\limits_{-T}^{T}e^{it(\lambda_1-\lambda_2)}(1-\frac{|t|}{T})dt]d\lambda_1 d\lambda_2.$$

Оскільки функція

$$\Phi_T(\lambda)=\int\limits_{-T}^{T}e^{it\lambda}(1-\frac{|t|}{T})dt=\frac{1}{T}\left(\frac{\sin(T\lambda/2)}{\lambda/2}\right)^2, \lambda\in\mathbb{R},$$

приймає дійсні значення та є парною, і пов'язана з ядром Фейєра співвідношенням: $\Phi_T(\lambda)=2\pi\Phi_T(\lambda), \lambda\in\mathbb{R}$, остаточно дістаємо

$$\mathrm{E}Z_{T,\Delta}(\tau_1)Z_{T,\Delta}(\tau_2)=\frac{2\pi}{c^2}\int\limits_{-\infty}^{\infty}\int\limits_{-\infty}^{\infty}[e^{i(\tau_1-\tau_2)\lambda_2}|H^*(\lambda_2)|^2+e^{i(\tau_1\lambda_1+\tau_2\lambda_2)}H^*(\lambda_1)H^*(\lambda_2)]\times$$

$$\times\Phi_T(\lambda_2-\lambda_1)f_\Delta(\lambda_1)f_\Delta(\lambda_2)d\lambda_1 d\lambda_2.$$

Таким чином, лему 2.2 доведено. $\square$

Зазначимо, що кореляційна функція процесу $Z_{T,\Delta}$ є невід'ємно визначеною неперервною функцією на $\mathbb{R}\times\mathbb{R}$.

*Зауваження 2.2.* Ядра Фейєра $\Phi_T$, що зустрічаються у формулі (2.8), залежать лише від одного параметра усереднення $T>0$, і не залежать від параметра схеми серій $\Delta>0$.

*Зауваження 2.3.* Наступні властивості ядер Фейєра вірні при всіх $T>0$:

(i) $\Phi_T(\lambda)\geq 0$, і $\Phi_T(-\lambda)=\Phi_T(\lambda), \lambda\in\mathbb{R}$;

(ii) $\int\limits_{-\infty}^{\infty}\Phi_T(\lambda)d\lambda=\|\Phi_T\|_1=1$;

(iii) $\forall \varepsilon > 0 \quad \lim_{T \to \infty} \int_{|\lambda| > \varepsilon} \Phi_T(\lambda) d\lambda = 0,$

й показують, що функції $(\Phi_T, T > 0)$ є $\delta$ - видною сім'єю при $T \to \infty$.

**Лема 2.3.** *Нехай* $H \in L_2(\mathbb{R})$*, тоді для всіх* $T > 0, \Delta > 0,$ *та* $\tau_1, \tau_2 \in \mathbb{R}$*, має місце нерівність*

$$|\mathrm{E} Z_{T,\Delta}(\tau_1) Z_{T,\Delta}(\tau_2)| \leq \frac{4\pi (\sup_{\Delta > 0} \| f_\Delta \|_\infty)^2 \| H^* \|_2^2}{c^2}.$$

*Доведення.* Будемо оцінювати функцію, визначену формулою (2.8), розбивши її на два доданки наступним чином:

$$|\mathrm{E} Z_{T,\Delta}(\tau_1) Z_{T,\Delta}(\tau_2)| =$$

$$= |\frac{2\pi}{c^2} \int_{-\infty}^{\infty} \int_{-\infty}^{\infty} e^{i(\tau_1 - \tau_2)\lambda_2} |H^*(\lambda_2)|^2 \Phi_T(\lambda_2 - \lambda_1) f_\Delta(\lambda_1) f_\Delta(\lambda_2) d\lambda_1 d\lambda_2 +$$

$$+ \frac{2\pi}{c^2} \int_{-\infty}^{\infty} \int_{-\infty}^{\infty} e^{i(\tau_1 \lambda_1 + \tau_2 \lambda_2)} H^*(\lambda_1) H^*(\lambda_2) \Phi_T(\lambda_2 - \lambda_1) f_\Delta(\lambda_1) f_\Delta(\lambda_2) d\lambda_1 d\lambda_2 |$$

В силу умови (2.1б), звідси отримаємо

$$|\mathrm{E} Z_{T,\Delta}(\tau_1) Z_{T,\Delta}(\tau_2)| \leq$$

$$\leq \frac{2\pi}{c^2} (\sup_{\Delta > 0} \| f_\Delta \|_\infty)^2 \int_{-\infty}^{\infty} \int_{-\infty}^{\infty} |H^*(\lambda_2)|^2 \Phi_T(\lambda_2 - \lambda_1) d\lambda_1 d\lambda_2 +$$

$$+ \frac{2\pi}{c^2} (\sup_{\Delta > 0} \| f_\Delta \|_\infty)^2 \int_{-\infty}^{\infty} \int_{-\infty}^{\infty} |H^*(\lambda_1)| |H^*(\lambda_2)| \Phi_T(\lambda_2 - \lambda_1) d\lambda_1 d\lambda_2 =$$

$$= \frac{2\pi}{c^2} (\sup_{\Delta > 0} \| f_\Delta \|_\infty)^2 \| H^* \|_2^2 \| \Phi_T \|_1 +$$

$$+ \frac{2\pi}{c^2} (\sup_{\Delta > 0} \| f_\Delta \|_\infty)^2 \int_{-\infty}^{\infty} |H^*(\lambda_2)| (|H^*| * \Phi_T)(\lambda_2) d\lambda_2.$$

З нерівності Юнга (див. теорему А.2), застосованої до другого доданка, й того, що $\|\Phi_T\|_1 = 1$, в результаті дістанемо

$$|\mathrm{E} Z_{T,\Delta}(\tau_1) Z_{T,\Delta}(\tau_2)| \leq \frac{2\pi}{c^2}(\sup_{\Delta>0}\|f_\Delta\|_\infty)^2 \|H^*\|_2^2 + \frac{2\pi}{c^2}(\sup_{\Delta>0}\|f_\Delta\|_\infty)^2 \|H^*\|_2^2 \|\Phi_T\|_1 =$$

$$= \frac{4\pi}{c^2}(\sup_{\Delta>0}\|f_\Delta\|_\infty)^2 \|H^*\|_2^2.$$

Таким чином, лему 2.3 доведено. □

**Наслідок 2.1.** *Оцінка кореляційної функції процесу $Z_{T,\Delta}$, наведена у лемі 2.3, не залежить від параметрів $T > 0, \Delta > 0$, та значень $\tau_1, \tau_2 \in \mathbb{R}$. Тобто при $H \in L_2(\mathbb{R})$, насправді встановлено співвідношення:*

$$\sup_{T,\Delta>0}\sup_{\tau_1,\tau_2\in\mathbb{R}} |\mathrm{E} Z_{T,\Delta}(\tau_1) Z_{T,\Delta}(\tau_2)| \leq \frac{4\pi}{c^2}(\sup_{\Delta>0}\|f_\Delta\|_\infty)^2 \|H^*\|_2^2.$$

**Наслідок 2.2.** *Нехай виконуються умови леми 2.3, тоді для всіх $T > 0$, $\Delta > 0$, та $\tau \in \mathbb{R}$, має місце нерівність*

$$\mathrm{V}ar(H_{T,\Delta}(\tau)) = \mathrm{E}|H_{T,\Delta}(\tau) - \mathrm{E} H_{T,\Delta}(\tau)|^2 \leq \frac{4\pi(\sup_{\Delta>0}\|f_\Delta\|_\infty)^2 \|H^*\|_2^2}{c^2 T}.$$

### 2.2.2 Гранична поведінка кореляційної функції процесу $Z_{T,\Delta}$

У пункті 2.2.1 був знайдений вигляд кореляційної функції $C_{T,\Delta}$ процесу $Z_{T,\Delta}$ (див. формулу (2.8)). У даному пункті знаходиться границя цієї функції при прямуванні параметрів $T$ і $\Delta$ до нескінченності.

Далі запис $(T,\Delta) \to \infty$, означає, що одночасно $T \to \infty$ та $\Delta \to \infty$.

Покладемо для всіх $\tau_1, \tau_2 \in \mathbb{R}$

$$C_\infty(\tau_1,\tau_2) = \frac{1}{2\pi}\int_{-\infty}^{\infty}\left[e^{i(\tau_1-\tau_2)\lambda}|H^*(\lambda)|^2 + e^{i(\tau_1+\tau_2)\lambda}(H^*(\lambda))^2\right]d\lambda. \qquad (2.11)$$

Зазначимо, що функція $C_\infty = (C_\infty(\tau_1,\tau_2), \tau_1,\tau_2 \in \mathbb{R})$, коректно визначена та є неперервною функцією за умови $H \in L_2(\mathbb{R})$. Зокрема, цю функцію можна записати у наступному вигляді:

$$C_\infty(\tau_1,\tau_2) = \frac{1}{\pi}\int_0^\infty [\cos(\tau_1-\tau_2)\lambda \cdot |H^*(\lambda)|^2 +$$

$$+\cos(\tau_1+\tau_2)\lambda \cdot Re(H^*(\lambda))^2 - \sin(\tau_1+\tau_2)\lambda \cdot Im(H^*(\lambda))^2]d\lambda.$$

звідки видно, що вона дійснозначна.

*Зауваження 2.4.* Розглянемо частинні випадки зображень для функції $C_\infty$:

(I) Нехай $H$ - парна дійснозначна імпульсна перехідна функція. Тоді її перетворення Фур'є - Планшереля має вигляд

$$H^*(\lambda) = 2\int_0^\infty H(\tau)\cos(\tau\lambda)d\tau,$$

і є чисто дійсною функцією, тому $(H^*(\lambda))^2 = |H^*(\lambda)|^2, \lambda \in \mathbb{R}$. Враховуючи це, з формули (2.11), після деяких перетворень, дістанемо

$$C_\infty(\tau_1,\tau_2) = \frac{2}{\pi}\int_0^\infty |H^*(\lambda)|^2 \cos(\tau_1\lambda)\cos(\tau_2\lambda)d\lambda.$$

(II) Нехай $H$ - непарна дійснозначна імпульсна перехідна функція. Тоді її перетворення Фур'є - Планшереля має вигляд

$$H^*(\lambda) = -2i\int_0^\infty H(\tau)\sin(\tau\lambda)d\tau,$$

і є чисто уявною функцією, тому $(H^*(\lambda))^2 = -|H^*(\lambda)|^2, \lambda \in \mathbb{R}$. Враховуючи це, з формули (2.11), після деяких перетворень, дістанемо

$$C_\infty(\tau_1,\tau_2) = \frac{2}{\pi}\int_0^\infty |H^*(\lambda)|^2 \sin(\tau_1\lambda)\sin(\tau_2\lambda)d\lambda.$$

**Теорема 2.1.** *Нехай* $H \in L_2(\mathbb{R})$, *тоді для всіх* $\tau_1,\tau_2 \in \mathbb{R}$, *має місце рівність*

$$\lim_{(T,\Delta)\to\infty} \mathrm{E} Z_{T,\Delta}(\tau_1)Z_{T,\Delta}(\tau_2) = C_\infty(\tau_1,\tau_2).$$

*Доведення.* При доведенні теореми 2.1 розвиваються загальні підходи, запропоновані у роботах [60] та [57]. Для простоти викладок, розіб'ємо доведення на два етапи. На першому етапі ми доведемо, що

$$\lim_{(T,\Delta)\to\infty} C^{(1)}_{T,\Delta}(\tau_1,\tau_2) = C^{(1)}_\infty(\tau_1,\tau_2), \tag{2.12}$$

де

$$C^{(1)}_{T,\Delta}(\tau_1,\tau_2) = \frac{2\pi}{c^2}\int_{-\infty}^{\infty}\int_{-\infty}^{\infty} e^{i(\tau_1-\tau_2)\lambda_2}|H^*(\lambda_2)|^2 f_\Delta(\lambda_1)f_\Delta(\lambda_2)\Phi_T(\lambda_2-\lambda_1)d\lambda_1 d\lambda_2;$$

$$C^{(1)}_\infty(\tau_1,\tau_2) = \frac{1}{2\pi}\int_{-\infty}^{\infty} e^{i(\tau_1-\tau_2)\lambda}|H^*(\lambda)|^2 d\lambda,\ \tau_1,\tau_2 \in \mathbb{R}.$$

На другому етапі отримаємо наступну рівність:

$$\lim_{(T,\Delta)\to\infty} C^{(2)}_{T,\Delta}(\tau_1,\tau_2) = C^{(2)}_\infty(\tau_1,\tau_2), \tag{2.13}$$

де

$$C^{(2)}_{T,\Delta}(\tau_1,\tau_2) = \frac{2\pi}{c^2}\int_{-\infty}^{\infty}\int_{-\infty}^{\infty} e^{i(\tau_1\lambda_1+\tau_2\lambda_2)} H^*(\lambda_1)H^*(\lambda_2) f_\Delta(\lambda_1)f_\Delta(\lambda_2)\Phi_T(\lambda_2-\lambda_1)d\lambda_1 d\lambda_2;$$

$$C^{(2)}_\infty(\tau_1,\tau_2) = \frac{1}{2\pi}\int_{-\infty}^{\infty} e^{i(\tau_1+\tau_2)\lambda}(H^*(\lambda))^2 d\lambda,\ \tau_1,\tau_2 \in \mathbb{R}.$$

- *Етап 1.* Зафіксуємо довільне $b > 0$, і через $\mathbb{I}_{[-b/2,b/2]}$ позначимо індикатор сегмента $[-b/2, b/2]$. Далі, враховуючи, що $\int_{-\infty}^{\infty}\Phi_T(\lambda)d\lambda = 1$, для будь-яких

$T > 0, \Delta > 0,$ та $\tau_1, \tau_2 \in \mathbb{R},$ розглянемо різницю

$$C_{\infty}^{(1)}(\tau_1, \tau_2) - C_{T,\Delta}^{(1)}(\tau_1, \tau_2) = \frac{1}{2\pi}[d_1(b) + d_2(b,T,\Delta) + d_3(b,T,\Delta)],$$

де

$$d_1(b) = \int_{-\infty}^{\infty} e^{i(\tau_1-\tau_2)\lambda} |H^*(\lambda)|^2 [1 - \mathbb{I}_{[-b/2,b/2]}(\lambda)]d\lambda;$$

$$d_2(b,T,\Delta) = \int_{-\infty}^{\infty}\int_{-\infty}^{\infty} e^{i(\tau_1-\tau_2)\lambda_2} |H^*(\lambda_2)|^2 \mathbb{I}_{[-b/2,b/2]}(\lambda_2) \times$$

$$\times \left[1 - \left(\frac{2\pi}{c}\right)^2 f_{\Delta}(\lambda_1) f_{\Delta}(\lambda_2)\right] \Phi_T(\lambda_2 - \lambda_1) d\lambda_1 d\lambda_2;$$

$$d_3(b,T,\Delta) =$$

$$= \left(\frac{2\pi}{c}\right)^2 \int_{-\infty}^{\infty}\int_{-\infty}^{\infty} e^{i(\tau_1-\tau_2)\lambda_2} |H^*(\lambda_2)|^2 [\mathbb{I}_{[-b/2,b/2]}(\lambda_2) - 1] f_{\Delta}(\lambda_1) f_{\Delta}(\lambda_2) \Phi_T(\lambda_2 - \lambda_1) d\lambda_1 d\lambda_2.$$

Тепер розглянемо поведінку кожної з величин $d_k, k \in \{1,2,3\}$. Оскільки для будь-якого $b > 0$

$$|d_1(b)| = |\int_{|\lambda|>b/2} e^{i(\tau_1-\tau_2)\lambda} |H^*(\lambda)|^2 d\lambda| \leq \int_{|\lambda|>b/2} |H^*(\lambda)|^2 d\lambda,$$

та $|H^*|^2 \in L_1(\mathbb{R}),$ то одержимо

$$\lim_{b \to \infty} d_1(b) = 0. \tag{2.14}$$

Для будь-яких $b > 0$, та $T > 0, \Delta > 0,$ має місце наступна оцінка для $d_2$:

$$|d_2(b,T,\Delta)| \leq B_1(T,\Delta) + B_2(T,\Delta),$$

де

$$B_1(T,\Delta) = |\iint_{D_b} e^{i(\tau_1-\tau_2)\lambda_2} |H^*(\lambda_2)|^2 \mathbb{I}_{[-b/2,b/2]}(\lambda_2) \times$$

$$\times \left[1 - \left(\frac{2\pi}{c}\right)^2 f_{\Delta}(\lambda_1) f_{\Delta}(\lambda_2)\right] \Phi_T(\lambda_2 - \lambda_1) d\lambda_1 d\lambda_2 |;$$

$$B_2(T,\Delta) = |\iint_{\mathbb{R}^2 \setminus D_b} e^{i(\tau_1-\tau_2)\lambda_2} |H^*(\lambda_2)|^2 \mathbb{I}_{[-b/2,b/2]}(\lambda_2) \times$$

$$\times \left[1 - \left(\frac{2\pi}{c}\right)^2 f_\Delta(\lambda_1) f_\Delta(\lambda_2)\right] \Phi_T(\lambda_2 - \lambda_1) d\lambda_1 d\lambda_2 |,$$

та $D_b = [-b,b] \times [-b,b]$. Оскільки $\|\Phi_T\|_1 = 1$, то перший з щойно наведених членів можна оцінити так

$$B_1(T,\Delta) \leq$$

$$\leq \iint_{D_b} |H^*(\lambda_2)|^2 \mathbb{I}_{[-b/2,b/2]}(\lambda_2) \left|1 - \left(\frac{2\pi}{c}\right)^2 f_\Delta(\lambda_1) f_\Delta(\lambda_2)\right| \Phi_T(\lambda_2 - \lambda_1) d\lambda_1 d\lambda_2 \leq$$

$$\leq \sup_{(\lambda_1,\lambda_2) \in D_b} \left|1 - \left(\frac{2\pi}{c}\right)^2 f_\Delta(\lambda_1) f_\Delta(\lambda_2)\right| \|H^*\|_2^2 \|\Phi_T\|_1 =$$

$$= \sup_{(\lambda_1,\lambda_2) \in D_b} \left|1 - \left(\frac{2\pi}{c}\right)^2 f_\Delta(\lambda_1) f_\Delta(\lambda_2)\right| \|H^*\|_2^2.$$

В силу умови (2.1б), можна записати

$$\sup_{(\lambda_1,\lambda_2) \in D_b} \left|1 - \left(\frac{2\pi}{c}\right)^2 f_\Delta(\lambda_1) f_\Delta(\lambda_2)\right| \leq$$

$$\leq \sup_{(\lambda_1,\lambda_2) \in D_b} \left[\left|1 - \left(\frac{2\pi}{c}\right) f_\Delta(\lambda_1)\right| + \left(\frac{2\pi}{c}\right) |f_\Delta(\lambda_1)| \left|1 - \left(\frac{2\pi}{c}\right) f_\Delta(\lambda_2)\right|\right] \leq$$

$$\leq \left[1 + \frac{2\pi}{c} \sup_{\Delta > 0} \|f_\Delta\|_\infty\right] \sup_{-b \leq \lambda \leq b} \left|1 - \frac{2\pi}{c} f_\Delta(\lambda)\right|, \qquad (2.15)$$

звідки, використовуючи умову (2.1г), в результаті одержуємо

$$\lim_{(T,\Delta) \to \infty} B_1(T,\Delta) = 0.$$

Для даного $b > 0$, розглянемо множину $\Pi(b/2) = \{(\lambda_1, \lambda_2) \in \mathbb{R}^2 : |\lambda_2 - \lambda_1| \leq b/2\}$. Оскільки

$$B_2(T,\Delta) \leq \iint_{\mathbb{R}^2 \setminus \Pi(b/2)} |H^*(\lambda_2)|^2 \, \mathbb{I}_{[-b/2,b/2]}(\lambda_2) \left|1 - \left(\frac{2\pi}{c}\right)^2 f_\Delta(\lambda_1) f_\Delta(\lambda_2)\right| \times$$

$$\times \Phi_T(\lambda_2 - \lambda_1) d\lambda_1 d\lambda_2 \leq$$

$$\leq \left[1 + \left(\frac{2\pi}{c} \sup_{\Delta>0} \|f_\Delta\|_\infty\right)^2\right] \int_{-b/2}^{b/2} |H^*(\lambda_2)|^2 d\lambda_2 \times \int_{|\lambda|>b/2} \Phi_T(\lambda) d\lambda \leq$$

$$\leq \left[1 + \left(\frac{2\pi}{c} \sup_{\Delta>0} \|f_\Delta\|_\infty\right)^2\right] \|H^*\|_2^2 \int_{|\lambda|>b/2} \Phi_T(\lambda) d\lambda,$$

й для будь-якого $b>0$, $\lim_{T \to \infty} \int_{|\lambda|>b/2} \Phi_T(\lambda) d\lambda = 0,$ то можна записати

$$\lim_{(T,\Delta) \to \infty} B_2(T,\Delta) = 0.$$

В силу того, що $|d_2| \leq B_1(T,\Delta) + B_2(T,\Delta)$, остаточно дістанемо

$$\lim_{(T,\Delta) \to \infty} d_2(T,\Delta) = 0. \tag{2.16}$$

Розглянемо величину $d_3$. Для будь-яких $b>0$, та $T>0, \Delta>0$, враховуючи, що $\|\Phi_T\|_1 = 1$, можна записати оцінку

$$|d_3(b,T,\Delta)| \leq$$

$$\leq \left(\frac{2\pi}{c}\right)^2 \int_{-\infty}^{\infty} \int_{-\infty}^{\infty} |H^*(\lambda_2)|^2 |\mathbb{I}_{[-b/2,b/2]}(\lambda_2) - 1| |f_\Delta(\lambda_1)| |f_\Delta(\lambda_2)| \Phi_T(\lambda_2 - \lambda_1) d\lambda_1 d\lambda_2 \leq$$

$$\leq \left(\frac{2\pi}{c} \sup_{\Delta>0} \|f_\Delta\|_\infty\right)^2 \int_{-\infty}^{\infty} \int_{-\infty}^{\infty} |H^*(\lambda_2)|^2 |\mathbb{I}_{[-b/2,b/2]}(\lambda_2) - 1| \Phi_T(\lambda_2 - \lambda_1) d\lambda_1 d\lambda_2 \leq$$

$$\leq \left(\frac{2\pi}{c} \sup_{\Delta>0} \|f_\Delta\|_\infty\right)^2 \|\Phi_T\|_1 \int_{|\lambda|>b/2} |H^*(\lambda)|^2 d\lambda =$$

$$= \left(\frac{2\pi}{c} \sup_{\Delta>0} \|f_\Delta\|_\infty\right)^2 \int_{|\lambda|>b/2} |H^*(\lambda)|^2 d\lambda,$$

яка є рівномірною по параметрах $T>0, \Delta>0$. Оскільки $|H^*|^2 \in L_1(\mathbb{R})$, звідси

випливає, що

$$\limsup_{b\to\infty, T,\Delta>0} |d_3(b,T,\Delta)| = 0. \tag{2.17}$$

З формул (2.14), (2.16) та (2.17) остаточно дістанемо

$$\limsup_{(T,\Delta)\to\infty} |C_\infty^{(1)}(\tau_1,\tau_2) - C_{T,\Delta}^{(1)}(\tau_1,\tau_2)| \le \frac{1}{2\pi}[\limsup_{b\to\infty}|d_1(b)| +$$

$$+\limsup_{b\to\infty}(\limsup_{(T,\Delta)\to\infty}|d_2(b,T,\Delta)|) + \limsup_{b\to\infty}(\limsup_{(T,\Delta)\to\infty}|d_3(b,T,\Delta)|)] = 0.$$

Таким чином, співвідношення (2.12) доведене.

• *Етап 2.* Розглянемо простір комплекснозначних неперервних фінітних функцій $C_f(\mathbb{R})$, визначених на $\mathbb{R}$. Це означає, що якщо $h \in C_f(\mathbb{R})$, тоді $h$ неперервна на $\mathbb{R}$ та існує таке додатне число $a_0(h)$, що $h(\lambda) = 0$ для $|\lambda| > a_0(h)$. Зазначимо, що будь-яка $h \in C_f(\mathbb{R})$ є рівномірно неперервною функцією.

В інтегральному зображенні кореляційної функції процесу $Z_{T,\Delta}$ (див. (2.8)) фігурують функції $(H^*(\lambda)e^{i\tau_k\lambda}, \lambda \in \mathbb{R}) \in L_2(\mathbb{R}), k \in \{1,2\}$. Оскільки простір $C_f(\mathbb{R})$ - всюди щільний у $L_2(\mathbb{R})$, то для будь-якого $\varepsilon > 0$ та всіх $\lambda \in \mathbb{R}$, можна вибрати такі функції $h_k^\varepsilon \in C_f(\mathbb{R}), k \in \{1,2\}$, що $\|H^*(\lambda)e^{i\tau_k\lambda} - h_k^\varepsilon(\lambda)\|_2 < \varepsilon$.

Для будь-яких $T > 0, \Delta > 0$, при всіх $\tau_1, \tau_2 \in \mathbb{R}$, та $\varepsilon > 0$, розглянемо різницю

$$C_\infty^{(2)}(\tau_1,\tau_2) - C_{T,\Delta}^{(2)}(\tau_1,\tau_2) = \frac{1}{2\pi}[d_1(\varepsilon) + d_2(\varepsilon,T,\Delta) + d_3(\varepsilon,T,\Delta)],$$

де

$$d_1(\varepsilon) = \int_{-\infty}^{\infty} [e^{i(\tau_1+\tau_2)\lambda}(H^*(\lambda))^2 - h_1^\varepsilon(\lambda)h_2^\varepsilon(\lambda)]d\lambda;$$

$$d_2(\varepsilon,T,\Delta) = \int_{-\infty}^{\infty} h_1^\varepsilon(\lambda) h_2^\varepsilon(\lambda) d\lambda -$$

$$- \int_{-\infty}^{\infty}\int_{-\infty}^{\infty} h_1^\varepsilon(\lambda_1) h_2^\varepsilon(\lambda_2) \left[\left(\frac{2\pi}{c}\right)^2 f_\Delta(\lambda_1) f_\Delta(\lambda_2)\right] \Phi_T(\lambda_2 - \lambda_1) d\lambda_1 d\lambda_2;$$

$$d_3(\varepsilon,T,\Delta) = \int_{-\infty}^{\infty}\int_{-\infty}^{\infty} [h_1^\varepsilon(\lambda_1) h_2^\varepsilon(\lambda_2) - e^{i(\tau_1\lambda_1+\tau_2\lambda_2)} H^*(\lambda_1) H^*(\lambda_2)] \times$$

$$\times \left[\left(\frac{2\pi}{c}\right)^2 f_\Delta(\lambda_1) f_\Delta(\lambda_2)\right] \Phi_T(\lambda_2 - \lambda_1) d\lambda_1 d\lambda_2.$$

Тепер розглянемо поведінку кожної з величин $d_k, k \in \{1,2,3\}$. З нерівності Коші-Буняковського випливає, що

$$|d_1(\varepsilon)| \leq \|e^{i\tau_1\lambda} H^*(\lambda)\|_2 \|e^{i\tau_2\lambda} H^*(\lambda) - h_2^\varepsilon(\lambda)\|_2 +$$

$$+ \|e^{i\tau_2\lambda} H^*(\lambda)\|_2 \|e^{i\tau_1\lambda} H^*(\lambda) - h_1^\varepsilon(\lambda)\|_2 +$$

$$+ \|e^{i\tau_1\lambda} H^*(\lambda) - h_1^\varepsilon(\lambda)\|_2 \|e^{i\tau_2\lambda} H^*(\lambda) - h_2^\varepsilon(\lambda)\|_2 < \varepsilon[2\|H^*\|_2 + \varepsilon],$$

звідки дістаємо

$$\lim_{\varepsilon \to 0} d_1(\varepsilon) = 0. \qquad (2.18)$$

З нерівності Юнга для згорток (див. теорему А.2), застосованої до доданка $d_3$, отримаємо оцінку

$$|d_3(\varepsilon,T,\Delta)| \leq \left(\frac{2\pi}{c}\right)^2 [\sup_{(\lambda_1,\lambda_2)\in\mathbb{R}^2} |f_\Delta(\lambda_1) f_\Delta(\lambda_1)|] \times$$

$$\times \int_{-\infty}^{\infty}\int_{-\infty}^{\infty} |h_1^\varepsilon(\lambda_1) h_2^\varepsilon(\lambda_2) - e^{i(\tau_1\lambda_1+\tau_2\lambda_2)} H^*(\lambda_1) H^*(\lambda_2)| \Phi_T(\lambda_2-\lambda_1) d\lambda_1 d\lambda_2$$

$$\leq \left(\frac{2\pi}{c} \sup_{\Delta>0} \|f_\Delta\|_\infty\right)^2 \|\Phi_T\|_1 [\|e^{i\tau_1\lambda_1} H^*(\lambda_1)\|_2 \|h_2^\varepsilon(\lambda_2) - e^{i\tau_2\lambda_2} H^*(\lambda_2)\|_2 +$$

$$+ \|e^{i\tau_2\lambda_2} H^*(\lambda_2)\|_2 \|h_1^\varepsilon(\lambda_1) - e^{i\tau_1\lambda_1} H^*(\lambda_1)\|_2 +$$

$$+ \| h_1^\varepsilon(\lambda_1) - e^{i\tau_1\lambda_1} H^*(\lambda_1) \|_2 \| h_2^\varepsilon(\lambda_2) - e^{i\tau_2\lambda_2} H^*(\lambda_2) \|_2 ] <$$

$$< \left( \frac{2\pi}{c} \sup_{\Delta>0} \| f_\Delta \|_\infty \right)^2 \varepsilon [ 2 \| H^* \|_2 + \varepsilon ],$$

яка є рівномірною по параметрах $T > 0, \Delta > 0$.

Таким чином, дістанемо рівність

$$\lim_{\varepsilon \to 0} \sup_{T, \Delta > 0} | d_3(\varepsilon, T, \Delta) | = 0. \qquad (2.19)$$

Тепер розглянемо величину $d_2$, й зауважимо, що при будь-яких $\varepsilon > 0$, та $T > 0, \Delta > 0$, має місце нерівність

$$| d_2(\varepsilon, T, \Delta) | \le E_1(\varepsilon, T) + E_2(\varepsilon, T, \Delta),$$

де

$$E_1(\varepsilon, T) = \left| \int_{-\infty}^{\infty} h_1^\varepsilon(\lambda) h_2^\varepsilon(\lambda) d\lambda - \int_{-\infty}^{\infty} \int_{-\infty}^{\infty} h_1^\varepsilon(\lambda_1) h_2^\varepsilon(\lambda_2) \Phi_T(\lambda_2 - \lambda_1) d\lambda_1 d\lambda_2 \right|;$$

$$E_2(\varepsilon, T, \Delta) =$$

$$= \left| \int_{-\infty}^{\infty} \int_{-\infty}^{\infty} h_1^\varepsilon(\lambda_1) h_2^\varepsilon(\lambda_2) \left[ 1 - \left( \frac{2\pi}{c} \right)^2 f_\Delta(\lambda_1) f_\Delta(\lambda_2) \right] \Phi_T(\lambda_2 - \lambda_1) d\lambda_1 d\lambda_2 \right|.$$

Для фіксованого $\varepsilon > 0$ розглянемо смугу $\Pi(\varepsilon) = \{ (\lambda_1, \lambda_2) \in \mathbb{R}^2 : |\lambda_2 - \lambda_1| \le \varepsilon \}$. Оскільки $\| \Phi_T \|_1 = 1$, то для першого доданка отримаємо

$$E_1(\varepsilon, T) \le \int\!\!\int_{\Pi(\varepsilon)} | h_1^\varepsilon(\lambda_1) | | h_2^\varepsilon(\lambda_1) - h_2^\varepsilon(\lambda_2) | \Phi_T(\lambda_2 - \lambda_1) d\lambda_1 d\lambda_2 +$$

$$+ \int\!\!\int_{\mathbb{R}^2 \setminus \Pi(\varepsilon)} | h_1^\varepsilon(\lambda_1) | | h_2^\varepsilon(\lambda_1) - h_2^\varepsilon(\lambda_2) | \Phi_T(\lambda_2 - \lambda_1) d\lambda_1 d\lambda_2 \le$$

$$\le \max_{(\lambda_1, \lambda_1) \in \Pi(\varepsilon)} | h_2^\varepsilon(\lambda_1) - h_2^\varepsilon(\lambda_2) | \| h_1^\varepsilon \|_1 \| \Phi_T \|_1 + 2 \| h_2^\varepsilon \|_\infty \| h_1^\varepsilon \|_1 \int_{|\lambda|>\varepsilon} \Phi_T(\lambda) d\lambda \le$$

$$\le \| h_1^\varepsilon \|_1 \left[ \max_{(\lambda_1, \lambda_1) \in \Pi(\varepsilon)} | h_2^\varepsilon(\lambda_1) - h_2^\varepsilon(\lambda_2) | + 2 \| h_2^\varepsilon \|_\infty \int_{|\lambda|>\varepsilon} \Phi_T(\lambda) d\lambda \right].$$

В силу того, що $h_2^\varepsilon$ - рівномірно неперервна функція, й для будь-якого $\varepsilon > 0$, $\lim\limits_{T\to\infty} \int\limits_{|\lambda|>\varepsilon} \Phi_T(\lambda)d\lambda = 0$, далі дістанемо

$$\limsup_{(T,\Delta)\to\infty} E_1(\varepsilon,T) = \limsup_{\varepsilon\to 0}(\limsup_{(T,\Delta)\to\infty} E_1(\varepsilon,T)) = 0.$$

Тепер розглянемо другий доданок у виразі оцінки для $d_2$. Покладемо $b = \max\{a_0(h_1^\varepsilon), a_0(h_2^\varepsilon)\}$. Так як $\{h_1^\varepsilon, h_2^\varepsilon\} \subset C_f(\mathbb{R})$ та $\||h_2^\varepsilon|*\Phi_T\|_2 \leq \|h_2^\varepsilon\|_2$ (див. теорему А.2), то для всіх $\varepsilon > 0$, та $T > 0, \Delta > 0$, маємо оцінку

$$E_2(\varepsilon,T,\Delta) \leq$$

$$\leq \iint\limits_{D_b} |h_1^\varepsilon(\lambda_1)||h_2^\varepsilon(\lambda_2)|\left|1-\left(\frac{2\pi}{c}\right)^2 f_\Delta(\lambda_1)f_\Delta(\lambda_2)\right|\Phi_T(\lambda_2-\lambda_1)d\lambda_1 d\lambda_2 +$$

$$+ \int\int\limits_{\mathbb{R}^2\setminus\Pi(b/2)} |h_1^\varepsilon(\lambda_1)||h_2^\varepsilon(\lambda_2)|\left|1-\left(\frac{2\pi}{c}\right)^2 f_\Delta(\lambda_1)f_\Delta(\lambda_2)\right|\Phi_T(\lambda_2-\lambda_1)d\lambda_1 d\lambda_2 \leq$$

$$\leq \sup_{(\lambda_1,\lambda_2)\in D_b}\left|1-\left(\frac{2\pi}{c}\right)^2 f_\Delta(\lambda_1)f_\Delta(\lambda_2)\right|\|h_1^\varepsilon\|_2\||h_2^\varepsilon|*\Phi_T\|_2 +$$

$$+ \left[1+\left(\frac{2\pi}{c}\sup_{\Delta>0}\|f_\Delta\|_\infty\right)^2\right]\|h_2^\varepsilon\|_\infty\|h_1^\varepsilon\|_1 \int\limits_{|\lambda|>b/2}\Phi_T(\lambda)d\lambda \leq$$

$$\leq \sup_{(\lambda_1,\lambda_2)\in D_b}\left|1-\left(\frac{2\pi}{c}\right)^2 f_\Delta(\lambda_1)f_\Delta(\lambda_2)\right|\|h_1^\varepsilon\|_2\|h_2^\varepsilon\|_2 +$$

$$+ \left[1+\left(\frac{2\pi}{c}\sup_{\Delta>0}\|f_\Delta\|_\infty\right)^2\right]\|h_2^\varepsilon\|_\infty\|h_1^\varepsilon\|_1 \int\limits_{|\lambda|>b/2}\Phi_T(\lambda)d\lambda.$$

З формули (2.15) та властивостей ядер Фейєра (див. зауваження 2.3) випливає співвідношення

$$\limsup_{(T,\Delta)\to\infty} E_2(\varepsilon,T,\Delta) = 0,$$

справедливе для будь-якого $\varepsilon > 0$.

В силу того, що $|d_2| \leq E_1(\varepsilon, T) + E_2(\varepsilon, T, \Delta)$ при будь-якому $\varepsilon > 0$, в результаті отримаємо

$$\limsup_{(T,\Delta) \to \infty} d_2(\varepsilon, T, \Delta) = 0. \qquad (2.20)$$

З формул (2.18) - (2.20), при всіх $\tau_1, \tau_2 \in \mathbb{R}$, остаточно дістанемо

$$\limsup_{(T,\Delta) \to \infty} |C_\infty^{(2)}(\tau_1, \tau_2) - C_{T,\Delta}^{(2)}(\tau_1, \tau_2)| \leq \frac{1}{2\pi}[\limsup_{\varepsilon \to 0}|d_1(\varepsilon)|$$

$$+ \limsup_{\varepsilon \to 0}(\limsup_{(T,\Delta) \to \infty} |d_2(\varepsilon, T, \Delta)|) + \limsup_{\varepsilon \to 0}(\limsup_{(T,\Delta) \to \infty} |d_3(\varepsilon, T, \Delta)|)] = 0.$$

Таким чином, співвідношення (2.13) доведене.

Об'єднуючи перший та другий етапи, отримаємо при всіх $\tau_1, \tau_2 \in \mathbb{R}$

$$\limsup_{(T,\Delta) \to \infty} C_{T,\Delta}(\tau_1, \tau_2) = \sum_{j=1}^{2} \limsup_{(T,\Delta) \to \infty} C_{T,\Delta}^{(j)}(\tau_1, \tau_2) = \sum_{j=1}^{2} C_\infty^{(j)}(\tau_1, \tau_2) = C_\infty(\tau_1, \tau_2).$$

Таким чином, теорему 2.1 доведено. $\square$

Порівняємо отриманий результат з відповідним твердженням у статті В. Булдигіна та Фу Лі [58].

*Зауваження 2.5.* У **лемі 4** [58] наведено ряд умов для встановлення граничної поведінки кореляційної функції процесу $Z_{T,\Delta}$, а саме, якщо:

1) $H \in L_2(\mathbb{R})$;

2) існує $p > 1$ таке, що $H^* \in L_p(\mathbb{R})$;

3) функція $H^*$ неперервна майже скрізь (відносно міри Лебега) на $\mathbb{R}$,

тоді для всіх $\tau_1, \tau_2 \in \mathbb{R}$: $\limsup\limits_{(T,\Delta) \to \infty} C_{T,\Delta}(\tau_1, \tau_2) = C_\infty(\tau_1, \tau_2).$

Зокрема, дана границя має місце, якщо $H \in L_1(\mathbb{R}) \bigcap L_2(\mathbb{R})$.

### 2.2.3 Асимптотична нормальність скінченновимірних розподілів процесу $Z_{T,\Delta}$

У цьому пункті досліджується асимптотична нормальність скінченновимірних розподілів процесу $Z_{T,\Delta}$ (див. (2.7)). Зазначимо, що ніяких додаткових умов на перехідну функцію $H$ та характер прямування параметрів $T,\Delta$ до безмежності не накладається.

Теорема 2.1 демонструє, що функція $C_\infty$, визначена формулою (2.11), є невід'ємно визначеною на $\mathbb{R}\times\mathbb{R}$. Таким чином, завжди можна стверджувати, що існує центрований дійснозначний гауссівський процес $Z = (Z(\tau), \tau \in \mathbb{R})$ з кореляційною функцією $C_\infty$, тобто

$$\mathrm{E}Z(\tau_1)Z(\tau_2) = C_\infty(\tau_1,\tau_2) = \qquad (2.21)$$

$$= \frac{1}{2\pi}\int_{-\infty}^{\infty}\left[e^{i(\tau_1-\tau_2)\lambda}|H^*(\lambda)|^2 + e^{i(\tau_1+\tau_2)\lambda}(H^*(\lambda))^2\right]d\lambda, \ \tau_1,\tau_2 \in \mathbb{R}.$$

Без втрати загальності, припускається, що процес $Z$ заданий на тому ж базовому ймовірнісному просторі $(\Omega, \mathfrak{F}, \mathrm{P})$, що й процес $Z_{T,\Delta}$.

Має місце наступне твердження про асимптотичну нормальність процесу $Z_{T,\Delta}$ при $(T,\Delta) \to \infty$.

**Теорема 2.2.** *Нехай $H \in L_2(\mathbb{R})$, тоді для всіх $m \in \mathbb{N}$ та $\tau_1,\tau_2,...,\tau_m \in \mathbb{R}$, має місце рівність*

$$\lim_{(T,\Delta)\to\infty} cum(Z_{T,\Delta}(\tau_j), j=1,...,m) = \begin{cases} 0, & m=1; \\ C_\infty(\tau_1,\tau_2), & m=2; \\ 0, & m\geq 3, \end{cases} \qquad (2.22)$$

*де $cum(Z_{T,\Delta}(\tau_j), j=1,...,m)$ - сумісний кумулянт набору випадкових величин $Z_{T,\Delta}(\tau_j), j=1,...,m.$*

*Зокрема, всі скінченновимірні розподіли процесу $(Z_{T,\Delta}(\tau), \tau \in \mathbb{R})$ слабко збігаються до відповідних скінченновимірних розподілів центрованого гауссівського процесу $(Z(\tau), \tau \in \mathbb{R})$.*

*Доведення.* Нехай $m \in \{1,2\}$. Формула (2.22) має місце, оскільки процес $Z_{T,\Delta}$ центрований, та задовольняються умови теореми 2.1.

Нехай $m \geq 3$. Доведення формули (2.22) проведемо у два етапи. На першому етапі покажемо, що $cum(Z_{T,\Delta}(\tau_j), j=1,...,m)$ можна подати у виді скінченної суми інтегралів з циклічним зачепленням ядер. Другий етап доведення присвячений перевірці умов теореми А.3, згідно з якою кожен доданок у сумі прямує до нуля при $(T, \Delta) \to \infty$. Для встановлення цих фактів використовується діаграмно-кумулянтний метод Д. Бриллінджера[7] та підходи, запропоновані В. Булдигіним, Ф. Уцетом та В. Зайцем [61], [62].

- *Етап 1*. З властивостей лінійності кумулянтів та їх інваріантності відносно зсувів, а також з означень процесу $Z_{T,\Delta}$ та оцінки $H_{T,\Delta}$ (див. (2.7) і (2.4) відповідно) випливає, що

$$cum(Z_{T,\Delta}(\tau_j), j=1,...,m) = \tag{2.23}$$

$$= \left(\frac{1}{c^2 T}\right)^{\frac{m}{2}} \int_0^T ... \int_0^T cum(Y_\Delta(t_j + \tau_j) X_\Delta(t_j), j=1,...,m) dt_1 ... dt_m.$$

Оскільки $X_\Delta$ та $Y_\Delta$ - центровані сумісно гауссівські процеси, то з теореми А.5, застосованої до кумулянта, що стоїть під знаком інтеграла у формулі (2.23), дістанемо

$$cum(Y_\Delta(t_j + \tau_j) X_\Delta(t_j), j=1,...,m) = \sum \prod_{p=1}^m cum(D_p^{(2)}), \tag{2.24}$$

де сума береться по всіх можливих невпорядкованих нерозкладних розбиттях двостовпцевої таблиці

$$D_{m \times 2} = \begin{bmatrix} Y_\Delta(t_1 + \tau_1) & X_\Delta(t_1) \\ Y_\Delta(t_2 + \tau_2) & X_\Delta(t_2) \\ \vdots & \vdots \\ Y_\Delta(t_m + \tau_m) & X_\Delta(t_m) \end{bmatrix}$$

на пари $\{D_1^{(2)},...,D_m^{(2)}\}$.

Так як порядок елементів у $\{D_1^{(2)},...,D_m^{(2)}\}$ не має значення, можна завжди припускати, що розбиття задовольняє умовам:

(1)     $D_p^{(2)} \bigcap D_q^{(2)} = \varnothing$ для $p \neq q$;

(2)     $D_{m \times 2} = D_1^{(2)} \bigcup ... \bigcup D_m^{(2)}$;

(3)     якщо $m \geq 3$, тоді для будь-якого $p = 1,...,m$ множина $D_p^{(2)}$ не співпадає з жодним із рядків $r_1,...,r_m$ таблиці 2; більш того, для будь-якого $1 \leq \nu < m$, об'єднання довільної кількості елементів розбиття $\{D_1^{(2)},...,D_m^{(2)}\}$ не повинно співпадати з об'єднанням будь-яких $\nu$ рядків таблиці $D_{m \times 2}$;

(4)     для будь-якого $p = 1,...,m-1$, існує єдиний рядок $\tilde{r}_p$ таблиці $D_{m \times 2}$ такий, що $\tilde{r}_p \subset D_p^{(2)} \bigcup D_{p+1}^{(2)}$;

(5)     $r_1 \subset D_1^{(2)} \bigcup D_m^{(2)}$.

З умов (1)-(5) випливає, що елементи розбиття $\{D_1^{(2)},...,D_m^{(2)}\}$ послідовно зачіплюються один за інший та "вичерпують" всю таблицю $D_{m \times 2}$; $D_m^{(2)}$ зачіплюється за $D_1^{(2)}$.

Далі, для заданого $D_p^{(2)}$, запишемо $D_p^{(2)} = D_{j,j}^{(2)}$, де $j$ та $j$ - номери тих двох рядків таблиці $D_{m \times 2}$, в яких містяться елементи множини $D_p^{(2)}$. Отже, довільне невпорядковане нерозкладне розбиття $\{D_1^{(2)},...,D_m^{(2)}\}$ можна записати

у виді

$$\vec{D}^{(2)} = \{D^{(2)}_{j_1,j_2}, D^{(2)}_{j_2,j_3},...,D^{(2)}_{j_{m-1},j_m}, D^{(2)}_{j_m,j_{m+1}}\},$$

де $(j_1,...,j_m)$ - впорядкований набір такий, що $j_1 = 1$, $(j_2, j_3,..., j_m)$ - перестановка чисел множини $\{2,3,...,m\}$, та $j_{m+1} = j_1$.

Аналіз структури $\vec{D}^{(2)}$ показує, що можна виділити три різні класи його елементів. До першого класу належать невпорядковані множини виду

$$D^{(2)}_{j,j} = \{X_\Delta(t_j), X_\Delta(t_j)\}, \ j \neq j.$$

Другий клас утворюють невпорядковані множини виду

$$D^{(2)}_{j,j} = \{Y_\Delta(t_j + \tau_j), Y_\Delta(t_j + \tau_j)\}, \ j \neq j.$$

Третій клас містить невпорядковані множини виду

$$D^{(2)}_{j,j} = \{Y_\Delta(t_j + \tau_j), X_\Delta(t_j)\}, \ j \neq j.$$

Відповідно позначимо класи $G_1(\vec{D}^{(2)}), G_2(\vec{D}^{(2)})$ та $G_3(\vec{D}^{(2)})$. Нехай $m_\nu(\vec{D}^{(2)}) = card G_\nu(\vec{D}^{(2)}), \nu = 1,2,3$, тоді для будь-якого $\vec{D}^{(2)}$

$$m_1(\vec{D}^{(2)}) = m_2(\vec{D}^{(2)}); m_1(\vec{D}^{(2)}) + m_2(\vec{D}^{(2)}) + m_3(\vec{D}^{(2)}) = m. \qquad (2.25)$$

Тепер знайдемо кумулянти від наборів величин, що належать кожному з трьох класів. Якщо $D^{(2)}_{j,j} \in G_1(\vec{D}^{(2)})$, то

$$cum(D^{(2)}_{j,j}) = K_\Delta(t_j - t_j) = \int_{-\infty}^{\infty} e^{i(t_j - t_j)\lambda_j} f_\Delta(\lambda_j) d\lambda_j;$$

якщо $D^{(2)}_{j,j} \in G_2(\vec{D}^{(2)})$, то

$$cum(D^{(2)}_{j,j}) = EY_\Delta(t_j + \tau_j) Y_\Delta(t_j + \tau_j) =$$
$$= \int_{-\infty}^{\infty} e^{i(t_j - t_j)\lambda_j} \cdot e^{i(\tau_j - \tau_j)\lambda_j} |H^*(\lambda_j)|^2 f_\Delta(\lambda_j) d\lambda_j;$$

якщо $D^{(2)}_{j,j} \in G_3(\vec{D}^{(2)})$, то

$$cum(D^{(2)}_{j,j}) = EY_\Delta(t_j + \tau_j)X_\Delta(t_j) = \int_{-\infty}^{\infty} e^{i(t_j-t_j)\lambda_j} \cdot e^{i\tau_j\lambda_j} H^*(\lambda_j) f_\Delta(\lambda_j) d\lambda_j.$$

З останніх формул випливає представлення для добутку кумулянтів з $\vec{D}^{(2)}$:

$$\prod_{p=1}^{m} cum(D^{(2)}_p) = \int \ldots \int_{\mathbb{R}^m} \left[ \prod_{k=1}^{m} e^{i(t_{j_{k+1}}-t_{j_k})\lambda_{j_k}} \right] \varphi_0(\vec{\lambda}, \vec{\tau}, \vec{D}^{(2)}) \times \qquad (2.26)$$

$$\times \left( \prod_{k=1}^{m} \varphi_{j_k}(\lambda_{j_{k+1}}, \Delta, \vec{D}^{(2)}) \right) d\lambda_{j_1} \ldots d\lambda_{j_m},$$

де $\vec{\lambda} = (\lambda_1, \ldots, \lambda_m)$; $\vec{\tau} = (\tau_1, \ldots, \tau_m)$; $(j_1, \ldots, j_m)$ - набір такий, що $j_{m+1} = j_1 = 1$, $(j_2, \ldots, j_m)$ - перестановка чисел множини $\{2, \ldots, m\}$.

Для кожного розбиття $\vec{D}^{(2)}$ функція $\varphi_0(\cdot, \vec{D}^{(2)})$ є добутком деяких з функцій $e^{i(\tau_{j_{k+1}}-\tau_{j_k})\lambda_{j_k}}, e^{i\tau_{j_k}\lambda_{j_k}}$ та індикаторних функцій $\mathbb{I}_\mathbb{R}(\lambda_{j_k}), k=1,\ldots,m$. Таким чином, має місце рівність

$$\sup_{\vec{D}^{(2)}} \sup_{\vec{\lambda},\vec{\tau}} |\varphi_0(\vec{\lambda}, \vec{\tau}, \vec{D}^{(2)})| = 1. \qquad (2.27)$$

Для будь-яких $\Delta > 0$, та розбиття $\vec{D}^{(2)}$ кожна функція з множини $F(\Delta, \vec{D}^{(2)}) = \{\varphi_1(\cdot, \Delta, \vec{D}^{(2)}), \ldots, \varphi_m(\cdot, \Delta, \vec{D}^{(2)})\}$ може збігатися лише з однією з наступних $f_\Delta, |H^*|^2 f_\Delta$ та $H^* f_\Delta$. Крім того,

$$card\{\varphi \in F(\Delta, \vec{D}^{(2)}) : \varphi = f_\Delta\} = m_1(\vec{D}^{(2)});$$

$$card\{\varphi \in F(\Delta, \vec{D}^{(2)}) : \varphi = |H^*|^2 f_\Delta\} = m_2(\vec{D}^{(2)}); \qquad (2.28)$$

$$card\{\varphi \in F(\Delta, \vec{D}^{(2)}) : \varphi = |H^*| f_\Delta\} = m_3(\vec{D}^{(2)}).$$

Перетворимо вираз $\prod_{k=1}^{m} e^{i(t_{j_k}-t_{j_{k+1}})\lambda_{j_k}} = \prod_{k=1}^{m} e^{it_{j_{k+1}}(\lambda_{j_{k+1}}-\lambda_{j_k})}$, де $\lambda_{m+1} = \lambda_1$, та

розглянемо $m$-кратний інтеграл

$$\int_0^T \cdots \int_0^T \left(\prod_{k=1}^m e^{it_{j_{k+1}}(\lambda_{j_{k+1}} - \lambda_{j_k})}\right) dt_{j_1} \ldots dt_{j_m} = \prod_{k=1}^m \left(\int_0^T e^{it_{j_{k+1}}(\lambda_{j_{k+1}} - \lambda_{j_k})} dt_{j_{k+1}}\right) =$$

$$= \prod_{k=1}^m \Phi_T(\lambda_{j_{k+1}} - \lambda_{j_k}),$$

де $\Phi_T(\lambda) = \int_0^T e^{it\lambda} dt = \dfrac{e^{it\lambda T} - 1}{i\lambda}, \lambda \in \mathbb{R}.$

Враховуючи останні перетворення, з формул (2.23) - (2.26) отримаємо

$$cum(Z_{T,\Delta}(\tau_j), j = 1, \ldots, m) = \qquad (2.29)$$

$$= \left(\frac{2\pi}{c^2}\right)^{\frac{m}{2}} \sum_{\vec{D}^{(2)}} \int \cdots \int_{\mathbb{R}^m} \left[\prod_{k=1}^m \Phi^{(T)}(\lambda_{k+1} - \lambda_k)\right] \varphi_0(\vec{\lambda}, \vec{\tau}, \vec{D}^{(2)}) \times$$

$$\times \left(\prod_{k=1}^m \varphi_k(\lambda_k, \Delta, \vec{D}^{(2)})\right) d\lambda_1 \ldots d\lambda_m,$$

де $\lambda_{m+1} = \lambda_1$, та

$$\Phi^{(T)}(\lambda) = \left(\frac{1}{2\pi T}\right)^{\frac{1}{2}} \Phi_T(\lambda) = \left(\frac{1}{2\pi T}\right)^{\frac{1}{2}} \frac{e^{iT\lambda} - 1}{i\lambda}, \lambda \in \mathbb{R}. \qquad (2.30)$$

Звідси, виключаючи $\varphi_0$ та залежність $\varphi_k, k = 1, \ldots, m,$ від $\Delta$ (див. (2.27), (2.1б)), отримаємо наступну оцінку для кумулянта:

$$|cum(Z_{T,\Delta}(\tau_j), j = 1, \ldots, m)| = \left(\frac{2\pi}{c^2}\right)^{\frac{m}{2}} (\sup_{\Delta > 0} \|f_\Delta\|_\infty)^m \times \qquad (2.31)$$

$$\times \sum_{\vec{D}^{(2)}} \int \cdots \int_{\mathbb{R}^m} |\prod_{k=1}^m \Phi^{(T)}(\lambda_{k+1} - \lambda_k) \overline{\varphi}_k(\lambda_k, \vec{D}^{(2)})| d\lambda_1 \ldots d\lambda_m,$$

де $\lambda_{m+1} = \lambda_1$, та для будь-якого $\vec{D}^{(2)}$

$$\overline{\varphi}_k = \begin{cases} \mathbb{I}_{\mathbb{R}}, & \text{якщо } \varphi_k(\cdot, \Delta, \overrightarrow{D}^{(2)}) = f_\Delta; \\ |H^*|^2, & \text{якщо } \varphi_k(\cdot, \Delta, \overrightarrow{D}^{(2)}) = |H^*|^2 f_\Delta; \\ |H^*|, & \text{якщо } \varphi_k(\cdot, \Delta, \overrightarrow{D}^{(2)}) = |H^*| f_\Delta. \end{cases} \quad (2.32)$$

Оцінкою (2.31) для кумулянта $cum(Z_{T,\Delta}(\tau_j), j=1,...,m)$ завершується перший етап доведення.

- *Етап 2.* Зафіксуємо $\overrightarrow{D}^{(2)}$. Для інтеграла з циклічним зачепленням ядер ($\lambda_{m+1} = \lambda_1$), що фігурує в (2.31), виду

$$I_m^{(T)}(\overrightarrow{D}^{(2)}) = \int...\int_{\mathbb{R}^m} |\prod_{k=1}^{m} \Phi^{(T)}(\lambda_{k+1} - \lambda_k) \overline{\varphi}_k(\lambda_k, \overrightarrow{D}^{(2)})| d\lambda_1...d\lambda_m,$$

далі будемо використовувати теорему А.3, переформульовану для нашого випадку наступним чином.

**Лема 2.5** *Нехай $m \in \mathbb{N} \setminus \{1, 2\}$. Якщо: (i) ядра $(\Phi^{(T)}, T > 0)$ задовольняють умові $\|\Phi^{(T)}\|_p \leq C_p [T]^{\frac{1}{2}-\frac{1}{p}}$ для всіх $p \in (1, \infty]$, де $C_p$ - деяка додатна стала, що не залежить від параметра $T$; (ii) серед функцій $\overline{\varphi}_k, k=1,...,m,$ існує $m_1 \geq 0$ функцій з простору $L_1(\mathbb{R})$, $m_\infty = m_1$ функцій з простору $L_\infty(\mathbb{R})$, та $m_2 = m - 2m_1$ функцій з простору $L_2(\mathbb{R})$. Тоді*

$$\lim_{T \to \infty} I_m^{(T)}(\overrightarrow{D}^{(2)}) = 0.$$

Перевіримо умови цієї леми в нашому випадку:

(i) Для ядер

$$\Phi^{(T)}(\lambda) = \left(\frac{1}{2\pi T}\right)^{\frac{1}{2}} \frac{e^{iT\lambda}-1}{i\lambda}, \lambda \in \mathbb{R},$$

при всіх $p \in (1, \infty]$, виконується умова:

$$\|\Phi^{(T)}\|_p = T^{\frac{1}{2}-\frac{1}{p}} \cdot \left[\frac{1}{\sqrt{2\pi}}\int_0^\infty \left|\frac{\sin\lambda}{\lambda}\right|^p d\lambda\right]^{\frac{1}{p}}.$$

(ii) Для вказаного розбиття $\vec{D}^{(2)}$ множина функцій $\vec{F}(\vec{D}^{(2)}) = \{\overline{\varphi}_1(\cdot,\vec{D}^{(2)}),...,\overline{\varphi}_m(\cdot,\vec{D}^{(2)})\}$ розбивається на класи:

$$M_\infty(\vec{D}^{(2)}) = \{\overline{\varphi} \in \vec{F}(\vec{D}^{(2)}) : \overline{\varphi} = \mathbb{I}_\mathbb{R}\};$$
$$M_1(\vec{D}^{(2)}) = \{\overline{\varphi} \in \vec{F}(\vec{D}^{(2)}) : \overline{\varphi} = |H^*|^2\}; \qquad (2.33)$$
$$M_2(\vec{D}^{(2)}) = \{\overline{\varphi} \in \vec{F}(\vec{D}^{(2)}) : \overline{\varphi} = |H^*|\};$$

З формул (2.25), (2.28), (2.32) та (2.33) випливає, що

$$card[M_1(\vec{D}^{(2)})] = card[M_\infty(\vec{D}^{(2)})];$$
$$card[M_1(\vec{D}^{(2)})] + card[M_2(\vec{D}^{(2)})] + card[M_\infty(\vec{D}^{(2)})] = m$$

Оскільки $H \in L_2(\mathbb{R})$, то мають місце включення

$$M_\infty(\vec{D}^{(2)}) \subset L_\infty(\mathbb{R}), \ M_1(\vec{D}^{(2)}) \subset L_1(\mathbb{R}), \ M_2(\vec{D}^{(2)}) \subset L_2(\mathbb{R}).$$

Отже, справджуються всі умови переформульованої теореми А.3, тому

$$\lim_{T\to\infty} I_m^{(T)}(\vec{D}^{(2)}) = 0. \qquad (2.34)$$

Оскільки в скінченній сумі з (2.31) виключена залежність від параметра $\Delta > 0$, й для кожного її доданка $I_m^{(T)}(\vec{D}^{(2)})$ виконується рівність (2.34), то

$$\lim_{(T,\Delta)\to\infty} cum(Z_{T,\Delta}(\tau_j), j=1,...,m) = 0.$$

Так як гауссівський розподіл однозначно визначається своїми кумулянтами, то з формули (2.22) в силу теореми А.9 випливає, що всі скінченновимірні розподіли процесу $(Z_{T,\Delta}(\tau), \tau \in \mathbb{R})$ слабко збігаються до відповідних скінченновимірних розподілів гауссівського процесу $(Z(\tau), \tau \in \mathbb{R})$ при $(T,\Delta) \to \infty$. Таким чином, теорему 2.2 доведено. $\square$

Твердження теореми 2.2 можна сформулювати за допомогою моментів.

**Теорема 2.3.** *Нехай $H \in L_2(\mathbb{R})$, тоді для всіх $m \in \mathbb{N}$ та $\tau_1, \tau_2, ..., \tau_m \in \mathbb{R}$, має місце рівність*

$$\lim_{(T,\Delta)\to\infty} \mathrm{E}[\prod_{j=1}^{m} Z_{T,\Delta}(\tau_j)] = \mathrm{E}[\prod_{j=1}^{m} Z(\tau_j)]. \qquad (2.35)$$

*Зокрема, всі скінченновимірні розподіли процесу $(Z_{T,\Delta}(\tau), \tau \in \mathbb{R})$ слабко збігаються до відповідних скінченновимірних розподілів центрованого гауссівського процесу $(Z(\tau), \tau \in \mathbb{R})$.*

*Доведення.* Для гауссівського центрованого процесу $Z$ з кореляційною функцією $C_\infty$ (див. (2.21)), при всіх $m \in \mathbb{N}$ та $\tau_1, \tau_2, ..., \tau_m \in \mathbb{R}$

$$\mathrm{cum}(Z(\tau_j), j=1,...,m) = \begin{cases} 0, & m \neq 2; \\ C_\infty(\tau_1, \tau_2), & m = 2. \end{cases}$$

З цього факту та теореми 2.2 випливає, що

$$\lim_{(T,\Delta)\to\infty} \mathrm{cum}(Z_{T,\Delta}(\tau_j), j=1,...,m) = \mathrm{cum}(Z(\tau_j), j=1,...,m).$$

Збіжність відповідних моментів (2.35) має місце з представлення моментів вектора $(Z(\tau_j), j=1,...,m)$ за допомогою кумулянтів:

$$\mathrm{E}[\prod_{k=1}^{m} Z(\tau_k)] = \sum \prod_{p=1}^{q} \mathrm{cum}(Z_{T,\Delta}(\tau_j) \mid j \in \nu_p),$$

де сума береться по всіх можливих невпорядкованих розбиттях $\{\nu_1, ..., \nu_q\}, q = 1, ..., m,$ множини чисел $\{1, 2, ..., m\}$, граничного співвідношення (2.22) та формули Іссерліса (див. теорему А.7).

Оскільки гауссівський розподіл однозначно визначається своїми моментами, то з формули (2.35) та теореми Маркова (див. теорему А.8)

випливає, що всі скінченновимірні розподіли процесу $(Z_{T,\Delta}(\tau), \tau \in \mathbb{R})$ слабко збігаються до відповідних скінченновимірних розподілів гауссівського процесу $(Z(\tau), \tau \in \mathbb{R})$ при $(T,\Delta) \to \infty$. Таким чином, теорему 2.3 доведено. $\square$

Порівняємо отриманий результат з відповідними твердженнями у статті В. Булдигіна та Фу Лі [59], де наведено ряд умов для встановлення асимптотичної нормальності скінченновимірних розподілів процесу $Z_{T,\Delta}$ залежно від обмеженості або необмеженості функції $H^*$.

*Зауваження 2.6.* У **теоремі 1** [59] стверджується, якщо:

1) $H \in L_2(\mathbb{R})$;

2) $H^* \in L_1(\mathbb{R}) \bigcap L_\infty(\mathbb{R})$;

3) функція $H^*$ неперервна майже скрізь на $\mathbb{R}$,

тоді для всіх $m \in \mathbb{N}$ та $\tau_1, \tau_2, ..., \tau_m \in \mathbb{R}$, має місце рівність

$$\lim_{(T,\Delta)\to\infty} \mathrm{E}[\prod_{j=1}^{m} Z_{T,\Delta}(\tau_j)] = \mathrm{E}[\prod_{j=1}^{m} Z(\tau_j)].$$

Зокрема, всі скінченновимірні розподіли процесу $(Z_{T,\Delta}(\tau), \tau \in \mathbb{R})$ слабко збігаються до відповідних скінченновимірних розподілів центрованого гауссівського процесу $(Z(\tau), \tau \in \mathbb{R})$.

Варто зазначити, що якщо $H \in L_1(\mathbb{R}) \bigcap L_2(\mathbb{R})$ і $H^* \in L_1(\mathbb{R})$, тоді твердження наведеної теореми 1 виконані.

У **теоремі 2** [59] стверджується, якщо:

1) $H \in L_2(\mathbb{R})$;

2) існує $p > 2$ таке, що $H^* \in L_{2p}(\mathbb{R})$;

3) функція $H^*$ неперервна майже скрізь на $\mathbb{R}$;

4) $T \to \infty$ і $\Delta \to \infty$ так, що для всіх $m \geq 3$

$$\frac{\|f_\Delta\|_2^{a(m)} \|f_\Delta\|_p^{b(m)}}{T^{(m-2)(p-2)/(2p)}} \to 0,$$

де $a(m) = \begin{cases} 1, & m \geq 4; \\ 0, & m = 3, \end{cases}$ $b(m) = \begin{cases} ent(\frac{m}{2}) - 1, & m \geq 4; \\ 1, & m = 3, \end{cases}$ та $ent(x)$ - ціла частина числа $x$.

Тоді для всіх $m \in \mathbb{N}$ та $\tau_1, \tau_2, ..., \tau_m \in \mathbb{R}$, має місце рівність

$$\lim_{(T,\Delta) \to \infty} \mathrm{E}[\prod_{j=1}^m Z_{T,\Delta}(\tau_j)] = \mathrm{E}[\prod_{j=1}^m Z(\tau_j)].$$

Зокрема, всі скінченновимірні розподіли процесу $(Z_{T,\Delta}(\tau), \tau \in \mathbb{R})$ слабко збігаються до відповідних скінченновимірних розподілів центрованого гауссівського процесу $(Z(\tau), \tau \in \mathbb{R})$.

### 2.2.4 Асимптотична нормальність розподілів процесу $Z_{T,\Delta}$ у просторі неперервних функцій

Після встановлення асимптотичної нормальності скінченновимірних розподілів процесу $Z_{T,\Delta}$ (теореми 2.2 і 2.3), природно поставити питання про асимптотичну нормальність цього процесу у просторі неперервних функцій.

Нехай $[a,b] \in \mathbb{R}$; $\mathrm{C}[a,b]$ - банахів простір неперервних дійснозначних функцій, визначених на сегменті $[a,b]$, з рівномірною нормою.

Нагадаємо деякі терміни, пов'язані з гауссівськими випадковими процесами (див. додаток А.7).

Нехай $S \in \mathbb{R}$, та $\rho(t,s), t,s \in S$, - деяка превдометрика. Нагадаємо, що псевдометрика задовольняє всім аксіомам метрики, за винятком того, що множина нулів $\{(t,s) \in S \times S : \rho(t,s) = 0\}$ може бути ширшою за діагональ $\{(t,s) \in S \times S : t = s\}$. Як завжди, $N_\rho(S,\varepsilon)$ - мінімальне число замкнених $\rho$-куль радіуса $\varepsilon > 0$, центри яких лежать всередині $S$, та які покривають $S$. Якщо не існує скінченного покриття для $S$, тоді $N_\rho(S,\varepsilon) = \infty$. Далі, $\mathcal{H}_\rho(S,\varepsilon) = \log N_\rho(S,\varepsilon)$ - метрична ентропія множини $S$ відносно $\rho$. Домовимось, що для довільного $\beta > 0$ нерівність

$$\int_{0+} \mathcal{H}_\rho^\beta(S,\varepsilon) d\varepsilon < \infty$$

означає наступне: для деякого (отже, й для всіх) $u > 0$ виконується нерівність

$$\int_0^u \mathcal{H}_\rho^\beta(S,\varepsilon) d\varepsilon < \infty.$$

Розглянемо функцію

$$q_H(\tau) = \left[ \int_{-\infty}^\infty \sin^2 \frac{\tau\lambda}{2} |H^*(\lambda)|^2 d\lambda \right]^{\frac{1}{2}}, \tau \in \mathbb{R}, \quad (2.36)$$

де $H^*$ - перетворення Фур'є - Планшереля функції $H$. Оскільки $H \in L_2(\mathbb{R})$, то $q_H$ коректно визначена та задає псевдометрики

$$\sigma(\tau_1, \tau_2) = q_H(|\tau_1 - \tau_2|); \quad \sqrt{\sigma}(\tau_1, \tau_2) = \sqrt{q_H(|\tau_1 - \tau_2|)}, \tau_1, \tau_2 \in \mathbb{R}.$$

Зауважимо, що якщо $H^*(\lambda) \neq 0$ на множині додатної лебегової міри, то $\sigma$ та $\sqrt{\sigma}$ є метриками.

Для будь-якого $\varepsilon > 0$ покладемо

$$\mathcal{H}_\sigma(\varepsilon) = \mathcal{H}_\sigma([0,1], \varepsilon); \quad \mathcal{H}_{\sqrt{\sigma}}(\varepsilon) = \mathcal{H}_{\sqrt{\sigma}}([0,1], \varepsilon).$$

Оскільки псевдометрики $\sigma$ та $\sqrt{\sigma}$ залежать лише від $|\tau_1 - \tau_2|$, то

$$\int_{0+} \mathcal{H}_\sigma^\beta(\varepsilon) d\varepsilon < \infty \Leftrightarrow \int_{0+} \mathcal{H}_\sigma^\beta([a,b],\varepsilon) d\varepsilon < \infty;$$

$$\int_{0+} \mathcal{H}_{\sqrt{\sigma}}(\varepsilon) d\varepsilon < \infty \Leftrightarrow \int_{0+} \mathcal{H}_{\sqrt{\sigma}}([a,b],\varepsilon) d\varepsilon < \infty$$

для будь-яких відрізка $[a,b] \in \mathbb{R}$, та числа $\beta > 0$.

Крім того, оскільки для всіх $\tau_1, \tau_2 \in \mathbb{R}$,

$$\sigma(\tau_1, \tau_2) \leq \left[ \max_{\tau_1, \tau_2 \in \mathbb{R}} \sigma(\tau_1, \tau_2) \right]^{\frac{1}{2}} \sqrt{\sigma(\tau_1, \tau_2)} \leq \| H^* \|_2^{\frac{1}{2}} \sqrt{\sigma(\tau_1, \tau_2)}, \qquad (2.37)$$

то зі співвідношення $\int\limits_{0+} \mathcal{H}_{\sqrt{\sigma}}(\varepsilon) < \infty,$ отримаємо $\int\limits_{0+} \mathcal{H}_{\sigma}(\varepsilon) < \infty,$ звідки, в свою чергу, випливає

$$\int\limits_{0+} \mathcal{H}_{\sigma}^{\frac{1}{2}}(\varepsilon) < \infty.$$

Далі запис

$$Z_{T,\Delta} \overset{C[a,b]}{\Rightarrow} Z$$

означає слабку збіжність процесу $Z_{T,\Delta}$ до процесу $Z$ у просторі неперервних функцій при $(T, \Delta) \to \infty$; тобто для будь-якого $C[a,b]$-неперервного дійснозначного функціоналу $F$ розподіли випадкової величини $F(Z_{T,\Delta})$ слабко збігаються до розподілу випадкової величини $F(Z)$. Будемо вважати, що

$$Z_{T,\Delta} = (Z_{T,\Delta}(\tau), \tau \in [a,b]);$$

$$Z = (Z(\tau), \tau \in [a,b]),$$

й ці процеси є сепарабельними. (Таке припущення не є обмежувальним в силу стохастичної неперервності вказаних процесів).

Для всіх $T > 0, \Delta > 0,$ визначимо сім'ю псевдометрик

$$\rho_{(T,\Delta)}(\tau_1, \tau_2) = \left( \mathrm{E} | Z_{T,\Delta}(\tau_2) - Z_{T,\Delta}(\tau_1) |^2 \right)^{\frac{1}{2}}, \ \tau_1, \tau_2 \in \mathbb{R}. \qquad (2.38)$$

**Лема 2.4** *Нехай $H \in L_2(\mathbb{R})$, тоді для всіх $T > 0, \Delta > 0,$ та $\tau_1, \tau_2 \in \mathbb{R}$, має місце нерівність*

$$\rho_{(T,\Delta)}(\tau_1, \tau_2) \leq \frac{4\sqrt{\pi}(\sup\limits_{\Delta > 0} \| f_\Delta \|_\infty)}{c} \| H^* \|_2^{\frac{1}{2}} \sqrt{\sigma(\tau_1, \tau_2)}, \qquad (2.39)$$

*де $c$ - стала з умови (2.1г). Крім того, псевдометрика $\rho_{(T,\Delta)}$ неперервна відносно псевдометрики $\sigma$.*

*Доведення.* З означення кореляційної функції процесу $Z_{T,\Delta}$ (див. формулу (2.8)) можна записати

$$E(Z_{T,\Delta}(\tau_2) - Z_{T,\Delta}(\tau_1))^2 = EZ_{T,\Delta}^2(\tau_2) - 2EZ_{T,\Delta}(\tau_1)Z_{T,\Delta}(\tau_2) + EZ_{T,\Delta}^2(\tau_1) =$$

$$= \frac{4\pi}{c^2} \int_{-\infty}^{\infty}\int_{-\infty}^{\infty} [1 - e^{i(\tau_1-\tau_2)\lambda_2}] |H^*(\lambda_2)|^2 \Phi_T(\lambda_2 - \lambda_1) f_\Delta(\lambda_1) f_\Delta(\lambda_2) d\lambda_1 d\lambda_2 +$$

$$+ \frac{2\pi}{c^2} \int_{-\infty}^{\infty}\int_{-\infty}^{\infty} [e^{i\tau_2(\lambda_1+\lambda_2)} - 2e^{i(\tau_1\lambda_1+\tau_2\lambda_2)} + e^{i\tau_1(\lambda_1+\lambda_2)}] H^*(\lambda_1) H^*(\lambda_2) \times$$

$$\times \Phi_T(\lambda_2 - \lambda_1) f_\Delta(\lambda_1) f_\Delta(\lambda_2) d\lambda_1 d\lambda_2 =$$

$$= \frac{4\pi}{c^2} \int_{-\infty}^{\infty}\int_{-\infty}^{\infty} [\sin\frac{(\tau_1-\tau_2)\lambda_2}{2}]^2 |H^*(\lambda_2)|^2 \Phi_T(\lambda_2 - \lambda_1) f_\Delta(\lambda_1) f_\Delta(\lambda_2) d\lambda_1 d\lambda_2 +$$

$$+ \frac{2\pi}{c^2} \int_{-\infty}^{\infty}\int_{-\infty}^{\infty} e^{i(\tau_1\lambda_1+\tau_2\lambda_2)} [e^{i\lambda_1(\tau_2-\tau_1)} - 2 + e^{i\lambda_2(\tau_1-\tau_2)}] H^*(\lambda_1) H^*(\lambda_2) \times$$

$$\times \Phi_T(\lambda_2 - \lambda_1) f_\Delta(\lambda_1) f_\Delta(\lambda_2) d\lambda_1 d\lambda_2.$$

Оскільки $|e^{i\lambda_1(\tau_2-\tau_1)} + e^{i\lambda_2(\tau_1-\tau_2)} - 2| \leq |e^{i\lambda_1(\tau_2-\tau_1)} - 1| + |e^{i\lambda_2(\tau_1-\tau_2)} - 1|$ й $|e^{i\tau\lambda} - 1| = 2|\sin\frac{\tau\lambda}{2}|$, то, оцінюючи попередню рівність із врахуванням умови (2.1б), парності ядра Фейєра та того, що $\|\Phi_T\|_1 = 1$, дістанемо

$$E|Z_{T,\Delta}(\tau_2) - Z_{T,\Delta}(\tau_1)|^2 \leq$$

$$\leq \frac{8\pi}{c^2} (\sup_{\Delta>0} \|f_\Delta\|_\infty)^2 \int_{-\infty}^{\infty}\int_{-\infty}^{\infty} [\sin\frac{(\tau_1-\tau_2)\lambda_2}{2}]^2 |H^*(\lambda_2)|^2 \Phi_T(\lambda_2 - \lambda_1) d\lambda_1 d\lambda_2 +$$

$$+\frac{4\pi}{c^2}(\sup_{\Delta>0}\|f_\Delta\|_\infty)^2\int_{-\infty}^{\infty}\int_{-\infty}^{\infty}[|\sin\frac{(\tau_2-\tau_1)\lambda_1}{2}|+|\sin\frac{(\tau_1-\tau_2)\lambda_2}{2}|]\times$$

$$\times|H^*(\lambda_1)||H^*(\lambda_2)|\Phi_T(\lambda_2-\lambda_1)d\lambda_1 d\lambda_2=$$

$$=\frac{8\pi}{c^2}(\sup_{\Delta>0}\|f_\Delta\|_\infty)^2\int_{-\infty}^{\infty}[\sin\frac{(\tau_1-\tau_2)\lambda_2}{2}]^2|H^*(\lambda_2)|^2[\int_{-\infty}^{\infty}\Phi_T(\lambda_2-\lambda_1)d\lambda_1]d\lambda_2+$$

$$+\frac{8\pi}{c^2}(\sup_{\Delta>0}\|f_\Delta\|_\infty)^2\int_{-\infty}^{\infty}\int_{-\infty}^{\infty}|\sin\frac{(\tau_1-\tau_2)\lambda_2}{2}||H^*(\lambda_1)||H^*(\lambda_2)|\Phi_T(\lambda_2-\lambda_1)d\lambda_1 d\lambda_2=$$

$$=\frac{8\pi}{c^2}(\sup_{\Delta>0}\|f_\Delta\|_\infty)^2\int_{-\infty}^{\infty}[\sin\frac{(\tau_1-\tau_2)\lambda_2}{2}]^2|H^*(\lambda_2)|^2\,d\lambda_2+$$

$$+\frac{8\pi}{c^2}(\sup_{\Delta>0}\|f_\Delta\|_\infty)^2\int_{-\infty}^{\infty}|\sin\frac{(\tau_1-\tau_2)\lambda_2}{2}||H^*(\lambda_2)|[\int_{-\infty}^{\infty}|H^*(\lambda_1)|\Phi_T(\lambda_2-\lambda_1)d\lambda_1]d\lambda_2=$$

$$=\frac{8\pi}{c^2}(\sup_{\Delta>0}\|f_\Delta\|_\infty)^2[\sigma(\tau_1,\tau_2)]^2+$$

$$+\frac{8\pi}{c^2}(\sup_{\Delta>0}\|f_\Delta\|_\infty)^2\int_{-\infty}^{\infty}|\sin\frac{(\tau_1-\tau_2)\lambda_2}{2}||H^*(\lambda_2)|[|H^*|*\Phi_T](\lambda_2)d\lambda_2.$$

Застосовуючи до другого доданку спочатку нерівність Коші - Буняковського, а потім нерівність Юнга для згортки (див. теорему А.2), отримаємо

$$\mathrm{E}|Z_{T,\Delta}(\tau_2)-Z_{T,\Delta}(\tau_1)|^2\leq$$

$$\leq\frac{8\pi}{c^2}(\sup_{\Delta>0}\|f_\Delta\|_\infty)^2[\sigma(\tau_1,\tau_2)]^2+$$

$$+\frac{8\pi}{c^2}(\sup_{\Delta>0}\|f_\Delta\|_\infty)^2[\int_{-\infty}^{\infty}[\sin\frac{(\tau_1-\tau_2)\lambda_2}{2}]^2|H^*(\lambda_2)|^2\,d\lambda_2]^{\frac{1}{2}}\||H^*|*\Phi_T\|_2=$$

$$=\frac{8\pi}{c^2}(\sup_{\Delta>0}\|f_\Delta\|_\infty)^2[\|H^*\|_2^{\frac{1}{2}}\sqrt{\sigma}(\tau_1,\tau_2)]^2+$$

$$+\frac{8\pi}{c^2}(\sup_{\Delta>0}\|f_\Delta\|_\infty)^2[\int_{-\infty}^{\infty}(\sin\frac{(\tau_1-\tau_2)\lambda_2}{2})^2|H^*(\lambda_2)|^2\,d\lambda_2]^{\frac{1}{2}}\|H^*\|_2=$$

$$= \frac{16\pi}{c^2} (\sup_{\Delta > 0} \| f_\Delta \|_\infty)^2 \| H^* \|_2 \, \sigma(\tau_1, \tau_2),$$

де використане означення функції $\sigma$ та формула зв'язку між псевдометриками $\sigma$ і $\sqrt{\sigma}$ (див. (2.37)).

Беручи квадратний корінь з обох частин отриманої нерівності, дістанемо для всіх $T > 0, \Delta > 0,$ та $\tau_1, \tau_2 \in \mathbb{R}$,

$$\left( \mathrm{E} | Z_{T,\Delta}(\tau_2) - Z_{T,\Delta}(\tau_1) |^2 \right)^{\frac{1}{2}} \leq \frac{4\sqrt{\pi}(\sup_{\Delta > 0} \| f_\Delta \|_\infty)}{c} \| H^* \|_2^{\frac{1}{2}} \sqrt{\sigma}(\tau_1, \tau_2).$$

З доведеної нерівності (2.39) та теореми Лебега про мажоровану збіжність, маємо $\rho_{(T,\Delta)}(\tau_1, \tau_2) \to 0$ при $\sqrt{\sigma}(\tau_1, \tau_2) \to 0$, тобто псевдометрика $\rho_{(T,\Delta)}$ неперервна відносно псевдометрики $\sqrt{\sigma}$. Збіжність псевдометрики $\sqrt{\sigma}$ до нуля еквівалентна збіжності псевдометрики $\sigma$ в силу означення. Таким чином, лему 2.4 доведено.

$\square$

**Наслідок 2.3.** *Оцінка для середньоквадратичних відхилень процесу $Z_{T,\Delta}$, наведена у лемі 2.4, не залежить від параметрів $T > 0, \Delta > 0$, тобто насправді встановлено співвідношення:*

$$\sup_{T,\Delta > 0} \left( \mathrm{E} | Z_{T,\Delta}(\tau_2) - Z_{T,\Delta}(\tau_1) |^2 \right)^{\frac{1}{2}} \leq \frac{4\sqrt{\pi}(\sup_{\Delta > 0} \| f_\Delta \|_\infty)}{c} \| H^* \|_2^{\frac{1}{2}} \sqrt{\sigma}(\tau_1, \tau_2), \ \tau_1, \tau_2 \in \mathbb{R}.$$

**Наслідок 2.4.** *З означення псевдометрики $\sigma$ та теореми Лебега про мажоровану збіжність, маємо $\sigma(\tau_1, \tau_2) \to 0$ при $|\tau_1 - \tau_2| \to 0$, звідки в силу означення (2.38) та нерівності (2.39) випливає, що процеси $Z_{T,\Delta}$ при всіх $T > 0, \Delta > 0,$ є неперервними у середньому квадратичному.*

У наступній теоремі наводяться достатні умови для неперервності майже напевно процесів $Z_{T,\Delta}$ та $Z$, а також слабкої збіжності $Z_{T,\Delta}$ до $Z$ при $(T,\Delta) \to \infty$ у просторі неперервних функцій.

**Теорема 2.4.** *Нехай $H \in L_2(\mathbb{R})$ та виконується нерівність*

$$\int_{0+} \mathcal{H}_{\sqrt{\sigma}}(\varepsilon)d\varepsilon < \infty, \qquad (2.40)$$

*тоді для будь-якого $[a,b] \subset \mathbb{R}$ мають місце наступні твердження:*

*(I)   $Z \in C[a,b]$ майже напевно;*

*(II)   $Z_{T,\Delta} \in C[a,b]$ майже напевно, $T>0, \Delta>0$;*

*(III)   $Z_{T,\Delta} \overset{C[a,b]}{\Rightarrow} Z$.*

*Зокрема, для будь-якого $x>0$*

$$\lim_{(T,\Delta)\to\infty} \mathrm{P}\left\{ \sup_{\tau \in [a,b]} |Z_{T,\Delta}(\tau)| > x \right\} = \mathrm{P}\left\{ \sup_{\tau \in [a,b]} |Z(\tau)| > x \right\}.$$

*Доведення.* (I) З теореми Дадлі про неперервність гауссівських процесів (див.теорему А.10) пункт (I) даної теореми має місце, якщо для будь-якого $[a,b] \subset \mathbb{R}$

$$\int_{0+} \mathcal{H}^{\frac{1}{2}}_{d_Z}([a,b],\varepsilon)d\varepsilon < \infty, \qquad (2.41)$$

де $d_Z(\tau_1,\tau_2) = [\mathrm{E}|Z(\tau_1)-Z(\tau_2)|^2]^{\frac{1}{2}}, \tau_1,\tau_2 \in \mathbb{R}$; а сам процес $Z$ є неперервним у середньому квадратичному.

Оскільки для всіх $\tau_1, \tau_2 \in \mathbb{R}$

$$d_Z^2(\tau_1,\tau_2) \le$$

$$\le \frac{2}{\pi}\int_{-\infty}^{\infty}(1-e^{i(\tau_1-\tau_2)\lambda})|H^*(\lambda)|^2\,d\lambda + \frac{1}{2\pi}\int_{-\infty}^{\infty}|e^{i(\tau_1-\tau_2)\lambda}-2+e^{i(\tau_2-\tau_1)\lambda}||H^*(\lambda)|^2\,d\lambda =$$

$$= \frac{4}{\pi}\int_{-\infty}^{\infty}[\sin\frac{(\tau_1-\tau_2)\lambda}{2}]^2|H^*(\lambda)|^2 d\lambda = \frac{4}{\pi}\sigma^2(\tau_1,\tau_2),$$

й з теореми Лебега про мажоровану збіжність, маємо $\sigma(\tau_1, \tau_2) \to 0$ при $d(\tau_1, \tau_2) = |\tau_1 - \tau_2| \to 0$, то з цієї оцінки випливає, що процес $Z$ є неперервним у середньому квадратичному.

Відповідно, нерівність (2.41) виконується, якщо $\int_{0+}^{1} \mathcal{H}_\sigma^{\frac{1}{2}}([a,b], \varepsilon) d\varepsilon < \infty$. Остання умова випливає з формули (2.40). Таким чином, пункт (I) доведено; тобто процес $Z \in C[a,b]$ майже напевно.

(II) В силу леми 2.4 при всіх $T > 0, \Delta > 0$ та $\tau_1, \tau_2 \in \mathbb{R}$, виконується оцінка

$$\rho_{T,\Delta}(\tau_1, \tau_2) = [\mathrm{E} |Z_{T,\Delta}(\tau_1) - Z_{T,\Delta}(\tau_2)|^2]^{\frac{1}{2}} \leq \frac{4\sqrt{\pi}(\sup_{\Delta > 0} \|f_\Delta\|_\infty)}{c} \|H^*\|_2^{\frac{1}{2}} \sqrt{\sigma(\tau_1, \tau_2)},$$

тому з умови (2.40) для будь-якого $[a,b] \subset \mathbb{R}$ має місце $\int_{0+}^{1} \mathcal{H}_{\rho_{T,\Delta}}([a,b], \varepsilon) d\varepsilon < \infty$, звідки, в свою чергу, випливає

$$\int_{0+}^{1} \mathcal{H}_{\rho_{T,\Delta}}^{\frac{1}{2}}([a,b], \varepsilon) d\varepsilon < \infty.$$

Крім цього, з наслідку 2.4 відомо, що процес $Z_{T,\Delta}$ є неперервним у середньому квадратичному, тому за теоремою Дадлі (див.теорему А.10) маємо: $Z_{T,\Delta} \in \mathrm{C}[a,b]$ майже напевно при всіх $T > 0, \Delta > 0$. Таким чином, пункт (II) доведено.

(III) У теоремі 2.2 або теоремі 2.3 встановлено, що за умови $H \in L_2(\mathbb{R})$:

(A) всі скінченновимірні розподіли процесу $Z_{T,\Delta}$ прямують при $(T, \Delta) \to \infty$ до відповідних скінченновимірних розподілів процесу $Z$.

Оскільки $Z_{T,\Delta}$ - сім'я (по $T > 0, \Delta > 0$) квадратично-гауссівських процесів, то з теореми А.12 випливає, що

$$\sup_{T,\Delta>0}\sup_{\tau_1,\tau_2\in[a,b]} E\exp\{|\frac{Z_{T,\Delta}(\tau_1)-Z_{T,\Delta}(\tau_2)}{\sqrt{8}\rho_{T,\Delta}(\tau_1,\tau_2)}|\}<\infty. \qquad (2.42)$$

Згідно з наслідком 2.3 маємо, що псевдометрика при всіх $\tau_1,\tau_2\in\mathbb{R}$

$$\rho_\infty(\tau_1,\tau_2)=\sup_{T,\Delta>0}\rho_{T,\Delta}(\tau_1,\tau_2)\leq\frac{4\sqrt{\pi}(\sup_{\Delta>0}\|f_\Delta\|_\infty)}{c}\|H^*\|_2^{\frac{1}{2}}\sqrt{\sigma(\tau_1,\tau_2)}, \qquad (2.43)$$

мажорується псевдометрикою $\sqrt{\sigma}$, яка є неперервною відносно рівномірної метрики $d(\tau_1,\tau_2)=|\tau_1-\tau_2|$. З теореми Лебега про мажоровану збіжність випливає, що псевдометрика $\rho_\infty(\tau_1,\tau_2)$ неперервна відносно метрики $d(\tau_1,\tau_2)$ при всіх $\tau_1,\tau_2\in\mathbb{R}$.

З умови (2.40) та нерівності (2.39) випливає, що для будь-якого $[a,b]\subset\mathbb{R}$

$$\lim_{u\downarrow 0}\sup_{T,\Delta>0}\int_0^u \mathcal{H}_{\rho_{T,\Delta}}([a,b],\varepsilon)d\varepsilon=0. \qquad (2.44)$$

Оскільки сім'я $Z_{T,\Delta}\in C[a,b]$ майже напевно при $T>0,\Delta>0$, і виконуються співвідношення (2.42)-(2.44), то з теореми А.11, випливає, що для будь-якого $[a,b]\subset\mathbb{R}$

(Б) $\quad \forall h>0 \quad \lim_{\delta\downarrow 0}\sup_{T,\Delta>0}P\{\sup_{\substack{\tau_1,\tau_2\in[a,b]\\|\tau_1-\tau_2|<\delta}}|Z_{T,\Delta}(\tau_1)-Z_{T,\Delta}(\tau_2)|>h\}=0.$

Оскільки виконуються умови (А) та (Б), то з теореми Прохорова про слабку збіжність процесів у просторі неперервних функцій (див. теорему А.13), випливає $Z_{T,\Delta}\overset{C[a,b]}{\Rightarrow}Z$. Пункт (III) доведено.

Оскільки функція розподілу супремума модуля центрованого неперервного майже напевно гауссівського процесу, заданого на відрізку числової осі, є неперервною функцією для всіх $x>0$ [34], то з пункту (III) випливає

$$\lim_{(T,\Delta)\to\infty} P\left\{\sup_{\tau\in[a,b]}|Z_{T,\Delta}(\tau)|>x\right\} = P\left\{\sup_{\tau\in[a,b]}|Z(\tau)|>x\right\}.$$

Отже, теорема 2.4 доведена. $\square$

*Зауваження 2.7.* Твердження (I) теореми 2.4 виконується і при слабшій умові, ніж (2.40), а саме, якщо є збіжним інтеграл Дадлі

$$\int_{0+}^{1}\mathcal{H}_{\sigma}^{\frac{1}{2}}(\varepsilon)d\varepsilon<\infty.$$

Зокрема, у роботі [69] встановлено, що вказана умова має місце, якщо існує таке $\beta>0$, що

$$\int_{-\infty}^{\infty}|H^{*}(\lambda)|^{2}\ln^{1+\beta}(1+|\lambda|)d\lambda<\infty.$$

Зазначимо, що другу умову перевірити значно простіше, оскільки вона задається у термінах функцій, а не їх ентропійних характеристик.

У наступному твердженні формулюється просте обмеження на функцію $H^{*}$, при якому виконується нерівність (2.40).

*Зауваження 2.8* Умова (2.40) виконується (див. приклад 3.2.8 [14]), якщо існують такі числа $K>0, \beta>0$ і $\tau_{0}>0$, що для всіх $\tau\in(0,\tau_{0})$ виконується нерівність

$$q_{H}(\tau)\leq\frac{K}{\ln|\tau|^{4+\beta}},$$

в силу якої при $|\tau_{1}-\tau_{2}|<\tau_{0}$

$$\frac{c^{2}}{16\pi(\sup_{\Delta>0}\|f_{\Delta}\|_{\infty})\|H^{*}\|_{2}}E|Z_{T,\Delta}(\tau_{1})-Z_{T,\Delta}(\tau_{2})|^{2}\leq$$

$$\leq(\int_{-\infty}^{\infty}[\sin\frac{(\tau_{1}-\tau_{2})\lambda}{2}]^{2}|H^{*}(\lambda)|^{2}d\lambda)^{\frac{1}{2}}\leq\frac{K}{\ln|\tau_{1}-\tau_{2}|^{4+\beta}}.$$

Відомо [31], що ця нерівність виконується, якщо

$$\int_{-\infty}^{\infty} |H^*(\lambda)|^2 \ln^{4+\beta}(1+|\lambda|)d\lambda < \infty. \qquad (2.45)$$

Таким чином, з умови (2.45) при деякому $\beta > 0$, випливає нерівність (2.40).

Теорема 2.4 дає підґрунтя для побудови довірчих функціональних інтервалів граничного процесу. Зокрема, має місце наступне твердження.

*Зауваження 2.9.* Нехай $d_Z(t,s) = [\mathrm{E}|Z(t)-Z(s)|^2]^{\frac{1}{2}}, t,s \in \mathbb{R}$, - середньоквадратичне відхилення граничного центрованого гауссівського процесу $Z$; $\varepsilon_0 = \sup\limits_{t,s \in [a,b]} d_Z(t,s)$. Так як виконується умова Дадлі

$$I(\varepsilon_0) = \frac{1}{\sqrt{2}} \int_0^{\varepsilon_0} \mathcal{H}_{d_Z}^{\frac{1}{2}}([a,b],\varepsilon)d\varepsilon < \infty,$$

то при всіх $x \geq 8I(\varepsilon_0)$, має місце наступна оцінка "хвоста" розподілу супремума процесу $Z$ (див. приклад 3.4.1 [14])

$$\mathrm{P}\left\{\sup_{\tau \in [a,b]} |Z(\tau)| > x\right\} \leq 2\exp\{-\frac{1}{2\varepsilon_0^2}(x - \sqrt{8xI(\varepsilon_0)})^2\}.$$

## 2.3 Асимптотична незсуненість та конзистентність оцінки $H_{T,\Delta}$

Підрозділ присвячений вивченню деяких характеристик якості оцінки $H_{T,\Delta}$, а саме – асимптотичній незсуненості та конзистентності у середньому квадратичному. Для встановлення основних результатів використовуються позначення та факти з підрозділу 2.2.

**Асимптотична незсуненість оцінки.** Нагадаємо, що оцінка $H_{T,\Delta}$ є зсуненою, тобто, взагалі кажучи, $\mathrm{E}H_{T,\Delta}(\tau) \neq H(\tau), \tau \in \mathbb{R}$ (див. підрозділ 2.1). Далі нас цікавлять умови, за яких похибка оцінювання виду $H_{T,\Delta}(\tau) - H(\tau), \tau \in \mathbb{R}$, прямує до нуля при збіжності параметра $\Delta$ до безмежності.

Покладемо

$$\hat{v}_\Delta(\tau) = [\mathbf{E} H_{T,\Delta}(\tau) - H(\tau)], \tau \in \mathbb{R}. \qquad (2.46)$$

Анулювати невипадкову функцію $\hat{v}_\Delta$ при $\Delta \to \infty$ можна за рахунок накладення додаткових умов на порядок гладкості перехідної функції $H$ та асимптотичних властивостей кореляційних функцій $K_\Delta$. В силу формули (2.5) та рівності

$$\int_{-\infty}^{\infty} K_\Delta(t) dt = 2\pi f_\Delta(0),$$

має місце наступне представлення для $\hat{v}_\Delta$, справедливе при всіх $\tau \in \mathbb{R}$:

$$\hat{v}_\Delta(\tau) = \frac{1}{c} \int_{-\infty}^{\infty} K_\Delta(\tau - s) H(s) ds = \frac{1}{c} \int_{-\infty}^{\infty} K_\Delta(s) H(\tau - s) ds =$$

$$= \frac{1}{c} \int_{-\infty}^{\infty} K_\Delta(s) [H(\tau - s) - H(\tau)] ds + (\frac{2\pi f_\Delta(0)}{c} - 1) H(\tau). \qquad (2.47)$$

Формула (2.47) якраз і "диктує" ряд додаткових умов.

(А) Нехай $\alpha \in (0,1]$. Говоритимемо, що функція $H$ є локально ліпшицевою на $[a,b] \subseteq \mathbb{R}$, з показником $\alpha$, якщо існують такі сталі $\delta > 0$ та $M > 0$, що

$$\forall t, s \in [a,b] \ \exists \delta > 0 : |t - s| < \delta \Rightarrow |H(t) - H(s)| < M |t - s|^\alpha;$$

і, відповідно, будемо позначати: $H \in Lip_\alpha[a,b]$.

*Зауваження 2.10.* Якщо при деякому $\alpha \in (0,1]$ $H \in Lip_\alpha[a,b]$, тоді $H$ рівномірно неперервна на $[a,b]$.

(Б) Припустимо, що при заданому $\alpha \in (0,1]$ виконуються балансні умови:

$$\forall \delta > 0 \quad \lim_{\Delta \to \infty} \int_{\delta}^{\infty} K_{\Delta}(t) dt = 0; \tag{2.48а}$$

$$\forall \delta > 0 \quad \lim_{\Delta \to \infty} \int_{\delta}^{\infty} K_{\Delta}^{2}(t) dt = 0; \tag{2.48б}$$

$$\exists \delta > 0 \quad \lim_{\Delta \to \infty} \int_{-\delta}^{\delta} |K_{\Delta}(t)||t|^{\alpha} dt = 0. \tag{2.48в}$$

Умови (2.1д) та (2.48а)-(2.48в) показують, в якому саме сенсі слід розуміти $\delta$-видність сім'ї кореляційних функцій $K_{\Delta}$ вхідних процесів при $\Delta \to \infty$.

Після введення умов (А) та (Б) має місце наступне твердження.

**Теорема 2.5.** *Нехай задане $\alpha \in (0,1]$; $H \in Lip_{\alpha}(\mathbb{R}) \cap L_{2}(\mathbb{R})$ та виконуються умови (2.48а) - (2.48в), тоді:*

*(I) для будь-якого $\tau \in \mathbb{R}$   $\lim_{\Delta \to \infty} \hat{v}_{\Delta}(\tau) = 0$;*

*(II) для будь-якого $[a,b] \subset \mathbb{R}$   $\lim_{\Delta \to \infty} \sup_{\tau \in [a,b]} |\hat{v}_{\Delta}(\tau)| = 0.$*

*Доведення.* (I)   Нехай $\delta = \min\{\delta_{1}, \delta_{2}\}$, де $\delta_{1}$ - величина, що фігурує в означенні локальної ліпшивості функції $H$ (див. умову (А)); $\delta_{2}$ - величина, що фігурує в умові (2.48в). З представлення (2.47), нерівності Коші - Буняковського та парності кореляційної функції $K_{\Delta}$ випливає, що

$$|\hat{v}_{\Delta}(\tau)| \leq \frac{1}{c} | \int_{|s|<\delta} K_{\Delta}(s)[H(\tau-s) - H(\tau)] ds | +$$

$$+ \frac{1}{c} | \int_{|s| \geq \delta} K_{\Delta}(s)[H(\tau-s) - H(\tau)] ds | + | \frac{2\pi f_{\Delta}(0)}{c} - 1 || H(\tau) | \leq$$

$$\leq \frac{M}{c} | \int_{|s|<\delta} |K_{\Delta}(s)||s|^{\alpha} ds | + \frac{1}{c} | \int_{|s| \geq \delta} K_{\Delta}(s) H(\tau-s) ds | +$$

$$+\frac{1}{c}|\int\limits_{|s|\geq\delta}|K_\Delta(s)|ds\,||H(\tau)|+|\frac{2\pi f_\Delta(0)}{c}-1||H(\tau)|\leq$$

$$\leq\frac{M}{c}\int\limits_{-\delta}^{\delta}|K_\Delta(s)||s|^\alpha\,ds+\frac{2}{c}\|H\|_2\,[\int\limits_{\delta}^{\infty}K_\Delta^2(s)ds]^{\frac{1}{2}}+$$

$$+\frac{2}{c}|\int\limits_{\delta}^{\infty}|K_\Delta(s)|ds\,||H(\tau)|+|\frac{2\pi f_\Delta(0)}{c}-1||H(\tau)|.$$

З цієї нерівності, умови (2.1г) та балансних умов (Б) (див. формули (2.48а) - (2.48в)) випливає, що для будь-якого $\tau\in\mathbb{R}$

$$\lim_{\Delta\to\infty}|\hat{v}_\Delta(\tau)|=0.$$

Таким чином, пункт (I) леми 2.5 доведено.

(II) Оскільки $H$ неперервна на будь-якому $[a,b]\subset\mathbb{R}$ (див. зауваження 2.10), то за теоремою Вейєрштасса вона обмежена на цьому відрізку, тобто

$$\sup_{\tau\in[a,b]}|H(\tau)|<\infty.$$

З цього факту, нерівності

$$\sup_{\tau\in[a,b]}|\hat{v}_\Delta(\tau)|\leq\frac{M}{c}\int\limits_{-\delta}^{\delta}|K_\Delta(s)||s|^\alpha\,ds+\frac{2}{c}\|H\|_2\,[\int\limits_{\delta}^{\infty}K_\Delta^2(s)ds]^{\frac{1}{2}}+$$

$$+\frac{2}{c}|\int\limits_{\delta}^{\infty}|K_\Delta(s)|ds|\sup_{\tau\in[a,b]}|H(\tau)|+|\frac{2\pi f_\Delta(0)}{c}-1|\sup_{\tau\in[a,b]}|H(\tau)|,$$

та балансних умов (Б) (див. формули (2.48а) - (2.48в)) випливає, що

$$\lim_{\Delta\to\infty}\sup_{\tau\in[a,b]}|\hat{v}_\Delta(\tau)|=0.$$

Таким чином, пункт (II) леми 2.5 доведено. Лема 2.5 доведена повністю. □

Наведемо кілька прикладів функцій, які справджують умовам (Б).

*Приклад 2.2.* Нехай задане $\alpha\in(0,1]$. Спектральні щільності

$$f_\Delta = (\frac{c}{2\pi}\exp\left(-\frac{\lambda^2}{\Delta}\right), \lambda \in \mathbb{R})$$

та відповідні їм кореляційні функції

$$K_\Delta = (\frac{c}{2}\sqrt{\frac{\Delta}{\pi}}\exp(-\frac{\Delta t^2}{4}), t \in \mathbb{R}),$$

процесів $X_\Delta, \Delta > 0$, задовольняють умовам (2.48а) - (2.48в) при $\Delta \to \infty$.

*Доведення.* Використовуючи нерівність $\int\limits_a^\infty e^{-x^2}dx \leq \sqrt{\frac{\pi}{8}}\exp(-\frac{a^2}{2}), a > 0$, з вигляду кореляційної функції для будь-яких $\Delta > 0$ та $\delta > 0$ неважко обчислити

$$\int\limits_\delta^\infty K_\Delta(t)dt = \frac{c}{\sqrt{\pi}}\int\limits_{\frac{\delta\sqrt{\Delta}}{2}}^\infty e^{-x^2}dx \leq \frac{c}{\sqrt{8}}\exp(-\frac{\delta^2\Delta}{8});$$

$$\int\limits_\delta^\infty K_\Delta^2(t)dt = \frac{c^2\sqrt{\Delta}}{\pi\sqrt{8}}\int\limits_{\delta\sqrt{\frac{\Delta}{2}}}^\infty e^{-x^2}dx \leq \frac{c^2\sqrt{\Delta}}{8\sqrt{\pi}}\exp(-\frac{\delta^2\Delta}{4});$$

Оскільки експонента "забиває" будь-яку степеневу функцію, то при $\Delta \to \infty$, виконуються умови (2.48а) та (2.48б). Далі, так як для будь-якого $\delta > 0$

$$\int\limits_{-\delta}^\delta |K_\Delta(t)||t|^\alpha\,dt = \frac{1}{(\sqrt{\Delta})^\alpha}\frac{c}{2\sqrt{\pi}}\int\limits_{-\delta\sqrt{\Delta}}^{\delta\sqrt{\Delta}}\exp(-\frac{x^2}{4})|x|^\alpha\,dx \leq$$

$$\leq \frac{1}{(\sqrt{\Delta})^\alpha}\times\frac{c}{\sqrt{\pi}}\int\limits_0^\infty\exp(-\frac{x^2}{4})x^\alpha dx,$$

де невласний інтеграл збігається до сталої, незалежної від $\Delta > 0$. З цієї нерівності випливає умова (2.48в). □

*Приклад 2.3.* Нехай задане $\alpha \in (0,1]$. Спектральні щільності

$$f_\Delta = (\frac{c}{2\pi}\frac{\Delta}{\Delta + \lambda^2}, \lambda \in \mathbb{R})$$

та відповідні їм кореляційні функції

$$K_\Delta = (\frac{c\sqrt{\Delta}}{2}\exp(-\sqrt{\Delta}|t|), t \in \mathbb{R}),$$

процесів $X_\Delta, \Delta > 0$, задовольняють умовам (2.48а) - (2.48в) при $\Delta \to \infty$.

*Доведення.* З вигляду кореляційної функції для будь-яких $\Delta > 0$ та $\delta > 0$ неважко обчислити

$$\int_\delta^\infty K_\Delta(t)dt = \frac{c}{2}\exp(-\delta\sqrt{\Delta});$$

$$\int_\delta^\infty K_\Delta^2(t)dt = \frac{c^2\sqrt{\Delta}}{8}\exp(-2\delta\sqrt{\Delta}).$$

Оскільки експонента "забиває" будь-яку степеневу функцію, то при $\Delta \to \infty$ виконуються умови (2.48а) та (2.48б). Далі, так як для будь-якого $\delta > 0$

$$\int_{-\delta}^{\delta} |K_\Delta(t)||t|^\alpha\, dt = \frac{1}{(\sqrt{\Delta})^\alpha}\int_{-\delta\sqrt{\Delta}}^{\delta\sqrt{\Delta}} \frac{c}{2\pi}\exp(-|x|)|x|^\alpha\, dx \le$$

$$\le \frac{1}{(\sqrt{\Delta})^\alpha} \times \frac{c}{\pi}\int_0^\infty \exp(-x)x^\alpha dx,$$

де невласний інтеграл збігається до сталої, незалежної від $\Delta > 0$. З цієї нерівності випливає умова (2.48в). □

*Приклад 2.4.* Нехай задане $\alpha \in (0,1)$. Спектральні щільності

$$f_\Delta = (\frac{c}{2\pi}\exp\left(-\frac{|\lambda|}{\sqrt{\Delta}}\right), \lambda \in \mathbb{R})$$

та відповідні їм кореляційні функції

$$K_\Delta = (\frac{c}{\pi}\frac{\sqrt{\Delta}}{1+\Delta t^2}, t \in \mathbb{R}),$$

процесів $X_\Delta, \Delta > 0$, задовольняють умовам (2.48а) - (2.48в) при $\Delta \to \infty$.

*Доведення.* З вигляду кореляційної функції для будь-яких $\Delta > 0$ та $\delta > 0$ неважко обчислити

$$\int_\delta^\infty K_\Delta(t)dt = \frac{c}{\pi}[\frac{\pi}{2} - \arctan(\sqrt{\Delta}\delta)];$$

$$\int_\delta^\infty K_\Delta^2(t)dt = \frac{c^2\sqrt{\Delta}}{2\pi^2}[\arctan(\sqrt{\Delta}t) + \frac{\sqrt{\Delta}t}{1+\Delta t^2}]|_\delta^\infty =$$

$$= \frac{c^2\sqrt{\Delta}}{2\pi^2}[\lim_{N \to \infty}\arctan(\frac{\sqrt{\Delta}(N-\delta)}{1+\Delta N\delta}) - \frac{\sqrt{\Delta}\delta}{1+\Delta\delta^2}] =$$

$$= \frac{c^2\sqrt{\Delta}}{2\pi^2}[\arctan(\frac{1}{\sqrt{\Delta}\delta}) - \frac{\sqrt{\Delta}\delta}{1+\Delta\delta^2}] =$$

$$= \frac{c^2\sqrt{\Delta}}{2\pi^2}[\frac{1}{\sqrt{\Delta}\delta} - \frac{1}{3(\sqrt{\Delta}\delta)^3} - \frac{\sqrt{\Delta}\delta}{1+\Delta\delta^2} + o(\frac{1}{(\sqrt{\Delta}\delta)^5})] =$$

$$= \frac{1}{\Delta} \times \frac{c^2}{6\delta^3\pi^2}[\frac{2\Delta\delta^2-1}{1+\Delta\delta^2} + o(\frac{1}{(\sqrt{\Delta}\delta)^2})].$$

З цих рівностей видно, що при $\Delta \to \infty$, виконуються умови (2.48а) та (2.48б). Далі, так як для будь-якого $\delta > 0$

$$\int_{-\delta}^\delta |K_\Delta(t)||t|^\alpha dt = \frac{1}{(\sqrt{\Delta})^\alpha}\frac{c}{\pi}\int_{-\delta\sqrt{\Delta}}^{\delta\sqrt{\Delta}}\frac{|x|^\alpha}{1+x^2}dx \leq \frac{1}{(\sqrt{\Delta})^\alpha}\times\frac{2c}{\pi}\int_0^\infty\frac{x^\alpha}{1+x^2}dx,$$

де невласний інтеграл збігається до сталої, незалежної від $\Delta > 0$. З цієї нерівності випливає умова (2.48в). □

*Приклад 2.5.* Нехай задане $\alpha \in (0,1)$. Спектральні щільності

$$f_\Delta = (\frac{c}{2\pi}(1 - \frac{|\lambda|}{\sqrt{\Delta}})\mathbb{I}_{[-\sqrt{\Delta},\sqrt{\Delta}]}(\lambda), \lambda \in \mathbb{R})$$

та відповідні їм кореляційні функції

$$K_\Delta = (\frac{c}{2\pi\sqrt{\Delta}}(\frac{\sin(\frac{\sqrt{\Delta}t}{2})}{\frac{t}{2}})^2, t \in \mathbb{R}),$$

процесів $X_\Delta, \Delta > 0$, задовольняють умовам (2.48а) - (2.48в) при $\Delta \to \infty$.

*Доведення.* Використовуючи нерівність $\sin x \leq x, x \geq 0$, для кореляційної функції для будь-яких $\Delta > 0$ та $\delta > 0$ неважко обчислити

$$\int_\delta^\infty K_\Delta(t)dt \leq \frac{2c}{\pi\delta\sqrt{\Delta}};$$

$$\int_\delta^\infty K_\Delta^2(t)dt \leq \frac{4c^2}{3\pi^2\delta^3\Delta}.$$

З цих оцінок видно, що при $\Delta \to \infty$, виконуються умови (2.48а) та (2.48б). Далі, так як для будь-якого $\delta > 0$

$$\int_{-\delta}^{\delta} |K_\Delta(t)||t|^\alpha \, dt = \frac{1}{(\sqrt{\Delta})^\alpha} \frac{c}{2\pi} \int_{-\delta\sqrt{\Delta}}^{\delta\sqrt{\Delta}} (\frac{\sin\frac{x}{2}}{\frac{x}{2}})^2 |x|^\alpha dx =$$

$$= \frac{1}{(\sqrt{\Delta})^\alpha} \times \frac{c}{\pi} \int_0^\infty (\frac{\sin\frac{x}{2}}{\frac{x}{2}})^2 x^\alpha dx,$$

де невласний інтеграл збігається до сталої, незалежної від $\Delta > 0$. З цієї нерівності випливає умова (2.48в). □

**Конзистентність оцінки.** Нагадаємо, що оцінка, нормована $\sqrt{T}$ та центрована своїм середнім значенням, мала вид

$$Z_{T,\Delta}(\tau) = \sqrt{T}[H_{T,\Delta}(\tau) - \mathrm{E}H_{T,\Delta}(\tau)], \tau \in \mathbb{R},$$

і детально досліджувалась у підрозділі 2.2. З наступного представлення

$$\mathrm{E}|H_{T,\Delta}(\tau) - H(\tau)|^2 = \mathrm{E}|H_{T,\Delta}(\tau) - \mathrm{E}H_{T,\Delta}(\tau)|^2 + |\mathrm{E}H_{T,\Delta}(\tau) - H(\tau)|^2 =$$

$$= \frac{\mathrm{E}|Z_{T,\Delta}(\tau)|^2}{T} + |\hat{v}_\Delta(\tau)|^2, \tau \in \mathbb{R}, \qquad (2.49)$$

випливає таке твердження.

**Теорема 2.6.** *Нехай задане $\alpha \in (0,1]$; $H \in Lip_\alpha(\mathbb{R}) \cap L_2(\mathbb{R})$ та виконуються умови (2.48а) - (2.48в), тоді для будь-якого $\tau \in \mathbb{R}$*

$$\lim_{(T,\Delta) \to \infty} \mathrm{E}|H_{T,\Delta}(\tau) - H(\tau)|^2 = 0.$$

*Тобто, $H_{T,\Delta}$ є конзистентною у середньому квадратичному в точці $\tau$.*

*Доведення.* Використовуючи лему 2.3, із зображення (2.49) маємо оцінку

$$\mathrm{E}|H_{T,\Delta}(\tau) - H(\tau)|^2 = \frac{\mathrm{E}|Z_{T,\Delta}(\tau)|^2}{T} + |\hat{v}_\Delta(\tau)|^2 \leq$$

$$\leq \frac{1}{T} \times \frac{4\pi}{c^2} (\sup_{\Delta > 0} \|f_\Delta\|_\infty)^2 \|H^*\|_2^2 + |\hat{v}_\Delta(\tau)|^2.$$

Звідси та пункту (I) леми 2.5, для будь-якої точки $\tau \in \mathbb{R}$ виконується

$$\lim_{(T,\Delta) \to \infty} \mathrm{E}|H_{T,\Delta}(\tau) - H(\tau)|^2 = 0.$$

Таким чином, лема 2.6 доведена повністю. □

Зазначимо, що в роботах В. Булдигіна та Фу Лі [58, 59] питання цього підрозділу не розглядались.

## 2.4 Асимптотична поведінка похибки оцінювання $W_{T,\Delta}$

Підрозділ присвячений дослідженню асимптотичних властивостей нормованої похибки оцінювання $W_{T,\Delta}$. Для встановлення основних результатів використовуються позначення та факти з підрозділів 2.2 та 2.3 .

### 2.4.1 Анулювання функції $V_{T,\Delta}$

У пункті 2.2.3 показано, що умови $H \in L_2(\mathbb{R})$ достатньо, щоб при довільному прямуванні $T \to \infty, \Delta \to \infty$, оцінка $H_{T,\Delta}$, центрована своїм середнім значенням та нормована $\sqrt{T}$, була асимптотично нормальною. Асимптотична незсуненість та конзистентність дають лише часткову характеристику якості оцінки $H_{T,\Delta}(\tau)$ (див. підрозділ 2.3). Насправді, інтерес полягає у дослідженні асимптотичної поведінки нормованої похибки оцінювання $\sqrt{T}[H_{T,\Delta}(\tau) - H(\tau)], \tau \in \mathbb{R}$. У даному пункті наведені умови, при яких $\sqrt{T}[\mathrm{E}H_{T,\Delta}(\tau) - H(\tau)] \to 0$ при балансному прямуванні $T \to \infty, \Delta \to \infty$.

Покладемо
$$V_{T,\Delta}(\tau) = \sqrt{T}[\mathrm{E}H_{T,\Delta}(\tau) - H(\tau)], \tau \in \mathbb{R}. \qquad (2.50)$$

Анулювати невипадкову функцію $V_{T,\Delta}$ можна за рахунок підбору характеру сумісного прямування параметрів $T, \Delta$ до безмежності, та накладення додаткових умов на порядок локальної гладкості функції $H$.

Враховуючи представлення (2.46), наведене у підрозділі 2.3, для $V_{T,\Delta}$ має місце рівність при всіх $\tau \in \mathbb{R}$:

$$V_{T,\Delta}(\tau) = \sqrt{T}\hat{v}_\Delta(\tau) =$$
$$= \frac{\sqrt{T}}{c} \int_{-\infty}^{\infty} K_\Delta(s)[H(\tau - s) - H(\tau)]ds + \sqrt{T}\left(\frac{2\pi f_\Delta(0)}{c} - 1\right)H(\tau), \qquad (2.51)$$

Формула (2.51) якраз і "диктує" ряд додаткових умов:

(А) Нехай при деякому $\alpha \in (0,1]$ $H \in Lip_\alpha(\mathbb{R})$;

(Б) Припустимо, що при заданому $\alpha \in (0,1]$ параметри $T \to \infty, \Delta \to \infty$ так, що виконуються балансні умови:

$$\sqrt{T}[1 - \frac{2\pi f_\Delta(0)}{c}] \to 0; \qquad (2.52а)$$

$$\forall \delta > 0 \quad \sqrt{T} \int_\delta^\infty K_\Delta(t) dt \to 0; \qquad (2.52б)$$

$$\forall \delta > 0 \quad \sqrt{T}[\int_\delta^\infty K_\Delta^2(t) dt]^{\frac{1}{2}} \to 0; \qquad (2.52в)$$

$$\exists \delta > 0 \quad \sqrt{T} \int_{-\delta}^\delta |K_\Delta(t)||t|^\alpha dt \to 0. \qquad (2.52г)$$

Після введення умов (А) та (Б), має місце твердження.

**Лема 2.5.** *Нехай задане $\alpha \in (0,1]$; $H \in Lip_\alpha(\mathbb{R}) \cap L_2(\mathbb{R})$ та $T \to \infty, \Delta \to \infty$ так, що виконуються умови (2.52а) - (2.52г), тоді:*

*(I) для будь-якого $\tau \in \mathbb{R}$ $\quad V_{T,\Delta}(\tau) \to 0$;*

*(II) для будь-якого $[a,b] \subset \mathbb{R}$ $\quad \sup_{\tau \in [a,b]} |V_{T,\Delta}(\tau)| \to 0.$*

*Доведення.* Доведення цієї леми повторює основні прийоми доведення тверджень 2.5, з врахуванням множника $\sqrt{T}$ у балансних умовах (2.52а) -(2.52г). □

Зазначимо, що для встановлення відповідного твердження в роботі В. Булдигіна та Фу Лі [59] для функції $H$ розглядалась множина точок локальної ліпшивості з показником $\alpha \in (0,1]$, і такі ж балансні умови (Б).

Наведемо декілька прикладів функцій, які задовольняють балансним умовам (Б). Твердження наводяться без доведень, бо відповідні функції були раніше детально розглянуті (див. приклади 2.1-2.5).

*Приклад 2.6.* Нехай задане $\alpha \in (0,1]$. Спектральні щільності

$$f_\Delta = (\frac{c}{2\pi}\exp\left(-\frac{\lambda^2}{\Delta}\right), \lambda \in \mathbb{R})$$

та відповідні їм кореляційні функції

$$K_\Delta = (\frac{c}{2}\sqrt{\frac{\Delta}{\pi}}\exp(-\frac{\Delta t^2}{4}), t \in \mathbb{R}),$$

збурюючих процесів $X_\Delta, \Delta > 0$, задовольняють умовам (2.1а) - (2.1д) та (2.52а) - (2.52г), якщо $T \to \infty, \Delta \to \infty$ так, що має місце співвідношення

$$T\Delta^{-\alpha} \to 0. \tag{2.53}$$

*Приклад 2.7.* Нехай задане $\alpha \in (0,1]$. Спектральні щільності

$$f_\Delta = (\frac{c}{2\pi}\frac{\Delta}{\Delta + \lambda^2}, \lambda \in \mathbb{R})$$

та відповідні їм кореляційні функції

$$K_\Delta = (c\sqrt{\Delta}\exp(-\sqrt{\Delta}|t|), t \in \mathbb{R}),$$

збурюючих процесів $X_\Delta, \Delta > 0$, задовольняють умовам (2.1а) - (2.1д) та (2.52а) - (2.52г), якщо $T \to \infty, \Delta \to \infty$ так, що має місце співвідношення

$$T\Delta^{-\alpha} \to 0. \tag{2.53}$$

*Приклад 2.8.* Нехай задане $\alpha \in (0,1)$. Спектральні щільності

$$f_\Delta = (\frac{c}{2\pi}\exp\left(-\frac{|\lambda|}{\sqrt{\Delta}}\right), \lambda \in \mathbb{R})$$

та відповідні їм кореляційні функції

$$K_\Delta = (\frac{c}{\pi}\frac{\sqrt{\Delta}}{1 + \Delta t^2}, t \in \mathbb{R}),$$

збурюючих процесів $X_\Delta, \Delta > 0$, задовольняють умовам (2.1а) - (2.1д) та (2.52а) - (2.52г), якщо $T \to \infty, \Delta \to \infty$ так, що має місце співвідношення

$$T\Delta^{-\alpha} \to 0. \tag{2.53}$$

*Приклад 2.9.* Нехай задане $\alpha \in (0,1)$. Спектральні щільності

$$f_\Delta = (\frac{c}{2\pi}(1 - \frac{|\lambda|}{\sqrt{\Delta}})\mathbb{I}_{[-\sqrt{\Delta},\sqrt{\Delta}]}(\lambda), \lambda \in \mathbb{R})$$

та відповідні їм кореляційні функції

$$K_\Delta = (\frac{c}{2\pi\sqrt{\Delta}}(\frac{\sin(\frac{\sqrt{\Delta}t}{2})}{\frac{t}{2}})^2, t \in \mathbb{R}),$$

збурюючих процесів $X_\Delta, \Delta > 0$, задовольняють умовам (2.1а) - (2.1д) та (2.52а) - (2.52г), якщо $T \to \infty, \Delta \to \infty$ так, що має місце співвідношення

$$T\Delta^{-\alpha} \to 0. \tag{2.53}$$

### 2.4.2 Асимптотична нормальність скінченновимірних розподілів похибки $W_{T,\Delta}$

У цьому пункті розглядається асимптотична нормальність скінченновимірних розподілів нормованої похибки оцінки імпульсної перехідної функції. Покладемо

$$W_{T,\Delta}(\tau) = \sqrt{T}[H_{T,\Delta}(\tau) - H(\tau)], \tau \in \mathbb{R}. \tag{2.55}$$

Для дослідження асимптотичної поведінки процесу $W_{T,\Delta} = (W_{T,\Delta}(\tau), \tau \in \mathbb{R})$, його зручно представити у вигляді суми

$$W_{T,\Delta} = Z_{T,\Delta} + V_{T,\Delta}, \tag{2.56}$$

де
$$Z_{T,\Delta}(\tau) = \sqrt{T}[H_{T,\Delta}(\tau) - \mathrm{E}H_{T,\Delta}(\tau)], \tau \in \mathbb{R};$$

$$V_{T,\Delta}(\tau) = \sqrt{T}[\mathrm{E}H_{T,\Delta}(\tau) - H(\tau)], \tau \in \mathbb{R}.$$

У пунктах 2.2.1 - 2.2.3 розглядались властивості випадкового процесу $Z_{T,\Delta}$, а саме, доведена його асимптотична нормальність при довільному прямуванні параметрів $T, \Delta$ до безмежності. У пункті 2.4.1 знайдені умови, які дозволяють анулювати вклад невипадкової функції $V_{T,\Delta}$ у виразі для $W_{T,\Delta}$ за рахунок додаткових умов на гладкість перехідної функції $H$ та балансних умов на кореляційні функції процесів $X_\Delta$ і прямування параметрів $T, \Delta$ до безмежності. Суть подальшої роботи - об'єднати ці характеристики.

Розглянемо асимптотичну поведінку кореляційної функції процесу $W_{T,\Delta}$ при прямуванні параметрів $T, \Delta$ до безмежності.

**Теорема 2.7.** *Нехай для деякого $\alpha \in (0,1]$ $H \in Lip_\alpha(\mathbb{R}) \cap L_2(\mathbb{R})$. Якщо $T \to \infty, \Delta \to \infty$ так, що виконуються умови (2.52а) - (2.52г), тоді для всіх $\tau_1, \tau_2 \in \mathbb{R}$, має місце співвідношення*

$$\mathrm{E}W_{T,\Delta}(\tau_1)W_{T,\Delta}(\tau_2) \to C_\infty(\tau_1, \tau_2) = \qquad (2.57)$$
$$= \frac{1}{2\pi}\int_{-\infty}^{\infty}[e^{i(\tau_1-\tau_2)\lambda}|H^*(\lambda)|^2 + e^{i(\tau_1+\tau_2)\lambda}(H^*(\lambda))^2]d\lambda.$$

*Доведення.* Твердження теореми 2.7 випливає з представлення

$$\mathrm{E}W_{T,\Delta}(\tau_1)W_{T,\Delta}(\tau_2) = \mathrm{E}Z_{T,\Delta}(\tau_1)Z_{T,\Delta}(\tau_2) + V_{T,\Delta}(\tau_1)V_{T,\Delta}(\tau_2),$$

теореми 2.1 та леми 2.5, пункт (I). Таким чином, теорему 2.7 доведено. □

З теореми 2.7 випливає, що гранична кореляційна функція процесу $W_{T,\Delta}$ співпадає з кореляційною функцією деякого гауссівського центрованого процесу $Z = (Z(\tau), \tau \in \mathbb{R})$, що зустрічається у підрозділах 2.2.3 - 2.2.4. Зокрема, має місце таке твердження про асимптотичну поведінку розподілів $W_{T,\Delta}$ при прямуванні параметрів $T, \Delta$ до безмежності.

**Теорема 2.8.** *Нехай для деякого* $\alpha \in (0,1]$ $H \in Lip_\alpha(\mathbb{R}) \cap L_2(\mathbb{R})$. *Якщо* $T \to \infty, \Delta \to \infty$ *так, що виконуються умови (2.52а) - (2.52г), тоді для всіх* $m \in \mathbb{N}$ *та* $\tau_1, \tau_2, ..., \tau_m \in \mathbb{R}$, *має місце співвідношення*

$$cum(W_{T,\Delta}(\tau_j), j=1,...,m) \to \begin{cases} 0, & m = 1; \\ C_\infty(\tau_1, \tau_2), & m = 2; \\ 0, & m \geq 3, \end{cases} \quad (2.58)$$

*де* $cum(W_{T,\Delta}(\tau_j), j=1,...,m)$ - *сумісний кумулянт набору випадкових величин* $W_{T,\Delta}(\tau_j), j=1,...,m.$

*Зокрема, всі скінченновимірні розподіли процесу* $(W_{T,\Delta}(\tau), \tau \in \mathbb{R})$ *слабко збігаються до відповідних скінченновимірних розподілів центрованого гауссівського процесу* $(Z(\tau), \tau \in \mathbb{R})$ *при вказаному характері прямування* $T$ *і* $\Delta$ *до нескінченності.*

*Доведення.* З представлення процесу $W_{T,\Delta}$ (див. формулу (2.56)), та властивості інваріантності кумулянтів порядку $m \geq 2$ відносно зсувів [7], маємо

$$cum(W_{T,\Delta}(\tau_j), j=1,...,m) = \begin{cases} V_{T,\Delta}(\tau_1), & m = 1; \\ cum(Z_{T,\Delta}(\tau_j), j=1,...,m), & m \geq 2. \end{cases} \quad (2.59)$$

В силу леми 2.5, пункт (I), з формули (2.59) при $m=1$ дістанемо

$$cum(W_{T,\Delta}(\tau_1)) \to 0$$

при вказаному характері прямування $T$ і $\Delta$ до нескінченності.

При $m \geq 2$ формула (2.58) випливає з рівності (2.59) та теореми 2.2, справедливої при довільному характері прямування $T \to \infty, \Delta \to \infty$. Отже, формула (2.58) доведена.

Оскільки гауссівський розподіл однозначно визначається своїми кумулянтами, то з формули (2.58) в силу теореми А.9 випливає, що всі скінченновимірні розподіли процесу $(W_{T,\Delta}(\tau), \tau \in \mathbb{R})$ слабко збігаються до відповідних скінченновимірних розподілів гауссівського процесу $(Z(\tau), \tau \in \mathbb{R})$ при вказаному характері прямування $T$ і $\Delta$ до нескінченності. Таким чином, теорему 2.8 доведено. $\square$

**Теорема 2.9.** *Нехай для деякого $\alpha \in (0,1]$ $H \in Lip_\alpha(\mathbb{R}) \cap L_2(\mathbb{R})$. Якщо $T \to \infty, \Delta \to \infty$ так, що виконуються умови (2.52а) - (2.52г), тоді для всіх $m \in \mathbb{N}$ та $\tau_1, \tau_2, ..., \tau_m \in \mathbb{R}$, має місце співвідношення*

$$\mathrm{E}[\prod_{j=1}^{m} W_{T,\Delta}(\tau_j)] \to \mathrm{E}[\prod_{j=1}^{m} Z(\tau_j)]. \tag{2.60}$$

*Зокрема, всі скінченновимірні розподіли процесу $(W_{T,\Delta}(\tau), \tau \in \mathbb{R})$ слабко збігаються до відповідних скінченновимірних розподілів центрованого гауссівського процесу $(Z(\tau), \tau \in \mathbb{R})$ при вказаному характері прямування $T$ і $\Delta$ до нескінченності.*

*Доведення.* З представлення процесу $W_{T,\Delta}$ (див. (2.56)) випливає, що

$$\mathrm{E}\left[\prod_{j=1}^{m} W_{T,\Delta}(\tau_j)\right] = \mathrm{E} \sum_{k_1=0, k_2=0, ..., k_m=0}^{1} \left[\prod_{j=1}^{m} Z_{T,\Delta}^{k_j}(\tau_j) V_{T,\Delta}^{1-k_j}(\tau_j)\right] =$$

$$= \sum_{k_1=0, k_2=0, ..., k_m=0}^{1} \left[\mathrm{E}\prod_{j=1}^{m} Z_{T,\Delta}^{k_j}(\tau_j)\right] \prod_{j=1}^{m} V_{T,\Delta}^{1-k_j}(\tau_j).$$

При вказаному характері прямування $T$ і $\Delta$ до нескінченності, з останньої формули, застосовуючи лему 2.5, пункт (I), та теорему 2.3, дістанемо, що

$$E\left[\prod_{j=1}^{m} W_{T,\Delta}(\tau_j)\right] \to E\left[\prod_{j=1}^{m} Z_{T,\Delta}(\tau_j)\right] \to E\left[\prod_{j=1}^{m} Z(\tau_j)\right]$$

для будь-яких $m \in \mathbb{N}$ та $\tau_1, \tau_2, ..., \tau_m \in \mathbb{R}$. Отже, формула (2.60) доведена.

Оскільки гауссівський розподіл однозначно визначається своїми моментами, то з формули (2.60) та теореми Маркова (див. теорему А.8) випливає, що всі скінченновимірні розподіли процесу $(W_{T,\Delta}(\tau), \tau \in \mathbb{R})$ слабко збігаються до відповідних скінченновимірних розподілів гауссівського процесу $(Z(\tau), \tau \in \mathbb{R})$ при вказаному характері прямування $T$ і $\Delta$ до нескінченності. Таким чином, теорему 2.9 доведено. □

### 2.4.3 Асимптотична нормальність розподілів процесу $W_{T,\Delta}$ у просторі неперервних функцій

Після встановлення асимптотичної нормальності скінченновимірних розподілів процесу $W_{T,\Delta}$ (теореми 2.8 і 2.9), природно поставити питання про асимптотичну нормальність цього процесу у просторі неперервних функцій. У цьому пункті використовуються термінологія та позначення пункту 2.2.4.

Будемо вважати, що гауссівський процес

$$W_{T,\Delta} = (W_{T,\Delta}(\tau), \tau \in [a,b])$$

є сепарабельним. (Таке припущення не є обмежувальним в силу стохастичної неперервності цього процесу). Має місце наступне твердження про асимптотичну нормальність процесу $W_{T,\Delta}$ у просторі неперервних функцій.

**Теорема 2.10.** *Нехай для деякого $\alpha \in (0,1]$ $H \in Lip_{\alpha}(\mathbb{R}) \bigcap L_2(\mathbb{R})$, та виконується нерівність*

$$\int_{0+}\mathcal{H}_{\sqrt{\sigma}}(\varepsilon)d\varepsilon<\infty, \qquad (2.40)$$

*тоді для будь-якого $[a,b]\subset\mathbb{R}$ мають місце наступні твердження:*

*(I)    $Z\in C[a,b]$ майже напевно;*

*(II)   $W_{T,\Delta}\in C[a,b]$ майже напевно, $T>0, \Delta>0$;*

*Крім того, якщо $T\to\infty, \Delta\to\infty$ так, що виконуються умови (2.52а) - (2.52г), тоді*

*(III)  $W_{T,\Delta}\overset{C[a,b]}{\Rightarrow}Z$.*

*Зокрема, при вказаному характері прямування $T$ і $\Delta$ до нескінченності, для будь-якого $x>0$*

$$\mathrm{P}\left\{\sup_{\tau\in[a,b]}\left|W_{T,\Delta}(\tau)\right|>x\right\}\to\mathrm{P}\left\{\sup_{\tau\in[a,b]}\left|Z(\tau)\right|>x\right\}.$$

*Доведення.* (I)    Використовуючи формулу (2.40), пункт (I) теореми 2.10 був доведений у пункті (I) теореми 2.4.

(II) Функція $(\mathrm{E}H_{T,\Delta}(\tau), \tau\in\mathbb{R})$ - невипадкова і неперервна на $\mathbb{R}$, як нормована (на $\frac{1}{c}$) сумісна кореляційна функція процесів $X_\Delta$ та $Y_\Delta$, кожен з яких є неперервним у середньому квадратичному. Згідно із зауваженням 2.10, перехідна функція $H$ є неперервною на $\mathbb{R}$. Тому $V_{T,\Delta}$, як нормована різниця цих функцій, є неперервною на $\mathbb{R}$ функцією.

Пункт (II) теореми 2.10 випливає з представлення процесу $W_{T,\Delta}=Z_{T,\Delta}+V_{T,\Delta}$ (див. (2.56)) та пункту (II) теореми 2.4.

(III)   Пункт (III) теореми 2.10 випливає з з представлення процесу $W_{T,\Delta}$ (див. формулу (2.56)), пункту (II) теореми 2.4, а також твердження (II) леми 2.5 та теореми А.14.                                                                                        □

## 2.5 Приклади функцій $H$ та процесів $X_\Delta$, які задовольняють схемі оцінювання

У підрозділі наводяться приклади імпульсних перехідних функцій $H \in L_2(\mathbb{R})$ та збурюючих процесів $(X_\Delta, \Delta > 0)$, які гарантують асимптотичну нормальність відповідних оцінки та її похибки.

*Приклад 2.10.* Розглянемо імпульсну перехідну функцію виду

$$H(t) = \frac{t}{1+t^2}, t \in \mathbb{R},$$

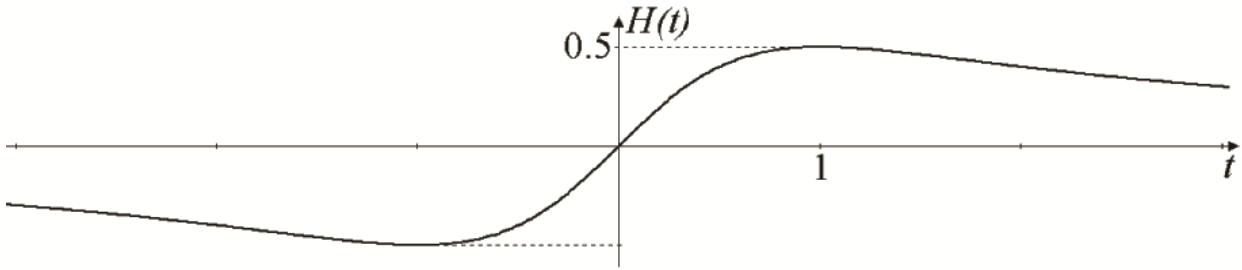

**Рис. 1.**

та дослідимо її властивості:

(I)   $H \notin L_1(\mathbb{R})$;

(II)   $H \in L_2(\mathbb{R})$,

оскільки $\|H\|_2^2 = 2\int\limits_0^\infty \frac{t^2}{(1+t^2)^2} dt \asymp 2\int\limits_0^\infty \frac{1}{(1+t)^2} dt < \infty$.

(III)   $H \in Lip_\alpha(\mathbb{R})$, при всіх $\alpha \in (0,1]$.

Дійсно, нехай $s < t$. Тоді при $\alpha = 1$,

$$|H(t) - H(s)| \leq \sup_{z \in \mathbb{R}} \left|\frac{1-z^2}{(1+z^2)^2}\right| |t-s| = 3|t-s|.$$

При будь-якому $\alpha \in (0,1)$, з нерівності Гельдера дістанемо

$$|H(t)-H(s)|=\left|\int_s^t H'(z)dz\right|=\left|\int_{-\infty}^{\infty}\frac{1-z^2}{(1+z^2)^2}\mathbb{I}_{[s,t]}(z)dz\right|\le$$

$$\le\left(\int_{-\infty}^{\infty}\left|\frac{1-z^2}{(1+z^2)^2}\right|^p dz\right)^{\frac{1}{p}}|t-s|^{\frac{1}{q}}=C(\alpha)|t-s|^\alpha,$$

де $\frac{1}{p}+\frac{1}{q}=1$, та $\alpha=\frac{1}{q}$, $C(\alpha)=\left\|\frac{1-z^2}{(1+z^2)^2}\right\|_{\frac{1}{1-\alpha}}<\infty$.

(IV) Перетворення Фур'є - Планшереля [22]

$$H^*(\lambda)=-2i\int_0^{\infty}\frac{t}{1+t^2}\sin(\lambda t)dt=\begin{cases}-\operatorname{sign}(\lambda)\pi i e^{-|\lambda|}, & \lambda\ne 0;\\ 0, & \lambda=0.\end{cases}$$

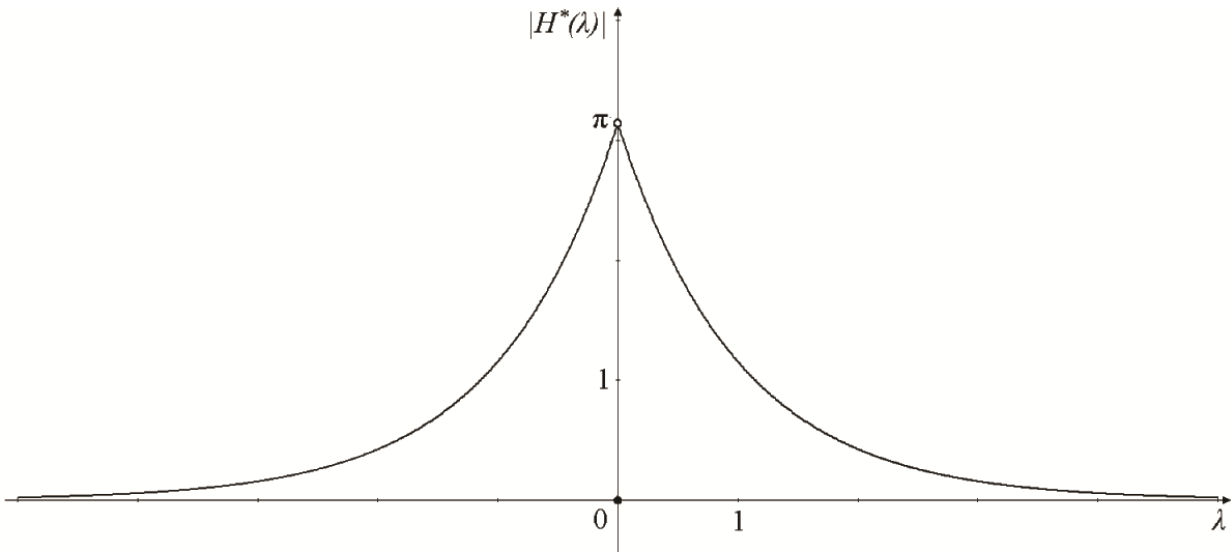

**Рис. 2.**

(V) Для будь-якого $\beta>0$

$$\int_{-\infty}^{\infty}|H^*(\lambda)|^2\ln^{4+\beta}(1+|\lambda|)d\lambda=2\pi^2\int_0^{\infty}e^{-2\lambda}\ln^{4+\beta}(1+\lambda)d\lambda<\infty,$$

оскільки при $\lambda>0$: $\ln^{4+\beta}(1+\lambda)<\lambda^{4+\beta}$, та експонента $e^{-2\lambda}$ спадає швидше, ніж будь-яка степенева функція.

Таким чином, для всіх $\alpha\in(0,1]$ $H\in Lip_\alpha(\mathbb{R})\cap L_2(\mathbb{R})$, а також в силу зауваження 2.8 з умови (V) випливає збіжність ентропійного інтеграла

$$\int_{0+}\mathcal{H}_{\sqrt{\sigma}}(\varepsilon)d\varepsilon<\infty. \tag{2.40}$$

Якщо в якості спектральних щільностей $f_\Delta$ процесів $X_\Delta$ розглядати функції з прикладів 2.2- 2.5 й врахувати, що $T\to\infty, \Delta\to\infty$ так, що

$$T\Delta^{-\alpha}\to 0, \tag{2.53}$$

то має місце функціональна теорема 2.10.

Зазначимо (див. зауваження 2.6), що оскільки $H^* \in L_1(\mathbb{R})\bigcap L_\infty(\mathbb{R})$, та є неперервною функцією майже скрізь на $\mathbb{R}$, то дану функцію $H$ можна оцінювати, використовуючи результати робіт В. Булдигіна та Фу Лі [58, 59]. □

*Приклад 2.11.* Розглянемо імпульсну перехідну функцію виду

$$H(t) = \frac{1-\cos(\mu t)}{\pi t}, t \in \mathbb{R}, \ \mu > 0,$$

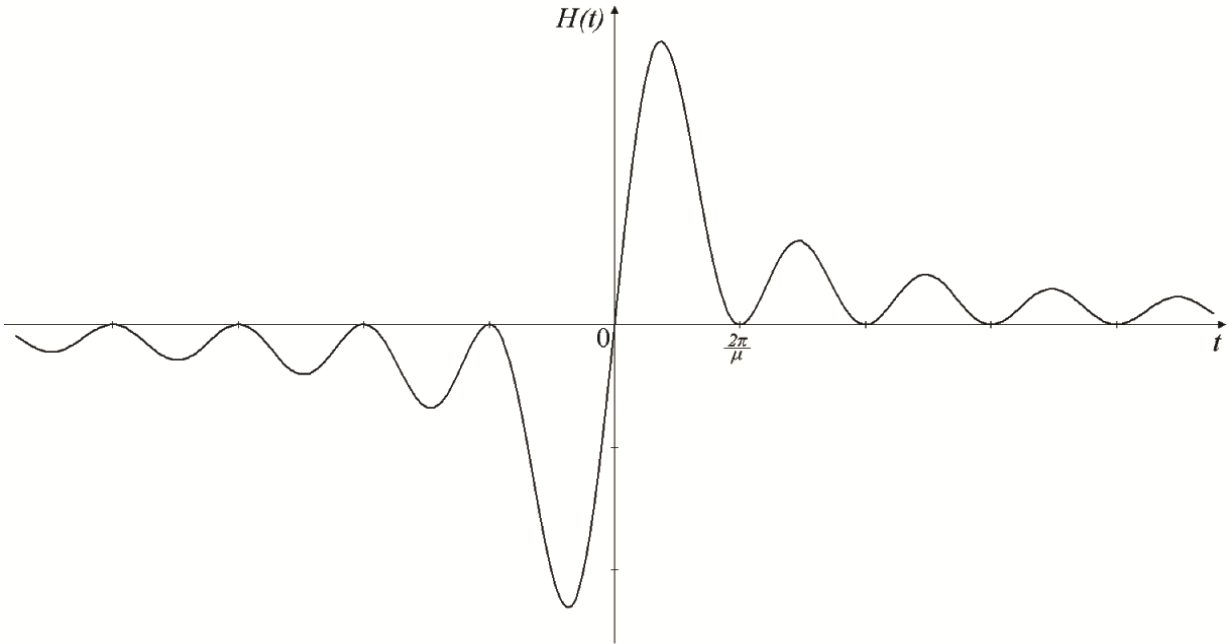

**Рис. 3.**

та дослідимо її властивості:

(I) $\lim\limits_{t\to 0}H(t) = \lim\limits_{t\to 0}\frac{\mu\sin(\mu t)}{\pi} = 0$, тобто $t = 0$ - точка усувного розриву;

(II) $H \notin L_1(\mathbb{R})$;

(III) $H \in L_2(\mathbb{R})$;

(IV) $H \in Lip_\alpha(\mathbb{R})$, при всіх $\alpha \in (0,1]$.

Дійсно, нехай $s < t$. Тоді при $\alpha = 1$,

$$|H(t) - H(s)| \leq \sup_{z \in \mathbb{R}} \left| \frac{\mu z \sin(\mu z) - 1 + \cos(\mu z)}{\pi z^2} \right| |t - s| = \frac{\mu^2}{2\pi} |t - s|.$$

При будь-якому $\alpha \in (0,1)$, з нерівності Гельдера дістанемо

$$|H(t) - H(s)| = \left| \int_s^t H'(z) dz \right| = \left| \int_{-\infty}^{\infty} \left( \frac{\mu z \sin(\mu z) - 1 + \cos(\mu z)}{\pi z^2} \right) \mathbb{I}_{[s,t]}(z) dz \right| \leq$$

$$\leq \frac{1}{\pi} \left( \mu \left\| \frac{\sin(\mu z)}{z} \right\|_p + \left\| \frac{1 - \cos(\mu z)}{z^2} \right\|_p \right) |t - s|^{\frac{1}{q}} = C(\alpha) |t - s|^\alpha,$$

де $\frac{1}{p} + \frac{1}{q} = 1$, та $\alpha = \frac{1}{q}$, $C(\alpha) = \frac{1}{\pi} \left( \mu \left\| \frac{\sin(\mu z)}{z} \right\|_{\frac{1}{1-\alpha}} + \left\| \frac{1 - \cos(\mu z)}{z^2} \right\|_{\frac{1}{1-\alpha}} \right) < \infty.$

(V) Перетворення Фур'є - Планшереля

$$H^*(\lambda) = i\left[\mathbb{I}_{[-\mu,0]}(\lambda) - \mathbb{I}_{[0,\mu]}(\lambda)\right] = -i\, sign(\lambda), |\lambda| \leq \mu.$$

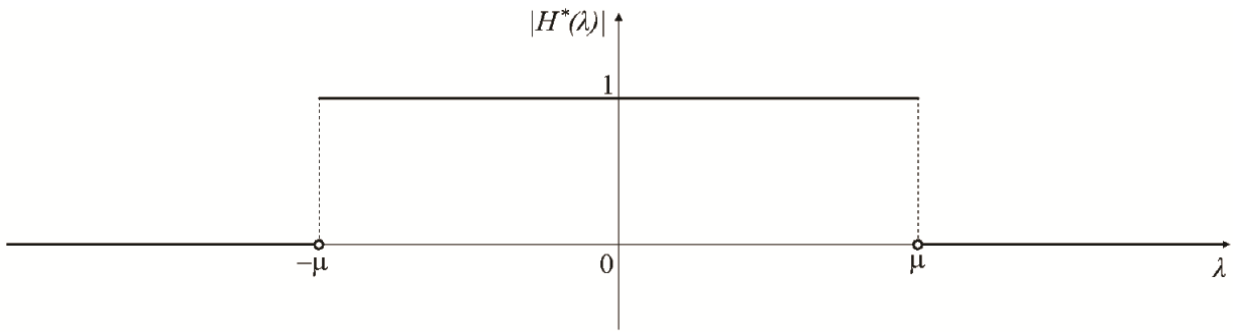

**Рис. 4.**

Дійсно, оскільки $H^* \in L_1(\mathbb{R})$, то з формули оберненого перетворення Фур'є [30], бачимо

$$\frac{1}{2\pi} \int_{-\infty}^{\infty} e^{it\lambda} H^*(\lambda) d\lambda = \frac{i}{2\pi} \left( \int_{-\mu}^{0} - \int_{0}^{\mu} \right) e^{it\lambda} d\lambda = \frac{1 - \cos(\mu t)}{\pi t} = H(t)$$

(VI) Для будь-якого $\beta > 0$

$$\int_{-\infty}^{\infty} |H^*(\lambda)|^2 \ln^{4+\beta}(1+|\lambda|) d\lambda = \left(\int_{-\mu}^{0} + \int_{0}^{\mu}\right) \ln^{4+\beta}(1+|\lambda|) d\lambda =$$

$$= 2\int_{0}^{\mu} \ln^{4+\beta}(1+\lambda) d\lambda < 2\int_{0}^{\mu} (\lambda)^{4+\beta} d\lambda = \frac{2\mu^{5+\beta}}{5+\beta} < \infty.$$

Таким чином, для всіх $\alpha \in (0,1]$ $H \in Lip_\alpha(\mathbb{R}) \cap L_2(\mathbb{R})$, а також в силу зауваження 2.8 з умови (VI) випливає збіжність ентропійного інтеграла

$$\int_{0+} \mathcal{H}_{\sqrt{\sigma}}(\varepsilon) d\varepsilon < \infty. \qquad (2.40)$$

Якщо в якості спектральних щільностей $f_\Delta$ процесів $X_\Delta$ розглядати функції з прикладів 2.2- 2.5 й врахувати, що $T \to \infty, \Delta \to \infty$ так, що

$$T\Delta^{-\alpha} \to 0, \qquad (2.53)$$

то має місце функціональна теорема 2.10.

Зазначимо (див. зауваження 2.6), що оскільки $H^* \in L_1(\mathbb{R}) \bigcap L_\infty(\mathbb{R})$, та є неперервною функцією майже скрізь на $\mathbb{R}$, то дану функцію $H$ можна оцінювати, використовуючи результати робіт В. Булдигіна та Фу Лі [58, 59]. □

*Приклад 2.12*. Розглянемо імпульсну перехідну функцію виду

$$H(t) = \frac{\sin(\mu t)}{\pi t}, t \in \mathbb{R}, \ \mu > 0,$$

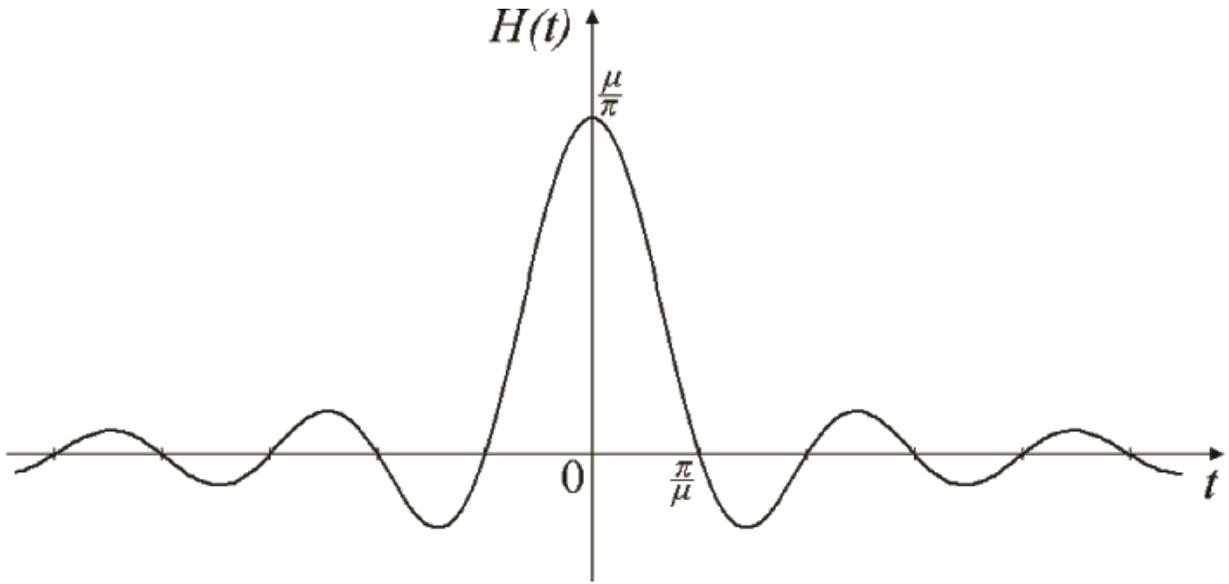

**Рис. 5.**

та дослідимо її властивості:

(I) $\lim\limits_{t \to 0} H(t) = \dfrac{\mu}{\pi}$, тобто $t = 0$ - точка усувного розриву;

(II) $H \notin L_1(\mathbb{R})$;

(III) $H \in L_2(\mathbb{R})$;

(IV) $H \in Lip_\alpha(\mathbb{R})$, при всіх $\alpha \in (0,1]$.

Дійсно, нехай $s < t$. Тоді при $\alpha = 1$,

$$|H(t) - H(s)| \leq \sup_{z \in \mathbb{R}} \left| \frac{\mu z \cos(\mu z) - \sin(\mu z)}{\pi z^2} \right| |t-s| = \frac{\mu^3}{\pi} |t-s|.$$

При будь-якому $\alpha \in (0,1)$, з нерівності Гельдера дістанемо

$$|H(t)-H(s)| = \left| \int_s^t H'(z) dz \right| = \left| \int_{-\infty}^{\infty} \left( \frac{\mu z \cos(\mu z) - \sin(\mu z)}{\pi z^2} \right) \mathbb{I}_{[s,t]}(z) dz \right| \leq$$

$$\leq \frac{2^{\frac{1}{p}}}{\pi} \left( \mu \left\| \frac{\cos(\mu z)}{z} \mathbb{I}_{[1,\infty]}(z) \right\|_p + \left\| \frac{\sin(\mu z)}{z^2} \mathbb{I}_{[1,\infty]}(z) \right\|_p + \frac{\mu^3}{\pi} \right) |t-s|^{\frac{1}{q}} = C(\alpha) |t-s|^\alpha,$$

де $\dfrac{1}{p} + \dfrac{1}{q} = 1$, та $\alpha = \dfrac{1}{q}$,

$$C(\alpha) = \frac{2^{1-\alpha}}{\pi}\left(\mu\left\|\frac{\cos(\mu z)}{z}\mathbb{I}_{[1,\infty]}(z)\right\|_{\frac{1}{1-\alpha}} + \left\|\frac{\sin(\mu z)}{z^2}\mathbb{I}_{[1,\infty]}(z)\right\|_{\frac{1}{1-\alpha}} + \frac{\mu^3}{\pi}\right) < \infty.$$

(V)  Перетворення Фур'є - Планшереля

$$H^*(\lambda) = \mathbb{I}_{[-\mu,\mu]}(\lambda), \lambda \in \mathbb{R}.$$

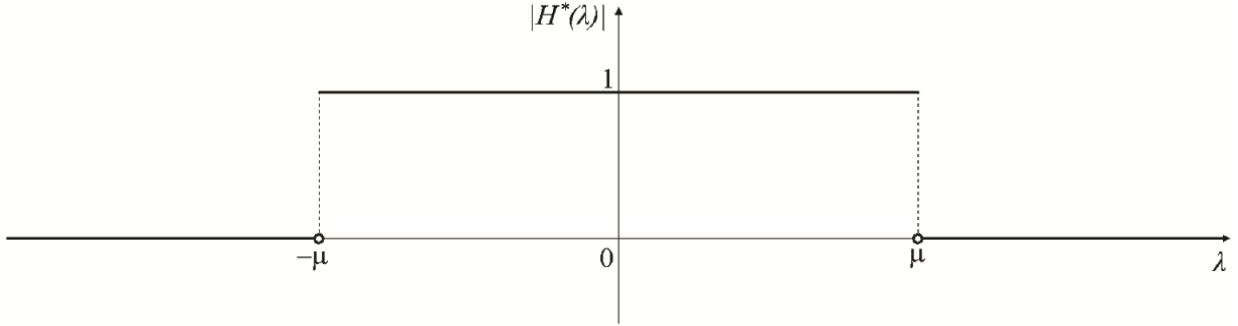

**Рис. 6.**

Дійсно, оскільки $H^* \in L_1(\mathbb{R})$, то з формули оберненого перетворення Фур'є [30], бачимо

$$\frac{1}{2\pi}\int_{-\infty}^{\infty} e^{it\lambda} H^*(\lambda)d\lambda = \frac{1}{2\pi}\int_{-\mu}^{\mu} e^{it\lambda}d\lambda = \frac{\sin(\mu t)}{\pi t} = H(t).$$

Зазначимо, що $H^*(\lambda) = \mathbb{I}_{[-\mu,\mu]}(\lambda)$, називається *фільтром низьких частот*; він пропускає лише гармонічні складові процесів з частотами $\lambda$, для яких $-\mu \leq \lambda \leq \mu, \mu > 0$ (див.,[20]).

(VI)  Для будь-якого $\beta > 0$

$$\int_{-\infty}^{\infty} |H^*(\lambda)|^2 \ln^{4+\beta}(1+|\lambda|)d\lambda = \int_{-\mu}^{\mu} \ln^{4+\beta}(1+|\lambda|)d\lambda =$$

$$= 2\int_0^{\mu} \ln^{4+\beta}(1+\lambda)d\lambda < 2\int_0^{\mu} (\lambda)^{4+\beta}d\lambda = \frac{2\mu^{5+\beta}}{5+\beta} < \infty.$$

Таким чином, для всіх $\alpha \in (0,1]$ $H \in Lip_\alpha(\mathbb{R}) \cap L_2(\mathbb{R})$, а також в силу зауваження 2.8 з умови (VI) випливає збіжність ентропійного інтеграла

$$\int_{0+}\mathcal{H}_{\sqrt{\sigma}}(\varepsilon)d\varepsilon<\infty. \qquad (2.40)$$

Якщо в якості спектральних щільностей $f_\Delta$ процесів $X_\Delta$ розглядати функції з прикладів 2.2- 2.5 й врахувати, що $T\to\infty, \Delta\to\infty$ так, що

$$T\Delta^{-\alpha}\to 0, \qquad (2.53)$$

то має місце функціональна теорема 2.10.

Зазначимо (див. зауваження 2.6), що оскільки $H^* \in L_1(\mathbb{R})\bigcap L_\infty(\mathbb{R})$, та є неперервною функцією майже скрізь на $\mathbb{R}$, то дану функцію $H$ можна оцінювати, використовуючи результати робіт В. Булдигіна та Фу Лі [58, 59]. □

*Приклад 2.13.* Розглянемо імпульсну перехідну функцію виду

$$H(t)=\frac{\sin(\nu t)-\sin(\mu t)}{\pi t}, t\in\mathbb{R},\ \nu>\mu>0,$$

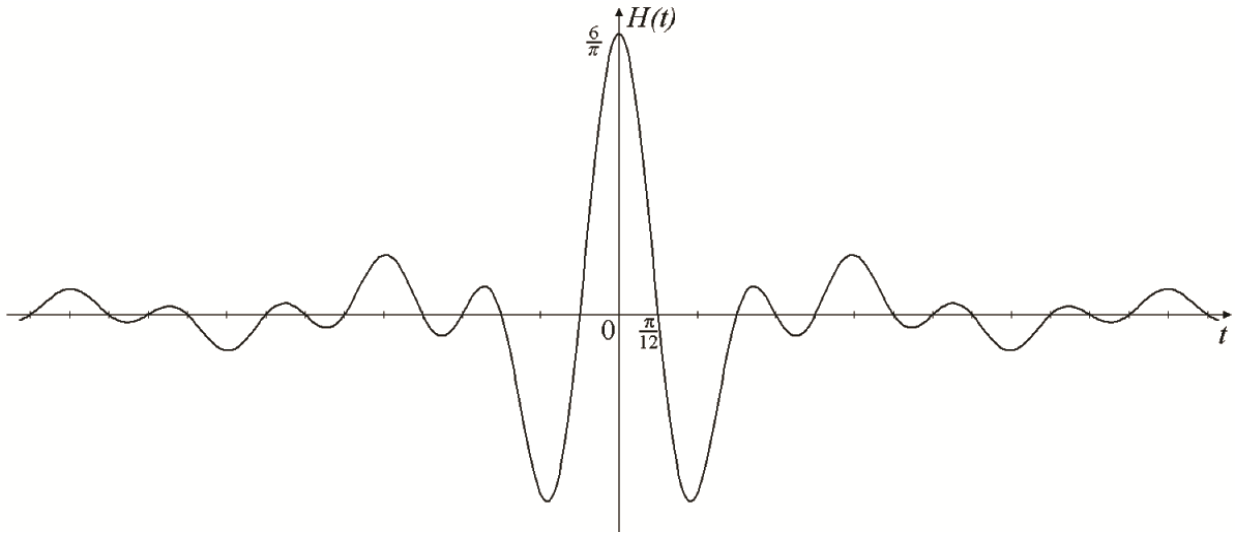

**Рис. 7.** Приклад функції $H$ з $\nu=9$ та $\mu=3$.

та дослідимо її властивості:

(I) $\lim\limits_{t\to 0}H(t)=\dfrac{\nu-\mu}{\pi}$, тобто $t=0$ - точка усувного розриву;

(II) $H\notin L_1(\mathbb{R})$;

(III) $H \in L_2(\mathbb{R})$;

(IV) $H \in Lip_\alpha(\mathbb{R})$, при всіх $\alpha \in (0,1]$.

Цей факт випливає з того, що дану функцію можна розглядати як різницю двох функцій з попереднього прикладу 2.12.

(V) Перетворення Фур'є - Планшереля
$$H^*(\lambda) = \mathbb{I}_{(-\nu,-\mu)}(\lambda) + \mathbb{I}_{(\mu,\nu)}(\lambda), \lambda \in \mathbb{R}.$$

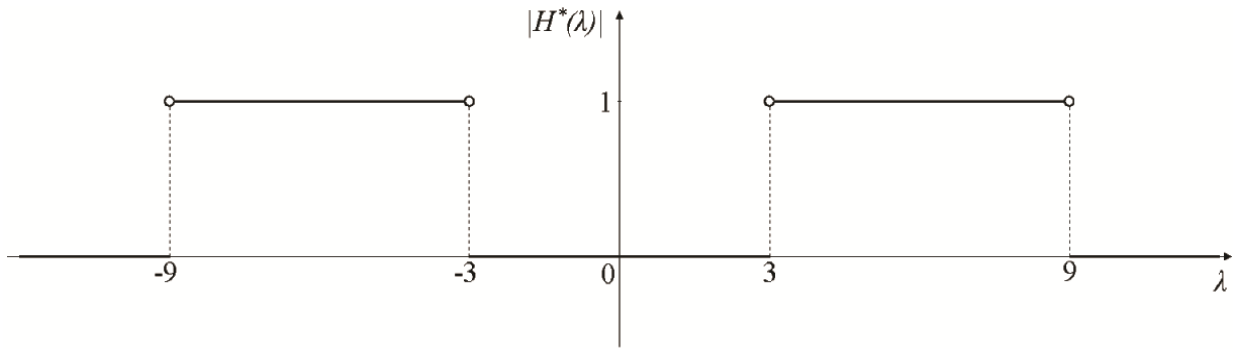

**Рис. 8.** Модуль перетворення Фур'є - Планшереля функції $H$ з $\nu = 9$ та $\mu = 3$.

Дійсно, оскільки $H^* \in L_1(\mathbb{R})$, то з формули оберненого перетворення Фур'є [30], бачимо

$$\frac{1}{2\pi}\int_{-\infty}^{\infty} e^{it\lambda} H^*(\lambda) d\lambda = \frac{1}{2\pi}\left(\int_{-\nu}^{-\mu} + \int_{\mu}^{\nu}\right) e^{it\lambda} d\lambda = \frac{\sin(\nu t) - \sin(\mu t)}{\pi t} = H(t).$$

Зазначимо, що $H^*(\lambda) = \mathbb{I}_{(-\nu,-\mu)}(\lambda) + \mathbb{I}_{(\mu,\nu)}(\lambda)$, називається *смуговим фільтром*; він пропускає (не змінюючи їх) всі гармонічні складові процесів з частотами $\lambda$, для яких $\mu \leq |\lambda| \leq \nu, \mu > 0$ (див., [20]).

(VI) Для будь-якого $\beta > 0$

$$\int_{-\infty}^{\infty} |H^*(\lambda)|^2 \ln^{4+\beta}(1+|\lambda|) d\lambda = \left(\int_{-\nu}^{\mu} + \int_{\mu}^{\nu}\right) \ln^{4+\beta}(1+|\lambda|) d\lambda =$$

$$= 2\int_\mu^\nu \ln^{4+\beta}(1+\lambda)d\lambda < 2\int_\mu^\nu (\lambda)^{4+\beta}d\lambda = \frac{2(\nu^{5+\beta}-\mu^{5+\beta})}{5+\beta} < \infty.$$

Таким чином, для всіх $\alpha \in (0,1]$ $H \in Lip_\alpha(\mathbb{R}) \cap L_2(\mathbb{R})$, а також в силу зауваження 2.8 з умови (VI) випливає збіжність ентропійного інтеграла

$$\int_{0+} \mathcal{H}_{\sqrt{\sigma}}(\varepsilon)d\varepsilon < \infty. \qquad (2.40)$$

Якщо в якості спектральних щільностей $f_\Delta$ процесів $X_\Delta$ розглядати функції з прикладів 2.2- 2.5 й врахувати, що $T \to \infty, \Delta \to \infty$ так, що

$$T\Delta^{-\alpha} \to 0, \qquad (2.53)$$

то має місце функціональна теорема 2.10.

Зазначимо (див. зауваження 2.6), що оскільки $H^* \in L_1(\mathbb{R}) \bigcap L_\infty(\mathbb{R})$, та є неперервною функцією майже скрізь на $\mathbb{R}$, то дану функцію $H$ можна оцінювати, використовуючи результати робіт В. Булдигіна та Фу Лі [58, 59]. □

Згідно з попереднім прикладом 2.13 має місце наступне узагальнення.

*Приклад 2.14.* Розглянемо імпульсну перехідну функцію виду

$$H(t) = \sum_{k=1}^{n} \frac{\sin(\nu_k t) - \sin(\mu_k t)}{\pi t}, t \in \mathbb{R},$$

де $0 < \mu_k < \nu_k \leq \mu_{k+1}, k = 1,...,n,$ та дослідимо її властивості:

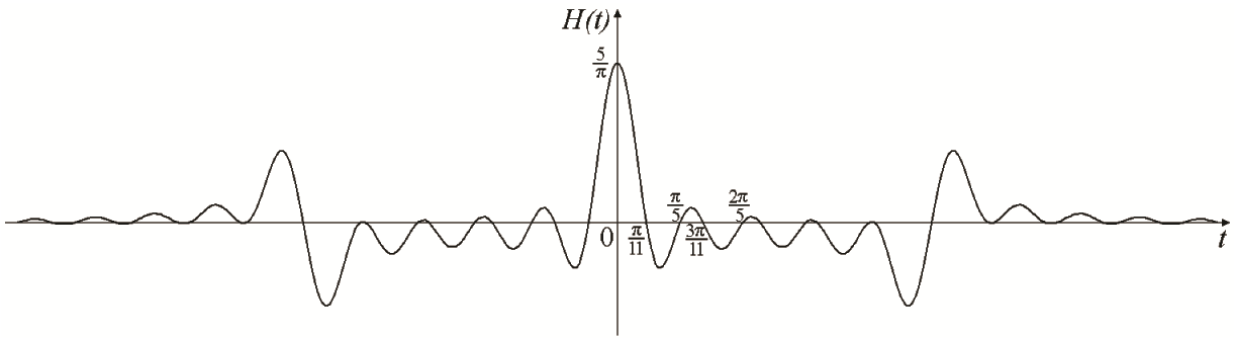

**Рис. 9.** Приклад функції $H$ з $n = 5$ та $\nu_i = 2i$, $\mu_i = 2i-1$, $i = \overline{1,5}$.

(I) $\lim\limits_{t \to 0} H(t) = \sum\limits_{k=1}^{n} \dfrac{\nu_k - \mu_k}{\pi}$, тобто $t = 0$ - точка усувного розриву;

(II) $H \notin L_1(\mathbb{R})$;

(III) $H \in L_2(\mathbb{R})$;

(IV) $H \in Lip_\alpha(\mathbb{R})$, при всіх $\alpha \in (0,1]$.

Цей факт випливає з того, що дану функцію можна розглядати як скінченну суму функцій з попереднього прикладу 2.13.

(V) Перетворення Фур'є - Планшереля

$$H^*(\lambda) = \sum_{k=1}^{n} \left[ \mathbb{I}_{(-\nu_k, -\mu_k)}(\lambda) + \mathbb{I}_{(\mu_k, \nu_k)}(\lambda) \right], \lambda \in \mathbb{R}.$$

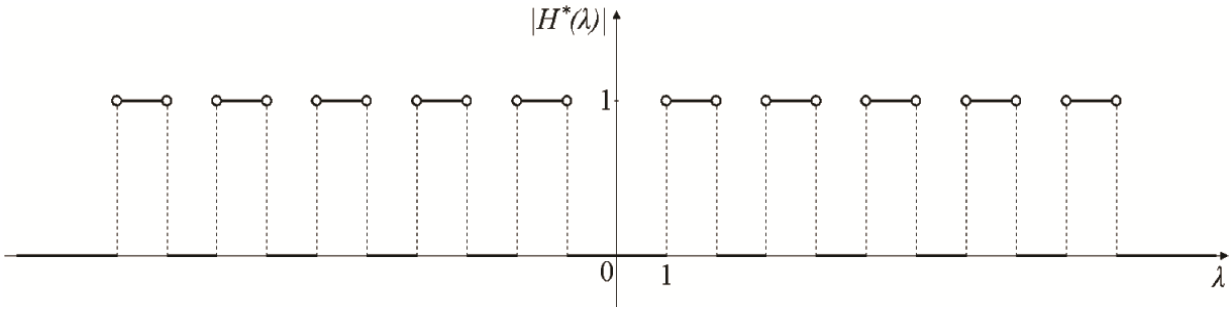

**Рис. 10.** Модуль перетворення Фур'є - Планшереля функції $H$ з $n = 5$ та $\nu_i = 2i$, $\mu_i = 2i - 1$, $i = \overline{1,5}$.

(VI) Для будь-якого $\beta > 0$

$$\int\limits_{-\infty}^{\infty} |H^*(\lambda)|^2 \ln^{4+\beta}(1+|\lambda|)d\lambda = 2\sum_{k=1}^{n} \int\limits_{\mu_k}^{\nu_k} \ln^{4+\beta}(1+\lambda)d\lambda < \dfrac{2}{5+\beta} \sum_{k=1}^{n} \left( \nu_k^{5+\beta} - \mu_k^{5+\beta} \right) < \infty.$$

Таким чином, для всіх $\alpha \in (0,1]$ $H \in Lip_\alpha(\mathbb{R}) \cap L_2(\mathbb{R})$, а також в силу зауваження 2.8 з умови (VI) випливає збіжність ентропійного інтеграла

$$\int_{0+} \mathcal{H}_{\sqrt{\sigma}}(\varepsilon)d\varepsilon < \infty. \tag{2.40}$$

Якщо в якості спектральних щільностей $f_\Delta$ процесів $X_\Delta$ розглядати функції з прикладів 2.2- 2.5 й врахувати, що $T \to \infty, \Delta \to \infty$ так, що

$$T\Delta^{-\alpha} \to 0, \qquad (2.53)$$

то має місце функціональна теорема 2.10.

Зазначимо (див. зауваження 2.6), що оскільки $H^* \in L_1(\mathbb{R}) \bigcap L_\infty(\mathbb{R})$, та є неперервною функцією майже скрізь на $\mathbb{R}$, то дану функцію $H$ можна оцінювати, використовуючи результати робіт В. Булдигіна та Фу Лі [58, 59]. □

У прикладах 2.10 - 2.14 перетворення Фур'є - Планшереля $H^*$ імпульсної перехідної функції $H$ безпосередньо обчислювалось в елементарних функціях. Тепер розглянемо ряд прикладів, коли для $H^*(\lambda)$ можна знайти лише асимптотику при відповідних значеннях $\lambda \in \mathbb{R}$. Для різноманіття вибірки прикладів розглянемо фізично здійснимі системи (див. додаток А.4).

*Приклад 2.15.* Розглянемо імпульсну перехідну функцію виду

$$H(t) = \begin{cases} \dfrac{1}{1+t}, & t \geq 0; \\ 0, & t < 0, \end{cases}$$

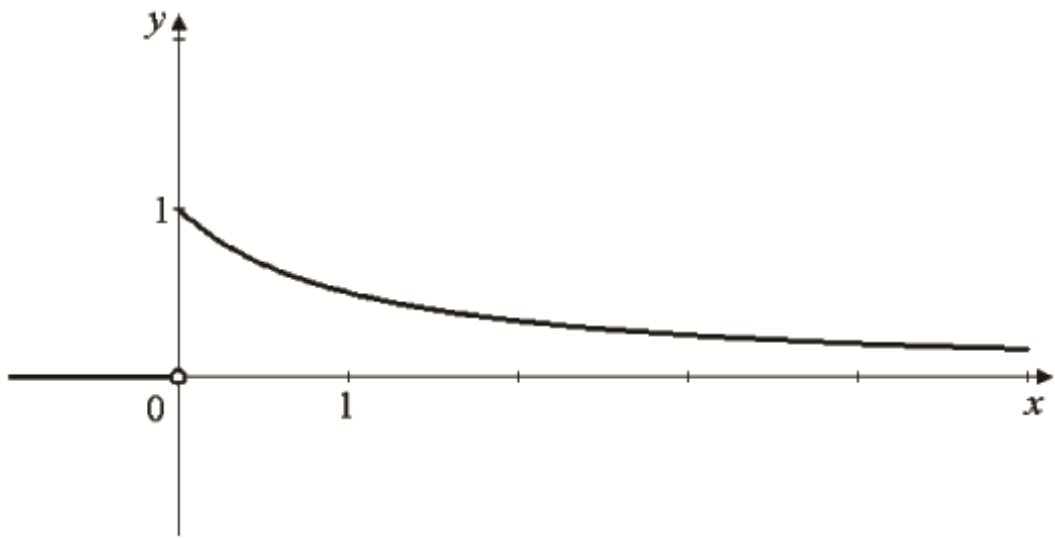

**Рис. 11.**

та дослідимо її властивості:

(I)   $H \notin L_1(\mathbb{R})$;

(II)  $H \in L_2(\mathbb{R})$;

(III) $H \in Lip_\alpha[0,\infty)$, при всіх $\alpha \in (0,1]$.

Дійсно, нехай $0 \leq s < t$. Тоді при $\alpha = 1$, $|H(t) - H(s)| \leq |t - s|$.

При будь-якому $\alpha \in (0,1)$, з нерівності Гельдера дістанемо

$$|H(t) - H(s)| = \left|\int_s^t H'(z)dz\right| = \left|\int_0^\infty \left(\frac{1}{(1+z)^2}\right)\mathbb{I}_{[s,t]}(z)dz\right| \leq$$

$$\leq \left\|\frac{1}{(1+z)^2}\mathbb{I}_{[0,\infty]}(z)\right\|_p |t-s|^{\frac{1}{q}} = C(\alpha)|t-s|^\alpha,$$

де $\dfrac{1}{p} + \dfrac{1}{q} = 1$, та $\alpha = \dfrac{1}{q}$, $C(\alpha) = \left(\dfrac{1-\alpha}{1+\alpha}\right)^{1-\alpha} < \infty.$

(IV)  Перетворення Фур'є - Планшереля

$$H^*(\lambda) = \begin{cases} e^{i\lambda}\int_1^\infty \dfrac{e^{-i\lambda t}}{t}dt, & \lambda \neq 0, \\ +\infty, & \lambda = 0; \end{cases}$$

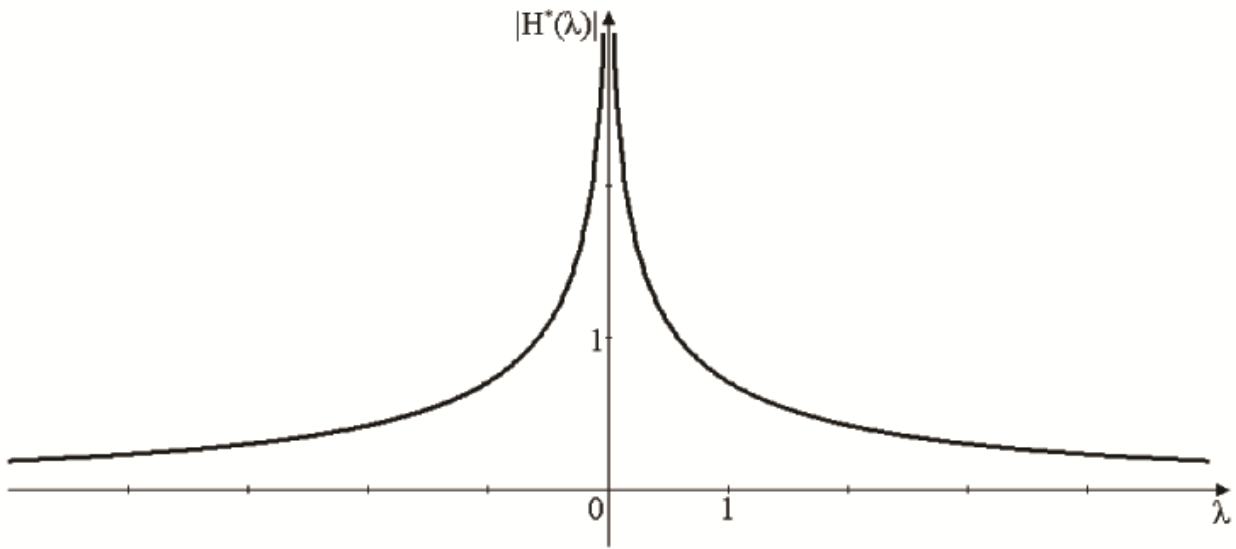

**Рис. 12.**

Розглянувши окремо випадки $\lambda > 0$ і $\lambda < 0$, маємо загальне зображення

$$H^*(\lambda) = e^{i\lambda} \int\limits_{|\lambda|}^{\infty} \frac{e^{-it \cdot sign(\lambda)}}{t} dt, \lambda \neq 0.$$

Знайдемо поведінку $H^*$ в околі особливих точок. Так, при $|\lambda| \to 0$

$$|H^*(\lambda)| = \left| e^{i\lambda} \int\limits_{|\lambda|}^{\infty} \frac{e^{-it \cdot sign(\lambda)}}{t} dt \right| | \int\limits_{|\lambda|}^{1} \frac{1}{t} dt + \left| \int\limits_{1}^{\infty} \frac{e^{-it}}{t} dt \right| = O(|\ln|\lambda||).$$

Далі, так як при $\lambda \neq 0$

$$\int\limits_{1}^{\infty} \frac{\cos(\lambda t)}{t} dt = -\frac{\sin(\lambda t)}{\lambda} + \frac{1}{\lambda} \int\limits_{1}^{\infty} \frac{\sin(\lambda t)}{t^2} dt;$$

$$\int\limits_{1}^{\infty} \frac{\sin(\lambda t)}{t} dt = \frac{\cos(\lambda t)}{\lambda} - \frac{1}{\lambda} \int\limits_{1}^{\infty} \frac{\cos(\lambda t)}{t^2} dt,$$

то при $|\lambda| \to \infty$

$$|H^*(\lambda)| = O\left( \frac{1}{|\lambda|} \right).$$

(V)    Для будь-якого $\beta > 0$

$$\int\limits_{-\infty}^{\infty} |H^*(\lambda)|^2 \ln^{4+\beta}(1+|\lambda|) d\lambda = 2\left( \int\limits_{0}^{1} + \int\limits_{1}^{\infty} \right) |H^*(\lambda)|^2 \ln^{4+\beta}(1+\lambda) d\lambda \asymp$$

$$\asymp \int\limits_{0}^{1} (\ln(\lambda))^2 \ln^{4+\beta}(1+\lambda) d\lambda + \int\limits_{1}^{\infty} \left( \frac{1}{\lambda} \right)^2 \ln^{4+\beta}(1+\lambda) d\lambda \asymp$$

$$\asymp \int\limits_{0}^{1} \lambda^{4+\beta} \ln^2(\lambda) d\lambda + (-1)^{5+\beta} \int\limits_{1}^{\infty} \ln^{4+\beta}\left( \frac{1}{\lambda} \right) d\left( \frac{1}{\lambda} \right) =$$

$$= \int\limits_{0}^{1} \lambda^{4+\beta} \ln^2(\lambda) d\lambda + (-1)^{4+\beta} \int\limits_{0}^{1} \ln^{4+\beta}(\lambda) d\lambda < \infty.$$

Таким чином, для всіх $\alpha \in (0,1]$ $H \in Lip_\alpha[0,\infty) \cap L_2(\mathbb{R})$, а також в силу зауваження 2.8 з умови (V) випливає збіжність ентропійного інтеграла для будь-якого $[a,b] \in [0,\infty)$

$$\int_{0+}\mathcal{H}_{\sqrt{\sigma}}(\varepsilon)d\varepsilon<\infty. \qquad (2.40)$$

Якщо в якості спектральних щільностей $f_\Delta$ процесів $X_\Delta$ розглядати функції з прикладів 2.2- 2.5 й врахувати, що $T\to\infty, \Delta\to\infty$ так, що

$$T\Delta^{-\alpha}\to 0, \qquad (2.53)$$

то має місце функціональна теорема 2.10.

Зазначимо (див. зауваження 2.6), що оскільки $H^* \in L_p(\mathbb{R}), p>1$, та є неперервною функцією майже скрізь на $\mathbb{R}$, то дану функцію $H$ можна оцінювати, використовуючи результати робіт В. Булдигіна та Фу Лі [58, 59]. □

*Приклад 2.16.* Розглянемо імпульсну перехідну функцію виду

$$H(t) = \begin{cases} \dfrac{1}{(1+t)^\gamma}, & t\geq 0; \\ 0, & t<0, \end{cases}$$

де $\gamma \in \left(\dfrac{1}{2},1\right)$, та дослідимо її властивості:

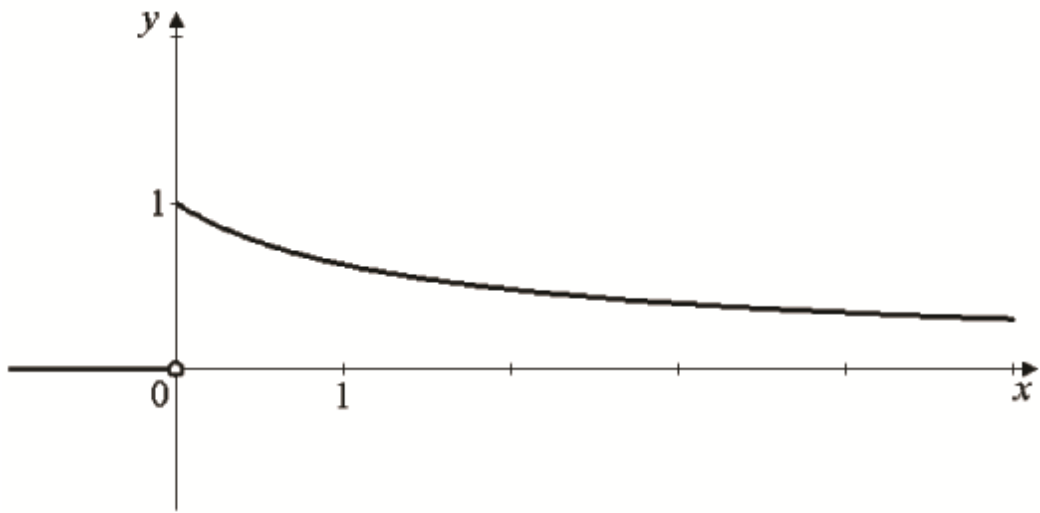

**Рис. 13.**

(I)   $H \notin L_1(\mathbb{R})$;

(II) $H \in L_2(\mathbb{R})$;

(III) $H \in Lip_\alpha[0,\infty)$, при всіх $\alpha \in (0,1]$.

Дійсно, нехай $0 \leq s < t$. Тоді при $\alpha = 1$, $|H(t) - H(s)| \leq \gamma |t-s|$.

При будь-якому $\alpha \in (0,1)$, з нерівності Гельдера дістанемо

$$|H(t) - H(s)| = \left|\int_s^t H'(z)dz\right| = \left|\int_0^\infty \left(\frac{-\gamma}{(1+z)^{\gamma+1}}\right)\mathbb{I}_{[s,t]}(z)dz\right| \leq$$

$$\leq \gamma \left\|\frac{1}{(1+z)^{\gamma+1}}\mathbb{I}_{[0,\infty]}(z)\right\|_p |t-s|^{\frac{1}{q}} = C(\alpha)|t-s|^\alpha,$$

де $\frac{1}{p} + \frac{1}{q} = 1$, та $\alpha = \frac{1}{q}$, $C(\alpha) = \gamma\left(\frac{1-\alpha}{\alpha+\gamma}\right)^{1-\alpha}$.

(IV) Перетворення Фур'є - Планшереля

$$H^*(\lambda) = \begin{cases} e^{i\lambda}\int_1^\infty \frac{e^{-i\lambda t}}{t^\gamma}dt, & \lambda \neq 0, \\ +\infty, & \lambda = 0; \end{cases}$$

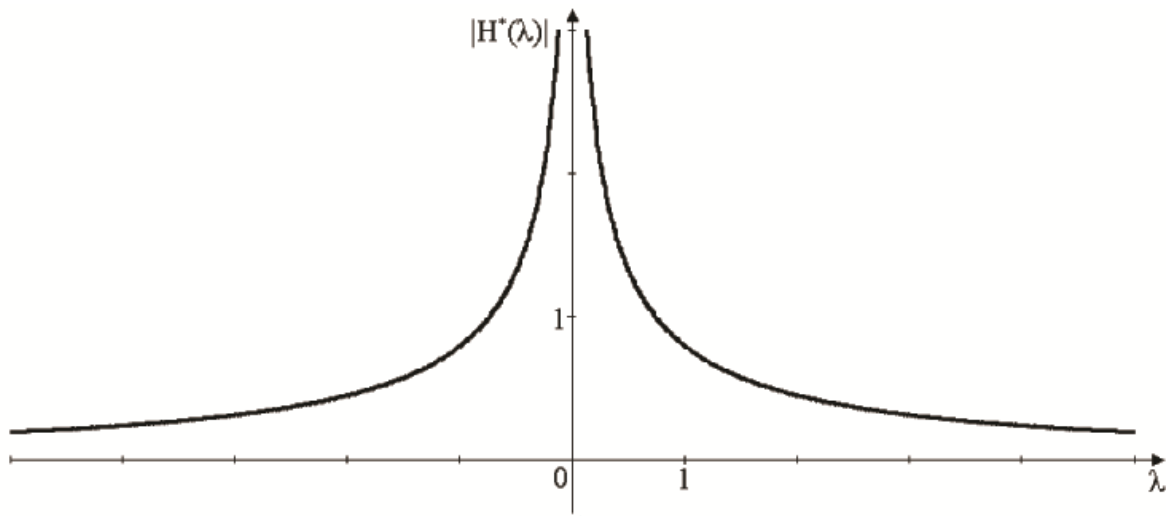

**Рис. 14.**

Розглянувши окремо випадки $\lambda > 0$ і $\lambda < 0$, маємо загальне зображення

$$H^*(\lambda) = \frac{e^{i\lambda}}{|\lambda|^{1-\gamma}} \int\limits_{|\lambda|}^{\infty} \frac{e^{-it \cdot sign(\lambda)}}{t^{\gamma}} dt, \lambda \neq 0.$$

Знайдемо поведінку $H^*$ в околі особливих точок. Використовуючи ейлерові інтеграли [37] з умовою $b > 0$

$$\int\limits_{0}^{\infty} \frac{\cos(bx)}{x^s} dx = \frac{\pi b^{s-1}}{2\Gamma(s)\cos\left(\frac{s\pi}{2}\right)}, \ 0 < s < 1;$$

$$\int\limits_{0}^{\infty} \frac{\sin(bx)}{x^s} dx = \frac{\pi b^{s-1}}{2\Gamma(s)\sin\left(\frac{s\pi}{2}\right)}, \ 0 < s < 2,$$

запишемо асимптотику при $|\lambda| \to 0$

$$|H^*(\lambda)| = \frac{1}{|\lambda|^{1-\gamma}} \left|\int\limits_{|\lambda|}^{\infty} \frac{e^{-it \cdot sign(\lambda)}}{t^{\gamma}} dt\right| \sim \frac{1}{|\lambda|^{1-\gamma}} \left|\int\limits_{0}^{\infty} \frac{e^{-it \cdot sign(\lambda)}}{t^{\gamma}} dt\right| =$$

$$= \frac{1}{|\lambda|^{1-\gamma}} \left|\int\limits_{0}^{\infty} \frac{\cos(t)}{t^{\gamma}} dt - i sign(\lambda) \int\limits_{0}^{\infty} \frac{\sin(t)}{t^{\gamma}} dt\right| =$$

$$= \frac{1}{|\lambda|^{1-\gamma}} \left|\frac{\pi}{2\Gamma(\gamma)\cos\left(\frac{\gamma\pi}{2}\right)} - i sign(\lambda) \frac{\pi}{2\Gamma(\gamma)\sin\left(\frac{\gamma\pi}{2}\right)}\right| = \frac{1}{|\lambda|^{1-\gamma}} \frac{\pi}{\Gamma(\gamma)\sin(\gamma\pi)}.$$

Далі, при $|\lambda| \to \infty$

$$|H^*(\lambda)| = \left|e^{i\lambda} \int\limits_{1}^{\infty} \frac{e^{-i\lambda t}}{t^{\gamma}} dt\right| = \frac{1}{|\lambda|} \left|-e^{i\lambda} + \gamma \int\limits_{1}^{\infty} \frac{e^{-i\lambda t}}{t^{\gamma+1}} dt\right| = O\left(\frac{1}{|\lambda|}\right).$$

(V) Для будь-якого $\beta > 0$

$$\int\limits_{-\infty}^{\infty} |H^*(\lambda)|^2 \ln^{4+\beta}(1+|\lambda|) d\lambda = 2\left(\int\limits_{0}^{1} + \int\limits_{1}^{\infty}\right) |H^*(\lambda)|^2 \ln^{4+\beta}(1+\lambda) d\lambda \asymp$$

$$\asymp \int_0^1 \left(\frac{1}{|\lambda|^{1-\gamma}}\right)^2 \ln^{4+\beta}(1+\lambda)d\lambda + \int_1^\infty \left(\frac{1}{\lambda}\right)^2 \ln^{4+\beta}(1+\lambda)d\lambda \asymp$$

$$\asymp \int_0^1 \lambda^{4+\beta+2\gamma} d\lambda + (-1)^{5+\beta} \int_1^\infty \ln^{4+\beta}\left(\frac{1}{\lambda}\right) d\left(\frac{1}{\lambda}\right) = \frac{1}{3+2\gamma+\beta} + (-1)^{4+\beta} \int_0^1 \ln^{4+\beta}(\lambda)d\lambda < \infty.$$

Таким чином, для всіх $\alpha \in (0,1]$ $H \in Lip_\alpha[0,\infty) \cap L_2(\mathbb{R})$, а також в силу зауваження 2.8 з умови (V) випливає збіжність ентропійного інтеграла для будь-якого $[a,b] \in [0,\infty)$

$$\int_{0+} \mathcal{H}_{\sqrt{\sigma}}(\varepsilon) d\varepsilon < \infty. \quad (2.40)$$

Якщо в якості спектральних щільностей $f_\Delta$ процесів $X_\Delta$ розглядати функції з прикладів 2.2- 2.5 й врахувати, що $T \to \infty, \Delta \to \infty$ так, що

$$T\Delta^{-\alpha} \to 0, \quad (2.53)$$

то має місце функціональна теорема 2.10.

Зазначимо, що $H^*$ є неперервною функцією майже скрізь на $\mathbb{R}$, але $H^* \not\in L_1(\mathbb{R}) \bigcap L_\infty(\mathbb{R})$, й при $\gamma \in \left(\frac{1}{2}, \frac{3}{4}\right]$ $H^* \not\in L_{2p}(\mathbb{R})$ при деякому $p > 2$. Тому, згідно з зауваженням 2.6, **дану функцію $H$ не можна оцінювати, використовуючи результати робіт В. Булдигіна та Фу Лі** [58, 59]. □

*Приклад 2.17.* Розглянемо імпульсну перехідну функцію виду

$$H(t) = \begin{cases} \dfrac{\cos(\mu t)}{(1+t)^\gamma}, & t \geq 0; \\ 0, & t < 0, \end{cases}$$

де $\mu > 0$ та $\gamma \in \left(\dfrac{1}{2}, 1\right)$, та дослідимо її властивості:

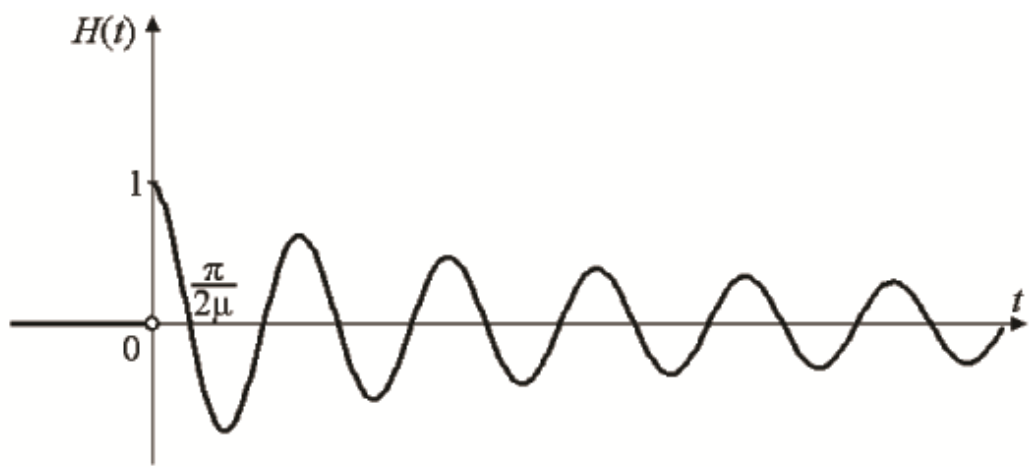

**Рис. 15.**

(I) $H \notin L_1(\mathbb{R})$;

(II) $H \in L_2(\mathbb{R})$;

(III) $H \in Lip_\alpha[0,\infty)$, при всіх $\alpha \in \left(\dfrac{1}{2}, 1\right]$.

Дійсно, нехай $0 \le s < t$. Тоді при $\alpha = 1$,

$$|H(t) - H(s)| \le \sup_{z \in [0,\infty)} \left| \frac{\mu \sin(\mu z)}{(1+z)^\gamma} + \frac{\gamma \cos(\mu z)}{(1+z)^{\gamma+1}} \right| \cdot |t-s| = (\mu + \gamma)|t-s|.$$

При будь-якому $\alpha \in (\dfrac{1}{2}, 1)$, з нерівності Гельдера дістанемо

$$|H(t) - H(s)| = \left| \int_s^t H'(z) dz \right| = \left| \int_0^\infty \left( \frac{\mu \sin(\mu z)}{(1+z)^\gamma} + \frac{\gamma \cos(\mu z)}{(1+z)^{\gamma+1}} \right) \mathbb{I}_{[s,t]}(z) dz \right| \le$$

$$\le \left( \mu \left\| \frac{\sin(\mu z)}{(1+z)^\gamma} \mathbb{I}_{[0,\infty]}(z) \right\|_p + \gamma \left\| \frac{\cos(\mu z)}{(1+z)^{\gamma+1}} \mathbb{I}_{[0,\infty]}(z) \right\|_p \right) |t-s|^{\frac{1}{q}} \le C(\alpha)|t-s|^\alpha,$$

де $\dfrac{1}{p} + \dfrac{1}{q} = 1$, та $\alpha = \dfrac{1}{q}$, $C(\alpha) = \mu \left( \dfrac{1-\alpha}{\alpha + \gamma - 1} \right)^{1-\alpha} + \gamma \left( \dfrac{1-\alpha}{\alpha + \gamma} \right)^{1-\alpha}$.

(IV) Перетворення Фур'є - Планшереля

$$H^*(\lambda) = \begin{cases} \dfrac{1}{2}\left[ e^{i(\lambda-\mu)} \int\limits_1^\infty \dfrac{e^{-i(\lambda-\mu)t}}{t^\gamma} dt + e^{i(\lambda+\mu)} \int\limits_1^\infty \dfrac{e^{-i(\lambda+\mu)t}}{t^\gamma} dt \right], & \lambda \neq \pm\mu, \\ +\infty, & \lambda = \pm\mu. \end{cases}$$

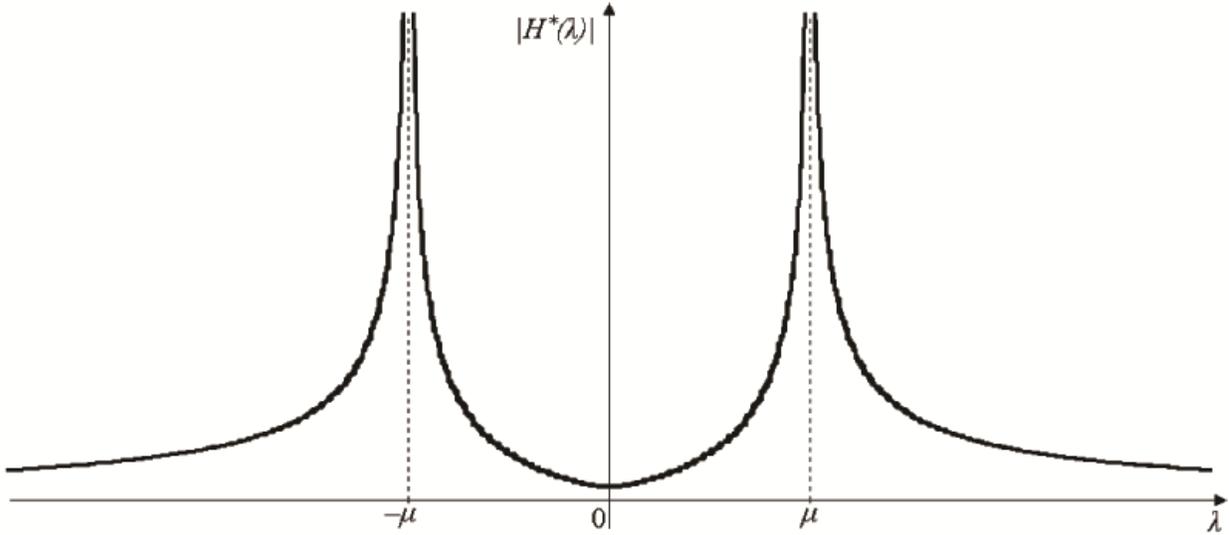

**Рис. 16.**

Розглянувши околи особливих точок $\lambda = \pm\mu$, маємо загальне зображення

$$H^*(\lambda) = \dfrac{1}{2}\left[ \dfrac{e^{i(\lambda-\mu)}}{|\lambda-\mu|^{1-\gamma}} \int\limits_{|\lambda-\mu|}^\infty \dfrac{e^{-it\cdot sign(\lambda-\mu)}}{t^\gamma} dt + \dfrac{e^{i(\lambda+\mu)}}{|\lambda+\mu|^{1-\gamma}} \int\limits_{|\lambda+\mu|}^\infty \dfrac{e^{-it\cdot sign(\lambda+\mu)}}{t^\gamma} dt \right], \lambda \neq \pm\mu.$$

Знайдемо поведінку $H^*$ в околі особливих точок. Використовуючи ейлерові інтеграли [37], запишемо асимптотику при $\lambda \to \mu$

$$|H^*(\lambda)| = \dfrac{1}{2}\left| \dfrac{e^{i(\lambda-\mu)}}{|\lambda-\mu|^{1-\gamma}} \int\limits_{|\lambda-\mu|}^\infty \dfrac{e^{-it\cdot sign(\lambda-\mu)}}{t^\gamma} dt + \dfrac{e^{i(\lambda+\mu)}}{|\lambda+\mu|^{1-\gamma}} \int\limits_{|\lambda+\mu|}^\infty \dfrac{e^{-it\cdot sign(\lambda+\mu)}}{t^\gamma} dt \right| \sim$$

$$\sim \dfrac{1}{2|\lambda-\mu|^{1-\gamma}} \left| \int\limits_0^\infty \dfrac{e^{-it\cdot sign(\lambda-\mu)}}{t^\gamma} dt \right| = \dfrac{1}{2|\lambda-\mu|^{1-\gamma}} \left| \int\limits_0^\infty \dfrac{\cos(t)}{t^\gamma} dt - i\,sign(\lambda-\mu) \int\limits_0^\infty \dfrac{\sin(t)}{t^\gamma} dt \right| =$$

$$= \dfrac{1}{2|\lambda-\mu|^{1-\gamma}} \left| \dfrac{\pi}{2\Gamma(\gamma)\cos\left(\dfrac{\gamma\pi}{2}\right)} - i\,sign(\lambda-\mu) \dfrac{\pi}{2\Gamma(\gamma)\sin\left(\dfrac{\gamma\pi}{2}\right)} \right| =$$

$$= \frac{1}{|\lambda - \mu|^{1-\gamma}} \frac{\pi}{2\Gamma(\gamma)\sin(\gamma\pi)}.$$

Аналогічно, при $\lambda \to -\mu$, маємо асимптотику

$$|H^*(\lambda)| \sim \frac{1}{2|\lambda + \mu|^{1-\gamma}} \left| \int_0^\infty \frac{e^{-it \cdot sign(\lambda+\mu)}}{t^\gamma} dt \right| =$$

$$= \frac{1}{2|\lambda+\mu|^{1-\gamma}} \left| \frac{\pi}{2\Gamma(\gamma)\cos\left(\frac{\gamma\pi}{2}\right)} - isign(\lambda+\mu)\frac{\pi}{2\Gamma(\gamma)\sin\left(\frac{\gamma\pi}{2}\right)} \right| =$$

$$= \frac{1}{|\lambda+\mu|^{1-\gamma}} \frac{\pi}{2\Gamma(\gamma)\sin(\gamma\pi)}.$$

Таким чином,

$$|H^*(\lambda)| \sim \frac{1}{||\lambda|-\mu|^{1-\gamma}} \frac{\pi}{2\Gamma(\gamma)\sin(\gamma\pi)}, \ |\lambda| \to \mu.$$

Далі, так як при $\lambda \neq \pm\mu$

$$\int_1^\infty \frac{\cos((\lambda \pm \mu)t)}{t^\gamma} dt = -\frac{\sin((\lambda \pm \mu)t)}{\lambda \pm \mu} + \frac{\gamma}{\lambda \pm \mu} \int_1^\infty \frac{\sin((\lambda \pm \mu)t)}{t^{\gamma+1}} dt;$$

$$\int_1^\infty \frac{\sin((\lambda \pm \mu)t)}{t^\gamma} dt = \frac{\cos((\lambda \pm \mu)t)}{\lambda \pm \mu} - \frac{\gamma}{\lambda \pm \mu} \int_1^\infty \frac{\cos((\lambda \pm \mu)t)}{t^{\gamma+1}} dt,$$

то при $||\lambda|-\mu| \to \infty$

$$|H^*(\lambda)| = O\left(\frac{1}{||\lambda|-\mu|}\right).$$

(V)     Для будь-якого $\beta > 0$

$$\int_{-\infty}^\infty |H^*(\lambda)|^2 \ln^{4+\beta}(1+|\lambda|)d\lambda = 2\left(\int_0^{\mu+1} + \int_{\mu+1}^\infty\right)|H^*(\lambda)|^2 \ln^{4+\beta}(1+\lambda)d\lambda \asymp$$

$$\asymp \int_0^{\mu+1} \left(\frac{1}{(\lambda-\mu)^{2-2\gamma}}\right)\ln^{4+\beta}(1+\lambda)d\lambda + \int_{\mu+1}^\infty \left(\frac{1}{\lambda-\mu}\right)^2 \ln^{4+\beta}(1+\lambda)d\lambda \leq$$

$$\leq \ln^{4+\beta}(2+\mu) \int_0^{\mu+1} \frac{1}{(\lambda-\mu)^{2-2\gamma}} d\lambda + (-1)^{5+\beta} \int_{\mu+1}^{\infty} \ln^{4+\beta}\left(\frac{1}{\lambda}\right) d\left(\frac{1}{\lambda}\right) =$$

$$= \ln^{4+\beta}(2+\mu) \frac{1+\mu^{2\gamma-1}}{2\gamma-1} + (-1)^{4+\beta} \int_0^{\frac{1}{1+\mu}} \ln^{4+\beta}(\lambda) d\lambda < \infty.$$

Таким чином, для всіх $\alpha \in \left(\frac{1}{2}, 1\right]$ $H \in Lip_\alpha[0,\infty) \cap L_2(\mathbb{R})$, а також в силу зауваження 2.8 з умови (VI) випливає збіжність ентропійного інтеграла для будь-якого $[a,b] \in [0,\infty)$

$$\int_{0+} \mathcal{H}_{\sqrt{\sigma}}(\varepsilon) d\varepsilon < \infty. \tag{2.40}$$

Якщо в якості спектральних щільностей $f_\Delta$ процесів $X_\Delta$ розглядати функції з прикладів 2.2- 2.5 й врахувати, що $T \to \infty, \Delta \to \infty$ так, що

$$T\Delta^{-\alpha} \to 0, \tag{2.53}$$

то має місце функціональна теорема 2.10.

Зазначимо, що $H^*$ є неперервною функцією майже скрізь на $\mathbb{R}$, але $H^* \not\in L_1(\mathbb{R}) \bigcap L_\infty(\mathbb{R})$, й при $\gamma \in \left(\frac{1}{2}, \frac{3}{4}\right]$ $H^* \not\in L_{2p}(\mathbb{R})$ при деякому $p > 2$. Тому, згідно з зауваженням 2.6, *дану функцію H не можна оцінювати, використовуючи результати робіт В. Булдигіна та Фу Лі* [58, 59]. □

*Приклад 2.18.* Розглянемо імпульсну перехідну функцію виду

$$H(t) = \begin{cases} \dfrac{\sin(\mu t)}{(1+t)^\gamma}, & t \geq 0; \\ 0, & t < 0, \end{cases}$$

де $\mu > 0$ та $\gamma \in \left(\dfrac{1}{2}, 1\right)$, та дослідимо її властивості:

(I)   $H \notin L_1(\mathbb{R})$;

(II)   $H \in L_2(\mathbb{R})$;

(III)   $H \in Lip_\alpha[0,\infty)$, при всіх $\alpha \in \left(\dfrac{1}{2}, 1\right]$.

Дійсно, нехай $0 \leq s < t$. Тоді при $\alpha = 1$,

$$|H(t) - H(s)| \leq \sup_{z \in [0,\infty)} \left| \frac{\mu \cos(\mu z)}{(1+z)^\gamma} - \frac{\gamma \sin(\mu z)}{(1+z)^{\gamma+1}} \right| \cdot |t-s| = (\mu + \gamma)|t-s|.$$

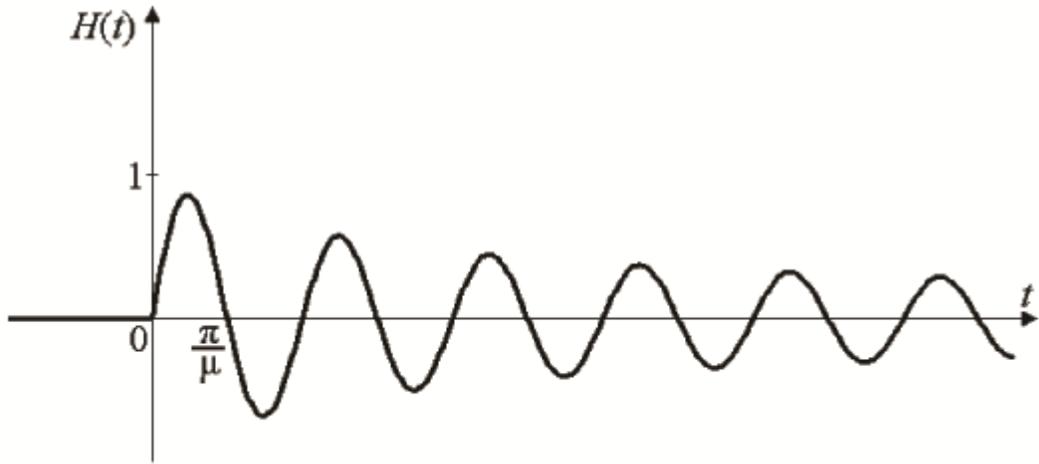

**Рис. 17.**

При будь-якому $\alpha \in \left(\dfrac{1}{2}, 1\right)$, з нерівності Гельдера дістанемо

$$|H(t) - H(s)| = \left| \int_s^t H'(z) dz \right| = \left| \int_0^\infty \left( \frac{\mu \cos(\mu z)}{(1+z)^\gamma} - \frac{\gamma \sin(\mu z)}{(1+z)^{\gamma+1}} \right) \mathbb{I}_{[s,t]}(z) dz \right| \leq$$

$$\leq \left( \mu \left\| \frac{\cos(\mu z)}{(1+z)^\gamma} \mathbb{I}_{[0,\infty]}(z) \right\|_p + \gamma \left\| \frac{\sin(\mu z)}{(1+z)^{\gamma+1}} \mathbb{I}_{[0,\infty]}(z) \right\|_p \right) |t-s|^{\frac{1}{q}} \leq C(\alpha) |t-s|^\alpha,$$

де $\dfrac{1}{p} + \dfrac{1}{q} = 1$, та $\alpha = \dfrac{1}{q}$, $C(\alpha) = \mu \left( \dfrac{1-\alpha}{\alpha + \gamma} \right)^{1-\alpha} + \gamma \left( \dfrac{1-\alpha}{\alpha + \gamma - 1} \right)^{1-\alpha}$.

(IV)   Перетворення Фур'є - Планшереля

$$H^*(\lambda) = \begin{cases} \dfrac{1}{2i}[e^{i(\lambda-\mu)}\int\limits_1^\infty \dfrac{e^{-i(\lambda-\mu)t}}{t^\gamma}dt - e^{i(\lambda+\mu)}\int\limits_1^\infty \dfrac{e^{-i(\lambda+\mu)t}}{t^\gamma}dt], & \lambda \neq \pm\mu, \\ \pm\infty, & \lambda = \pm\mu. \end{cases}$$

Розглянувши околи особливих точок $\lambda = \pm\mu$, маємо загальне зображення

$$H^*(\lambda) = \dfrac{1}{2i}\left[\dfrac{e^{i(\lambda-\mu)}}{|\lambda-\mu|^{1-\gamma}}\int\limits_{|\lambda-\mu|}^\infty \dfrac{e^{-it\cdot sign(\lambda-\mu)}}{t^\gamma}dt - \dfrac{e^{i(\lambda+\mu)}}{|\lambda+\mu|^{1-\gamma}}\int\limits_{|\lambda+\mu|}^\infty \dfrac{e^{-it\cdot sign(\lambda+\mu)}}{t^\gamma}dt\right], \lambda \neq \pm\mu.$$

Знайдемо поведінку $H^*$ в околі особливих точок. Використовуючи ейлерові інтеграли [37], запишемо асимптотику при $\lambda \to \mu$

$$|H^*(\lambda)| = \dfrac{1}{2}\left|\dfrac{e^{i(\lambda-\mu)}}{|\lambda-\mu|^{1-\gamma}}\int\limits_{|\lambda-\mu|}^\infty \dfrac{e^{-it\cdot sign(\lambda-\mu)}}{t^\gamma}dt - \dfrac{e^{i(\lambda+\mu)}}{|\lambda+\mu|^{1-\gamma}}\int\limits_{|\lambda+\mu|}^\infty \dfrac{e^{-it\cdot sign(\lambda+\mu)}}{t^\gamma}dt\right| \sim$$

$$\sim \dfrac{1}{2|\lambda-\mu|^{1-\gamma}}\left|\int\limits_0^\infty \dfrac{e^{-it\cdot sign(\lambda-\mu)}}{t^\gamma}dt\right| = \dfrac{1}{|\lambda-\mu|^{1-\gamma}}\dfrac{\pi}{2\Gamma(\gamma)\sin(\gamma\pi)}.$$

Аналогічно, при $\lambda \to -\mu$, маємо асимптотику

$$|H^*(\lambda)| \sim \dfrac{1}{2|\lambda+\mu|^{1-\gamma}}\left|\int\limits_0^\infty \dfrac{e^{-it\cdot sign(\lambda+\mu)}}{t^\gamma}dt\right| = \dfrac{1}{|\lambda+\mu|^{1-\gamma}}\dfrac{\pi}{2\Gamma(\gamma)\sin(\gamma\pi)}.$$

Таким чином,

$$|H^*(\lambda)| \sim \dfrac{1}{||\lambda|-\mu|^{1-\gamma}}\dfrac{\pi}{2\Gamma(\gamma)\sin(\gamma\pi)}, \ |\lambda| \to \mu.$$

Далі всі результати дублюють викладки попереднього прикладу. Тобто, при $||\lambda|-\mu| \to \infty$

$$|H^*(\lambda)| = O\left(\dfrac{1}{||\lambda|-\mu|}\right).$$

(V) Для будь-якого $\beta > 0$

$$\int\limits_{-\infty}^\infty |H^*(\lambda)|^2 \ln^{4+\beta}(1+|\lambda|)d\lambda = 2\left(\int\limits_0^{\mu+1} + \int\limits_{\mu+1}^\infty\right)|H^*(\lambda)|^2 \ln^{4+\beta}(1+\lambda)d\lambda < \infty.$$

Таким чином, для всіх $\alpha \in \left(\frac{1}{2},1\right]$ $H \in Lip_\alpha[0,\infty) \cap L_2(\mathbb{R})$, а також в силу зауваження 2.8 з умови (VI) випливає збіжність ентропійного інтеграла для будь-якого $[a,b] \in [0,\infty)$

$$\int_{0+} \mathcal{H}_{\sqrt{\sigma}}(\varepsilon) d\varepsilon < \infty. \tag{2.40}$$

Якщо в якості спектральних щільностей $f_\Delta$ процесів $X_\Delta$ розглядати функції з прикладів 2.2- 2.5 й врахувати, що $T \to \infty, \Delta \to \infty$ так, що

$$T\Delta^{-\alpha} \to 0, \tag{2.53}$$

то має місце функціональна теорема 2.10.

Зазначимо, що $H^*$ є неперервною функцією майже скрізь на $\mathbb{R}$, але $H^* \not\in L_1(\mathbb{R}) \bigcap L_\infty(\mathbb{R})$, й при $\gamma \in \left(\frac{1}{2}, \frac{3}{4}\right]$ $H^* \not\in L_{2p}(\mathbb{R})$ при деякому $p > 2$. Тому, згідно з зауваженням 2.6, ***дану функцію $H$ не можна оцінювати, використовуючи результати робіт В. Булдигіна та Фу Лі*** [58, 59]. □

*Зауваження 2.11.* Доведена в роботі функціональна теорема 2.10 дозволяє розглядати імпульсні перехідні функції у вигляді рядів

$$H(t) = \sum_{j=1}^\infty a_j H_j(t),$$

де $H_j(t) \in \left\{ \frac{1}{(1+|t|)^{\gamma_j}}, \frac{\sin(\mu_j t)}{(1+|t|)^{\gamma_j}}, \frac{\cos(\mu_j t)}{(1+|t|)^{\gamma_j}} \right\}$, та узгоджено вибираються коефіцієнти $a_j > 0, \mu_j > 0, \gamma_j \in \left(\frac{1}{2},1\right)$.

При цьому перетворення Фур'є-Планшереля $H^*$ матиме зліченне число точок розриву другого роду; тобто така система міститиме резонанс на множині лебегової міри нуль.

## 2.6 Висновки

У розділі 2 знайдено: умови асимптотичної незсуненості та конзистентності інтегральних корелограмних оцінок для імпульсної перехідної функції; умови асимптотичної нормальності відповідної оцінки та похибки нестійкої однорідної лінійної системи. Досліджено збіжність розподілів оцінки та похибки як в сенсі збіжності скінченновимірних розподілів, так і у сенсі слабкої збіжності у просторі неперервних функцій. За допомогою застосування діаграмно-кумулянтного методу та нерівностей для інтегралів з циклічним зачепленням ядер, при дослідженні асимптотичної поведінки розподілів інтегральної корелограмної оцінки та її похибки вдалось поліпшити відомі результати у такій схемі оцінювання. Зокрема, вдалось розширити клас перехідних функцій відносно попередніх результатів, для яких можливе застосування даної теорії оцінювання (див., приклади 2.16-2.18).

Результати одержані при мінімальних умовах на порядок інтегрованості імпульсної перехідної функції однорідної лінійної системи, і, в деякому сенсі, ставлять крапку на питання про асимптотичну нормальність інтегральної корелограмної оцінки та її похибки. Ці результати застосовано далі для корелограмного оцінювання імпульсних перехідних функцій нестійких однорідних лінійних систем з внутрішнім шумом.

Результати цього розділу були опубліковані у роботах [123], [125], [127]; і доповідалися на "Міжуніверситетській науковій конференції з математики та фізики для студентів та молодих вчених" (2009) [128], конференції "Фрактали та сучасна математика" (2009) [130], "XV Міжнародній науковій конференції ім. акад. М. Кравчука" (2014) [133].

# РОЗДІЛ 3

# АСИМПТОТИЧНІ ВЛАСТИВОСТІ КОРЕЛОГРАМНИХ ОЦІНОК ПЕРЕХІДНИХ ФУНКЦІЙ ЛІНІЙНИХ ОДНОРІДНИХ СИСТЕМ З ШУМОМ

Розділ 3 узагальнює результати розділу 2 в задачі корелограмного оцінювання невідомої дійснозначної перехідної функції $H$ лінійної системи з внутрішнім шумом та умовою $H \in L_2(\mathbb{R})$. За такого припущення на порядок інтегрованості перехідної функції, під розгляд потрапляють нестійкі системи з резонансними особливостями. Припускається, що на вхід однорідної лінійної системи подається сім'я стаціонарних центрованих гауссівських процесів близьких, в деякому сенсі, до білого шуму. На виході до реакції системи домішується внутрішній шум, який є стаціонарним центрованим гауссівським процесом, ортогональним до вхідного сигналу.

Використовуючи схему однієї вибірки, досліджуються умови асимптотичної незсуненості та консистентності корелограмних оцінок інтегрального типу, а також умови асимптотичної нормальності відповідних оцінок та похибки оцінювання, як у сенсі збіжності скінченновимірних розподілів, так і в сенсі слабкої збіжності відповідних розподілів у просторі неперервних функцій. Зокрема, наводяться приклади сім'ї збурюючих процесів та характеристик внутрішнього шуму, а також імпульсних перехідних функцій $H \in L_2(\mathbb{R})$, які задовольняють розв'язуваній задачі.

## 3.1 Оцінка імпульсної перехідної функції

У підрозділі визначається статистична оцінка імпульсної перехідної функції однорідної лінійної системи з внутрішнім шумом. Ця оцінка – відповідним чином нормована сумісна інтегральна корелограма між деяким гауссівським випадковим процесом, що збурює систему, та реакцією системи на це збурення. Цей підрозділ повторює термінологію підрозділу 2.1, розділ 2.

**Процеси, що збурюють систему.** Нехай $X_\Delta = (X_\Delta(t), t \in \mathbb{R})$, $\Delta > 0$, - сім'я вимірних сепарабельних стаціонарних центрованих дійснозначних гауссівських процесів, що збурюють лінійну систему з внутрішнім шумом (див. додаток А.4). Спектральні щільності $f_\Delta = (f_\Delta(\lambda), \lambda \in \mathbb{R})$, $\Delta > 0$, процесів $X_\Delta$ задовольняють умови:

$$f_\Delta(\lambda) = f_\Delta(-\lambda), \lambda \in \mathbb{R}; \quad (3.1а)$$

$$\sup_{\Delta > 0} \| f_\Delta \|_\infty < \infty; \quad (3.1б)$$

$$f_\Delta \in L_1(\mathbb{R}); \quad (3.1в)$$

$$\exists c \in (0,\infty) \ \forall a \in (0,\infty): \lim_{\Delta \to \infty} \sup_{-a \leq \lambda \leq a} \left| f_\Delta(\lambda) - \frac{c}{2\pi} \right| = 0; \quad (3.1г)$$

$$K_\Delta \in L_1(\mathbb{R}), \quad (3.1д)$$

де $K_\Delta(t) = \mathrm{E} X_\Delta(s+t) X_\Delta(s) = \int_{-\infty}^{\infty} e^{i\lambda t} f_\Delta(\lambda) d\lambda, \lambda \in \mathbb{R}$, - кореляційна функція процесу $X_\Delta$.

Далі припускається, що умови (3.1а) - (3.1д) завжди виконані. Умова (3.1г) показує, в якому саме сенсі слід розуміти "близькість" сім'ї процесів $X_\Delta, \Delta > 0$, до гауссівського білого шуму при $\Delta \to \infty$.

Згідно з умовами (3.1б) та (3.1в) маємо, що $f_\Delta \in L_2(\mathbb{R})$, тому з теореми

Фур'є - Планшереля (див. додаток А.1) випливає, що $K_\Delta \in L_2(\mathbb{R})$. Крім того, з умови (3.1в) маємо, що $K_\Delta(t), t \in \mathbb{R}$, - неперервна функція, тому процес $X_\Delta$ є неперервним у середньому квадратичному.

**Вигляд оцінки для перехідної функції.** Нехай $H \in L_2(\mathbb{R})$, - перехідна функція однорідної лінійної системи з шумом (див. додаток А.4). Відгук такої системи на вхідний процес $X_\Delta$ описується випадковим процесом

$$Y_\Delta(t) = \int\limits_{-\infty}^{\infty} H(s) X_\Delta(t-s) ds + U(t), \ t \in \mathbb{R}.$$

Припустимо, що внутрішній шум системи $U = (U(t), t \in \mathbb{R})$, - це вимірний сепарабельний стаціонарний центрований дійснозначний гауссівський процес, який є ортогональним до $X_\Delta$; тобто $\mathrm{E} X_\Delta(s) U(t) = 0, t, s \in \mathbb{R}$. Спектральна щільність $g = (g(\lambda), \lambda \in \mathbb{R})$, процесу $U$ задовольняє умовам:

$$g(\lambda) = g(-\lambda), \lambda \in \mathbb{R}; \tag{3.2а}$$

$$g \in L_1(\mathbb{R}); \tag{3.2б}$$

які далі припускаються завжди виконаними.

Зазначимо, що з умови (3.2б) випливає, що кореляційна функція $K_U(t) = \int\limits_{-\infty}^{\infty} e^{it\lambda} g(\lambda) d\lambda, t \in \mathbb{R}$, процесу $U$ є неперервною, тому процес $U$ - неперервний у середньому квадратичному.

Відповідь про коректність визначення процесу $Y_\Delta$ дає таке твердження.

**Лема 3.1** *Нехай $K_\Delta \in L_1(\mathbb{R})$, $H \in L_2(\mathbb{R})$ та $g \in L_1(\mathbb{R})$, тоді $Y_\Delta$ коректно визначений і є стаціонарним центрованим неперервним у середньому квадратичному гауссівським випадковим процесом зі спектральною щільністю*

$$\phi_\Delta(\lambda) = |H^*(\lambda)|^2 f_\Delta(\lambda) + g(\lambda), \lambda \in \mathbb{R}, \tag{3.3}$$

*де $H^*$ - перетворення Фур'є - Планшереля функції $H$. Крім того, процеси $X_\Delta$ і $Y_\Delta$ є сумісно стаціонарними та сумісно гауссівськими випадковими процесами.*

*Доведення.* Інтеграл, який задає випадковий процес $Y_\Delta$, коректно визначений як відповідна границя у середньому квадратичному, тоді й лише тоді, коли існує вираз

$$\int_{-\infty}^{\infty}\int_{-\infty}^{\infty} K_\Delta(t-s)H(s)H(t)dsdt + K_U(t-s). \qquad (3.4)$$

Існування та обмеженість першого доданка були доведені раніше у лемі 2.1 (див. підрозділ 2.1, розділ 2). Оскільки $K_U(t-s) \leq K_U(0) = \|g\|_1 < \infty$, то другий доданок існує й обмежений. Отже, існує вираз (3.4), а разом з тим й процес $Y_\Delta$. Інші властивості процесу $Y_\Delta$ мають місце в силу його означення (див. додаток А.4). Таким чином, лему 3.1 доведено. □

*Зауваження 3.1.* Вираз (3.4) існує, наприклад, якщо система з внутрішнім шумом є стійкою, тобто $H \in L_1(\mathbb{R})$. При цьому умову (3.1д) можна відкинути.

**Означення 3.1.** Оцінку для $H$ будемо шукати у виді інтегральної сумісної корелограми

$$H_{T,\Delta}(\tau) = \frac{1}{cT}\int_0^T Y_\Delta(t+\tau)X_\Delta(t)dt, \tau \in \mathbb{R}, \qquad (3.5)$$

де $c$ - стала з умови (3.1г), $\Delta$ - параметр схеми серій, й $T$ - довжина інтервалу усереднення $[0,T]$. Для побудови оцінки процеси $X_\Delta$ та $Y_\Delta$ мають спостерігатись на всій дійсній осі. Інтеграл, що фігурує у визначенні оцінки $H_{T,\Delta}$, слід розуміти як середньоквадратичний інтеграл Рімана.

З формули (3.5), для довільних $T>0, \Delta>0,$ при всіх $\tau \in \mathbb{R},$ маємо

$$\mathrm{E} H_{T,\Delta}(\tau) = \frac{1}{c} \mathrm{E} Y_\Delta(t+\tau) X_\Delta(t) = \frac{1}{c} \int_{-\infty}^{\infty} K_\Delta(\tau-s) H(s) ds. \qquad (3.6)$$

Взагалі кажучи, $\mathrm{E} H_{T,\Delta}(\tau) \neq H(\tau)$, тобто оцінка (3.5) є зсуненою.

**Вигляд похибки оцінювання.** Як і в розділі 2, для подальших досліджень додатково розглянемо *нормовану похибку корелограмного оцінювання* виду

$$W_{T,\Delta}(\tau) = \sqrt{T}[H_{T,\Delta}(\tau) - H(\tau)], \tau \in \mathbb{R}, \qquad (3.7)$$

та поставимо питання про її асимптотичні властивості при прямуванні параметрів $T, \Delta$ до безмежності. Для цього процес $W_{T,\Delta} = (W_{T,\Delta}(\tau), \tau \in \mathbb{R})$ розщеплюється на суму

$$W_{T,\Delta} = Z_{T,\Delta} + V_{T,\Delta},$$

де доданки справа визначаються далі, й для них спрацьовують такі факти:

- Вибір послідовності збурюючих процесів $(X_\Delta, \Delta > 0)$, апроксимуючої гауссівський білий шум, - є специфічним: *умови (3.1а) - (3.1д) є достатніми для нормалізації процесу $Z_{T,\Delta}$ при довільній збіжності параметрів $T, \Delta$ до безмежності.*

- Знищення вкладу невипадкової функції $V_{T,\Delta}$ у виразі для $W_{T,\Delta}$ можливе за рахунок використання $\delta$-видності сім'ї $K_\Delta$, підбору характеру сумісного прямування параметрів $T, \Delta$ до безмежності, та накладення умов на порядок локальної гладкості $H$.

## 3.2 Асимптотична поведінка процесу $Z_{T,\Delta}$

Даний підрозділ присвячений дослідженню асимптотичних властивостей нормованої оцінки $H_{T,\Delta}$, центрованої своїм середнім значенням.

Основні результати цього підрозділу встановлюються аналогічним чином з підрозділом 2.2, розділ 2.

### 3.2.1 Кореляційна функція процесу $Z_{T,\Delta}$

У цьому пункті встановлюється вид кореляційної функції емпіричного процесу $Z_{T,\Delta}$, пов'язаного з оцінкою перехідної функції. Досліджуються властивості цієї кореляційної функції. Нагадаємо, що оцінка $H_{T,\Delta}$ та функції $f_\Delta$ й $g$ визначені у підрозділі 3.1.

Покладемо

$$Z_{T,\Delta}(\tau) = \sqrt{T}[H_{T,\Delta}(\tau) - \mathrm{E}H_{T,\Delta}(\tau)], \tau \in \mathbb{R}, \tag{3.8}$$

і зазначимо, що визначений процес є центрованим $L_2$-процесом.

**Лема 3.2** *Нехай $H \in L_2(\mathbb{R})$ та $g \in L_1(\mathbb{R})$, тоді для всіх $T > 0, \Delta > 0$, та $\tau_1, \tau_2 \in \mathbb{R}$, має місце рівність*

$$\mathrm{E}Z_{T,\Delta}(\tau_1)Z_{T,\Delta}(\tau_2) = C_{T,\Delta}(\tau_1,\tau_2) =$$
$$= \frac{2\pi}{c^2} \int_{-\infty}^{\infty}\int_{-\infty}^{\infty} [e^{i(\tau_1-\tau_2)\lambda_2}(|H^*(\lambda_2)|^2 f_\Delta(\lambda_2) + g(\lambda_2)) + \tag{3.9}$$
$$+ e^{i(\tau_1\lambda_1+\tau_2\lambda_2)}H^*(\lambda_1)H^*(\lambda_2)f_\Delta(\lambda_2)]\Phi_T(\lambda_2-\lambda_1)f_\Delta(\lambda_1)d\lambda_1 d\lambda_2;$$

*де $\Phi_T$ - ядро Фейєра:*

$$\Phi_T(\lambda) = \frac{1}{2\pi T}\left(\frac{\sin(T\lambda/2)}{\lambda/2}\right)^2, \lambda \in \mathbb{R};$$

$c$ - стала з умови (3.1г), та $H^*$ - перетворення Фур'є-Планшереля функції $H$ у просторі $L_2(\mathbb{R})$.

*Доведення.* Оскільки кумулянт другого порядку - лінійна функція своїх аргументів, інваріантна відносно зсувів [7], то

$$\mathrm{E}Z_{T,\Delta}(\tau_1)Z_{T,\Delta}(\tau_2) = (\sqrt{T})^2 \mathrm{cov}(H_{T,\Delta}(\tau_1)H_{T,\Delta}(\tau_2)) =$$

$$= T \cdot (\frac{1}{cT})^2 \mathrm{cov}(\int_0^T Y_\Delta(u+\tau_1)X_\Delta(u)du, \int_0^T Y_\Delta(v+\tau_2)X_\Delta(v)dv) =$$

$$= \frac{1}{c^2 T} \int_0^T \int_0^T \mathrm{cov}(Y_\Delta(u+\tau_1)X_\Delta(u), Y_\Delta(v+\tau_2)X_\Delta(v))dudv. \qquad (3.10)$$

Для підінтегральної функції - кумулянта від набору, що є добутками центрованих випадкових величин, - застосуємо теорему Леонова - Ширяєва – Бриллінджера (див. теорему А.5) та лему 3.1:

$$\mathrm{cov}(Y_\Delta(u+\tau_1)X_\Delta(u), Y_\Delta(v+\tau_2)X_\Delta(v)) =$$

$$= \mathrm{cov}(Y_\Delta(u+\tau_1), Y_\Delta(v+\tau_2)) \cdot \mathrm{cov}(X_\Delta(v), X_\Delta(u)) +$$

$$+ \mathrm{cov}(Y_\Delta(u+\tau_1), X_\Delta(v)) \cdot \mathrm{cov}(Y_\Delta(v+\tau_2), X_\Delta(u)) =$$

$$= K_{Y_\Delta}(u-v+\tau_1-\tau_2) \cdot K_{X_\Delta}(v-u) + K_{Y_\Delta,X_\Delta}(u-v+\tau_1) \cdot K_{Y_\Delta,X_\Delta}(v-u+\tau_2),$$

де $K_{Y_\Delta}$ та $K_{X_\Delta}$ - відповідно кореляційні функції процесів $Y_\Delta$ та $X_\Delta$, а $K_{Y_\Delta,X_\Delta}$ - сумісна кореляційна функція процесів $Y_\Delta$ та $X_\Delta$.

Після заміни змінних та деяких обчислень, з формули (3.10) дістанемо:

$$\mathrm{E}Z_{T,\Delta}(\tau_1)Z_{T,\Delta}(\tau_2) = \qquad (3.11)$$

$$= \frac{1}{c^2} \int_{-T}^{T} [K_{Y_\Delta}(t+\tau_1-\tau_2) \cdot K_{X_\Delta}(-t) + K_{Y_\Delta,X_\Delta}(t+\tau_1) \cdot K_{Y_\Delta,X_\Delta}(-t+\tau_2)](1-\frac{|t|}{T})dt.$$

Підставивши наступні спектральні представлення:

$$K_{Y_\Delta}(t) = \int_{-\infty}^{\infty} e^{it\lambda}(|H^*(\lambda)|^2 f_\Delta(\lambda) + g(\lambda))d\lambda;$$

$$K_{X_\Delta}(t) = \int_{-\infty}^{\infty} e^{it\lambda} f_\Delta(\lambda)d\lambda;$$

$$K_{Y_\Delta, X_\Delta}(t) = \int_{-\infty}^{\infty} e^{it\lambda} H^*(\lambda) f_\Delta(\lambda)d\lambda,$$

у формулу (3.11), отримаємо вираз:

$$\mathrm{E} Z_{T,\Delta}(\tau_1) Z_{T,\Delta}(\tau_2) = \frac{1}{c^2} \int_{-T}^{T} (\int_{-\infty}^{\infty} e^{i(-t)\lambda_1} f_\Delta(\lambda_1) d\lambda_1) \times$$

$$\times (\int_{-\infty}^{\infty} e^{i(t+\tau_1-\tau_2)\lambda_2}(|H^*(\lambda_2)|^2 f_\Delta(\lambda_2) + g(\lambda_2))d\lambda_2)(1 - \frac{|t|}{T})dt +$$

$$+ \frac{1}{c^2} \int_{-T}^{T} (\int_{-\infty}^{\infty} e^{i(t+\tau_1)\lambda_1} H^*(\lambda_1) f_\Delta(\lambda_1) d\lambda_1) \cdot (\int_{-\infty}^{\infty} e^{i(-t+\tau_2)\lambda_2} H^*(\lambda_2) f_\Delta(\lambda_2) d\lambda_2) \times$$

$$\times (1 - \frac{|t|}{T})dt = \frac{1}{c^2} \int_{-\infty}^{\infty}\int_{-\infty}^{\infty} e^{i(\tau_1-\tau_2)\lambda_2}(|H^*(\lambda_2)|^2 f_\Delta(\lambda_2) + g(\lambda_2)) f_\Delta(\lambda_1) \times$$

$$\times [\int_{-T}^{T} e^{it(\lambda_2-\lambda_1)}(1 - \frac{|t|}{T})dt] d\lambda_1 d\lambda_2 +$$

$$+ \frac{1}{c^2} \int_{-\infty}^{\infty}\int_{-\infty}^{\infty} e^{i(\tau_1\lambda_1+\tau_2\lambda_2)} H^*(\lambda_1) H^*(\lambda_2) f_\Delta(\lambda_1) f_\Delta(\lambda_2) [\int_{-T}^{T} e^{it(\lambda_1-\lambda_2)}(1 - \frac{|t|}{T})dt] d\lambda_1 d\lambda_2.$$

Оскільки ядро Фейєра

$$\Phi_T(\lambda) = \frac{1}{2\pi} \int_{-T}^{T} e^{it\lambda}(1 - \frac{|t|}{T})dt = \frac{1}{2\pi T}\left(\frac{\sin(T\lambda/2)}{\lambda/2}\right)^2, \lambda \in \mathbb{R},$$

є дійснозначною парною функцією (див. зауваження 2.3, розділ 2), то в результаті дістанемо

$$EZ_{T,\Delta}(\tau_1)Z_{T,\Delta}(\tau_2) = \frac{2\pi}{c^2}\int_{-\infty}^{\infty}\int_{-\infty}^{\infty}[e^{i(\tau_1-\tau_2)\lambda_2}(|H^*(\lambda_2)|^2 f_\Delta(\lambda_2) + g(\lambda_2)) +$$

$$+e^{i(\tau_1\lambda_1+\tau_2\lambda_2)}H^*(\lambda_1)H^*(\lambda_2)f_\Delta(\lambda_2)]\Phi_T(\lambda_2-\lambda_1)f_\Delta(\lambda_1)d\lambda_1 d\lambda_2.$$

Таким чином, лему 3.2 доведено. □

Зазначимо, що кореляційна функція процесу $Z_{T,\Delta}$ є невід'ємно визначеною неперервною функцією на $\mathbb{R}\times\mathbb{R}$.

*Зауваження 3.2.* Ядра Фейєра $\Phi_T$, що зустрічаються у формулі (3.9), залежать лише від одного параметра усереднення $T>0$, і не залежать від параметра схеми серій $\Delta>0$.

**Лема 3.3** *Нехай $H \in L_2(\mathbb{R})$ та $g \in L_1(\mathbb{R})$, тоді для всіх $T>0, \Delta>0$, та $\tau_1,\tau_2 \in \mathbb{R}$, має місце нерівність*

$$|EZ_{T,\Delta}(\tau_1)Z_{T,\Delta}(\tau_2)| \leq \frac{2\pi}{c^2}(\sup_{\Delta>0}\|f_\Delta\|_\infty)[2(\sup_{\Delta>0}\|f_\Delta\|_\infty)\|H^*\|_2^2 + \|g\|_1].$$

*Доведення.* Будемо оцінювати функцію, визначену формулою (3.9), розбивши її на три доданки наступним чином:

$$|EZ_{T,\Delta}(\tau_1)Z_{T,\Delta}(\tau_2)|=$$

$$=|\frac{2\pi}{c^2}\int_{-\infty}^{\infty}\int_{-\infty}^{\infty}e^{i(\tau_1-\tau_2)\lambda_2}|H^*(\lambda_2)|^2\Phi_T(\lambda_2-\lambda_1)f_\Delta(\lambda_1)f_\Delta(\lambda_2)d\lambda_1 d\lambda_2 +$$

$$+\frac{2\pi}{c^2}\int_{-\infty}^{\infty}\int_{-\infty}^{\infty}e^{i(\tau_1-\tau_2)\lambda_2}g(\lambda_2)\Phi_T(\lambda_2-\lambda_1)f_\Delta(\lambda_1)d\lambda_1 d\lambda_2 +$$

$$+\frac{2\pi}{c^2}\int_{-\infty}^{\infty}\int_{-\infty}^{\infty}e^{i(\tau_1\lambda_1+\tau_2\lambda_2)}H^*(\lambda_1)H^*(\lambda_2)\Phi_T(\lambda_2-\lambda_1)f_\Delta(\lambda_1)f_\Delta(\lambda_2)d\lambda_1 d\lambda_2|.$$

В силу умови (3.1б), звідси отримаємо

$$|EZ_{T,\Delta}(\tau_1)Z_{T,\Delta}(\tau_2)| \leq \frac{2\pi}{c^2}(\sup_{\Delta>0}\|f_\Delta\|_\infty)^2 \int_{-\infty}^{\infty}\int_{-\infty}^{\infty}|H^*(\lambda_2)|^2\Phi_T(\lambda_2-\lambda_1)d\lambda_1 d\lambda_2 +$$

$$+\frac{2\pi}{c^2}(\sup_{\Delta>0}\|f_\Delta\|_\infty)\int_{-\infty}^{\infty}\int_{-\infty}^{\infty}g(\lambda_2)\Phi_T(\lambda_2-\lambda_1)d\lambda_1 d\lambda_2 +$$

$$+\frac{2\pi}{c^2}(\sup_{\Delta>0}\|f_\Delta\|_\infty)^2 \int_{-\infty}^{\infty}\int_{-\infty}^{\infty}|H^*(\lambda_1)||H^*(\lambda_2)|\Phi_T(\lambda_2-\lambda_1)d\lambda_1 d\lambda_2 =$$

$$=\frac{2\pi}{c^2}(\sup_{\Delta>0}\|f_\Delta\|_\infty)^2\|H^*\|_2^2\|\Phi_T\|_1 + \frac{2\pi}{c^2}(\sup_{\Delta>0}\|f_\Delta\|_\infty)\|g\|_1\|\Phi_T\|_1 +$$

$$+\frac{2\pi}{c^2}(\sup_{\Delta>0}\|f_\Delta\|_\infty)^2 \int_{-\infty}^{\infty}|H^*(\lambda_2)|(|H^*|*\Phi_T)(\lambda_2)d\lambda_2.$$

З нерівності Юнга (див. теорему А.2), застосованої до третього доданка, й того, що $\|\Phi_T\|_1 = 1$, в результаті дістанемо

$$|EZ_{T,\Delta}(\tau_1)Z_{T,\Delta}(\tau_2)| \leq \frac{2\pi}{c^2}(\sup_{\Delta>0}\|f_\Delta\|_\infty)^2\|H^*\|_2^2 + \frac{2\pi}{c^2}(\sup_{\Delta>0}\|f_\Delta\|_\infty)\|g\|_1 +$$

$$+\frac{2\pi}{c^2}(\sup_{\Delta>0}\|f_\Delta\|_\infty)^2\|H^*\|_2^2\|\Phi_T\|_1 =$$

$$=\frac{4\pi}{c^2}(\sup_{\Delta>0}\|f_\Delta\|_\infty)^2\|H^*\|_2^2 + \frac{2\pi}{c^2}(\sup_{\Delta>0}\|f_\Delta\|_\infty)\|g\|_1.$$

Таким чином, лему 3.3 доведено. □

**Наслідок 3.1** *Оцінка кореляційної функції процесу $Z_{T,\Delta}$, наведена у лемі 3.3, не залежить від параметрів $T>0, \Delta>0$, та значень $\tau_1, \tau_2 \in \mathbb{R}$. Тобто при $H \in L_2(\mathbb{R})$ та $g \in L_1(\mathbb{R})$, насправді встановлено співвідношення:*

$$\sup_{T,\Delta>0}\sup_{\tau_1,\tau_2\in\mathbb{R}}|EZ_{T,\Delta}(\tau_1)Z_{T,\Delta}(\tau_2)| \leq \frac{2\pi}{c^2}(\sup_{\Delta>0}\|f_\Delta\|_\infty)[2(\sup_{\Delta>0}\|f_\Delta\|_\infty)\|H^*\|_2^2 + \|g\|_1].$$

**Наслідок 3.2** *Нехай виконуються умови леми 3.3, тоді для всіх $T > 0, \Delta > 0$, та $\tau \in \mathbb{R}$, має місце нерівність*

$$\mathrm{Var}(H_{T,\Delta}(\tau)) = E\,|\,H_{T,\Delta}(\tau) - \mathrm{E}H_{T,\Delta}(\tau)\,|^2 \leq$$

$$\leq \frac{2\pi}{c^2 T}(\sup_{\Delta>0}\|\,f_\Delta\,\|_\infty)[2(\sup_{\Delta>0}\|\,f_\Delta\,\|_\infty)\|\,H^*\,\|_2^2 + \|\,g\,\|_1].$$

### 3.2.2 Гранична поведінка кореляційної функції процесу $Z_{T,\Delta}$

У пункті 3.2.1 був знайдений вигляд кореляційної функції $C_{T,\Delta}$ процесу $Z_{T,\Delta}$ (див. формулу (3.9)). У даному пункті знаходиться границя цієї функції при довільному прямуванні параметрів $T$ і $\Delta$ до нескінченності.

Далі запис $(T,\Delta) \to \infty$, означає, що одночасно $T \to \infty$ та $\Delta \to \infty$.

Покладемо для всіх $\tau_1, \tau_2 \in \mathbb{R}$

$$C_\infty(\tau_1,\tau_2) = \frac{1}{2\pi}\int_{-\infty}^{\infty}\left[e^{i(\tau_1-\tau_2)\lambda}(|H^*(\lambda)|^2 + \frac{2\pi}{c}g(\lambda)) + e^{i(\tau_1+\tau_2)\lambda}(H^*(\lambda))^2\right]d\lambda. \quad (3.12)$$

Зазначимо, що функція $C_\infty = (C_\infty(\tau_1,\tau_2), \tau_1,\tau_2 \in \mathbb{R})$, коректно визначена та є неперервною функцією за умови $H \in L_2(\mathbb{R})$ і $g \in L_1(\mathbb{R})$. Цю функцію можна записати у наступному вигляді:

$$C_\infty(\tau_1,\tau_2) = \frac{1}{\pi}\int_0^\infty[\cos(\tau_1-\tau_2)\lambda \cdot (|H^*(\lambda)|^2 + \frac{2\pi}{c}g(\lambda)) +$$

$$+\cos(\tau_1+\tau_2)\lambda \cdot Re(H^*(\lambda))^2 - \sin(\tau_1+\tau_2)\lambda \cdot Im(H^*(\lambda))^2]d\lambda.$$

звідки видно, що вона дійснозначна.

*Зауваження 3.3.* Розглянемо деякі частинні випадки функцій $C_\infty$:

(I) Нехай $H$ - парна дійснозначна імпульсна перехідна функція. Тоді з формули (3.12), після деяких перетворень, дістанемо

$$C_{\infty}(\tau_1,\tau_2) = \frac{2}{\pi}\int_0^{\infty}(|H^*(\lambda)|^2 + \frac{2\pi}{c}g(\lambda))\cos(\tau_1\lambda)\cos(\tau_2\lambda)d\lambda.$$

(II) Нехай $H$ - непарна дійснозначна імпульсна перехідна функція. Тоді з формули (3.12), після деяких перетворень, дістанемо

$$C_{\infty}(\tau_1,\tau_2) = \frac{2}{\pi}\int_0^{\infty}(|H^*(\lambda)|^2\sin(\tau_1\lambda)\sin(\tau_2\lambda) + \frac{2\pi}{c}g(\lambda)\cos(\tau_1\lambda)\cos(\tau_2\lambda))d\lambda.$$

**Теорема 3.1.** *Нехай $H \in L_2(\mathbb{R})$ та $g \in L_1(\mathbb{R})$, тоді для всіх $\tau_1,\tau_2 \in \mathbb{R}$, має місце рівність*

$$\lim_{(T,\Delta)\to\infty} \mathrm{E}Z_{T,\Delta}(\tau_1)Z_{T,\Delta}(\tau_2) = C_{\infty}(\tau_1,\tau_2).$$

*Доведення.* Для простоти викладок, розіб'ємо доведення на три етапи. На першому етапі ми покажемо, що

$$\lim_{(T,\Delta)\to\infty} C^{(1)}_{T,\Delta}(\tau_1,\tau_2) = C^{(1)}_{\infty}(\tau_1,\tau_2), \qquad (3.13)$$

де

$$C^{(1)}_{T,\Delta}(\tau_1,\tau_2) = \frac{2\pi}{c^2}\int_{-\infty}^{\infty}\int_{-\infty}^{\infty}e^{i(\tau_1-\tau_2)\lambda_2}|H^*(\lambda_2)|^2 f_{\Delta}(\lambda_1)f_{\Delta}(\lambda_2)\Phi_T(\lambda_2-\lambda_1)d\lambda_1d\lambda_2;$$

$$C^{(1)}_{\infty}(\tau_1,\tau_2) = \frac{1}{2\pi}\int_{-\infty}^{\infty}e^{i(\tau_1-\tau_2)\lambda}|H^*(\lambda)|^2d\lambda, \ \tau_1,\tau_2 \in \mathbb{R}.$$

На другому етапі отримаємо наступну рівність:

$$\lim_{(T,\Delta)\to\infty} C^{(2)}_{T,\Delta}(\tau_1,\tau_2) = C^{(2)}_{\infty}(\tau_1,\tau_2), \qquad (3.14)$$

де

$$C^{(2)}_{T,\Delta}(\tau_1,\tau_2) = \frac{2\pi}{c^2}\int_{-\infty}^{\infty}\int_{-\infty}^{\infty}e^{i(\tau_1\lambda_1+\tau_2\lambda_2)}H^*(\lambda_1)H^*(\lambda_2)f_{\Delta}(\lambda_1)f_{\Delta}(\lambda_2)\Phi_T(\lambda_2-\lambda_1)d\lambda_1d\lambda_2;$$

$$C^{(2)}_{\infty}(\tau_1,\tau_2) = \frac{1}{2\pi}\int_{-\infty}^{\infty}e^{i(\tau_1+\tau_2)\lambda}(H^*(\lambda))^2d\lambda, \ \tau_1,\tau_2 \in \mathbb{R}.$$

На третьому етапі ми доведемо, що

$$\lim_{(T,\Delta)\to\infty} C^{(3)}_{T,\Delta}(\tau_1,\tau_2) = C^{(3)}_{\infty}(\tau_1,\tau_2), \qquad (3.15)$$

де

$$C^{(3)}_{T,\Delta}(\tau_1,\tau_2) = \frac{2\pi}{c^2}\int_{-\infty}^{\infty}\int_{-\infty}^{\infty} e^{i(\tau_1-\tau_2)\lambda_2} f_\Delta(\lambda_1) g(\lambda_2) \Phi_T(\lambda_2-\lambda_1) d\lambda_1 d\lambda_2;$$

$$C^{(3)}_{\infty}(\tau_1,\tau_2) = \frac{1}{c}\int_{-\infty}^{\infty} e^{i(\tau_1-\tau_2)\lambda} g(\lambda) d\lambda, \ \ \tau_1,\tau_2 \in \mathbb{R}.$$

- *Етапи 1 - 2*. Формули (3.13) та (3.14) доводяться так же, як і формули (2.12) та (2.13) відповідно у теоремі 2.1, розділ 2.

- *Етап 3*. Для будь-яких $T>0, \Delta>0$ та $\tau_1,\tau_2 \in \mathbb{R}$, розглянемо різницю

$$C^{(3)}_{\infty}(\tau_1,\tau_2) - C^{(3)}_{T,\Delta}(\tau_1,\tau_2) =$$

$$= \frac{2\pi}{c}\int_{-\infty}^{\infty}\int_{-\infty}^{\infty} e^{i(\tau_1-\tau_2)\lambda_2}\left[\frac{c}{2\pi}-f_\Delta(\lambda_1)\right] g(\lambda_2)\Phi_T(\lambda_2-\lambda_1) d\lambda_1 d\lambda_2,$$

де у представленні використана парність ядер Фейєра, а також $\|\Phi_T\|_1 = 1$ (див. зауваження 2.3, розділ 2).

Визначимо множину $P(b,b) = \{(\lambda_1,\lambda_2)\in\mathbb{R}^2 : |\lambda_2-\lambda_1|\leq b\}$ для довільного $b>0$. Оскільки

$$J_1(b,T,\Delta) = \left|\int\int_{P(b,b)} e^{i(\tau_1-\tau_2)\lambda_2}\left[\frac{c}{2\pi}-f_\Delta(\lambda_1)\right] g(\lambda_2)\Phi_T(\lambda_2-\lambda_1) d\lambda_1 d\lambda_2\right| \leq$$

$$\leq \sup_{-b\leq\lambda\leq b}\left|\frac{c}{2\pi}-f_\Delta(\lambda)\right| \|g\|_1 \|\Phi_T\|_1 = \sup_{-b\leq\lambda\leq b}\left|\frac{c}{2\pi}-f_\Delta(\lambda)\right| \|g\|_1,$$

то з умов (3.1г) та (3.2б) отримаємо, що для будь-якого $b>0$

$$\limsup_{(T,\Delta)\to\infty} J_1(b,T,\Delta) = 0. \qquad (3.16)$$

Для будь-яких $T>0, \Delta>0$ та $\tau_1,\tau_2 \in \mathbb{R}$, і зафіксованого довільного

$b > 0,$ розглянемо нерівність

$$J_2(b,T,\Delta) = \left|\int\int_{\mathbb{R}^2\setminus P(b,b)} e^{i(\tau_1-\tau_2)\lambda_2}\left[\frac{c}{2\pi} - f_\Delta(\lambda_1)\right]g(\lambda_2)\Phi_T(\lambda_2-\lambda_1)d\lambda_1 d\lambda_2\right| \leq$$

$$\leq \left[\frac{c}{2\pi} + \sup_{\Delta>0}\|f_\Delta\|_\infty\right]\|g\|_1 \int_{|\lambda|>b}\Phi_T(\lambda)d\lambda.$$

В силу умов (3.1г) та (3.2б), оскільки для довільного $b > 0$, $\lim_{T\to\infty}\int_{|\lambda|>b}\Phi_T(\lambda)d\lambda = 0$

, то звідси дістанемо

$$\limsup_{(T,\Delta)\to\infty}J_2(b,T,\Delta) = \limsup_{b\to\infty}(\limsup_{(T,\Delta)\to\infty}J_2(b,T,\Delta)) = 0. \tag{3.17}$$

З формул (3.16) і (3.17), отримаємо для всіх $\tau_1, \tau_2 \in \mathbb{R}$

$$\limsup_{(T,\Delta)\to\infty}|C_\infty^{(3)}(\tau_1,\tau_2) - C_{T,\Delta}^{(3)}(\tau_1,\tau_2)| \leq$$

$$\leq \frac{2\pi}{c^2}[\limsup_{b\to\infty}(\limsup_{(T,\Delta)\to\infty}J_1(b,T,\Delta)) + \limsup_{b\to\infty}(\limsup_{(T,\Delta)\to\infty}J_1(b,T,\Delta))] = 0.$$

Таким чином, співвідношення (3.15) доведене.

В результаті, об'єднуючи всі три етапи, отримаємо при всіх $\tau_1, \tau_2 \in \mathbb{R}$

$$\limsup_{(T,\Delta)\to\infty}C_{T,\Delta}(\tau_1,\tau_2) = \sum_{j=1}^{3}\limsup_{(T,\Delta)\to\infty}C_{T,\Delta}^{(j)}(\tau_1,\tau_2) = \sum_{j=1}^{3}C_\infty^{(j)}(\tau_1,\tau_2) = C_\infty(\tau_1,\tau_2).$$

Таким чином, теорему 3.1 доведено. □

### 3.2.3 Асимптотична нормальність скінченновимірних розподілів процесу $Z_{T,\Delta}$

У цьому пункті досліджується асимптотична нормальність скінченновимірних розподілів процесу $Z_{T,\Delta}$ (див. (3.8)). Зазначимо, що ніяких додаткових умов на перехідну функцію $H$ та характер прямування параметрів $T, \Delta$ до безмежності не накладається.

Теорема 3.1 демонструє, що функція $C_\infty$, визначена формулою (3.12), є невід'ємно визначеною на $\mathbb{R} \times \mathbb{R}$. Таким чином, завжди можна стверджувати, що існує центрований дійснозначний гауссівський процес $Z = (Z(\tau), \tau \in \mathbb{R})$ з кореляційною функцією $C_\infty$, тобто

$$\mathrm{E} Z(\tau_1) Z(\tau_2) = C_\infty(\tau_1, \tau_2) = \qquad (3.18)$$

$$= \frac{1}{2\pi} \int_{-\infty}^{\infty} \left[ e^{i(\tau_1 - \tau_2)\lambda}(|H^*(\lambda)|^2 + \frac{2\pi}{c} g(\lambda)) + e^{i(\tau_1 + \tau_2)\lambda}(H^*(\lambda))^2 \right] d\lambda, \ \tau_1, \tau_2 \in \mathbb{R}.$$

Без втрати загальності, припускається, що процес $Z$ заданий на тому ж базовому ймовірнісному просторі $(\Omega, \mathfrak{F}, \mathrm{P})$, що й процес $Z_{T,\Delta}$.

Має місце наступне твердження про асимптотичну нормальність $Z_{T,\Delta}$.

**Теорема 3.2.** *Нехай $H \in L_2(\mathbb{R})$ та $g \in L_1(\mathbb{R})$, тоді для всіх $m \in \mathbb{N}$ та $\tau_1, \tau_2, ..., \tau_m \in \mathbb{R}$, має місце рівність*

$$\lim_{(T,\Delta) \to \infty} cum(Z_{T,\Delta}(\tau_j), j=1,...,m) = \begin{cases} 0, & m = 1; \\ C_\infty(\tau_1, \tau_2), & m = 2; \\ 0, & m \geq 3, \end{cases} \qquad (3.19)$$

*де $cum(Z_{T,\Delta}(\tau_j), j = 1,...,m)$ - сумісний кумулянт набору випадкових величин $Z_{T,\Delta}(\tau_j), j = 1,...,m$.*

*Зокрема, всі скінченновимірні розподіли процесу $(Z_{T,\Delta}(\tau), \tau \in \mathbb{R})$ слабко збігаються до відповідних скінченновимірних розподілів центрованого гауссівського процесу $(Z(\tau), \tau \in \mathbb{R})$.*

*Доведення.* Доведення аналогічне до доведення теореми 2.2, розділ 2.

Нехай $m \in \{1,2\}$. Формула (3.19) має місце, оскільки процес $Z_{T,\Delta}$ центрований, та задовольняються умови теореми 3.1.

Нехай $m \geq 3$. Доведення формули (3.19) проведемо у два етапи. На першому етапі покажемо, що $cum(Z_{T,\Delta}(\tau_j), j=1,...,m)$ можна подати у виді скінченної суми інтегралів з циклічним зачепленням ядер. Другий етап присвячений перевірці умов теореми А.3, згідно з якою кожен доданок у сумі прямує до нуля при $(T,\Delta) \to \infty$. Як і раніше, для встановлення цих фактів будемо використовувати діаграмно-кумулянтний метод Д. Бриллінджера [7] та підходи робіт В. Булдигіна, Ф. Уцета та В. Зайця [61, 62].

- *Етап 1.* З властивостей лінійності кумулянтів та їх інваріантності відносно зсувів, а також з означень процесу $Z_{T,\Delta}$ та оцінки $H_{T,\Delta}$ (див. (3.8) і (3.5) відповідно) випливає, що

$$cum(Z_{T,\Delta}(\tau_j), j=1,...,m) = \qquad (3.20)$$

$$= \left(\frac{1}{c^2 T}\right)^{\frac{m}{2}} \int_0^T ... \int_0^T cum(Y_\Delta(t_j + \tau_j) X_\Delta(t_j), j=1,...,m) dt_1 ... dt_m.$$

Оскільки $X_\Delta$ та $Y_\Delta$ - центровані сумісно гауссівські процеси, то з теореми А.5, застосованої до кумулянта, що стоїть під знаком інтеграла у формулі (3.20), дістанемо

$$cum(Y_\Delta(t_j + \tau_j) X_\Delta(t_j), j=1,...,m) = \sum \prod_{p=1}^m cum(D_p^{(2)}), \qquad (3.21)$$

де сума береться по всіх можливих невпорядкованих нерозкладних розбиттях двостовпцевої таблиці

$$D_{m\times 2} = \begin{bmatrix} Y_\Delta(t_1 + \tau_1) & X_\Delta(t_1) \\ Y_\Delta(t_2 + \tau_2) & X_\Delta(t_2) \\ \vdots & \vdots \\ Y_\Delta(t_m + \tau_m) & X_\Delta(t_m) \end{bmatrix}$$

на пари $\{D_1^{(2)},...,D_m^{(2)}\}$.

Так як порядок елементів у $\{D_1^{(2)},...,D_m^{(2)}\}$ не має значення, можна завжди припускати, що розбиття задовольняє умовам:

(1)   $D_p^{(2)} \bigcap D_q^{(2)} = \varnothing$ для $p \neq q$;

(2)   $D_{m\times 2} = D_1^{(2)} \bigcup ... \bigcup D_m^{(2)}$;

(3)   якщо $m \geq 3$, тоді для будь-якого $p = 1,...,m$ множина $D_p^{(2)}$ не співпадає з жодним із рядків $r_1,...,r_m$ таблиці $D_{m\times 2}$; більш того, для будь-якого $1 \leq \nu < m$, об'єднання довільної кількості елементів розбиття $\{D_1^{(2)},...,D_m^{(2)}\}$ не повинно співпадати з об'єднанням будь-яких $\nu$ рядків таблиці $D_{m\times 2}$;

(4)   для будь-якого $p = 1,...,m-1$, існує єдиний рядок $\tilde{r}_p$ таблиці $D_{m\times 2}$ такий, що $\tilde{r}_p \subset D_p^{(2)} \bigcup D_{p+1}^{(2)}$;

(5)   $r_1 \subset D_1^{(2)} \bigcup D_m^{(2)}$.

З умов (1)-(5) випливає, що елементи розбиття $\{D_1^{(2)},...,D_m^{(2)}\}$ послідовно зачіплюються один за іншим та "вичерпують" всю таблицю $D_{m\times 2}$; $D_m^{(2)}$ зачіплюється за $D_1^{(2)}$.

Далі, для заданого $D_p^{(2)}$, запишемо $D_p^{(2)} = D_{j,j}^{(2)}$, де $j$ та $j$ - номери тих двох рядків таблиці $D_{m\times 2}$, в яких містяться елементи множини $D_p^{(2)}$. Отже, довільне невпорядковане нерозкладне розбиття $\{D_1^{(2)},...,D_m^{(2)}\}$ можна записати як

$$\vec{D}^{(2)} = \{D^{(2)}_{j_1,j_2}, D^{(2)}_{j_2,j_3},...,D^{(2)}_{j_{m-1},j_m}, D^{(2)}_{j_m,j_{m+1}}\},$$

де $(j_1,...,j_m)$ - впорядкований набір такий, що $j_1 = 1$, $(j_2, j_3,...,j_m)$ - перестановка чисел множини $\{2,3,...,m\}$, та $j_{m+1} = j_1$.

Аналіз структури $\vec{D}^{(2)}$ показує, що можна виділити три різні класи його елементів. До першого класу належать невпорядковані множини виду

$$D^{(2)}_{j,j} = \{X_\Delta(t_j), X_\Delta(t_j)\}, \ j \neq j.$$

Другий клас утворюють невпорядковані множини виду

$$D^{(2)}_{j,j} = \{Y_\Delta(t_j + \tau_j), Y_\Delta(t_j + \tau_j)\}, \ j \neq j.$$

Третій клас містить невпорядковані множини виду

$$D^{(2)}_{j,j} = \{Y_\Delta(t_j + \tau_j), X_\Delta(t_j)\}, \ j \neq j.$$

Відповідно позначимо класи $G_1(\vec{D}^{(2)}), G_2(\vec{D}^{(2)})$ та $G_3(\vec{D}^{(2)})$. Нехай $m_\nu(\vec{D}^{(2)}) = card G_\nu(\vec{D}^{(2)}), \nu = 1,2,3$, тоді для будь-якого $\vec{D}^{(2)}$

$$m_1(\vec{D}^{(2)}) = m_2(\vec{D}^{(2)}); m_1(\vec{D}^{(2)}) + m_2(\vec{D}^{(2)}) + m_3(\vec{D}^{(2)}) = m. \quad (3.22)$$

Тепер знайдемо кумулянти від наборів величин, що належать кожному з трьох класів. Якщо $D^{(2)}_{j,j} \in G_1(\vec{D}^{(2)})$, то

$$cum(D^{(2)}_{j,j}) = K_\Delta(t_j - t_j) = \int_{-\infty}^{\infty} e^{i(t_j - t_j)\lambda_j} f_\Delta(\lambda_j) d\lambda_j;$$

якщо $D^{(2)}_{j,j} \in G_2(\vec{D}^{(2)})$, то

$$cum(D^{(2)}_{j,j}) = EY_\Delta(t_j + \tau_j) Y_\Delta(t_j + \tau_j) =$$
$$= \int_{-\infty}^{\infty} e^{i(t_j - t_j)\lambda_j} \cdot e^{i(\tau_j - \tau_j)\lambda_j} (|H^*(\lambda_j)|^2 f_\Delta(\lambda_j) + g(\lambda_j)) d\lambda_j;$$

якщо $D^{(2)}_{j,j} \in G_3(\vec{D}^{(2)})$, то

$$cum(D_{j,j}^{(2)}) = EY_\Delta(t_j + \tau_j)X_\Delta(t_j) = \int_{-\infty}^{\infty} e^{i(t_j - t_j)\lambda_j} \cdot e^{i\tau_j \lambda_j} H^*(\lambda_j) f_\Delta(\lambda_j) d\lambda_j.$$

З останніх формул випливає представлення для добутку кумулянтів з $\vec{D}^{(2)}$:

$$\prod_{p=1}^{m} cum(D_p^{(2)}) = \int \ldots \int_{\mathbb{R}^m} \left[ \prod_{k=1}^{m} e^{i(t_{j_{k+1}} - t_{j_k})\lambda_{j_k}} \right] \varphi_0(\vec{\lambda}, \vec{\tau}, \vec{D}^{(2)}) \times \qquad (3.23)$$

$$\times \left( \prod_{k=1}^{m} \varphi_{j_k}(\lambda_{j_{k+1}}, \Delta, \vec{D}^{(2)}) \right) d\lambda_{j_1} \ldots d\lambda_{j_m},$$

де $\vec{\lambda} = (\lambda_1, \ldots, \lambda_m)$; $\vec{\tau} = (\tau_1, \ldots, \tau_m)$; $(j_1, \ldots, j_m)$ - набір такий, що $j_{m+1} = j_1 = 1$, $(j_2, \ldots, j_m)$ - перестановка чисел множини $\{2, \ldots, m\}$.

Для кожного розбиття $\vec{D}^{(2)}$ функція $\varphi_0(\cdot, \vec{D}^{(2)})$ є добутком деяких з функцій $e^{i(\tau_{j_{k+1}} - \tau_{j_k})\lambda_{j_k}}, e^{i\tau_{j_k}\lambda_{j_k}}$ та індикаторних функцій $\mathbb{I}_\mathbb{R}(\lambda_{j_k}), k = 1, \ldots, m$. Таким чином, має місце рівність

$$\sup_{\vec{D}^{(2)}} \sup_{\vec{\lambda}, \vec{\tau}} | \varphi_0(\vec{\lambda}, \vec{\tau}, \vec{D}^{(2)}) | = 1. \qquad (3.24)$$

Для будь-яких $\Delta > 0$, та розбиття $\vec{D}^{(2)}$ кожна функція з множини $F(\Delta, \vec{D}^{(2)}) = \{\varphi_1(\cdot, \Delta, \vec{D}^{(2)}), \ldots, \varphi_m(\cdot, \Delta, \vec{D}^{(2)})\}$ може збігатися лише з однією з наступних $f_\Delta, |H^*|^2 f_\Delta + g$ та $H^* f_\Delta$. Крім того,

$$card\{\varphi \in F(\Delta, \vec{D}^{(2)}) : \varphi = f_\Delta\} = m_1(\vec{D}^{(2)});$$

$$card\{\varphi \in F(\Delta, \vec{D}^{(2)}) : \varphi = |H^*|^2 f_\Delta + g\} = m_2(\vec{D}^{(2)}); \qquad (3.25)$$

$$card\{\varphi \in F(\Delta, \vec{D}^{(2)}) : \varphi = |H^*| f_\Delta\} = m_3(\vec{D}^{(2)}).$$

Перетворимо вираз $\prod_{k=1}^{m} e^{i(t_{j_k} - t_{j_{k+1}})\lambda_{j_k}} = \prod_{k=1}^{m} e^{it_{j_{k+1}}(\lambda_{j_{k+1}} - \lambda_{j_k})}$, де $\lambda_{m+1} = \lambda_1$, та розглянемо $m$-кратний інтеграл

$$\int_0^T \cdots \int_0^T \left(\prod_{k=1}^m e^{it_{j_{k+1}}(\lambda_{j_{k+1}} - \lambda_{j_k})}\right) dt_{j_1} \ldots dt_{j_m} = \prod_{k=1}^m \left(\int_0^T e^{it_{j_{k+1}}(\lambda_{j_{k+1}} - \lambda_{j_k})} dt_{j_{k+1}}\right) =$$

$$= \prod_{k=1}^m \Phi_T(\lambda_{j_{k+1}} - \lambda_{j_k}),$$

де $\Phi_T(\lambda) = \int_0^T e^{it\lambda} dt = \dfrac{e^{it\lambda T} - 1}{i\lambda}, \lambda \in \mathbb{R}.$

Враховуючи останні перетворення, з формул (3.20) - (3.23) отримаємо

$$\mathrm{cum}(Z_{T,\Delta}(\tau_j), j = 1, \ldots, m) = \qquad (3.26)$$

$$= \left(\frac{2\pi}{c^2}\right)^{\frac{m}{2}} \sum_{\vec{D}^{(2)}} \int \cdots \int_{\mathbb{R}^m} \left[\prod_{k=1}^m \Phi^{(T)}(\lambda_{k+1} - \lambda_k)\right] \varphi_0(\vec{\lambda}, \vec{\tau}, \vec{D}^{(2)}) \times$$

$$\times \left(\prod_{k=1}^m \varphi_k(\lambda_k, \Delta, \vec{D}^{(2)})\right) d\lambda_1 \ldots d\lambda_m,$$

де $\lambda_{m+1} = \lambda_1$, та

$$\Phi^{(T)}(\lambda) = \left(\frac{1}{2\pi T}\right)^{\frac{1}{2}} \Phi_T(\lambda) = \left(\frac{1}{2\pi T}\right)^{\frac{1}{2}} \frac{e^{iT\lambda} - 1}{i\lambda}, \lambda \in \mathbb{R}. \qquad (3.27)$$

Звідси, виключаючи $\varphi_0$ та залежність $\varphi_k, k = 1, \ldots, m,$ від $\Delta$ (див. (3.24), (3.16)), отримаємо наступну оцінку для кумулянта:

$$|\mathrm{cum}(Z_{T,\Delta}(\tau_j), j = 1, \ldots, m)| = \left(\frac{2\pi}{c^2}\right)^{\frac{m}{2}} (\max\{\sup_{\Delta > 0} \|f_\Delta\|_\infty; 1\})^m \times \qquad (3.28)$$

$$\times \sum_{\vec{D}^{(2)}} \int \cdots \int_{\mathbb{R}^m} |\prod_{k=1}^m \Phi^{(T)}(\lambda_{k+1} - \lambda_k) \overline{\varphi}_k(\lambda_k, \vec{D}^{(2)})| d\lambda_1 \ldots d\lambda_m,$$

де $\lambda_{m+1} = \lambda_1,$ та для будь-якого $\vec{D}^{(2)}$

$$\overline{\varphi}_k = \begin{cases} \mathbb{I}_{\mathbb{R}}, & \text{якщо } \varphi_k(\cdot, \Delta, \overrightarrow{D}^{(2)}) = f_\Delta; \\ |H^*|^2 + g, & \text{якщо } \varphi_k(\cdot, \Delta, \overrightarrow{D}^{(2)}) = |H^*|^2 f_\Delta + g; \\ |H^*|, & \text{якщо } \varphi_k(\cdot, \Delta, \overrightarrow{D}^{(2)}) = |H^*| f_\Delta. \end{cases} \qquad (3.29)$$

Оцінкою (3.28) для кумулянта $cum(Z_{T,\Delta}(\tau_j), j = 1,...,m)$ завершується перший етап доведення.

- *Етап 2*. Зафіксуємо $\overrightarrow{D}^{(2)}$. Для інтеграла з циклічним зачепленням ядер ($\lambda_{m+1} = \lambda_1$), що фігурує в (3.28), виду

$$I_m^{(T)}(\overrightarrow{D}^{(2)}) = \int ... \int_{\mathbb{R}^m} |\prod_{k=1}^{m} \Phi^{(T)}(\lambda_{k+1} - \lambda_k) \overline{\varphi}_k(\lambda_k, \overrightarrow{D}^{(2)})| d\lambda_1 ... d\lambda_m,$$

далі будемо використовувати теорему А.3, адаптовану до нашого випадку.

(i) Для ядер

$$\Phi^{(T)}(\lambda) = \left(\frac{1}{2\pi T}\right)^{\frac{1}{2}} \frac{e^{iT\lambda} - 1}{i\lambda}, \lambda \in \mathbb{R},$$

при всіх $p \in (1, \infty]$, виконується умова:

$$\|\Phi^{(T)}\|_p = T^{\frac{1}{2} - \frac{1}{p}} \cdot \left[\frac{1}{\sqrt{2\pi}} \int_0^\infty |\frac{\sin \lambda}{\lambda}|^p d\lambda\right]^{\frac{1}{p}}.$$

(ii) Множина функцій $\overline{F}(\overrightarrow{D}^{(2)}) = \{\overline{\varphi}_1(\cdot, \overrightarrow{D}^{(2)}),...,\overline{\varphi}_m(\cdot, \overrightarrow{D}^{(2)})\}$ для вказаного розбиття $\overrightarrow{D}^{(2)}$ розбивається на три класи:

$$M_\infty(\overrightarrow{D}^{(2)}) = \{\overline{\varphi} \in \overline{F}(\overrightarrow{D}^{(2)}) : \overline{\varphi} = \mathbb{I}_{\mathbb{R}}\};$$
$$M_1(\overrightarrow{D}^{(2)}) = \{\overline{\varphi} \in \overline{F}(\overrightarrow{D}^{(2)}) : \overline{\varphi} = |H^*|^2 + g\}; \qquad (3.30)$$
$$M_2(\overrightarrow{D}^{(2)}) = \{\overline{\varphi} \in \overline{F}(\overrightarrow{D}^{(2)}) : \overline{\varphi} = |H^*|\};$$

З формул (3.23), (3.25), (3.29) та (3.30) випливає, що

$$card[M_1(\overrightarrow{D}^{(2)})] = card[M_\infty(\overrightarrow{D}^{(2)})];$$

$$card[M_1(\vec{D}^{(2)})] + card[M_2(\vec{D}^{(2)})] + card[M_\infty(\vec{D}^{(2)})] = m$$

Оскільки $H \in L_2(\mathbb{R})$ і $g \in L_1(\mathbb{R})$, то мають місце включення

$$M_\infty(\vec{D}^{(2)}) \subset L_\infty(\mathbb{R}), \ M_1(\vec{D}^{(2)}) \subset L_1(\mathbb{R}), \ M_2(\vec{D}^{(2)}) \subset L_2(\mathbb{R}).$$

Тобто для даного $\vec{D}^{(2)}$, серед функцій $\overline{\varphi}_k(\cdot, \vec{D}^{(2)}), k = 1,...,m$, існує $m_2(\vec{D}^{(2)}) \geq 0$ функцій з простору $L_1(\mathbb{R})$, $m_1(\vec{D}^{(2)}) = m_2(\vec{D}^{(2)})$ функцій з простору $L_\infty(\mathbb{R})$, та $m_3(\vec{D}^{(2)}) = m - 2m_1(\vec{D}^{(2)})$ функцій з простору $L_2(\mathbb{R})$.

Отже, справджуються всі умови теореми А.3, тому

$$\lim_{T \to \infty} I_m^{(T)}(\vec{D}^{(2)}) = 0. \tag{3.31}$$

Оскільки в скінченній сумі з (3.28) виключена залежність від параметра $\Delta > 0$, й для кожного її доданка $I_m^{(T)}(\vec{D}^{(2)})$ виконується рівність (3.31), то

$$\lim_{(T,\Delta) \to \infty} cum(Z_{T,\Delta}(\tau_j), j = 1,...,m) = 0.$$

Так як гауссівський розподіл однозначно визначається своїми кумулянтами, то з формули (3.19) в силу теореми А.9 випливає, що всі скінченновимірні розподіли процесу $(Z_{T,\Delta}(\tau), \tau \in \mathbb{R})$ слабко збігаються до відповідних скінченновимірних розподілів гауссівського процесу $(Z(\tau), \tau \in \mathbb{R})$ при $(T, \Delta) \to \infty$. Таким чином, теорему 3.2 доведено. □

Твердження теореми 3.2 можна сформулювати за допомогою моментів.

**Теорема 3.3.** *Нехай $H \in L_2(\mathbb{R})$ та $g \in L_1(\mathbb{R})$, тоді для всіх $m \in \mathbb{N}$ та $\tau_1, \tau_2,..., \tau_m \in \mathbb{R}$, має місце рівність*

$$\lim_{(T,\Delta) \to \infty} \mathrm{E}[\prod_{j=1}^m Z_{T,\Delta}(\tau_j)] = \mathrm{E}[\prod_{j=1}^m Z(\tau_j)]. \tag{3.32}$$

*Зокрема, всі скінченновимірні розподіли процесу $(Z_{T,\Delta}(\tau), \tau \in \mathbb{R})$ слабко збігаються до відповідних скінченновимірних розподілів центрованого гауссівського процесу $(Z(\tau), \tau \in \mathbb{R})$.*

*Доведення.* Теорема 3.3 доводиться так же, як і теорема 2.3, розділ 2. □

### 3.2.4 Асимптотична нормальність розподілів процесу $Z_{T,\Delta}$ у просторі неперервних функцій

Після встановлення асимптотичної нормальності скінченновимірних розподілів процесу $Z_{T,\Delta}$ (теореми 3.2 і 3.3), природно поставити питання про асимптотичну нормальність цього процесу у просторі неперервних функцій.

Цей пункт використовує термінологію, пов'язану з гауссівськими випадковими процесами, та спирається на основні поняття й властивості псевдометрик, метричної масивності, ентропійних інтегралів, раніше розглянутих у пункті 2.2.4, розділ 2.

По аналогії з розділом 2, тепер розглянемо функцію

$$q_{H,g}(\tau) = \left[\int_{-\infty}^{\infty} \sin^2\frac{\tau\lambda}{2}(|H^*(\lambda)|^2 + g(\lambda))d\lambda\right]^{\frac{1}{2}}, \tau \in \mathbb{R}, \quad (3.33)$$

де $H^*$ - перетворення Фур'є - Планшереля функції $H$.

Оскільки $H \in L_2(\mathbb{R})$ і $g \in L_1(\mathbb{R})$, то $q_H$ коректно визначена та задає псевдометрики

$$\sigma(\tau_1, \tau_2) = q_{H,g}(|\tau_1 - \tau_2|);$$

$$\sqrt{\sigma(\tau_1, \tau_2)} = \sqrt{q_{H,g}(|\tau_1 - \tau_2|)}, \tau_1, \tau_2 \in \mathbb{R}.$$

Зауважимо, що якщо одночасно $H^*(\lambda) \neq 0$ та $g(\lambda) \neq 0$, на множині додатної лебегової міри, то $\sigma$ та $\sqrt{\sigma}$ є метриками.

Для будь-якого $\varepsilon > 0$ покладемо

$$\mathcal{H}_\sigma(\varepsilon) = \mathcal{H}_\sigma([0,1], \varepsilon);$$

$$\mathcal{H}_{\sqrt{\sigma}}(\varepsilon) = \mathcal{H}_{\sqrt{\sigma}}([0,1], \varepsilon).$$

Оскільки псевдометрики $\sigma$ та $\sqrt{\sigma}$ залежать лише від $|\tau_1 - \tau_2|$, то

$$\int_{0+} \mathcal{H}_\sigma^\beta(\varepsilon) d\varepsilon < \infty \Leftrightarrow \int_{0+} \mathcal{H}_\sigma^\beta([a,b], \varepsilon) d\varepsilon < \infty;$$

$$\int_{0+} \mathcal{H}_{\sqrt{\sigma}}(\varepsilon) d\varepsilon < \infty \Leftrightarrow \int_{0+} \mathcal{H}_{\sqrt{\sigma}}([a,b], \varepsilon) d\varepsilon < \infty$$

для будь-яких сегмента $[a,b] \in \mathbb{R}$, та числа $\beta > 0$.

Крім того, оскільки для всіх $\tau_1, \tau_2 \in \mathbb{R}$,

$$\sigma(\tau_1, \tau_2) \leq \left[\max_{\tau_1, \tau_2 \in \mathbb{R}} \sigma(\tau_1, \tau_2)\right]^{\frac{1}{2}} \sqrt{\sigma(\tau_1, \tau_2)} \leq (\|H^*\|_2^2 + \|g\|_1)^{\frac{1}{4}} \sqrt{\sigma(\tau_1, \tau_2)}, \quad (3.34)$$

то зі співвідношення $\int_{0+} \mathcal{H}_{\sqrt{\sigma}}(\varepsilon) < \infty$, отримаємо $\int_{0+} \mathcal{H}_\sigma(\varepsilon) < \infty$, звідки, в свою чергу, випливає

$$\int_{0+} \mathcal{H}_\sigma^{\frac{1}{2}}(\varepsilon) < \infty.$$

Далі запис

$$Z_{T,\Delta} \stackrel{C[a,b]}{\Rightarrow} Z$$

означає слабку збіжність процесу $Z_{T,\Delta}$ до процесу $Z$ у просторі неперервних функцій при $(T, \Delta) \to \infty$. Будемо вважати, що

$$Z_{T,\Delta} = (Z_{T,\Delta}(\tau), \tau \in [a,b])$$

$$Z = (Z(\tau), \tau \in [a,b]),$$

й ці процеси є сепарабельними. (Таке припущення не є обмежувальним в силу стохастичної неперервності вказаних процесів).

Для всіх $T > 0, \Delta > 0,$ визначимо сім'ю псевдометрик

$$\rho_{(T,\Delta)}(\tau_1, \tau_2) = \left(\mathrm{E} \, | Z_{T,\Delta}(\tau_2) - Z_{T,\Delta}(\tau_1) |^2 \right)^{\frac{1}{2}}, \ \tau_1, \tau_2 \in \mathbb{R}. \qquad (3.35)$$

**Лема 3.4** *Нехай $H \in L_2(\mathbb{R})$ та $g \in L_1(\mathbb{R})$, тоді для всіх $T > 0, \Delta > 0,$ та $\tau_1, \tau_2 \in \mathbb{R}$, має місце нерівність*

$$\rho_{(T,\Delta)}(\tau_1, \tau_2) \leq \frac{4\sqrt{\pi} (\max\{\sup_{\Delta > 0} \| f_\Delta \|_\infty; 1\})}{c} [\| H^* \|_2^2 + \| g \|_1]^{\frac{1}{4}} \sqrt{\sigma(\tau_1, \tau_2)}, \quad (3.36)$$

*де $c$ - стала з умови (3.1г). Крім того, псевдометрика $\rho_{(T,\Delta)}$ неперервна відносно псевдометрики $\sigma$.*

*Доведення.* З означення кореляційної функції процесу $Z_{T,\Delta}$ (див. формулу (3.9)) можна записати

$$\mathrm{E}(Z_{T,\Delta}(\tau_2) - Z_{T,\Delta}(\tau_1))^2 = \mathrm{E} Z_{T,\Delta}^2(\tau_2) - 2\mathrm{E} Z_{T,\Delta}(\tau_1) Z_{T,\Delta}(\tau_2) + \mathrm{E} Z_{T,\Delta}^2(\tau_1) =$$

$$= \frac{4\pi}{c^2} \int_{-\infty}^{\infty} \int_{-\infty}^{\infty} [1 - e^{i(\tau_1 - \tau_2)\lambda_2}](|H^*(\lambda_2)|^2 f_\Delta(\lambda_2) + g(\lambda_2)) \times$$

$$\times \Phi_T(\lambda_2 - \lambda_1) f_\Delta(\lambda_1) d\lambda_1 d\lambda_2 +$$

$$+ \frac{2\pi}{c^2} \int_{-\infty}^{\infty} \int_{-\infty}^{\infty} [e^{i\tau_2(\lambda_1 + \lambda_2)} - 2e^{i(\tau_1 \lambda_1 + \tau_2 \lambda_2)} + e^{i\tau_1(\lambda_1 + \lambda_2)}] H^*(\lambda_1) H^*(\lambda_2) \times$$

$$\times \Phi_T(\lambda_2 - \lambda_1) f_\Delta(\lambda_1) f_\Delta(\lambda_2) d\lambda_1 d\lambda_2 =$$

$$= \frac{4\pi}{c^2} \int_{-\infty}^{\infty} \int_{-\infty}^{\infty} [\sin \frac{(\tau_1 - \tau_2)\lambda_2}{2}]^2 (|H^*(\lambda_2)|^2 f_\Delta(\lambda_2) + g(\lambda_2)) \times$$

$$\times \Phi_T(\lambda_2 - \lambda_1) f_\Delta(\lambda_1) d\lambda_1 d\lambda_2 +$$

$$+\frac{2\pi}{c^2}\int\limits_{-\infty}^{\infty}\int\limits_{-\infty}^{\infty}e^{i(\tau_1\lambda_1+\tau_2\lambda_2)}[e^{i\lambda_1(\tau_2-\tau_1)}-2+e^{i\lambda_2(\tau_1-\tau_2)}]H^*(\lambda_1)H^*(\lambda_2)\times$$

$$\times\Phi_T(\lambda_2-\lambda_1)f_\Delta(\lambda_1)f_\Delta(\lambda_2)d\lambda_1 d\lambda_2.$$

Оскільки $|e^{i\lambda_1(\tau_2-\tau_1)}+e^{i\lambda_2(\tau_1-\tau_2)}-2|\leq|e^{i\lambda_1(\tau_2-\tau_1)}-1|+|e^{i\lambda_2(\tau_1-\tau_2)}-1|$ й $|e^{i\tau\lambda}-1|=2|\sin\frac{\tau\lambda}{2}|$, то, оцінюючи попередню рівність із врахуванням умови (3.1б), парності ядра Фейєра та того, що $\|\Phi_T\|_1=1$, дістанемо

$$\mathrm{E}|Z_{T,\Delta}(\tau_2)-Z_{T,\Delta}(\tau_1)|^2\leq$$

$$\leq\frac{8\pi}{c^2}(\max\{\sup_{\Delta>0}\|f_\Delta\|_\infty;1\})^2\int\limits_{-\infty}^{\infty}\int\limits_{-\infty}^{\infty}[\sin\frac{(\tau_1-\tau_2)\lambda_2}{2}]^2(|H^*(\lambda_2)|^2+g(\lambda_2))\times$$

$$\times\Phi_T(\lambda_2-\lambda_1)d\lambda_1 d\lambda_2+$$

$$+\frac{4\pi}{c^2}(\sup_{\Delta>0}\|f_\Delta\|_\infty)^2\int\limits_{-\infty}^{\infty}\int\limits_{-\infty}^{\infty}[|\sin\frac{(\tau_2-\tau_1)\lambda_1}{2}|+|\sin\frac{(\tau_1-\tau_2)\lambda_2}{2}|]\times$$

$$\times|H^*(\lambda_1)\|H^*(\lambda_2)|\Phi_T(\lambda_2-\lambda_1)d\lambda_1 d\lambda_2=$$

$$=\frac{8\pi}{c^2}(\max\{\sup_{\Delta>0}\|f_\Delta\|_\infty;1\})^2\int\limits_{-\infty}^{\infty}[\sin\frac{(\tau_1-\tau_2)\lambda_2}{2}]^2(|H^*(\lambda_2)|^2+g(\lambda_2))\times$$

$$\times[\int\limits_{-\infty}^{\infty}\Phi_T(\lambda_2-\lambda_1)d\lambda_1]d\lambda_2+$$

$$+\frac{8\pi}{c^2}(\sup_{\Delta>0}\|f_\Delta\|_\infty)^2\int\limits_{-\infty}^{\infty}\int\limits_{-\infty}^{\infty}|\sin\frac{(\tau_1-\tau_2)\lambda_2}{2}\|H^*(\lambda_1)\|H^*(\lambda_2)|\Phi_T(\lambda_2-\lambda_1)d\lambda_1 d\lambda_2=$$

$$=\frac{8\pi}{c^2}(\max\{\sup_{\Delta>0}\|f_\Delta\|_\infty;1\})^2\int\limits_{-\infty}^{\infty}[\sin\frac{(\tau_1-\tau_2)\lambda_2}{2}]^2(|H^*(\lambda_2)|^2+g(\lambda_2))d\lambda_2+$$

$$+\frac{8\pi}{c^2}(\sup_{\Delta>0}\|f_\Delta\|_\infty)^2\int\limits_{-\infty}^{\infty}|\sin\frac{(\tau_1-\tau_2)\lambda_2}{2}\|H^*(\lambda_2)|[\int\limits_{-\infty}^{\infty}|H^*(\lambda_1)|\Phi_T(\lambda_2-\lambda_1)d\lambda_1]d\lambda_2=$$

$$= \frac{8\pi}{c^2}(\max\{\sup_{\Delta>0}\|f_\Delta\|_\infty;1\})^2[\sigma(\tau_1,\tau_2)]^2 +$$

$$+ \frac{8\pi}{c^2}(\sup_{\Delta>0}\|f_\Delta\|_\infty)^2 \int\limits_{-\infty}^{\infty} |\sin\frac{(\tau_1-\tau_2)\lambda_2}{2}| |H^*(\lambda_2)|[|H^*|*\Phi_T](\lambda_2)d\lambda_2.$$

Застосовуючи до другого доданку спочатку нерівність Коші - Буняковського, потім нерівність Юнга для згортки (теорема А.2) й той факт, що для невід'ємної функції $g \in L_1(\mathbb{R})$ інтеграл $0 \le \int\limits_{-\infty}^{\infty}[\sin\frac{(\tau_1-\tau_2)\lambda}{2}]^2 g(\lambda)d\lambda \le +\infty,$ отримаємо

$$\mathrm{E}\,|Z_{T,\Delta}(\tau_2) - Z_{T,\Delta}(\tau_1)|^2 \le$$

$$\le \frac{8\pi}{c^2}(\max\{\sup_{\Delta>0}\|f_\Delta\|_\infty;1\})^2[\sigma(\tau_1,\tau_2)]^2 +$$

$$+ \frac{8\pi}{c^2}(\sup_{\Delta>0}\|f_\Delta\|_\infty)^2 [\int\limits_{-\infty}^{\infty}[\sin\frac{(\tau_1-\tau_2)\lambda_2}{2}]^2 |H^*(\lambda_2)|^2 d\lambda_2]^{\frac{1}{2}} \||H^*|*\Phi_T\|_2 \le$$

$$\le \frac{8\pi}{c^2}(\max\{\sup_{\Delta>0}\|f_\Delta\|_\infty;1\})^2[(\|H^*\|_2^2 + \|g\|_1)^{\frac{1}{4}}\sqrt{\sigma(\tau_1,\tau_2)}]^2 +$$

$$+ \frac{8\pi}{c^2}(\sup_{\Delta>0}\|f_\Delta\|_\infty)^2 [\int\limits_{-\infty}^{\infty}[\sin\frac{(\tau_1-\tau_2)\lambda_2}{2}]^2 |H^*(\lambda_2)|^2 d\lambda_2]^{\frac{1}{2}} \|H^*\|_2 \|\Phi_T\|_1 \le$$

$$\le \frac{8\pi}{c^2}(\max\{\sup_{\Delta>0}\|f_\Delta\|_\infty;1\})^2[(\|H^*\|_2^2 + \|g\|_1)^{\frac{1}{4}}\sqrt{\sigma(\tau_1,\tau_2)}]^2 +$$

$$+ \frac{8\pi}{c^2}(\sup_{\Delta>0}\|f_\Delta\|_\infty)^2 [\|H^*\|_2^2 + \|g\|_1]^{\frac{1}{2}} [\int\limits_{-\infty}^{\infty}[\sin\frac{(\tau_1-\tau_2)\lambda}{2}]^2 (|H^*(\lambda)|^2 + g(\lambda))d\lambda]^{\frac{1}{2}} \le$$

$$\le \frac{16\pi}{c^2}(\max\{\sup_{\Delta>0}\|f_\Delta\|_\infty;1\})^2[\|H^*\|_2^2 + \|g\|_1]^{\frac{1}{2}} \sigma(\tau_1,\tau_2),$$

де використане означення функції $\sigma$ та формула зв'язку між псевдометриками $\sigma$ і $\sqrt{\sigma}$ (див. нерівність (3.34)).

Беручи квадратний корінь з обох частин отриманої нерівності,

дістанемо для всіх $T>0, \Delta>0,$ та $\tau_1, \tau_2 \in \mathbb{R}$,

$$\rho_{(T,\Delta)}(\tau_1,\tau_2) = \left(\mathrm{E}\,|\,Z_{T,\Delta}(\tau_2) - Z_{T,\Delta}(\tau_1)|^2\right)^{\frac{1}{2}} \leq$$

$$\leq \frac{4\sqrt{\pi}(\max\{\sup_{\Delta>0}\|f_\Delta\|_\infty;1\})}{c}[\|H^*\|_2^2 + \|g\|_1]^{\frac{1}{4}}\sqrt{\sigma}(\tau_1,\tau_2).$$

З доведеної нерівності (3.36) та теореми Лебега про мажоровану збіжність, маємо $\rho_{(T,\Delta)}(\tau_1,\tau_2) \to 0$ при $\sqrt{\sigma}(\tau_1,\tau_2) \to 0$, тобто псевдометрика $\rho_{(T,\Delta)}$ неперервна відносно псевдометрики $\sqrt{\sigma}$. Збіжність псевдометрики $\sqrt{\sigma}$ до нуля еквівалентна збіжності псевдометрики $\sigma$ в силу означення. Таким чином, лему 3.4 доведено. $\square$

**Наслідок 3.3.** *Оцінка для середньоквадратичних відхилень процесу $Z_{T,\Delta}$, наведена у лемі 3.4, не залежить від параметрів $T>0, \Delta>0$, тобто насправді встановлено співвідношення:*

$$\sup_{T,\Delta>0}\left(\mathrm{E}\,|\,Z_{T,\Delta}(\tau_2) - Z_{T,\Delta}(\tau_1)|^2\right)^{\frac{1}{2}} \leq$$

$$\leq \frac{4\sqrt{\pi}(\max\{\sup_{\Delta>0}\|f_\Delta\|_\infty;1\})}{c}[\|H^*\|_2^2 + \|g\|_1]^{\frac{1}{4}}\sqrt{\sigma}(\tau_1,\tau_2),\ \tau_1,\tau_2\in\mathbb{R}.$$

**Наслідок 3.4.** *З означення псевдометрики $\sigma$ та теореми Лебега про мажоровану збіжність, маємо $\sigma(\tau_1,\tau_2) \to 0$ при $|\tau_1 - \tau_2| \to 0$, звідки в силу означення (3.35) та нерівності (3.36) випливає, що процеси $Z_{T,\Delta}$ при всіх $T>0, \Delta>0$, є неперервними у середньому квадратичному.*

У наступній теоремі наводяться достатні умови для неперервності майже напевно процесів $Z_{T,\Delta}$ та $Z$, а також слабкої збіжності $Z_{T,\Delta}$ до $Z$ при $(T,\Delta) \to \infty$ у просторі неперервних функцій.

**Теорема 3.4.** *Нехай* $H \in L_2(\mathbb{R})$, $g \in L_1(\mathbb{R})$, *та виконується нерівність*

$$\int\limits_{0+} \mathcal{H}_{\sqrt{\sigma}}(\varepsilon)d\varepsilon < \infty, \tag{3.37}$$

*тоді для будь-якого* $[a,b] \subset \mathbb{R}$ *мають місце наступні твердження:*

*(I)* $Z \in C[a,b]$ *майже напевно;*

*(II)* $Z_{T,\Delta} \in C[a,b]$ *майже напевно,* $T > 0, \Delta > 0$;

*(III)* $Z_{T,\Delta} \overset{C[a,b]}{\Rightarrow} Z$.

*Зокрема, для будь-якого* $x > 0$

$$\lim_{(T,\Delta) \to \infty} \mathrm{P}\left\{\sup_{\tau \in [a,b]} |Z_{T,\Delta}(\tau)| > x\right\} = \mathrm{P}\left\{\sup_{\tau \in [a,b]} |Z(\tau)| > x\right\}.$$

*Доведення.* (I) З теореми Дадлі про неперервність гауссівських процесів (див. теорему А.10) пункт (I) даної теореми має місце, якщо для будь-якого $[a,b] \subset \mathbb{R}$

$$\int\limits_{0+} \mathcal{H}_{d_Z}^{\frac{1}{2}}([a,b],\varepsilon)d\varepsilon < \infty, \tag{3.38}$$

де $d_Z(\tau_1,\tau_2) = [\mathrm{E}\,|Z(\tau_1) - Z(\tau_2)|^2]^{\frac{1}{2}}, \tau_1, \tau_2 \in \mathbb{R}$; а сам процес $Z$ є неперервним у середньому квадратичному.

Оскільки для невід'ємної функції $g \in L_1(\mathbb{R})$ інтеграл $0 \leq \int\limits_{-\infty}^{\infty}[\sin\frac{(\tau_1-\tau_2)\lambda}{2}]^2 g(\lambda)d\lambda \leq +\infty$, то для всіх $\tau_1, \tau_2 \in \mathbb{R}$, має місце оцінка

$$d_Z^2(\tau_1,\tau_2) \leq \frac{2}{\pi}\int\limits_{-\infty}^{\infty} |1 - e^{i(\tau_1-\tau_2)\lambda}|(|H^*(\lambda)|^2 + \frac{2\pi}{c}g(\lambda))d\lambda +$$

$$+\frac{1}{2\pi}\int_{-\infty}^{\infty}|e^{i(\tau_1-\tau_2)\lambda}-2+e^{i(\tau_2-\tau_1)\lambda}|\|H^*(\lambda)|^2\,d\lambda=$$

$$=\frac{2}{\pi}(\max\{\frac{2\pi}{c};1\})\int_{-\infty}^{\infty}[\sin\frac{(\tau_1-\tau_2)\lambda}{2}]^2(|H^*(\lambda)|^2+g(\lambda))d\lambda+$$

$$+\frac{1}{\pi}\int_{-\infty}^{\infty}[\sin\frac{(\tau_1-\tau_2)\lambda}{2}]^2|H^*(\lambda)|^2\,d\lambda\le$$

$$\le(\max\{\frac{4}{\pi};\frac{8}{c}\})\int_{-\infty}^{\infty}[\sin\frac{(\tau_1-\tau_2)\lambda}{2}]^2(|H^*(\lambda)|^2+g(\lambda))d\lambda=$$

$$=(\max\{\frac{4}{\pi};\frac{8}{c}\})\sigma^2(\tau_1,\tau_2).$$

З теореми Лебега про мажоровану збіжність випливає, що $\sigma(\tau_1,\tau_2)\to 0$ при $d(\tau_1,\tau_2)=|\tau_1-\tau_2|\to 0$. З цього факту та встановленої оцінки маємо, що процес $Z$ є неперервним у середньому квадратичному.

Відповідно, нерівність (3.38) виконується, якщо $\int_{0+}\mathcal{H}_\sigma^{\frac{1}{2}}([a,b],\varepsilon)d\varepsilon<\infty$. Остання умова випливає з формули (3.37). Таким чином, пункт (I) доведено; тобто процес $Z\in C[a,b]$ майже напевно.

(II) В силу леми 3.4 при всіх $T>0, \Delta>0,$ та $\tau_1,\tau_2\in\mathbb{R}$, виконується оцінка

$$\rho_{T,\Delta}(\tau_1,\tau_2)=[\mathrm{E}|Z_{T,\Delta}(\tau_1)-Z_{T,\Delta}(\tau_2)|^2]^{\frac{1}{2}}$$

$$\le\frac{4\sqrt{\pi}(\max\{\sup_{\Delta>0}\|f_\Delta\|_\infty;1\})}{c}[\|H^*\|_2^2+\|g\|_1]^{\frac{1}{4}}\sqrt{\sigma(\tau_1,\tau_2)},$$

тому з умови (3.37) для будь-якого $[a,b]\subset\mathbb{R}$ має місце $\int_{0+}\mathcal{H}_{\rho_{T,\Delta}}([a,b],\varepsilon)d\varepsilon<\infty$, звідки, в свою чергу, випливає $\int_{0+}\mathcal{H}_{\rho_{T,\Delta}}^{\frac{1}{2}}([a,b],\varepsilon)d\varepsilon<\infty.$

Крім цього, з наслідку 3.4 відомо, що процес $Z_{T,\Delta}$ є неперервним у середньому квадратичному, тому за теоремою Дадлі (див.теорему A.10) маємо: $Z_{T,\Delta} \in C[a,b]$ майже напевно при всіх $T>0, \Delta>0$. Таким чином, пункт (II) доведено.

(III) У теоремах 3.2 або 3.3 встановлено, що за умов $H \in L_2(\mathbb{R})$ та $g \in L_1(\mathbb{R})$ :

(A) всі скінченновимірні розподіли процесу $Z_{T,\Delta}$ прямують при $(T,\Delta) \to \infty$ до відповідних скінченновимірних розподілів процесу $Z$.

Оскільки $Z_{T,\Delta}$ - сім'я (по $T>0, \Delta>0$) квадратично-гауссівських процесів, то з теореми A.12 випливає, що

$$\sup_{T,\Delta>0} \sup_{\tau_1,\tau_2 \in [a,b]} E\exp\{|\frac{Z_{T,\Delta}(\tau_1) - Z_{T,\Delta}(\tau_2)}{\sqrt{8}\rho_{T,\Delta}(\tau_1,\tau_2)}|\} < \infty. \tag{3.39}$$

Згідно з наслідком 3.7 маємо, що при всіх $\tau_1,\tau_2 \in \mathbb{R}$, псевдометрика

$$\rho_\infty(\tau_1,\tau_2) = \sup_{T,\Delta>0} \rho_{T,\Delta}(\tau_1,\tau_2) \leq \tag{3.40}$$

$$\leq \frac{4\sqrt{\pi}(\max\{\sup_{\Delta>0}\|f_\Delta\|_\infty; 1\})}{c}[\|H^*\|_2^2 + \|g\|_1]^{\frac{1}{4}}\sqrt{\sigma}(\tau_1,\tau_2),$$

мажорується псевдометрикою $\sqrt{\sigma}$, яка є неперервною відносно рівномірної метрики $d(\tau_1,\tau_2) = |\tau_1 - \tau_2|$. З теореми Лебега про мажоровану збіжність випливає, що псевдометрика $\rho_\infty(\tau_1,\tau_2)$ неперервна відносно метрики $d(\tau_1,\tau_2)$ при всіх $\tau_1,\tau_2 \in \mathbb{R}$.

З умови (3.37) та нерівності (3.36) випливає, що для будь-якого $[a,b] \subset \mathbb{R}$

$$\lim_{u \downarrow 0} \sup_{T,\Delta>0} \int_0^u \mathcal{H}_{\rho_{T,\Delta}}([a,b],\varepsilon)d\varepsilon = 0. \tag{3.41}$$

Оскільки сім'я $Z_{T,\Delta} \in C[a,b]$ майже напевно, $T>0, \Delta>0$, і виконуються співвідношення (3.39)-(3.41), то з теореми А.11, випливає, що для будь-якого $[a,b] \subset \mathbb{R}$

(Б) $\quad \forall h > 0 \quad \lim_{\delta \downarrow 0} \sup_{T,\Delta > 0} P\{ \sup_{\substack{\tau_1,\tau_2 \in [a,b] \\ |\tau_1-\tau_2|<\delta}} |Z_{T,\Delta}(\tau_1) - Z_{T,\Delta}(\tau_2)| > h \} = 0.$

Так як виконуються умови (А) та (Б), то з теореми Прохорова про слабку збіжність процесів у просторі неперервних функцій (див. теорему А.13), випливає $Z_{T,\Delta} \overset{C[a,b]}{\Rightarrow} Z$. Пункт (III) доведено.

Оскільки функція розподілу супремума модуля центрованого неперервного майже напевно гауссівського процесу, заданого на відрізку числової осі, є неперервною функцією для всіх $x > 0$ [34], то з пункту (III) випливає

$$\lim_{(T,\Delta)\to\infty} P\left\{ \sup_{\tau \in [a,b]} |Z_{T,\Delta}(\tau)| > x \right\} = P\left\{ \sup_{\tau \in [a,b]} |Z(\tau)| > x \right\}.$$

Отже, теорема 3.4 доведена. $\square$

*Зауваження 3.4.* Твердження (I) теореми 3.4 виконується і при слабшій умові, ніж (3.37), а саме, якщо є збіжним інтеграл Дадлі:

$$\int_{0+}^{1} \mathcal{H}_\sigma^{\frac{1}{2}}(\varepsilon) d\varepsilon < \infty.$$

У роботі [69] встановлено, що вказана умова має місце, якщо існує таке $\beta > 0$, що

$$\int_{-\infty}^{\infty} (|H^*(\lambda)|^2 + g(\lambda)) \ln^{1+\beta}(1+|\lambda|) d\lambda < \infty.$$

Зазначимо, що другу умову перевірити значно простіше, оскільки вона задається у термінах функцій, а не їх ентропійних характеристик.

У наступному твердженні формулюється просте обмеження на функції $H^*$ та $g$, при яких виконується нерівність (3.37).

*Зауваження 3.5.* Умова (3.37) виконується (див.,[31]), якщо існує таке $\beta > 0$,

що
$$\int_{-\infty}^{\infty}(|H^*(\lambda)|^2 + g(\lambda))\ln^{4+\beta}(1+|\lambda|)d\lambda < \infty. \qquad (3.42)$$

Теорема 3.4 дає підґрунтя для побудови довірчих функціональних інтервалів граничного процесу. Наприклад, має місце наступне твердження.

*Зауваження 3.6.* Нехай $d_Z(t,s) = [\mathrm{E}|Z(t) - Z(s)|^2]^{\frac{1}{2}}, t, s \in \mathbb{R}$, - середньоквадратичне відхилення граничного центрованого гауссівського процесу $Z$; $\varepsilon_0 = \sup_{t,s\in[a,b]} d_Z(t,s)$. Так як виконується умова Дадлі

$$I(\varepsilon_0) = \frac{1}{\sqrt{2}}\int_0^{\varepsilon_0} \mathcal{H}_{d_Z}^{\frac{1}{2}}([a,b],\varepsilon)d\varepsilon < \infty,$$

то при всіх $x \geq 8I(\varepsilon_0)$, має місце наступна оцінка "хвоста" розподілу супремума процесу $Z$ (див. приклад 3.4.1 [14])

$$\mathrm{P}\left\{\sup_{\tau\in[a,b]}|Z(\tau)| > x\right\} \leq 2\exp\{-\frac{1}{2\varepsilon_0^2}(x - \sqrt{8xI(\varepsilon_0)})^2\}.$$

## 3.3 Асимптотична незсуненість та конзистентність оцінки $H_{T,\Delta}$

Підрозділ присвячений вивченню деяких характеристик якості оцінки $H_{T,\Delta}$, а саме асимптотичній незсуненості та конзистентності у середньому квадратичному. Для встановлення основних результатів використовуються позначення та факти з підрозділів 3.2 та 2.3, розділ 2.

**Асимптотична незсуненість оцінки.** Нагадаємо, що оцінка $H_{T,\Delta}$ є зсуненою, тобто, взагалі кажучи, $\mathrm{E}H_{T,\Delta}(\tau) \neq H(\tau), \tau \in \mathbb{R}$ (див. підрозділ 3.1). Далі нас цікавлять умови, за яких прямує до нуля похибка оцінювання виду $H_{T,\Delta}(\tau) - H(\tau), \tau \in \mathbb{R}$, при збіжності параметра $\Delta$ до безмежності.

Покладемо

$$\hat{v}_\Delta(\tau) = [\mathrm{E} H_{T,\Delta}(\tau) - H(\tau)], \tau \in \mathbb{R}. \qquad (3.43)$$

Згідно з викладками підрозділу 2.3, розділ 2, анулювати невипадкову функцію $\hat{v}_\Delta$ при $\Delta \to \infty$ можна за рахунок накладення додаткових умов на порядок гладкості перехідної функції $H$ та асимптотичних властивостей кореляційних функцій $K_\Delta$. А саме:

(А) Нехай при деякому $\alpha \in (0,1]$ $H \in Lip_\alpha(\mathbb{R})$.

(Б) Припустимо, що при заданому $\alpha \in (0,1]$ виконуються балансні умови:

$$\forall \delta > 0 \quad \lim_{\Delta \to \infty} \int_\delta^\infty K_\Delta(t) dt = 0; \qquad (3.44а)$$

$$\forall \delta > 0 \quad \lim_{\Delta \to \infty} \int_\delta^\infty K_\Delta^2(t) dt = 0; \qquad (3.44б)$$

$$\exists \delta > 0 \quad \lim_{\Delta \to \infty} \int_{-\delta}^\delta |K_\Delta(t)| |t|^\alpha dt = 0. \qquad (3.44в)$$

Умови (3.1д) та (3.44а)-(3.44в) показують, в якому саме сенсі слід розуміти $\delta$-видність сім'ї кореляційних функцій $K_\Delta$ вхідних процесів при $\Delta \to \infty$.

Має місце наступне твердження, яке дублює лему 2.5, розділ 2, й тому наводиться без доведення.

**Теорема 3.5.** *Нехай задане* $\alpha \in (0,1]$, $H \in Lip_\alpha(\mathbb{R}) \cap L_2(\mathbb{R})$ *та виконуються умови (3.44а) - (3.44в), тоді:*

*(I) для будь-якого* $\tau \in \mathbb{R}$ $\quad \lim_{\Delta \to \infty} \hat{v}_\Delta(\tau) = 0;$

*(II) для будь-якого* $[a,b] \subset \mathbb{R}$ $\quad \lim_{\Delta \to \infty} \sup_{\tau \in [a,b]} |\hat{v}_\Delta(\tau)| = 0.$

**Консистентність оцінки.** Нагадаємо, що оцінка, нормована $\sqrt{T}$ та центрована своїм середнім значенням, мала вид

$$Z_{T,\Delta}(\tau) = \sqrt{T}[H_{T,\Delta}(\tau) - \mathrm{E}H_{T,\Delta}(\tau)], \tau \in \mathbb{R},$$

і детально досліджувалась у підрозділі 3.2.

З наступного представлення

$$\mathrm{E}|H_{T,\Delta}(\tau) - H(\tau)|^2 = \mathrm{E}|H_{T,\Delta}(\tau) - \mathrm{E}H_{T,\Delta}(\tau)|^2 + |\mathrm{E}H_{T,\Delta}(\tau) - H(\tau)|^2 =$$

$$= \frac{\mathrm{E}|Z_{T,\Delta}(\tau)|^2}{T} + |\hat{v}_\Delta(\tau)|^2, \tau \in \mathbb{R}, \qquad (3.45)$$

випливає таке твердження.

**Теорема 3.6.** *Нехай задане* $\alpha \in (0,1]$ $H \in Lip_\alpha(\mathbb{R}) \cap L_2(\mathbb{R})$; $g \in L_1(\mathbb{R})$, *та виконуються умови (3.44а) - (3.44в), тоді для будь-якого* $\tau \in \mathbb{R}$

$$\lim_{(T,\Delta)\to\infty} \mathrm{E}|H_{T,\Delta}(\tau) - H(\tau)|^2 = 0.$$

*Тобто, $H_{T,\Delta}$ є консистентною у середньому квадратичному в точці $\tau$.*

*Доведення.* Використовуючи лему 3.3, із зображення (3.45) маємо оцінку

$$\mathrm{E}|H_{T,\Delta}(\tau) - H(\tau)|^2 = \frac{\mathrm{E}|Z_{T,\Delta}(\tau)|^2}{T} + |\hat{v}_\Delta(\tau)|^2 \leq$$

$$\leq \frac{1}{T} \times \frac{2\pi}{c^2}(\sup_{\Delta>0}\|f_\Delta\|_\infty)[2(\sup_{\Delta>0}\|f_\Delta\|_\infty)\|H^*\|_2^2 + \|g\|_1] + |\hat{v}_\Delta(\tau)|^2.$$

Звідси та пункту (I) теореми 3.5, для будь-якої точки $\tau \in \mathbb{R}$ виконується

$$\lim_{(T,\Delta)\to\infty} \mathrm{E}|H_{T,\Delta}(\tau) - H(\tau)|^2 = 0.$$

Таким чином, теорем 3.6 доведена повністю. □

## 3.4 Асимптотична поведінка похибки оцінювання $W_{T,\Delta}$

Підрозділ присвячений дослідженню асимптотичних властивостей нормованої похибки оцінювання $W_{T,\Delta}$. Для встановлення основних результатів використовуються позначення та факти, наведені у підрозділі 3.2.

### 3.4.1 Анулювання функції $V_{T,\Delta}$

У пункті 3.2.3 показано, що умов $H \in L_2(\mathbb{R})$ та $g \in L_1(\mathbb{R})$, достатньо, щоб при довільному прямуванні параметрів $T \to \infty, \Delta \to \infty$, оцінка $H_{T,\Delta}$, центрована своїм середнім значенням та нормована $\sqrt{T}$, була асимптотично нормальною. Цей факт дає лише часткову характеристику якості оцінки $H_{T,\Delta}$, оскільки $H_{T,\Delta}$ є зсуненою (див. підрозділ 3.1). Насправді, інтерес полягає у дослідженні асимптотичної поведінки похибки оцінювання $H_{T,\Delta}(\tau) - H(\tau), \tau \in \mathbb{R}$. У даному пункті наведені умови, при яких $\sqrt{T}[\mathrm{E}H_{T,\Delta}(\tau) - H(\tau)] \to 0$ при балансному прямуванні $T \to \infty, \Delta \to \infty$.

Всі викладки цього пункту в точності повторюють основні факти, раніше встановлені у пункті 2.4.1, розділ 2.

Покладемо

$$V_{T,\Delta}(\tau) = \sqrt{T}[\mathrm{E}H_{T,\Delta}(\tau) - H(\tau)], \tau \in \mathbb{R}. \qquad (3.46)$$

Накладемо додаткові умови на порядок локальної гладкості перехідної функції $H$ та характер сумісного прямування параметрів $T, \Delta$ до безмежності:

(А) Нехай при деякому $\alpha \in (0,1]$ $H \in Lip_\alpha(\mathbb{R})$;

(Б) Припустимо, що при заданому $\alpha \in (0,1]$ параметри $T \to \infty, \Delta \to \infty$ так, що виконуються балансні умови:

$$\sqrt{T}[1-\frac{2\pi f_\Delta(0)}{c}]\to 0; \qquad (3.47а)$$

$$\forall \delta>0 \quad \sqrt{T}\int_\delta^\infty K_\Delta(t)dt\to 0; \qquad (3.47б)$$

$$\forall \delta>0 \quad T\int_\delta^\infty K_\Delta^2(t)dt\to 0; \qquad (3.47в)$$

$$\exists \delta>0 \quad \sqrt{T}\int_{-\delta}^\delta |K_\Delta(t)||t|^\alpha dt\to 0. \qquad (3.47г)$$

Наступне твердження дублює лему 2.5, розділ 2; в ньому наводяться асимптотичні властивості функції $V_{T,\Delta}$ після введення умов (А) та (Б).

**Лема 3.5.** *Нехай $\alpha \in (0,1]$; $H \in Lip_\alpha(\mathbb{R}) \cap L_2(\mathbb{R})$ та $T\to\infty, \Delta\to\infty$ так, що виконуються умови (3.47а) - (3.47г), тоді:*

*(I) для будь-якого $\tau \in \mathbb{R}$ $\quad V_{T,\Delta}(\tau)\to 0$;*

*(II) для будь-якого $[a,b]\subset \mathbb{R}$ $\sup\limits_{\tau\in[a,b]}|V_{T,\Delta}(\tau)|\to 0.$*

Для наглядності результатів, нагадаємо приклади збурюючих процесів, які задовольняють балансним умовам (Б).

*Приклад 3.1.* Нехай задане $\alpha \in (0,1]$. Спектральні щільності

$$f_\Delta = (\frac{c}{2\pi}\exp\left(-\frac{\lambda^2}{\Delta}\right), \lambda \in \mathbb{R})$$

та відповідні їм кореляційні функції

$$K_\Delta = (\frac{c}{2}\sqrt{\frac{\Delta}{\pi}}\exp(-\frac{\Delta t^2}{4}), t \in \mathbb{R}),$$

збурюючих процесів $X_\Delta, \Delta > 0$, задовольняють умовам (3.1а) - (3.1д) та (3.47а) - (3.47г), якщо $T\to\infty, \Delta\to\infty$ так, що має місце співвідношення

$$T\Delta^{-\alpha}\to 0. \qquad (3.48)$$

*Приклад 3.2.* Нехай задане $\alpha \in (0,1]$. Спектральні щільності

$$f_\Delta = (\frac{c}{2\pi}\frac{\Delta}{\Delta + \lambda^2}, \lambda \in \mathbb{R})$$

та відповідні їм кореляційні функції

$$K_\Delta = (c\sqrt{\Delta}\exp(-\sqrt{\Delta}|t|), t \in \mathbb{R}),$$

збурюючих процесів $X_\Delta, \Delta > 0$, задовольняють умовам (3.1а) - (3.1д) та (3.47а) - (3.47г), якщо $T \to \infty, \Delta \to \infty$ так, що має місце співвідношення

$$T\Delta^{-\alpha} \to 0. \qquad (3.48)$$

*Приклад 3.3.* Нехай задане $\alpha \in (0,1)$. Спектральні щільності

$$f_\Delta = (\frac{c}{2\pi}\exp\left(-\frac{|\lambda|}{\sqrt{\Delta}}\right), \lambda \in \mathbb{R})$$

та відповідні їм кореляційні функції

$$K_\Delta = (\frac{c}{\pi}\frac{\sqrt{\Delta}}{1+\Delta t^2}, t \in \mathbb{R}),$$

збурюючих процесів $X_\Delta, \Delta > 0$, задовольняють умовам (3.1а) - (3.1д) та (3.47а) - (3.47г), якщо $T \to \infty, \Delta \to \infty$ так, що має місце співвідношення

$$T\Delta^{-\alpha} \to 0. \qquad (3.48)$$

*Приклад 3.4.* Нехай задане $\alpha \in (0,1)$. Спектральні щільності

$$f_\Delta = (\frac{c}{2\pi}(1-\frac{|\lambda|}{\sqrt{\Delta}})\mathbb{I}_{[-\sqrt{\Delta},\sqrt{\Delta}]}(\lambda), \lambda \in \mathbb{R})$$

та відповідні їм кореляційні функції

$$K_\Delta = (\frac{c}{2\pi\sqrt{\Delta}}(\frac{\sin(\frac{\sqrt{\Delta}t}{2})}{\frac{t}{2}})^2, t \in \mathbb{R}),$$

збурюючих процесів $X_\Delta, \Delta > 0$, задовольняють умовам (3.1а) - (3.1д) та (3.47а) - (3.47г), якщо $T \to \infty, \Delta \to \infty$ так, що має місце співвідношення

$$T\Delta^{-\alpha} \to 0. \qquad (3.48)$$

### 3.4.2 Асимптотична нормальність скінченновимірних розподілів похибки $W_{T,\Delta}$

У цьому пункті розглядається асимптотична нормальність скінченновимірних розподілів нормованої похибки оцінки імпульсної перехідної функції. Покладемо

$$W_{T,\Delta}(\tau) = \sqrt{T}[H_{T,\Delta}(\tau) - H(\tau)], \tau \in \mathbb{R}. \qquad (3.49)$$

Для дослідження асимптотичної поведінки процесу $W_{T,\Delta} = (W_{T,\Delta}(\tau), \tau \in \mathbb{R})$, його зручно представити у вигляді суми

$$W_{T,\Delta} = Z_{T,\Delta} + V_{T,\Delta}, \qquad (3.50)$$

де

$$Z_{T,\Delta}(\tau) = \sqrt{T}[H_{T,\Delta}(\tau) - \mathrm{E}H_{T,\Delta}(\tau)], \tau \in \mathbb{R};$$

$$V_{T,\Delta}(\tau) = \sqrt{T}[\mathrm{E}H_{T,\Delta}(\tau) - H(\tau)], \tau \in \mathbb{R}.$$

У пунктах 3.2.1 - 3.2.3 доведена асимптотична нормальність випадкового процесу $Z_{T,\Delta}$ при довільному прямуванні параметрів $T, \Delta$ до безмежності. У пункті 3.4.1 знайдені умови, які дозволяють анулювати вклад невипадкової функції $V_{T,\Delta}$ у виразі для нормованої похибки оцінювання $W_{T,\Delta}$ за рахунок додаткових умов на гладкість функції $H$, балансних умов на кореляційні функції процесів $(X_\Delta, \Delta > 0)$ та прямування параметрів $T, \Delta$ до безмежності. Суть подальшої роботи - об'єднати ці характеристики.

Розглянемо асимптотичну поведінку кореляційної функції процесу $W_{T,\Delta}$ при прямуванні параметрів $T, \Delta$ до безмежності.

**Теорема 3.7** *Нехай для деякого $\alpha \in (0,1]$ $H \in Lip_\alpha(\mathbb{R}) \cap L_2(\mathbb{R})$, та $g \in L_1(\mathbb{R})$. Якщо $T \to \infty, \Delta \to \infty$ так, що виконуються умови (3.47а) - (3.47г), тоді для всіх $\tau_1, \tau_2 \in \mathbb{R}$, має місце співвідношення*

$$EW_{T,\Delta}(\tau_1)W_{T,\Delta}(\tau_2) \to C_\infty(\tau_1,\tau_2) = \frac{1}{2\pi}\int_{-\infty}^{\infty}[e^{i(\tau_1-\tau_2)\lambda}|H^*(\lambda)|^2 + e^{i(\tau_1+\tau_2)\lambda}(H^*(\lambda))^2]d\lambda.$$

(3.51)

*Доведення.* Твердження теореми 3.7 випливає з представлення

$$EW_{T,\Delta}(\tau_1)W_{T,\Delta}(\tau_2) = EZ_{T,\Delta}(\tau_1)Z_{T,\Delta}(\tau_2) + V_{T,\Delta}(\tau_1)V_{T,\Delta}(\tau_2),$$

теореми 3.1 та леми 3.5, пункт (I). Таким чином, теорему 3.7 доведено. □

З теореми 3.7 випливає, що гранична кореляційна функція процесу $W_{T,\Delta}$ співпадає з кореляційною функцією деякого центрованого гауссівського процесу $Z = (Z(\tau), \tau \in \mathbb{R})$, що зустрічається у підрозділах 3.2.3 - 3.2.4. Зокрема, має місце таке твердження про асимптотичну поведінку розподілів $W_{T,\Delta}$ при прямуванні параметрів $T, \Delta$ до безмежності.

**Теорема 3.8.** *Нехай для деякого $\alpha \in (0,1]$ $H \in Lip_\alpha(\mathbb{R}) \cap L_2(\mathbb{R})$, та $g \in L_1(\mathbb{R})$. Якщо $T \to \infty, \Delta \to \infty$ так, що виконуються умови (3.47а) - (3.47г), тоді для всіх $m \in \mathbb{N}$ та $\tau_1, \tau_2, ..., \tau_m \in \mathbb{R}$, має місце співвідношення*

$$cum(W_{T,\Delta}(\tau_j), j=1,...,m) \to \begin{cases} 0, & m \neq 2; \\ C_\infty(\tau_1,\tau_2), & m = 2, \end{cases} \quad (3.52)$$

*де $cum(W_{T,\Delta}(\tau_j), j=1,...,m)$ - сумісний кумулянт набору випадкових величин $W_{T,\Delta}(\tau_j), j=1,...,m$.*

*Зокрема, всі скінченновимірні розподіли процесу $(W_{T,\Delta}(\tau), \tau \in \mathbb{R})$ слабко збігаються до відповідних скінченновимірних розподілів центрованого гауссівського процесу $(Z(\tau), \tau \in \mathbb{R})$ при вказаному характері прямування $T$ і $\Delta$ до нескінченності.*

*Доведення.* Теорема 3.8 доводиться так же, як і теорема 2.8, розділ 2. □

**Теорема 3.9.** *Нехай для деякого* $\alpha \in (0,1]$ $H \in Lip_\alpha(\mathbb{R}) \cap L_2(\mathbb{R})$, *та* $g \in L_1(\mathbb{R})$. *Якщо* $T \to \infty, \Delta \to \infty$ *так, що виконуються умови (3.47а) - (3.47г), тоді для всіх* $m \in \mathbb{N}$ *та* $\tau_1, \tau_2, ..., \tau_m \in \mathbb{R}$, *має місце співвідношення*

$$E[\prod_{j=1}^{m} W_{T,\Delta}(\tau_j)] \to E[\prod_{j=1}^{m} Z(\tau_j)]. \qquad (3.53)$$

*Зокрема, всі скінченновимірні розподіли процесу* $(W_{T,\Delta}(\tau), \tau \in \mathbb{R})$ *слабко збігаються до відповідних скінченновимірних розподілів центрованого гауссівського процесу* $(Z(\tau), \tau \in \mathbb{R})$ *при вказаному характері прямування* $T$ *і* $\Delta$ *до нескінченності.*

*Доведення.* Теорема 3.9 доводиться так же, як і теорема 2.9, розділ 2. □

### 3.4.3 Асимптотична нормальність розподілів процесу $W_{T,\Delta}$ у просторі неперервних функцій

Після встановлення асимптотичної нормальності скінченновимірних розподілів процесу $W_{T,\Delta}$ (теореми 3.8 і 3.9), природно поставити питання про асимптотичну нормальність цього процесу у просторі неперервних функцій. У цьому пункті використовуються термінологія та позначення пункту 3.2.4.

Будемо вважати, що центрований гауссівський процес

$$W_{T,\Delta} = (W_{T,\Delta}(\tau), \tau \in [a,b])$$

є сепарабельними. (Таке припущення не є обмежувальним в силу стохастичної неперервності цього процесу). Має місце наступне твердження про асимптотичну нормальність процесу $W_{T,\Delta}$ у просторі неперервних функцій.

**Теорема 3.10.** *Нехай для деякого* $\alpha \in (0,1]$ $H \in Lip_\alpha(\mathbb{R}) \bigcap L_2(\mathbb{R})$; $g \in L_1(\mathbb{R})$, *та виконується нерівність*

$$\int_{0+} \mathcal{H}_{\sqrt{\sigma}}(\varepsilon) d\varepsilon < \infty, \qquad (3.37)$$

*тоді для будь-якого $[a,b] \subset \mathbb{R}$ мають місце наступні твердження:*

*(I)     $Z \in C[a,b]$ майже напевно;*

*(II)    $W_{T,\Delta} \in C[a,b]$ майже напевно, $T > 0, \Delta > 0$;*

*Якщо $T \to \infty, \Delta \to \infty$ так, що виконуються умови (3.47а) - (3.47г), тоді*

*(III)   $W_{T,\Delta} \overset{C[a,b]}{\Rightarrow} Z$.*

*Зокрема, при вказаному характері прямування $T$ і $\Delta$ до нескінченності, для будь-якого $x > 0$*

$$P\left\{\sup_{\tau \in [a,b]} |W_{T,\Delta}(\tau)| > x\right\} \to P\left\{\sup_{\tau \in [a,b]} |Z(\tau)| > x\right\}.$$

*Доведення.* (I)   Використовуючи формулу (3.37), пункт (I) теореми 3.10 був доведений у пункті (I) теореми 3.4.

(II)   Функція $(EH_{T,\Delta}(\tau), \tau \in \mathbb{R})$ - невипадкова і неперервна на $\mathbb{R}$, як нормова-на (на $\dfrac{1}{c}$) сумісна кореляційна функція процесів $X_\Delta$ та $Y_\Delta$, кожен з яких є неперервним у середньому квадратичному. Згідно із зауваженням 2.10, розділ 2, перехідна функція $H$ є неперервною на $\mathbb{R}$. Тому $V_{T,\Delta}$, як нормована різни-ця цих функцій, є неперервною на $\mathbb{R}$ функцією.

Пункт (II) теореми 3.10 випливає з представлення процесу $W_{T,\Delta} = Z_{T,\Delta} + V_{T,\Delta}$ (див. (3.50)) та пункту (II) теореми 3.4.

(III)   Пункт (III) теореми 3.10 випливає з з представлення процесу $W_{T,\Delta}$ (див. формулу (3.50)), пункту (II) теореми 3.4, а також твердження (II) леми 3.5 та теореми А.14.                                          □

Наведемо приклади імпульсних перехідних функцій $H \in L_2(\mathbb{R})$ та характеристик збурюючих процесів $X_\Delta$ й процесу внутрішнього шуму $U$, які гарантують асимптотичну нормальність відповідної оцінки та похибки.

*Приклад 3.5.* В якості спектральної щільності $g$ процесу $U$ можна розглядати, наприклад, одну з наступних функцій:

(I) $\quad g(\lambda) = \mathbb{I}_{[-\mu,\mu]}(\lambda), \lambda \in \mathbb{R}$, де $\mu > 0$;

(II) $\quad g(\lambda) = e^{-\mu|\lambda|}, \lambda \in \mathbb{R}$, де $\mu > 0$;

(III) $\quad g(\lambda) = e^{-\mu\lambda^2}, \lambda \in \mathbb{R}$, де $\mu > 0$;

(IV) $\quad g(\lambda) = \dfrac{1}{2\pi\mu}\left(\dfrac{\sin(\mu\lambda/2)}{\lambda/2}\right)^2, \lambda \in \mathbb{R}$, де $\mu > 0$;

(V) $\quad g(\lambda) = \dfrac{\mu}{\pi}\dfrac{1}{\lambda^2 + \mu^2}, \lambda \in \mathbb{R}$, де $\mu > 0$.

Тоді $g$ - парна функція та $g \in L_1(\mathbb{R})$. Крім того, тоді при довільному $\beta > 0$

$$\int_{-\infty}^{\infty} g(\lambda)\ln^{4+\beta}(1+|\lambda|)d\lambda < \infty. \qquad (3.54)$$

В якості імпульсної перехідної функції $H$ можна розглядати, наприклад, одну з наступних функцій:

(I) $\quad H(t) = \dfrac{t}{1+t^2}, t \in \mathbb{R}$;

(II) $\quad H(t) = \dfrac{1-\cos(\mu t)}{\pi t}, t \in \mathbb{R}$, де $\mu > 0$;

(III) $\quad H(t) = \dfrac{\sin(\mu t)}{\pi t}, t \in \mathbb{R}$, де $\mu > 0$;

(IV) $\quad H(t) = \dfrac{\sin(\nu t) - \sin(\mu t)}{\pi t}, t \in \mathbb{R}$, де $\nu > \mu > 0,$;

(V) $\quad H(t) = \sum_{k=1}^{n}\dfrac{\sin(\nu_k t) - \sin(\mu_k t)}{\pi t}, t \in \mathbb{R}$, де $0 < \mu_k < \nu_k \leq \mu_{k+1}, k = 1,...,n$;

(VI) $\quad H(t) = \begin{cases} \dfrac{1}{1+t}, & t \geq 0; \\ 0, & t < 0, \end{cases}$

(VII)    $H(t) = \begin{cases} \dfrac{1}{(1+t)^\gamma}, & t \geq 0; \\ 0, & t < 0, \end{cases}$ де $\gamma \in (\dfrac{3}{4}, 1)$;

(VIII)    $H(t) = \begin{cases} \dfrac{\cos(\mu t)}{(1+t)^\gamma}, & t \geq 0; \\ 0, & t < 0, \end{cases}$ де $\mu > 0$ та $\gamma \in (\dfrac{3}{4}, 1)$.

(IX)    $H(t) = \begin{cases} \dfrac{\sin(\mu t)}{(1+t)^\gamma}, & t \geq 0; \\ 0, & t < 0, \end{cases}$ де $\mu > 0$ та $\gamma \in (\dfrac{3}{4}, 1)$.

Тоді:

у випадках (I) - (V) для всіх $\alpha \in (0,1]$  $H \in Lip_\alpha(\mathbb{R}) \cap L_2(\mathbb{R})$;

у випадку (VI) для всіх $\alpha \in (0,1]$  $H \in Lip_\alpha[0,\infty) \cap L_2(\mathbb{R})$;

у випадках (VII)-(IX) для всіх $\alpha \in (\dfrac{1}{2}, 1]$  $H \in Lip_\alpha[0,\infty) \cap L_2(\mathbb{R})$.

У прикладах 2.10-2.18, розділ 2, було показано, що при довільному $\beta > 0$

$$\int_{-\infty}^{\infty} |H^*(\lambda)|^2 \ln^{4+\beta}(1+|\lambda|) d\lambda < \infty, \tag{3.55}$$

де $H^*$ - перетворення Фур'є - Планшереля функції $H$. Таким чином, з формул (3.54) і (3.55) випливає збіжність інтеграла

$$\int_{-\infty}^{\infty} (|H^*(\lambda)|^2 + g(\lambda)) \ln^{4+\beta}(1+|\lambda|) d\lambda < \infty,$$

звідки в силу зауваження 3.5 випливає збіжність ентропійного інтеграла для будь-якого $[a,b] \subset \mathbb{R}$ (випадки (I) - (V)) або $[a,b] \subset [0,\infty)$ (випадки (VI) - (IX))

$$\int_{0+} \mathcal{H}_{\sqrt{\sigma}}(\varepsilon) d\varepsilon < \infty. \tag{3.37}$$

Нарешті, якщо в якості спектральних щільностей $f_\Delta$ сім'ї збурюючих процесів $X_\Delta, \Delta > 0$, розглядати функції з прикладів 3.1- 3.4, й врахувати, що $T \to \infty, \Delta \to \infty$ так, що

$$T\Delta^{-\alpha} \to 0, \tag{3.48}$$

то має місце функціональна теорема 3.10.

## 3.5 Висновки

У розділі 3 знайдено: умови асимптотичної незсуненості та конзистентності у середньому квадратичному інтегральної корелограмної оцінки імпульсної перехідної функції; умови асимптотичної нормальності відповідної оцінки та похибки для нестійкої однорідної лінійної системи з внутрішнім шумом. Досліджено збіжність розподілів оцінки та похибки як в сенсі збіжності скінченновимірних розподілів, так і у сенсі слабкої збіжності у просторі неперервних функцій. Всі результати розділу 3 узагальнюють результати попереднього розділу 2, а також покращують відомі результати у даній схемі оцінювання. Знайдено приклади збурюючих процесів, внутрішнього шуму та імпульсних перехідних функцій $H \in L_2(\mathbb{R})$, для яких спрацьовує розглядувана теорія оцінювання. Результати є застосовними до нестійких систем з резонансними особливостями (число яких має нульову лебегову міру). Результати одержані при мінімальних умовах на порядок інтегрованості перехідної функції нестійкої однорідної лінійної системи з внутрішнім шумом. Задача про асимптотичну нормальність інтегральної корелограмної оцінки та її похибки за умови, що внутрішній шум не залежить від параметра схеми серій, та ортогональний до сім'ї збурюючих процесів, розв'язана повністю.

Результати цього розділу були опубліковані у роботах [124, 126] і доповідалися на Міжнародній конференції «Stochastic analysis and random dynamics» (2009) [129], "XIII Міжнародній науковій конференції ім. акад. М. Кравчука" (2010) [131], Міжнародній конференції "Modern Stochastics: Theory and Applications, II" [132].

# ЗАГАЛЬНІ ВИСНОВКИ

Дисертаційна робота присвячена оцінюванню невідомої дійснозначної імпульсної перехідної функції неперервної однорідної лінійної системи без шуму та з шумом. Основне обмеження - інтегрованість у квадраті відповідної перехідної функції; таким чином, під розгляд потрапляють нестійкі системи з резонансними особливостями (число яких може бути зліченним), зокрема, системи типу Вольтерра.

Припускається, що на вхід системи подається сім'я (по $\Delta$) дійснозначних центрованих стаціонарних гауссівських процесів, "близьких" до білого шуму при $\Delta \to \infty$. Для систем з внутрішнім шумом, додатково припускається, що вхідні процеси ортогональні до внутрішнього шуму системи, який не залежить від $\Delta$, і є також дійснозначним центрованим стаціонарним гауссівським процесом. У роботі будується корелограмна оцінка інтегрального типу й встановлюються наступні факти:

- умови асимптотичної нормальності оцінки імпульсної перехідної функції як у сенсі збіжності скінченновимірних розподілів, так і в сенсі їх слабкої збіжності у просторі неперервних функцій;
- умови асимптотичної незсуненості та конзистентності у середньому квадратичному оцінки імпульсної перехідної функції;
- умови асимптотичної нормальності похибки оцінювання імпульсної перехідної функції як у сенсі збіжності скінченновимірних розподілів, так і в сенсі їх слабкої збіжності у просторі неперервних функцій;
- конкретні приклади збурюючих процесів, внутрішнього шуму системи та імпульсних перехідних функцій, які задовольняють схемі оцінювання.

Усі результати роботи є новими і покращують наявні розв'язки задачі оцінювання імпульсної перехідної функції за допомогою корелограм інтегрального типу. Вдосконалення вдається зробити за рахунок застосування кумулянтного методу аналізу розподілу та теорії багатовимірних згорткових інтегралів з циклічним зачепленням ядер, залежних від параметрів.

Робота носить теоретичний характер. Отримані результати можна застосовувати при оцінюванні або ідентифікації імпульсних перехідних функцій однорідних лінійних систем з внутрішнім шумом.

# Література

# ДОДАТОК А

# ПОПЕРЕДНІ ВІДОМОСТІ

У Додатку А вводяться основні поняття та формулюються відомі твердження, необхідні для загального тексту. Наведені властивості перетворення Фур'є - Планшереля; згорткова нерівність Юнга та деякі властивості інтегралів з циклічним зачепленням ядер; характеристики лінійних систем, що збурюються стаціонарними гауссівськими процесами; основні відомості про кумулянти та моменти гауссівських процесів; а також базові теореми про збіжність розподілів у просторі неперервних функцій.

## А.1 Перетворення Фур'є - Планшереля функцій з $L_2(\mathbb{R})$

*Перетворенням Фур'є - Планшереля* функції $\phi \in L_2(\mathbb{R})$ називається наступний інтеграл

$$\phi^*(\lambda) = \lim_{N \to +\infty} \int_{-N}^{N} e^{-i\lambda x}\phi(x)dx = \int_{-\infty}^{\infty} e^{-i\lambda x}\phi(x)dx, \lambda \in \mathbb{R}.$$

(збіжність середнього інтеграла розуміється в сенсі $L_2(\mathbb{R})$).

Згідно з теоремою Планшереля [30], $\phi^*$ існує, а також $\phi^* \in L_2(\mathbb{R})$ і $\|\phi^*\|_2 = \sqrt{2\pi}\|\phi\|_2$. Якщо $\phi \in L_1(\mathbb{R}) \bigcap L_2(\mathbb{R})$, то $\phi^*$ співпадає з перетворенням Фур'є у звичайному сенсі. Крім того, якщо $\phi_1, \phi_2 \in L_2(\mathbb{R})$, то

$$\int_{-\infty}^{\infty} \phi_1(x)\overline{\phi_2(x)}dx = \int_{-\infty}^{\infty} \phi_1^*(\lambda)\overline{\phi_2^*(\lambda)}d\lambda.$$

## А.2 Теорема Фубіні - Тонеллі

Нехай $(\Omega, \mathfrak{F}, P)$ - повний ймовірнісний простір.

Нижче наведена теорема Фубіні - Тонеллі [30], яка дає обґрунтування пронесенню математичного сподівання під знак інтегралу (і навпаки).

**Теорема А.1.** *Нехай $\xi(\omega,t):\Omega\times\mathbb{R}^m\to\mathbb{R}$ - вимірна випадкова функція. Тоді справедливі твердження:*

*(1) вибіркові функції $\xi(t)$ м.н. $B(\mathbb{R}^m)$ - вимірні по $t\in\mathbb{R}^m$;*

*(2) якщо для всіх $t\in\mathbb{R}^m$ існує $\mathrm{E}\xi(t)$, то функція $\mathrm{E}\xi(t)$ є $B(\mathbb{R}^m)$ - вимірною по $t\in\mathbb{R}^m$;*

*(3) нехай $A\subset B(\mathbb{R}^m)$. Випадкова функція $\xi(t)$ буде інтегрованою на $\Omega\times A$ тоді й лише тоді, коли хоча б один з повторних інтегралів*

$$\int\limits_A \mathrm{E}|\xi(t)|\,dt, \qquad \mathrm{E}\int\limits_A |\xi(t)|\,dt$$

*є скінченним, причому в цьому випадку подвійний інтеграл співпадає з повторними:*

$$\int\limits_A \mathrm{E}|\xi(t)|\,dt = \mathrm{E}\int\limits_A |\xi(t)|\,dt.$$

### А.3 Нерівності для згорток

Сформулюємо відому *нерівність Юнга* [70] для згортки двох функцій.

**Теорема А.2.** *Нехай $1\le a,b,c\le +\infty$, і $\dfrac{1}{a}+\dfrac{1}{b}=1+\dfrac{1}{c}$. Якщо $f\in L_a(\mathbb{R})$ і $\phi\in L_b(\mathbb{R})$, тоді $f*\phi\in L_c(\mathbb{R})$ та має місце нерівність*

$$\|f*\phi\|_c\le\|f\|_a\|\phi\|_b.$$

*Тут через $f*\phi$ позначено згортку функцій $f$ і $\phi$, тобто*

$$(f*\phi)(t)=\int\limits_{-\infty}^{\infty} f(t-s)\phi(s)\,ds.$$

Тепер розглянемо деякі властивості багатократних згорток із циклічним зачепленням аргументів.

Нехай $n \in \mathbb{N} \setminus \{1\}$; $K_j = (K_j(x), x \in \mathbb{R}^m)$ і $\phi_j = (\phi_j(x), x \in \mathbb{R}^m)$, $j = 1,...,n$, - вимірні функції. *Інтегралом з циклічним зачепленням ядер* $K_j, j = 1,...,n,$ називається інтеграл виду:

$$I_n = I(K_1,...,K_n;\phi_1,...,\phi_n) = \int \cdots \int_{\mathbb{R}^m \times \cdots \times \mathbb{R}^m} [\prod_{j=1}^{n} K_j(x_j - x_{j+1})\phi_j(x_j)]dx_1...dx_n,$$

де $x_{n+1} = x_1$; та збіжність якого розуміється в сенсі Лебега. Поряд з наведеним інтегралом будемо розглядати $\overline{I}_n = I(|K_1|,...,|K_n|;|\phi_1|,...,|\phi_n|)$. Зрозуміло: якщо $\overline{I}_n < \infty$, то $I_n$ коректно визначений та $|I_n| \leq \overline{I}_n$.

Припустимо, що ядра $K_1,...,K_n$, залежать від параметра $\theta$ з деякого параметричного простору $\Theta$, тобто $K_j = K_j^{\theta}, j = 1,...,n, \theta \in \Theta$.

Умови збіжності до нуля інтегралів $I_n^{\theta} = I(K_1^{\theta},...,K_n^{\theta};\phi_1,...,\phi_n)$ та $\overline{I}_n^{\theta} = I(|K_1^{\theta}|,...,|K_n^{\theta}|;|\phi_1|,...,|\phi_n|)$, залежних від параметра $\theta$ та функцій $\phi_j, j = 1,...,n,$ наводяться у наступному твердженні.

**Теорема А.3.** (Теорема 5.3 (В), [62]) *Нехай $n \in \mathbb{N} \setminus \{1, 2\}$. Припустимо:*

*(i) для ядер $K_j^{\theta}, j = 1,...,n,$ при всіх $p \in (1, \infty]$ існує така стала $C_p \geq 0$, що виконується мажорантна умова*

$$\max_j \| K_j^{\theta} \| \leq C_p [\sigma(\theta)]^{\frac{1}{2}-\frac{1}{p}},$$

*де $\sigma(\theta) > 0, \theta \in \Theta$;*

*(ii) серед функцій $\phi_1,...,\phi_n$, існує $n_1 \geq 0$ функцій з простору $L_1(\mathbb{R})$, $n_{\infty} = n_1$ функцій з простору $L_{\infty}(\mathbb{R})$, та $n_2 = n - 2n_1 \geq 0$ функцій з простору $L_2(\mathbb{R})$, тоді*

$$\lim_{\sigma(\theta) \to \infty} \overline{I}_n^{\theta} = 0.$$

### А.4   Лінійні системи. Стаціонарні процеси

Нехай задана однорідна лінійна система з імпульсною дійснозначною перехідною функцією $H = (H(\tau), \tau \in \mathbb{R})$. Це означає, що реакція системи на "допустимий" вхідний сигнал $x(t), t \in \mathbb{R}$, має вид

$$y(t) = \int_{-\infty}^{\infty} H(s)x(t-s)ds = \int_{-\infty}^{\infty} H(t-s)x(s)ds.$$

Система називається типу Вольтерра, якщо стосовно неї припускається фізична здійснимість, тобто функція $H$ задовольняє умові : $H(s) = 0, s < 0$. В цьому випадку відгук має вид

$$y(t) = \int_{-\infty}^{t} H(t-s)x(s)ds, t > 0.$$

Якщо $H \in L_1(\mathbb{R})$, то лінійна система називається *стійкою*, і її *частотною характеристикою* $H^*$ є звичайне перетворення Фур'є функції $H$ у просторі $L_1(\mathbb{R})$. Зазначимо, якщо умова $H \in L_1(\mathbb{R})$ не виконується, то лінійна система називається *нестійкою*.

Якщо $H \in L_2(\mathbb{R})$, то лінійна система є *нестійкою*, і її *частотною характеристикою* $H^*$ є перетворення Фур'є-Планшереля функції $H$ у просторі $L_2(\mathbb{R})$.

Зазначимо, що у дисертації розглядаються нестійкі однорідні лінійні системи з умовою $H \in L_2(\mathbb{R})$.

**Теорема А.4.** ([19, 31]) *Нехай $X = (X(t), t \in \mathbb{R})$ - дійснозначний, стаціонарний у широкому сенсі процес зі спектральною щільністю $f_X(\lambda), \lambda \in \mathbb{R}$; $\mathrm{E}X(t) = 0, t \in \mathbb{R}$;  $H$ - перехідна функція однорідної лінійної системи.*
*Покладемо*

$$Y(t) = \int_{-\infty}^{\infty} H(s)X(t-s)ds, t \in \mathbb{R}.$$

*Вказаний інтеграл визначений як відповідна границя у середньому квадратичному тоді й лише тоді, коли існує інтеграл*

$$\int_{-\infty}^{\infty}\int_{-\infty}^{\infty} K(t-s)H(s)H(t)dsdt,$$

*де $K(\tau) = \mathrm{E}X(t+\tau)X(t), \tau \in \mathbb{R}$, - кореляційна функція процесу $X$. При цьому спектральна щільність стаціонарного у широкому сенсі процесу $Y$ має вид*

$$f_Y(\lambda) = |H^*(\lambda)|^2 f_X(\lambda), \lambda \in \mathbb{R}.$$

*Крім того, процеси $X$ і $Y$ є сумісно стаціонарними у широкому сенсі. Якщо, додатково, процес $X$ є гауссівським, то процес $Y$ також є гауссівським, а процеси $X$ і $Y$ є сумісно гауссівськими.*

Також в роботі розглядається однорідна лінійна система з внутрішнім шумом та імпульсною перехідною функцією $H = (H(\tau), \tau \in \mathbb{R})$. Це означає, що реакція системи на "допустимий" вхідний сигнал $x(t), t \in \mathbb{R}$, має вид

$$y(t) = \int_{-\infty}^{\infty} H(s)x(t-s)ds + u(t), t \in \mathbb{R}.$$

де $u(t), t \in \mathbb{R}$, - випадковий процес, що описує внутрішній шум системи.

## А.5 Кумулянти та моменти гауссівських процесів

Твердження, наведені нижче, є безпосереднім наслідком відомої формули Леонова - Ширяєва - Бриллінджера [7], якщо врахувати, що у центрованих гауссівських векторів усі кумулянти вище другого порядку дорівнюють нулеві.

**Теорема А.5.** (Теорема 2.3.2, [7]) *Нехай $n \in \mathbb{N} \setminus \{1\}$; $\xi_1,...,\xi_n$ і $\eta_1,...,\eta_n$ - сумісно гауссівські випадкові величини; $\mathrm{E}\xi_j = 0$, $\mathrm{E}\eta_j = 0$, $j=1,...,n$. Тоді*

$$\mathrm{E}\prod_{j=1}^{n}(\xi_j\eta_j - \mathrm{E}\xi_j\eta_j) = \sum_{\{D_1,...,D_n\}}\prod_{j=1}^{n}\mathrm{cov}D_j,$$

*де сума береться по всіх невпорядкованих розбиттях $\{D_1,...,D_n\}$ таблиці*

$$2 = \begin{bmatrix} \xi_1 & \eta_1 \\ \xi_2 & \eta_2 \\ \vdots & \vdots \\ \xi_n & \eta_n \end{bmatrix}$$

*на неперетинні двоелементні підмножини $D_j, j=1,...,n$, які не співпадають з її рядками; $\mathrm{cov}D_j = \mathrm{E}\xi_k\eta_l$, якщо $D_j = \{\xi_k,\eta_l\}$.*

Сформулюємо дві характеризаційні властивості кумулянтів та моментів гауссівських центрованих процесів, взятих з монографій [7, 39]

**Теорема А.6.** *Нехай $(Z(\tau), \tau \in \mathbb{R})$ - гауссівський центрований процес з кореляційною функцією $C(\tau_1,\tau_2), \tau_1,\tau_2 \in \mathbb{R}$. Тоді для будь-яких $n \in \mathbb{N}$ та набору $\tau_1,...,\tau_n \in \mathbb{R}$, має місце рівність*

$$\mathrm{cum}(Z(\tau_1),...,Z(\tau_n)) = \begin{cases} 0, & n=1; \\ C(\tau_1,\tau_2), & n=2; \\ 0, & n \geq 3, \end{cases}$$

*де $\mathrm{cum}(Z(\tau_1),...,Z(\tau_n))$ - сумісний кумулянт набору випадкових величин $Z(\tau_j), j=1,...,n$.*

Наступне твердження у літературі відоме як "*формула Іссерліса*" (див., наприклад, [7, 39]).

**Теорема А.7.** *Нехай* $(Z(\tau), \tau \in \mathbb{R})$ *- гауссівський центрований процес з кореляційною функцією* $C(\tau_1, \tau_2), \tau_1, \tau_2 \in \mathbb{R}$. *Тоді для будь-яких* $n \in \mathbb{N}$ *та набору* $\tau_1, ..., \tau_n \in \mathbb{R}$, *має місце рівність*

$$\mathrm{E}\prod_{j=1}^{n} Z(\tau_j) = \begin{cases} 0, & n - \text{непарне число}; \\ \sum_{\{k_1,k_2\},...,\{k_{n-1},k_n\}} \prod_{\substack{j=1 \\ j-\text{непарне}}}^{n-1} C(\tau_{k_j}, \tau_{k_{j+1}}), & n - \text{парне число}, \end{cases}$$

*де в другому випадку сума береться по всіх невпорядкованих розбиттях множини* $\{1,...,n\}$ *на неперетинні двоелементні підмножини* $\{k_1, k_2\},...,\{k_{n-1}, k_n\}$.

### А.6 Слабка збіжність розподілів

Розглянемо два твердження про характеризацію багатовимірного розподілу за допомогою його моментів та кумулянтів (теорема Маркова та її наслідок) [4, 7].

**Теорема А.8.** *Нехай* $(Z(\tau), \tau \in \mathbb{R})$ *- деякий випадковий процес;* $(Z_\alpha(\tau), \tau \in \mathbb{R}), \alpha \in (0, \infty]$, *- сім'я (по* $\alpha$*) випадкових процесів. Якщо для будь-яких* $n \in \mathbb{N}$ *та набору* $\tau_1,...,\tau_n \in \mathbb{R}$, *має місце рівність*

$$\lim_{\alpha \to \infty} \mathrm{E}\prod_{j=1}^{n} Z_\alpha(\tau_j) = \mathrm{E}\prod_{j=1}^{n} Z(\tau_j),$$

*і якщо скінченновимірні розподіли процесу* $Z$ *однозначно визначаються його моментами, тоді* $Z_\alpha \Rightarrow Z$ *при* $\alpha \to \infty$.

**Теорема А.9.** ([4]) *Нехай* $(Z(\tau), \tau \in \mathbb{R})$ *- деякий випадковий процес;* $(Z_\alpha(\tau), \tau \in \mathbb{R}), \alpha \in (0, \infty]$, *- сім'я (по* $\alpha$*) випадкових процесів. Якщо для будь-яких* $n \in \mathbb{N}$ *та набору* $\tau_1,...,\tau_n \in \mathbb{R}$, *має місце рівність*

$$\lim_{\alpha \to \infty} cum(Z_\alpha(\tau_1),...,Z_\alpha(\tau_n)) = cum(Z(\tau_1),...,Z(\tau_n)),$$

*і якщо скінченновимірні розподіли процесу $Z$ однозначно визначаються його кумулянтами, тоді $Z_\alpha \Rightarrow Z$ при $\alpha \to \infty$.*

Оскільки гауссівський розподіл однозначно характеризується своїми моментами (або кумулянтами), то слабка збіжність розподілів випливає зі збіжності відповідних моментів (або кумулянтів).

**А.7      Слабка збіжність розподілів у просторі неперервних функцій**

Сформулюємо деякі терміни, пов'язані з гауссівськими випадковими процесами [14, 68, 34].

Нехай $S$ – параметрична множина – метричний компакт. Функція $\rho(t,s), t,s \in S$ називається *псевдометрикою* на $S$, якщо вона задовольняє всім аксіомам метрики, за винятком того, що множина $\{(t,s) \in S \times S : \rho(t,s) = 0\}$ може бути ширшою за діагональ $\{(t,s) \in S \times S : t = s\}$.

Позначимо через $N_\rho(S,\varepsilon)$ мінімальне число замкнених $\rho$-куль радіуса $\varepsilon > 0$, центри яких лежать всередині $S$, та які покривають $S$. Якщо не існує скінченного покриття для $S$, тоді $N_\rho(S,\varepsilon) = \infty$. Далі, як завжди, $\mathcal{H}_\rho(S,\varepsilon) = \log N_\rho(S,\varepsilon)$ - *метрична ентропія множини $S$ відносно $\rho$*. Для довільного $\beta > 0$ нерівність $\int_{0+} \mathcal{H}_\rho^\beta(S,\varepsilon) d\varepsilon < \infty$ означає, що для деякого (отже, й для всіх) $u > 0$ виконується нерівність $\int_0^u \mathcal{H}_\rho^\beta(S,\varepsilon) d\varepsilon < \infty$.

$C(S)$ – простір дійснозначних неперервних функцій, визначених на $S$, з рівномірною нормою.

Сформулюємо теорему Дадлі [69, 14] для неперервності майже напевно гауссівських процесів.

**Теорема А.10.** *Нехай $Z = (Z(t), t \in \mathbb{R})$ - центрований, неперервний в середньому квадратичному, гауссівський процес. Якщо для будь-якої компактної множини $S \subset \mathbb{R}$ його інтеграл Дадлі є збіжним, тобто*

$$\int\limits_{0+} \mathcal{H}_{d_Z}^{\frac{1}{2}}(S, \varepsilon) d\varepsilon < \infty,$$

*тоді $Z \in C(S)$ м. н. Тут $d_Z(\tau_1, \tau_2) = (\mathrm{E}\,|Z(\tau_1) - Z(\tau_2)|^2)^{\frac{1}{2}}, \tau_1, \tau_2 \in S$, - псевдометрика, породжена середньоквадратичними відхиленнями процесу $Z$.*

Сформулюємо декілька тверджень відносно слабкої збіжності випадкових процесів у просторі неперервних функцій.

**Теорема А.11.** (Лема 4.2.1, [14]) *Нехай $(S, d)$ - метричний компакт та $Z_\alpha = (Z_\alpha(s), s \in S), \alpha \in (0, \infty]$, - сім'я (по $\alpha$) м.н. вибірково неперервних випадкових процесів. Припустимо, що виконуються умови:*

*(i) для кожного $\alpha > 0$ існує така псевдометрика $\rho_\alpha$ на $S$, що*

$$\sup_{\alpha > 0} \sup_{\substack{s,t \in S \\ \rho_\alpha(s,t) \neq 0}} \mathrm{E} \exp\left\{\left|\frac{Z_\alpha(s) - Z_\alpha(t)}{\sqrt{8}\rho_\alpha(s,t)}\right|\right\} < \infty;$$

*(ii) псевдометрика $\rho_\infty(s,t) = \sup_{\alpha > 0} \rho_\alpha(s,t),\ s, t \in \mathbb{R}$, є неперервною відносно рівномірної метрики $d$ на $S$;*

*(iii) $\lim\limits_{u \downarrow 0} \sup\limits_{\alpha > 0} \int\limits_0^u \mathcal{H}_{\rho_\alpha}(S, \varepsilon) d\varepsilon = 0$.*

*Тоді для довільного $\varepsilon > 0$*

$$\lim\limits_{h \downarrow 0} \sup\limits_{\alpha > 0} \mathrm{P}\{\sup_{\substack{s,t \in S \\ d(s,t) < h}} |Z_\alpha(s) - Z_\alpha(t)| > \varepsilon\} = 0.$$

**Теорема А.12.** (Теорема 6.2.2, [14]) *Нехай $Z_\alpha = (Z_\alpha(s), s \in S), \alpha \in (0, \infty]$, - сім'я (по $\alpha$) квадратично-гауссівських випадкових процесів. Тоді справедлива оцінка*

$$\sup_{\alpha>0} \sup_{s,t \in S} \mathrm{E} \exp\left\{ \left| \frac{Z_\alpha(s) - Z_\alpha(t)}{\sqrt{8}\rho_\alpha(s,t)} \right| \right\} < \infty.$$

Зауважимо, що сумісні дискретні та інтегральні корелограми між сумісно стаціонарними центрованими гауссівськими процесами є крадратично-гауссівськими процесами [14].

Сформулюємо теорему Прохорова [5] про слабку збіжність процесів у просторі неперервних функцій.

**Теорема А.13.** *Нехай $Z_\alpha = (Z_\alpha(s), s \in S), \alpha \in (0, \infty]$, - сім'я (по $\alpha$) м.н. неперервних випадкових процесів на $S$. Якщо:*

*(i)* $Z_\alpha \Rightarrow Z$ *при* $\alpha \to \infty$;

*(ii)* $\forall \varepsilon > 0 \ \lim_{h \downarrow 0} \sup_{\alpha > 0} \mathrm{P}\{ \sup_{\substack{s,t \in S \\ d(s,t)<h}} |Z_\alpha(s) - Z_\alpha(t)| > \varepsilon \} = 0,$

*тоді* $Z_\alpha \overset{C(S)}{\Rightarrow} Z$ *при* $\alpha \to \infty$.

**Теорема А.14** ([5]) *Нехай при* $\alpha \to \infty$ $Z_\alpha \overset{C(S)}{\Rightarrow} Z$, *й послідовність неперервних функцій $a_\alpha(s), s \in S$, рівномірно збігається до нуля на $S$. Тоді* $Z_\alpha + a_\alpha \overset{C(S)}{\Rightarrow} Z$ *при* $\alpha \to \infty$.

# ДОДАТОК Б

# СПИСОК ОСНОВНИХ ПОЗНАЧЕНЬ І СКОРОЧЕНЬ

$\mathbb{R}$ — простір дійсних чисел

$\mathbb{N}$ — множина натуральних чисел

$\mathbb{R}^m$ — простір дійсних векторів розмірності $m$, зі стандартною нормою та скалярним добутком

$\mathbb{I}_A$ — індикатор множини $A \subset \mathbb{R}^m$

$B(\mathbb{R}^m)$ — $\sigma$-алгебра борелівських множин з $\mathbb{R}^m$

$(\Omega, \mathfrak{F}, P)$ — основний повний ймовірнісний простір

м.н. — майже напевно

E — математичне сподівання

c$ov$ — коваріація

c$um$ — кумулянт (семіінваріант)

$L_p(\mathbb{R}), p \in [1, \infty)$, — простір комплекснозначних функцій $\phi = (\phi(x), x \in \mathbb{R})$, інтегрованих у $p$-му степені по мірі Лебега, з нормою $\|\phi\|_p = [\int_{-\infty}^{\infty} |\phi(x)|^p \, dx]^{\frac{1}{p}}$

$L_\infty(\mathbb{R})$ — простір комплекснозначних обмежених функцій $\phi = (\phi(x), x \in \mathbb{R})$, з нормою $\|\phi\|_\infty = \sup_{x \in \mathbb{R}} |\phi(x)|$

$Lip_\alpha(\mathbb{R}), \alpha \in [0,1)$, — простір дійснозначних функцій, які рівномірно на $\mathbb{R}$ задовольняють умові Ліпшиця з показником $\alpha$

C$[a,b]$ — простір дійснозначних неперервних функцій визначених на $[a,b] \subset \mathbb{R}$, з рівномірною нормою

$\xi_n \Rightarrow \xi$ — збіжність скінченновимірних розподілів $\xi_n$ до відповідних розподілів $\xi$ при заданому характері $n \to \infty$

$\xi_n \overset{C[a,b]}{\Rightarrow} \xi$ — слабка збіжність розподілів $\xi_n$ до $\xi$ у просторі $C[a,b]$ при заданому характері $n \to \infty$

$A \asymp B$ — A поводиться так само, як і B (стосується збіжності інтегралів)

202